\documentclass[final]{siamltex}
\usepackage{algorithm,algpseudocode}
\usepackage{latexsym, graphicx, epsfig, amsmath, amsfonts,amssymb,bm}
\usepackage{caption, subcaption}
\usepackage{epstopdf}
\usepackage{booktabs}
\usepackage{array}
\usepackage{color}
\usepackage{url}
\allowdisplaybreaks
\def\be{\begin{equation}}
\def\ee{\end{equation}}

\def\v{\mathbf{v}}

\title{Data-Driven Deep Learning of Partial Differential Equations in Modal Space}

\author{Kailiang Wu \and
	Dongbin
	Xiu\thanks{Department of Mathematics,
		The Ohio State University, Columbus, OH 43210, USA.
		{\tt wu.3423@osu.edu, xiu.16@osu.edu.}} 
	Funding: This work was partially supported by AFOSR FA9550-18-1-0102.
}

\begin{document}
\maketitle
\begin{abstract}
We present a framework for recovering/approximating unknown time-dependent partial differential equation (PDE) using its solution data. 
Instead of identifying the terms in the underlying PDE, we seek to approximate the
evolution operator of the underlying PDE numerically. 
The evolution operator of the PDE, defined in infinite-dimensional space, maps the solution from a current time to a future time and completely characterizes 
the solution evolution of the underlying unknown PDE.
Our recovery strategy relies on approximation of the evolution operator in a properly
defined modal space, i.e., generalized Fourier space, in order to reduce the problem to finite dimensions. The finite dimensional approximation is then accomplished by
training a deep neural network structure, which is based on residual  network (ResNet), using the given data. 
Error analysis is provided to illustrate the predictive accuracy of the proposed method. 
 A set of examples of different types of PDEs, including inviscid Burgers' equation that develops discontinuity in its solution, are presented to demonstrate the effectiveness of the
proposed method.
\end{abstract}

% main text
\section{Introduction} \label{sec:intro}

Recently there has been an ongoing research effort to develop
data-driven methods for discovering unknown physical laws.
%The goal is to
%use solution data to
%recover the governing equations for the solution.
Earlier attempts such as
\cite{bongard2007automated,schmidt2009distilling} used
symbolic regression to select the proper physical laws and
determine the underlying dynamical systems. More recent efforts tend
to cast the problem as an approximation problem. In this approach, 
the sought-after governing equation is treated as
an unknown target function relating the data of the state variables to their temporal
derivatives. Methods along this line of approach usually seek
exact recovery of the equations by using certain sparse approximation 
techniques (e.g., \cite{tibshirani1996regression}) from a large set of dictionaries; see, for example, \cite{brunton2016discovering}. 
Studies have  been conducted to deal with noises in data
\cite{brunton2016discovering, schaeffer2017sparse,kang2019ident}, corruptions in data
\cite{tran2017exact}, limited data \cite{schaeffer2017extracting}, 
partial differential equations \cite{rudy2017data,
	schaeffer2017learning}, etc.
Variations of the approaches have been developed in conjunction with other
methods such as
model selection approach
\cite{Mangan20170009}, Koopman theory \cite{brunton2017chaos}, 
Gaussian process regression \cite{raissi2017machine,RAISSI2018125},
and expectation-maximization approach \cite{nguyen2019like}, to name a
few. 
Methods using standard basis functions and without requiring exact
recovery were also  developed for dynamical systems
\cite{WuXiu_JCPEQ18}  and Hamiltonian systems \cite{WuQinXiu2019}.

There is a recent surge of interest in developing methods
using modern machine learning techniques, particularly deep neural networks.
The studies include recovery of
ordinary differential equations (ODEs) \cite{raissi2018multistep,qin2018data,rudy2018deep} and partial differential equations (PDEs)
\cite{long2017pde,raissi2017physics1,raissi2017physics2,raissi2018deep,long2018pde,sun2019neupde}. 
It was shown that
residual network
(ResNet) is particularly suitable for equation recovery, in the sense
that it can be an exact integrator \cite{qin2018data}. 
Neural networks have also been explored for other
aspects of scientific computing, including 
reduced order modeling \cite{HesthavenU_JCP18,doi:10.1063/1.5113494}, solution of
conservation laws \cite{RayHeasthaven_JCP18,wang2019learning}, multiphase flow simulation \cite{wang2019efficient}, 
high-dimensional PDEs 
\cite{han2018solving,KhooLuYing_2018}, %MardtPWN_Nature18,
uncertainty quantification \cite{ChanE_JCP18,TripathyB_JCP18, Zabaras_2018,karumuri2019simulator}, etc.

The focus of this paper is on the development of a general numerical framework for approximating/learning unknown time-dependent PDE.
Even though the topic has been explored in several recent articles, cf.,
\cite{long2017pde,raissi2017physics1,raissi2017physics2,raissi2018deep,long2018pde,sun2019neupde},
the existing studies are relatively limited, as they mostly focus on 
 learning certain types of PDE or identifying the exact terms in the PDE
 from a (large) dictionary of possible terms.
 The specific novelty of this paper is that the proposed method seeks to
 recover/approximate the evolution operator of the underlying unknown
 PDE and is applicable for general class of PDEs. 
The evolution operator completely characterizes the time
 evolution of the solution. Its recovery allows
 one to conduct prediction of the underlying PDE and is effectively equivalent to
 the recovery of the equation. This is an extension of the
 equation recovery work from \cite{qin2018data}, where the flow map of
 the underlying unknown dynamical system is the goal of
 recovery. Unlike the ODE systems considered in \cite{qin2018data},
 PDE systems, which is the focus of this paper, are of infinite dimension. In order to cope with infinite
 dimension, our method first reduces the problem into finite dimensions 
 by utilizing a properly chosen modal space, i.e., generalized Fourier
 space. The equation recovery task is then transformed into recovery
 of the generalized Fourier coefficients, which follow a finite
 dimensional dynamical system. The approximation of the finite
 dimensional evolution operator of the reduced system is then carried
 out by using deep neural network, particularly the residual network
 (ResNet) which has been shown to be particularly suitable for this
 task \cite{qin2018data}.
One of the advantages of the proposed method is that, by focusing on
evolution operator, it
eliminates the need for time derivatives data of the
state variables. Time derivative data, often required by many existing
methods, are difficult
to acquire in practice and susceptible to (additional) errors when computed
numerically. Moreover, the proposed method can cope with solution data
that are more sparsely or unevenly distributed in time.
%
%On the other hand, the proposed framework allows one to approximate the unknown PDEs using data that are coarsely distributed in time. Such data may introduce challenges for some existing methods that require time derivative data, which can be difficult to be extracted accurately.  
Since the proposed framework is rather general, we present several
examples of recovering different types of PDEs. These include linear
advection, linear diffusion, viscous and inviscid nonlinear Burgers'
equations. The inviscid Burgers' equation represents a relatively
challenging problem, as it develops shock over time. Our results show
that the proposed method is able to accurately capture the evolution operator
using  only smooth data during training.
%That is, our training data 
%contain only smooth solutions.
The reconstructed evolution
operator is then able to produce shock structure developed over
time during prediction. Most of our examples are in one
dimensional physical space, as this allows us to easily and thoroughly
examine the solutions and their numerical errors. Our last
example is the recovery of a two-dimensional advection-diffusion
equation. It demonstrates the applicability of the method to multiple
dimensional PDEs.

This paper is organized as follows. After the problem setup in
Section \ref{sec:setup}, we discuss evolution operator and its finite dimensional
representation in Section \ref{sec:finite}. The numerical approach for
learning the evolution operator is then
presented in Section
\ref{sec:method}, along with an error analysis for the predictive
accuracy. Numerical examples are then presented in 
Section \ref{sec:examples} to
demonstrate the properties of the proposed approach.

\section{Problem Setup} \label{sec:setup}

Let us consider a state variable $u(x,t)$, which is governed by 
an unknown autonomous time-dependent PDE system 
\be \label{govern}
\begin{cases}
u_t = {\mathcal L}(u), \quad &(x,t) \in \Omega \times \mathbb R^+,
\\
{\mathcal B} (u) = 0, \quad &(x,t) \in \partial \Omega \times \mathbb
R^+,\\
u(x,0) = u_0(x), \quad & x \in \bar{\Omega}, 
\end{cases}
\ee
where $t$ denotes the time, $x$ is the spatial variable, 
 $\Omega$ is the physical domain, and ${\mathcal L}$ and ${\mathcal B}$
stand for the operators in the equations and boundary
conditions, respectively. Our basic assumption is that the
operator ${\mathcal L}$ is unknown. In this paper, we assume the
boundary conditions are known and focus on 
learning the PDE in the interior of the domain.

We assume data about the solution $u(x,t)$ are available at
certain time instances, loosely called ``snapshots'' hereafter.
%To
%further simplify the discussion here and without loss of generality,
%we assume a certain spatial reconstruction method has been adopted,
%for example, polynomial interpolation, such
%that the solution data can be interpreted as a continuous function
%over $x$.
That is, we have data
\be \label{data}
w(x,t_j) = u(x,t_j) + \epsilon(x, t_j), \qquad j=1,\dots, S,
\ee
where $S\geq 1$ is the total number of snapshots of the solution field
and $\epsilon(x,t_j)$ stands for the noises/errors when the data
are acquired. (Certain reconstruction procedure may be involved in
order to have the snapshot data in the form of \eqref{data}. This,
however, is not the focus of this paper.)
%Note that this term may include the standard random
%noises from measurement, as well as the numerical errors in
%constructing the continuous field $w$ from data of $u$ at discrete
%locations in $x$.

Our goal is to accurately reconstruct the evolution/dynamics of the unknown
governing equation \eqref{govern} via the snapshot data
\eqref{data}. Once an accurate reconstruction is achieved, it can be
used to provide
predictions of the solution.

\section{Finite Dimensional Approximation} \label{sec:finite}

While many of the existing equation learning methods seek
to directly 
approximate or learn the specific form of the governing equations, 
we adopt  a different framework, which seeks to approximate 
evolution operator of the underlying equations. Such an approach
was presented and analyzed in \cite{qin2018data} for recovery of
ODEs. For the PDE recovery problem considered in this paper, our first task is to
reduce the problem from infinite dimension to finite dimension.

\subsection{Evolution Operator}

Without specifying the form of the governing equation, we
loosely assume that for any fixed $t\ge 0$, the solution $u(x,t)$ belongs to 
a Hilbert space $\mathbb V$, with the space norm denoted by $\| \cdot \|_{\mathbb V}$. 
Moreover, we assume the known boundary conditions are linear, i.e.,
the boundary operator $\mathcal B$ in \eqref{govern} is linear.
For many commonly-used boundary conditions, e.g., Dirichlet, Neumann,
or periodic boundary conditions, etc.,  this assumption holds true.

We restrict our attention to autonomous PDEs. Consequently, there
exists an evolution operator
\be \label{evolution_operator}
{\mathcal E}_\Delta: {\mathbb V} \to {\mathbb V}, \qquad {\mathcal E}_\Delta u(\cdot,t)=  u(\cdot,t+\Delta).
\ee
Note that only the time difference, or time lag, $\Delta$ is
relevant, as the time variable $t$ can be arbitrarily shifted.
The evolution operator completely determines the solution over
time. Once it is accurately approximated, one can iteratively apply
the approximate evolution operator to
conduct prediction of the system.

\subsection{Finite Dimensional Evolution Operator}

%As a linear subspace of the separable Hilbert space $L^2(\Omega)$, the space 
%$L^2_{\mathcal B}(\Omega)$ is also a separable Hilbert space. 
%Thus, $L^2_{\mathcal B}(\Omega)$ has a countable  (orthonormal) basis, denoted by $\{\phi_j(x)\}_{j=1}^{+\infty}$. 
%Any function $u \in L^2_{\mathcal B}(\Omega)$ can be uniquely written as
%$$
%u = \sum_{j=1}^{+\infty} \widehat v_j \phi_j,
%$$
%and if the basis $\{\phi_j(x)\}_{j=1}^{+\infty}$ is orthonormal, then the modal expansion coefficients can be expressed as 
%$$
%\widehat v_j = ( u,\phi_j )_{\mathbb V}, \qquad j=1,2,\dots, 
%$$
%which satisfy the Parseval's identity
%$$
%\sum_{n=1}^{+\infty}  \left| {\widehat v}_j \right| ^2 =  \| u \|_{L^2 (\Omega) }^2 < +\infty.
%$$
%This means the sequence $( \widehat v_j )_{j\in \mathbb N_+} $ belongs to the standard $l^2$ space. The mapping $\widehat \Pi $ 

To make the PDE learning problem tractable, we consider a finite dimensional
space ${\mathbb V}_n \subset \mathbb V$. 
%reduce the learning problem to finite dimensions, we represent the solution of the underlying PDEs in a finite-dimensional subspace of $L^2_{\mathcal B}(\Omega)$ with a suitable basis. 
%The basis selection is fairly straightforward, 
%as it can follow the
%standard discretization of a given PDE in the physical domain
%$\Omega$.  
%In practical computations, one should choose a basis to provide
%sufficient resolution to the solution \eqref{data} of the PDE. Basis
%from any standard discretization techniques can be applied, for
%example, global polynomials used in spectral methods or piecewise
%polynomials used in finite elements or finite difference methods. 
%For many commonly seen boundary conditions, e.g., linear and
%stationary Dirichlet, Neumann, or periodic boundary conditions, it is possible
%to define a set of computable basis functions to satisfy the boundary
%conditions. Subsequently, the finite-dimensional modal space becomes
Let
\be
{\bf \Phi}(x) = \left(\phi_1(x), \dots, \phi_n(x)\right)^\dagger, \qquad
n\geq 1,
\ee
be a basis of ${\mathbb V}_n$ and satisfy $\mathcal{B}(\phi_j) = 0$, $1\le j \le n$, i.e.,
%${\mathbb V}_n$ can be expressed as 
 \be \label{V0}
\mathbb V_n = \textrm{span}\{\phi_j:~ \mathcal{B}(\phi_j) = 0, j=1,\dots,n\}.
\ee
An approximation of $u(x,t)$ from $\mathbb V_n$ can be written as 
\be \label{appsolu}
 u_n(x, t) =\sum_{j=1}^n
v_j(t)\phi_j(x) = \langle {\bf v}(t), {\bf \Phi}(x)\rangle,
\ee
where the last equality is written in vector notation after defining
%${\bf \Phi}=(\phi_1,\dots, \phi_n)^\dagger$ and
${\bf v} = (v_1,\dots, v_n)^\dagger$.  
Let $\mathcal{P}_n: {\mathbb V} \to \mathbb V_n$ be a   
projection operator.  For any $t$, we define the projection of the exact solution as 
\be \label{proj}
\widehat u_n(x, t) := {\mathcal P}_n u(x,t)=\langle \widehat {\bf v}(t), {\bf \Phi}(x)\rangle. 
\ee

To approximate the (unknown) infinite dimensional evolution operator
${\mathcal E}_\Delta$ in the finite dimensional space $\mathbb V_n$,
we consider a finite dimensional evolution operator ${\mathcal E}_{\Delta,n}$, which evolves an approximate solution $u_n
\in \mathbb V_n$, i.e. 
\be \label{finite_evo}
{\mathcal E}_{\Delta,n}: \mathbb{V}_n \to \mathbb{V}_n,
\qquad
{\mathcal E}_{\Delta,n}  u_n(\cdot, t) =  u_n(\cdot, t+\Delta).
\ee

In practice, one may choose any suitable finite dimensional operator
${\mathcal E}_{\Delta,n}$, as long as it provides a good
approximation to 
${\mathcal E}_\Delta$, i.e., ${\mathcal E}_{\Delta,n}\approx
{\mathcal E}_\Delta$. 
In this paper, we mostly employ  the following finite dimensional evolution operator 
\be \label{finite_evo2}
\widetilde {\mathcal E}_{\Delta,n} := {\mathcal P}_n {\mathcal E}_\Delta,
\ee
such that
\be
\widetilde {\mathcal E}_{\Delta,n} {v} = {\mathcal P}_n {\mathcal E}_\Delta {v}, \qquad \forall { v}\in \mathbb V_n.
\ee 
%The approximation accuracy of ${\mathcal E}_{\Delta,n}$ defined in \eqref{finite_evo2} is shown in Proposition \ref{prop2}. 
%It indicates that the approximation error is related to the projection error, which may be reduced by choosing a enough large approximation space 
%${\mathbb V}_n$ such that the structure of the exact solution can be
%accurately resolved.
(Note that if one chooses $v = {\mathcal P}_n u$, then this operator
closely resembles the evolution operator of spectral Galerkin method
for solving an {\em known} PDE.)
For this specific choice of ${\mathcal E}_{\Delta,n}$, we have
the following  error bound.
\begin{proposition} \label{prop2}
	Assume the evolution operator ${\mathcal E}_\Delta$
	\eqref{evolution_operator} of the underlying PDE is
	bounded. Let
	$t_k = k\Delta$, $k=0,1,\dots$, and 
	let
	$\varepsilon^{\tt proj}(t_k) := \|   u(\cdot,t_k) - {\mathcal
          P}_n u(\cdot,t_k) \|_{\mathbb V} $ be the projection error
        of the exact solution at $t_k$.
        Consider the approximate evolution operator $\widetilde {\mathcal E}_{\Delta,n}$ defined in \eqref{finite_evo2} and its
        corresponding 
        approximate solution:
	$$
	u_n (\cdot,t_k) = \widetilde {\mathcal E}_{\Delta,n}
        (\cdot,t_{k-1}), \quad k=1,2,\dots,  \qquad u_n (\cdot,0) =
        \widehat u_n(\cdot,0). 
	$$
        Then, the error in the approximate solution satisfies 
	\be \label{solu-err-bound-exactfinite}
	 \|  u_n (\cdot,t_k) -u  (\cdot, t_k) \|_{\mathbb V} \le  \sum_{j=0}^{k} \|  {\mathcal P}_n {\mathcal E}_\Delta \|^{k-j} 
	\varepsilon^{\tt proj}(t_j).
	\ee
\end{proposition}

\begin{proof}
	Let $e(t_k) := \|  u_n (\cdot,t_k) - \widehat u  (\cdot, t_k) \|_{\mathbb V}$. 
	For any $k\ge 1$, we have 
	\begin{align*}
	e(t_k) & = 
	\|  {\mathcal E}_{\Delta,n} u_n (\cdot,t_{k-1}) - {\mathcal P}_n u  (\cdot, t_k) \|_{\mathbb V} 
	\\
	& = \|  {\mathcal P}_n {\mathcal E}_\Delta u_n (\cdot,t_{k-1}) - {\mathcal P}_n {\mathcal E}_\Delta  u  (\cdot, t_{k-1}) \|_{\mathbb V} 
	\\
	& \le \|  {\mathcal P}_n {\mathcal E}_\Delta u_n (\cdot,t_{k-1}) - {\mathcal P}_n {\mathcal E}_\Delta  \widehat u_n  (\cdot, t_{k-1}) \|_{\mathbb V} 
	\\
	& \quad 
	+ \|  {\mathcal P}_n {\mathcal E}_\Delta \widehat u_n (\cdot,t_{k-1}) - {\mathcal P}_n {\mathcal E}_\Delta  u  (\cdot, t_{k-1}) \|_{\mathbb V} 
	\\
	& \le \|  {\mathcal P}_n {\mathcal E}_\Delta\|  e(t_{k-1}) +  \|  {\mathcal P}_n {\mathcal E}_\Delta\|  \varepsilon^{\tt proj}(t_{k-1}).
	\end{align*}
	By recursively applying the above inequality and using $e(t_0)
	= 0$, we obtain 
	$$
	 \|  u_n (\cdot,t_k) - \widehat u  (\cdot, t_k) \|_{\mathbb V} \le  \sum_{j=0}^{k-1} \|  {\mathcal P}_n {\mathcal E}_\Delta \|^{k-j} 
	\varepsilon^{\tt proj}(t_j).
	$$
	The estimate \eqref{solu-err-bound-exactfinite} is further obtained by using
	$$
	\|  u_n (\cdot,t_k) -u  (\cdot, t_k) \|_{\mathbb V} 
	\le  \|  u_n (\cdot, t_k) - \widehat u
	(\cdot, t_k) \|_{\mathbb V} + \|  \widehat u (\cdot, t_k) -u  (\cdot, t_k) \|_{\mathbb V}
	= e(t_k) + \varepsilon^{\tt proj}(t_k).
	$$
\end{proof}

We now define a  linear mapping
\be \label{Pi}
\Pi: \mathbb R^n \to \mathbb V_n, \qquad \Pi {\bf v} = \langle {\bf
  v}, {\bf \Phi} (x) \rangle  , \qquad {\bf v} \in \mathbb R^n,
\ee
which is a bijective mapping whose inverse exists. 
Subsequently, $\Pi: \mathbb R^n \to \mathbb V_n$ is 
an isomorphism. This mapping defines a unique correspondence between a
solution in $\mathbb V_n$ and its modal expansion coefficients in $\mathbb R^n$.

%We now assume there is an (unknown) autonomous governing
%equation system for $\widehat u_n(x, t)$ and define its corresponding finite
%dimensional evolution operator as
%We emphasize that in general the finite dimensional operator ${\mathcal E}_{\Delta,n}$ can only
%be an approximation to the true operator
%${\mathcal E}_\Delta$, i.e., ${\mathcal E}_{\Delta,n}\approx
%{\mathcal E}_\Delta$. This is because the exact representation of the
%infinite dimensional operator ${\mathcal E}_\Delta$ in the finite
%dimensional space typically results in nonlocal dependence and becomes
%non-autonomous.

\begin{proposition} \label{prop}
  Let ${\mathcal E}_{\Delta,n}$ be a finite dimensional evolution operator for $
  u_n\in \mathbb V_n$, as defined in \eqref{finite_evo}, and ${\bf
    v}\in \mathbb R^n$ be its coefficient vector as in \eqref{appsolu},
  then ${\bf v}$ follows an evolution operator
  \be \label{M_evo}
  {\mathcal M}_{\Delta,n}: \mathbb R^n\to \mathbb R^n, \qquad
   {\mathcal M}_{\Delta,n} {\bf v}(t) = {\bf v}(t+\Delta),
   \ee
   and
   \be \label{M}
   {\mathcal M}_{\Delta,n} = \Pi^{-1} {\mathcal E}_{\Delta,n} \Pi,
   \ee
   where $\Pi$ is the linear mapping defined in \eqref{Pi}.
   Furthermore, if ${\mathcal E}_{\Delta,n}$ is defined as in \eqref{finite_evo2}, then
     \be \label{M-our}
   {\mathcal M}_{\Delta,n} =  \Pi^{-1} {\mathcal P}_n {\mathcal E}_\Delta \Pi =: \widetilde {\mathcal M}_{\Delta,n}.
   \ee   
  \end{proposition}

  The proof is a trivial exercise of substituting \eqref{Pi} into
\eqref{finite_evo}.

Therefore, we have transformed the learning of the infinite
dimensional evolution operator ${\mathcal E}_\Delta$
\eqref{evolution_operator} for the true
solution $u\in  \mathbb V$ to the learning of its finite
dimensional approximation ${\mathcal E}_{\Delta,n}$ \eqref{finite_evo} for the
approximate solution $ u_n\in  \mathbb V_n$, which is equivalent to
the learning of the evolution operator $ {\mathcal M}_{\Delta,n}$
\eqref{M_evo} for its
expansion coefficient ${\bf
    v}\in \mathbb R^n$.

\section{Numerical Approach} \label{sec:method}

In this section, we discuss the detail of our PDE learning algorithm.
%We first present the general
%procedure and then discuss the specifics of data set construction and
%neural network structure for the task.
The general procedure consists of the following steps:
\begin{itemize}
  \item Choose a basis for the finite dimensional space
    $\mathbb V_n$ \eqref{V0} and a corresponding projection operator
    \eqref{proj}.
    \item Apply the projection operation and project the snapshot data
      \eqref{data} to $\mathbb V_n$ \eqref{V0} to obtain training data
      in modal space.
%    space. This is accomplished by applying the projection operator \eqref{p the modal expansion
  %  coefficient data for the snapshot training data;
    \item Choose an appropriate deep neural network structure to approximate
      the finite dimensional evolution operator ${\mathcal M}_{\Delta,n}$ \eqref{M_evo}
      and conduct the network training. 
      \item Conduct numerical prediction of the system by advancing
        the learned neural network model for the evolution operator.
      \end{itemize}
We remark that this is a fairly general procedure. One is certainly
not confined to using neural network. Other approximation methods can
also be applied to model the evolution operator. However, modern
neural networks, along with their advanced training algorithms, are
able to handle relatively high dimensional inputs. This feature makes
them more suitable to PDE learning.

\subsection{Basis Selection and Data Set Construction}
\label{sec:data}

The choice of basis functions is fairly straightforward -- any basis
suitable for spatial approximation of the solution data can be
used. These include piecewise polynomials, typically used in  finite
difference or finite elements methods, or orthogonal polynomials used
in spectral methods, etc. The basis should also be sufficiently fine
to resolve the structure of the true solution.

Once the basis functions are selected, one proceeds to employ a
suitable projection operator  $\mathcal{P}_n: {\mathbb V}
\to \mathbb V_n$ to represent the
solution in the finite dimensional form
\eqref{proj}.
This can be accomplished via piecewise
interpolation, as commonly used in finite elements and finite
difference methods, or orthogonal projection, which is often used in spectral
methods.

One of the key ingredients for learning the finite dimensional
evolution operator  ${\mathcal M}_{\Delta,n}$ \eqref{M_evo} is to 
acquire modal vector data in pairs, whose components are separated by
a time lag.
That is, let $J\geq
1$ be the total number of solution data pairs. Then, we define the $j$-th data pair as
\be \label{data-pair}
( {\bf v}_j(0), {\bf v}_j(\Delta_j)), \qquad j=1,\dots, J,
\ee
where $\Delta_j>0$ is the time lag. Note again that for the autonomous systems considered in this paper,
the time difference $\Delta_j$ is the only relavent variable.
%Therefore, one can arbitrarily sample an initial 
%solution and march forward only one step with $\Delta_j$ to obtain the solution data pair. 
Hereafter, we will assume, without loss of
generality, $\Delta_j = \Delta$ is a constant.

\subsubsection{Data Pairing via Snapshot Data Projection}
\label{sec:data0}

When solution snapshots are available in the form of \eqref{data}, we
first identify and create pairs of snapshots that
are separated by the  time lag $\Delta$. That is, we seek to have the
data arranged in the following form
\be \label{u_pairs}
(u^{(j)}(\cdot,t_{k_j}), u^{(j)}(\cdot,t_{k_j}+\Delta)), \qquad
j=1,\dots, J,
\ee
where $J\geq 1$ is the total number of pairs. Note that some (or, even all) of the pairs
may be
originated from the same ``trajectory''. That is, they are obtained
from the original data snapshots \eqref{data} with the same initial
condition. For more effective equation recovery, it is strongly preferred that
the data pairs are originated from different trajectories with a large
number of different initial conditions \cite{qin2018data}.
%In practice, the solution data \eqref{data} may be given and cannot be arbitrarily sampled by the practitioners. In this case, 
%direct projection of arbitrarily given solution data \eqref{data} 
%to the finite-dimensional space $\mathbb V_n$ may not correspond 
%to the operator \eqref{M-our} we define, but also provides a set of training data for learning other finite dimensional evolution operators. 
%Specifically, let $J\geq
%1$ be the total number of given solution data pairs. For $j=1,\dots,
%J$, denote the $j$-th data pair by 
%where $\Delta>0$ is the time lag between the two solution
%states and $t_{k_j}$ can be arbitrary.

Once the snapshot data pairs are constructed, we proceed to project
them onto the finite dimensional space $\mathbb V_n$ by applying the
projection operator \eqref{proj}. This then produces data pairs for
the modal expansion coefficients \eqref{data-pair}, where
%with the properly chosen
%basis functions and projection operator, and obtain their
%corresponding projection coefficients in the pairing \eqref{data-pair}, i.e., 
$$
{\bf v}_j(0) = \Pi^{-1} {\mathcal P}_n u^{(j)}(\cdot,t_{k_j}), \qquad 
{\bf v}_j(\Delta) = \Pi^{-1} {\mathcal P}_n u^{(j)}(\cdot,t_{k_j}+\Delta).
$$

Learning the finite dimensional evolution operator is then conducted
in the modal space in search for ${\mathcal M}_{\Delta,n}$
\eqref{M_evo}. This is 
accomplished by formally solving the following minimization problem:
\be \label{loss}
{\mathcal M}_{\Delta,n} = 
\mathop {\rm argmin} \limits_{ \mathcal{N}_\Delta: \mathbb R^n \to \mathbb R^n }  \frac{1}{J} 
\sum_{j=1}^J \|  \mathcal{N}_\Delta  {\bf v}_j(0) - {\bf v}_j(\Delta)
\|_2^2.
\ee
This corresponds to finding the operator ${\mathcal
  E}_{\Delta,n}$ \eqref{finite_evo} by formally solving
%
%Note that, in general, there may not exist an autonomous system governs the projection coefficients ${\bf v}(t)$, due to the nonlocal dependence. Therefore, one may not expect a finite dimensional evolution operator that exactly satisfies ${\bf v}_j(\Delta) = {\mathcal M}_{\Delta,n} {\bf v}_j(0) $. 
%In this case, the operators ${\mathcal E}_{\Delta,n}$ and ${\mathcal M}_{\Delta,n}$ can be formally defined as 
\be \label{data-E}
{\mathcal E}_{\Delta,n} :=\mathop {\rm argmin} \limits_{ E_\Delta: {\mathbb V}_n \to \mathbb V_n } 
 \frac{1}{J} 
\sum_{j=1}^J \|  E_\Delta \widehat u^{(j)}(\cdot,t_{k_j}) - \widehat u^{(j)}(\cdot,t_{k_j}+\Delta)   \|_{\mathbb V}^2,
%\\ \label{data-M}
%{\mathcal M}_{\Delta,n} := 
%\mathop {\rm argmin} \limits_{ M_\Delta: \mathbb R^n \to \mathbb R^n } \sum_{j=1}^J \|  M_\Delta  ( {\bf v}_j(0) )- {\bf v}_j(\Delta)  \|_2^2,
\ee
where $\widehat u^{(j)}  = {\mathcal P}_n u^{(j)}$.
%Here, note there is no general explicit form of ${\mathcal E}_{\Delta,n}$ and ${\mathcal M}_{\Delta,n}$.

\subsubsection{Data Pairing via Sampling in Modal Space}
\label{sec:data1}

When data collection procedure is in a controlled environment, it is then
possible to directly sample in
the modal space to generate the training data pairs. One example of
such a case is when the unknown PDE is controlled by a black-box
simulation software or a device. This would allow one to possibly
generate arbitrary ``initial'' conditions and then collect their
states at a later time.
%Let us assume that the solution data \eqref{data} can be collected by 
%the practitioners from an arbitrary initial 
%solution. This allows one to effectively sample the data pairs \eqref{data-pair}, as 
%the mapping between each data pair should correspond to the finite dimensional evolution operator \eqref{M-our} we define. 
%Because the evolution operator ${\mathcal E}_{\Delta,n}$ is defined on the finite dimensional space $\mathbb V_n$, it is reasonable to 
%collect the solution data starting from arbitrary initial solution
%$u(x,0) \in \mathbb V_n$. Then, the training data pairs
%\eqref{data-pair} are constructed in the following way:
In this case, the training data pairs \eqref{data-pair} can be
generated as follows.
\begin{itemize}
	\item Sample $J$ points $\{ {\bf v}_j(0)  \}_{j=1}^J$ over a
          domain $D \subset \mathbb R^n$, where 
	$D$ is a region in which one is interested in the solution
        behavior. The sampling can be conducted randomly or by using
        other proper sampling techniques.
        
	\item For each $j=1,\dots, J$, 
	construct initial solution 
	\begin{equation}\label{initialdata}
	u^{(j)}(x,0) = \Pi {\bf v}_j(0)  = \langle {\bf v}_j(0), {\bf \Phi} (x) \rangle.
      \end{equation}
        Then, march
	forward in time for the time lag $\Delta$, by using the
        underlying black-box simulation code or device, and obtain the
        solution snapshots
	$ u^{(j)}(x,\Delta)$.
	%The modal vector data ${\bf v}_j(\Delta)$ in
        %\eqref{data-pair} is then constructed by projecting the
        %solution data $u^{(j)}(x,\Delta) $ into the finite
        %dimensional space $\mathbb V_n$. That is
        Conduct the projection operation \eqref{proj} to obtain
	\begin{equation}\label{dataproj}
	{\bf v}_j(\Delta) = \Pi^{-1} {\mathcal P}_n u^{(j)}(x,\Delta).
	\end{equation}
      \end{itemize}

      The modal expansion coefficients data
      generated in this manner then satisfy, for each $j=1,\dots, J$,
      	%Combining \eqref{initialdata} and \eqref{dataproj} implies 
	$$
	{\bf v}_j(\Delta)
        %= \Pi^{-1} {\mathcal P}_n u^{(j)}(x,\Delta)
        = \Pi^{-1} {\mathcal P}_n {\mathcal E}_\Delta u^{(j)}(x,0)
        = \Pi^{-1} {\mathcal P}_n {\mathcal E}_\Delta \Pi {\bf v}_j(0)
        = \widetilde {\mathcal M}_{\Delta,n} {\bf v}_j(0),
	$$
        where $\widetilde {\mathcal M}_{\Delta,n}$ is  defined in
        \eqref{M-our} and corresponds to  the finite dimensional
        evolution operator defined in \eqref{finite_evo2}.
%	for noiseless data. When noises are involved, one has ${\bf v}_j(\Delta) 
%	\approx {\mathcal M}_{\Delta,n} {\bf v}_j(0)$. 
%	This indicates that the mapping between each data pair 
%	is indeed associated with 
%	the finite dimensional evolution operator $\mathcal M_\Delta$ defined in \eqref{M-our}, as desired.  

We remark that this procedure also  generates solution pairs
          along the way, i.e,
	$$
	( u^{(j)}(\cdot,0), u^{(j)}(\cdot,\Delta)  ), \qquad u^{(j)}(\cdot,0) \in \mathbb V_n, \qquad j=1,\dots, J.
	$$
        Note that here $u^{(j)}(\cdot,0) \in \mathbb V_n$ due to
        \eqref{initialdata} via the direct sampling in the modal space.
         These pairs are then different from
         those in \eqref{u_pairs}, whose components are not
         necessarily in $\mathbb V_n$.

\subsection{Neural Network Modeling of Evolution Operator}

Once the training data set \eqref{data-pair} becomes available, one
can proceed to learn the unknown governing equation. This is achieved
by learning the finite dimensional evolution operator relating the
solution coefficients in the data pairs.
As shown in Proposition \ref{prop}, finding the finite dimensional evolution
operator of the unknown PDE \eqref{finite_evo} is equivalent to find the evolution
operator \eqref{M} for the modal expansion coefficient
vectors.
Suppose the solution coefficient vectors formally follow an unknown
autonomous system,
%For a general finite dimensional evolution operator 
%${\mathcal M}_{\Delta,n}$ defined in \eqref{M_evo}, assume that  
%${\bf f}$ is the unknown right-hand-side of the corresponding equations governing the
%evolution of the approximate expansion coefficients ${\v}$, i.e.,
$d\v/dt = {\bf f}(\v)$, we then have, by using mean-value theorem,
\be \label{rela2}
\begin{split}
	\v(\Delta) &={\mathcal M}_{\Delta,n} ( \v(0)  ) 
	= \v(0) + \int_0^{\Delta}  {\bf f} ( \v(t) )   dt \\
	%= \v(0) + \Delta  \cdot  {\bf f} ( \v(\tau) )   \\
	& = \v(0) + \Delta \cdot  {\bf f} ( {\mathcal M}_{\tau,n} \v(0) ), \quad 0\leq \tau\leq \Delta,
\end{split}
\ee
where ${\mathcal M}_{\Delta,n}$ is the evolution operator defined in
\eqref{M} and $\tau$ depends on $\Delta$ and the form of the operator
$\mathcal L$ in \eqref{govern}.
%Let ${\bf f}$ be the unknown right-hand-side of the equations governing the
%evolution of the expansion coefficient ${\v}$, i.e., $d\v/dt = {\bf f}(\v)$. 
Furthermore, if ${\mathcal M}_{\Delta,n}$ takes the form of \eqref{M-our}, 
we then have
%\be
%\begin{split}
%	\v(\Delta) &={\mathcal M}_{\Delta,n} ( \v(0)  ) \\
%	&= \v(0) + \int_0^{\Delta}  {\bf f} ( \v(t) )   dt \\
%	&= \v(0) + \Delta  \cdot  {\bf f} ( \v(\tau) )   \\
%	& = \v(0) + \Delta \cdot  {\bf f} ( {\mathcal M}_\tau \v(0) ) , \quad 0\leq \tau\leq \Delta,
%\end{split}
%\ee
\be \label{rela1}
\begin{split}
	\v(\Delta) &=\widetilde {\mathcal M}_{\Delta,n} ( \v(0)  ) = \Pi^{-1} {\mathcal P}_n {\mathcal E}_\Delta \Pi  \v(0)   \\
	& = \Pi^{-1} {\mathcal P}_n {\mathcal E}_\Delta \langle {\bf v}(0), {\bf \Phi} (x) \rangle \\
	&= \Pi^{-1} {\mathcal P}_n  \left( \langle {\bf v}(0), {\bf \Phi} (x) \rangle + \int_0^\Delta {\mathcal L}    (  {\mathcal E}_t \Pi \v (0) ) dt \right) \\
	&=  \Pi^{-1} {\mathcal P}_n  \langle {\bf v}(0), {\bf \Phi} (x) \rangle + \Delta  \cdot  
	\Pi^{-1} {\mathcal P}_n     {\mathcal L}    (  {\mathcal E}_\tau \Pi \v (0) )    \\
	& = \v(0) + \Delta  \cdot  
	\Pi^{-1} {\mathcal P}_n     {\mathcal L}    (  {\mathcal E}_\tau \Pi \v (0) ), \quad 0\leq \tau\leq \Delta.
\end{split}
\ee
%where $\tau$ depends on $\Delta$ and the form of the operator $\mathcal L$ in \eqref{govern}. 
% One can see that, as long as $\Delta$ is reasonably small, the underlying evolution operator is close to identity. 
%Let $ {\mathcal F}_\Delta u(x,0):=  {\mathcal L} ( {\mathcal E}_\tau u(x,0) )$. Then 
%for the initial state satisfying $u(x,0) =  \Pi {\bf v}(0)$,  
%using \eqref{exactobj} gives 
%$$
%{\bf v} ( \Delta ) = {\mathcal M}_{\Delta,n}  {\bf v}(0) =  \Pi^{-1} {\mathcal P}_n {\mathcal E}_\Delta \Pi {\bf v}(0) = {\bf v}(0) + \Delta \cdot \Pi^{-1} {\mathcal P}_n {\mathcal F}_\Delta \Pi {\bf v}(0),
%$$
%which indicates that the reduced evolution operator ${\mathcal M}_{\Delta,n}$ is also close to identity. 

These relations suggest that, when the time lag $\Delta$ is small, it is
natural to adopt the residue network (ResNet) \cite{he2016deep}  to model the
evolution operator. This approach was presented and systematically
studied in \cite{qin2018data}, for learning of unknown dynamical systems
from data.
%Therefore, it is natural
%to use ResNet, which explicitly introduces the identity operator, to approximate the
%``residue'' of the reduced evolution operator. 
%
%By using the ResNet type neural network structure from
%\cite{qin2018data}, we seek to construct a nonlinear mapping between
%the expansion coefficients at two time instance. That is,
%The DNN structure used for learning the reduced evolution operator ${\mathcal M}_{\Delta,n}$ follows
%from the work of \cite{qin2018data}. It uses residual network (ResNet)
%\cite{he2016deep} as the basic building block.
%%%
The block ResNet structure for equation learning from
\cite{qin2018data} takes the following form
\be \label{Res}
{\bf v}(\Delta) = {\bm {\mathcal N}} ( {\bf v}(0);\Theta),
\ee
where ${\bm {\mathcal N}}$ denotes the nonlinear operator defined by
the underlying neural network with parameter set $\Theta$. For block
ResNet with $K\geq 1$ ResNet blocks, the
operator %${\bm {\mathcal N}}$ takes the following form
$$
 {\bm {\mathcal N}}  = \big( {\bf I} + {\bf N}(\bullet ;\Theta_{K-1}) \big) \circ \cdots \circ 
\big( {\bf I} + {\bf N}(\bullet ;\Theta_{0}) \big),
$$
with $\bf N$ is standard fully connected feedforward neural
network, and 
$\Theta=\{\Theta_i\}_{0\le i \le K-1} $ are the parameters (weights and biases)
in each ResNet block. (See \cite{qin2018data} for more details.) The
parameters are determined by minimizing the loss function
\eqref{loss}, i.e.,
\be \label{loss1}
L(\Theta) = \frac{1}{J}\sum_{j=1}^J \left\|  {\bm {\mathcal N}} ( {\v}_j(0); \Theta)   - \v_j(\Delta) \right\|_2^2,
\ee 
where $\| \cdot \|_2$ denotes vector 2-norm. 
%and $({\bf z}_0^{(j)}, \z^{(j)}_\Delta)$ are training data pairs in \eqref{set}. %hereafter. 
Let
$\Theta_*$ be the trained parameters, after satisfactory minimization
of the loss function.
A learned model is then constructed in the form of ${\bm {\mathcal N}} 
(\cdot; \Theta_*)$, which provides an
approximation to the evolution operator \eqref{M_evo}, 
\be \label{N_trained}
{\bm {\mathcal N}}(\cdot, \Theta_*) \approx {\mathcal M}_{\Delta,n}. 
\ee

The trained network model can then be used to provide prediction of
the system \eqref{govern}. For an arbitrary initial condition $u(x,0)\in \mathbb V $, we first conduct its projection \eqref{proj} and
obtain
\begin{equation}\label{model-initial}
\widehat{{\bf v}} (0) = \Pi^{-1}  \mathcal P_n u(\cdot,0),
\end{equation}
where $\Pi$ is the bijective mapping defined in \eqref{Pi}.
We then iteratively apply the neural network model \eqref{N_trained} to
obtain approximate solutions for the modal expansion coefficients at time instances $t_k = k\Delta$,
\begin{equation}\label{model-predict}
  \left\{
  \begin{split}
    &\widetilde {\v} (0) = \widehat{{\bf v}} (0), \\
%\begin{cases}
%\widetilde {\bf v} ( \Delta) = {\bm {\mathcal N}} ( {\bf v} (0) ; \Theta_*)
%\\
%\widetilde {\bf v} ( 2 \Delta) = {\bm {\mathcal N}} ( \widetilde {\bf v} ( \Delta) ; \Theta_*)
%\\
%\cdots
%\\
&\widetilde {\v} ( t_{k+1}) = {\bm {\mathcal N}} ( \widetilde {\v}
(t_k) ; \Theta_*), \qquad k=0,\dots.
\end{split}\right.
%\\
%\cdots
%\end{cases}
\end{equation}
The predicted solution fields $\widetilde u_n(x, t_k)$ are then obtained by
\begin{equation}\label{model-final}
\widetilde u_n(x, t_k) = \Pi \widetilde{\v}(t_k) = \langle
\widetilde {\bf v} (t_k), {\bf \Phi} (x) \rangle, \qquad k=1,\dots.
\end{equation}

\subsection{Error Analysis}\label{sec:analysis}

We now derive an error bound for the predicted solution
\eqref{model-final} of our neural network model, when
the finite dimensional evolution operator  \eqref{M-our} is the case, i.e.,
$ {{\mathcal M}}_{\Delta,n} =\widetilde {\mathcal M}_{\Delta,n} = \Pi^{-1} 
{\mathcal P}_n {\mathcal E}_\Delta \Pi $. 
For notational convenience, we assume that the basis functions
$\{\phi_j(x)\}_{j=1}^n$ of $\mathbb V_n$ are orthonormal. (Note that
non-orthogonal basis can always be orthogonalized via Gram-Schmidt procedure.)
Then, the bijective mapping $\Pi: \mathbb R^n \to \mathbb V_n$
defined in \eqref{Pi} is 
an isometric isomorphism and satisfies
$$
\|  \Pi {\bf v} \|_{\mathbb V} = \| {\bf v} \|_2, \qquad  {\bf v} \in \mathbb R^n.
$$
%According to \eqref{exactobj}, the training data \eqref{datapairs}, and the loss function \eqref{loss}, %our
%neural network is trained to learn the mapping from ${\bf v}(0)$ to ${\bf v} ( \Delta )$. Therefore, 
%the trained DNN operator ${\bm {\mathcal N}}(\cdot, \Theta_*)$ is an approximation to the ``reduced'' evolution operator ${\mathcal M}_{\Delta,n}$:
%$$
%{\bm {\mathcal N}}(\cdot, \Theta_*) \approx {\mathcal M}_{\Delta,n} = \Pi^{-1} {\mathcal P}_n {\mathcal E}_\Delta \Pi. 
%$$ 
Also, since it is well known that neural networks are universal
approximator for a general class of functions, we assume that the
training error in \eqref{N_trained} is bounded. We
then state the following result.

\begin{theorem}\label{thm:errorbound}
	Assume the evolution operator ${\mathcal E}_\Delta$
        \eqref{evolution_operator} of the underlying PDE is
        bounded. Also, assume the trained neural network model ${\bm
          {\mathcal N}}={\bm {\mathcal N}}(\cdot, \Theta_*)$
        is bounded, and its approximation error \eqref{N_trained} is
        bounded and denote
	$ \varepsilon^{\tt DNN} :=  \| {\bm {\mathcal N}} - \Pi^{-1}
        {\mathcal P}_n {\mathcal E}_\Delta \Pi \| < +\infty $. Let
        $t_k = k\Delta$, $k=0,1,\dots$, and let
	$\varepsilon^{\tt proj}(t_k) := \|   u(x,t_k) - {\mathcal P}_n u(x,t_k) \|_{\mathbb V} $ be the projection error of the exact solution at $t_k$. 
	Then the following error bounds hold: 
	\begin{equation}\label{coef-err-bound}
	\| \widetilde{\bf v} (t_k) - \widehat {\bf v}  (t_k) \|_2 \le \sum_{j=0}^{k-1} \| {\bm {\mathcal N}} \|^{k-1-j} \Big(  \varepsilon^{\tt DNN}   \| \widehat {\bf v} (t_j) \|_2 
	+ \varepsilon^{\tt proj}(t_j) \|  {\mathcal P}_n {\mathcal E}_\Delta \|   \Big),
	\end{equation}
	where $\widehat \v=\Pi^{-1}\widehat{u}$ is the modal expansion coefficient vector of the
        projected exact solution \eqref{proj} and $\widetilde{\v}$ is the
        coefficient vector predicted by the neural network model
        \eqref{model-predict}, and 
	\begin{equation}\label{solu-err-bound}
	\| \widetilde u_n (\cdot,t_k) -u  (\cdot, t_k) \|_{\mathbb V} \le \varepsilon^{\tt proj} (t_k) + \sum_{j=0}^{k-1} \| {\bm {\mathcal N}} \|^{k-1-j} \Big(  \varepsilon^{\tt DNN}   \| \widehat {\bf v} (t_j) \|_2 
	+ \varepsilon^{\tt proj}(t_j) \|  {\mathcal P}_n {\mathcal E}_\Delta \|   \Big),
      \end{equation}
      where $\widetilde{u}$ is the solution state predicted by the
      trained neural network model \eqref{model-final}.
\end{theorem}

\begin{proof}
	Let $e(t_k):= \| \widetilde{\bf v} (t_k) - \widehat {\bf v}  (t_k)
        \|_2$. Note that, for $k=1,\dots,$ 
	$$
	\widehat {\bf v} ( t_k ) = \Pi^{-1}  \widehat u_n(\cdot, t_k) = \Pi^{-1}
        {\mathcal P}_n u(\cdot, t_k) 
	= \Pi^{-1} {\mathcal P}_n {\mathcal E}_\Delta u(\cdot, t_{k-1}).  
	$$
	We then have 
	\begin{align*}
	e(t_k) & = \| \widetilde{\bf v} (t_k) - \Pi^{-1} {\mathcal P}_n {\mathcal E}_\Delta u(\cdot,t_{k-1}) \|_2 
	\\
	& \le \| \widetilde{\bf v} (t_k) - \Pi^{-1} {\mathcal P}_n
   {\mathcal E}_\Delta \widehat u_n (\cdot, t_{k-1}) \|_2 
	\\
	& \quad 
	+ \| \Pi^{-1} {\mathcal P}_n {\mathcal E}_\Delta \widehat u_n
   (\cdot, t_{k-1}) - \Pi^{-1} {\mathcal P}_n {\mathcal E}_\Delta u(\cdot, t_{k-1}) \|_2 
	\\
	&
	= \|  {\bm {\mathcal N}} ( \widetilde{\bf v} ( t_{k-1}); \Theta_*) - \Pi^{-1} {\mathcal P}_n {\mathcal E}_\Delta \Pi  {\bf v} ( t_{k-1}) \|_2 
	\\
	& \quad + \|  {\mathcal P}_n {\mathcal E}_\Delta \widehat u_n (\cdot,t_{k-1}) -  {\mathcal P}_n {\mathcal E}_\Delta u(\cdot,t_{k-1}) \|_{\mathbb V}
	\\
	& \le \|  {\bm {\mathcal N}} ( \widetilde{\bf v} ( t_{k-1}); \Theta_*) - {\bm {\mathcal N}} ( \widehat{\bf v} ( t_{k-1}); \Theta_*) \|_2 
	\\
	& \quad + \|  {\bm {\mathcal N}} ( \widehat {\bf v} ( t_{k-1}); \Theta_*) - \Pi^{-1} {\mathcal P}_n {\mathcal E}_\Delta \Pi  \widehat {\bf v} ( t_{k-1}) \|_2 
	\\
	& \quad + \|  {\mathcal P}_n {\mathcal E}_\Delta \| \|
   \widehat u_n (\cdot, t_{k-1}) -   u(\cdot, t_{k-1}) \|_{\mathbb V}
	\\
	& \le \| {\bm {\mathcal N}} \| e(t_{k-1}) + \varepsilon^{\tt DNN}  \|\widehat {\bf v} (t_{k-1}) \|_2
	+ \|  {\mathcal P}_n {\mathcal E}_\Delta \|  \varepsilon^{\tt proj}(t_{k-1}).
	\end{align*}
	By recursively applying the above inequality and using $e(t_0)
        = 0$, we obtain
	%\begin{equation}\label{error-ek}
	%e(t_k) \le \sum_{j=0}^{k-1} \| {\bm {\mathcal N}} \|^{k-1-j} \Big( \varepsilon^{\tt DNN}  \|\widehat {\bf v} (t_j) \|_2 
	%+ \varepsilon^{\tt proj}(t_j)  \|  {\mathcal P}_n {\mathcal E}_\Delta \| \Big),
      %\end{equation}
        %which is
        \eqref{coef-err-bound}. The proof is then complete by using
\begin{equation*}
	\begin{aligned}
	& \| \widetilde u_n (\cdot,t_k) -u  (\cdot, t_k) \|_{\mathbb V} 
	\\
	& \quad \le  \| \widetilde u_n (\cdot, t_k) - {\mathcal P}_n u
   (\cdot, t_k) \|_{\mathbb V} + \|  {\mathcal P}_n u (\cdot, t_k) -u  (\cdot, t_k) \|_{\mathbb V}
	\\
	& \quad 
	= \| \widetilde{\bf v} (t_k) -\widehat {\bf v}  (t_k) \|_2 + \varepsilon^{\tt proj}(t_k)
	= e(t_k) + \varepsilon^{\tt proj}(t_k).
	\end{aligned}
\end{equation*}
\end{proof}

Theorem \ref{thm:errorbound} indicates that the prediction error of the network model is affected by 
the approximation error of the neural network and the projection error, which is determined by the approximation space $\mathbb V_n$ and the regularity of the solution.

\section{Numerical Examples} \label{sec:examples}

In this section, we present numerical examples to verify the properties
of the proposed method. For benchmarking purpose, the true governing
PDEs are known in all of them examples. However, we use the true governing
equations only to generate training data, particularly by following
the procedure in Section \ref{sec:data1}. Our proposed learning method is
then applied to these synthetic data and produces a trained neural
netwok model for the underlying PDE. We will then use the neural
network model to conduct numerical predictions of the solution and
compare them against the reference solution produced by the governing
equations. Numerical errors will be reported, in term of relative
errors between the neural network prediction solutions and the
reference solutions.

The governing equations considered in this
section include: linear advection equation, linear diffusion equation, nonlinear viscous
Burgers' equation, nonlinear inviscid Burgers' equation, the last of
which produces shocks. We primarily focus on one dimension in physical
space with noiseless data, in order to conduct detailed and thorough examination of the
solution behavior. Data noises are introduced in one example to
demonstrate the applicability of the method.
%and do not introduce noises in the data. This is to allow us to
%conduct detailed and thorough examination of the solution behavior.
A
two-dimensional advection-diffusion problem is also presented
demonstrate the applicability of the method to multiple
dimensions.

In all examples, we employ global orthogonal polynomials as the basis
functions to define the finite dimensional space for training.
For benchmarking purpose, all training data are generated in modal
space, as described in Section \ref{sec:data1}. We also rather
arbitrarily impose a decay condition in certain examples to ensure the
higher modes are smaller, compared to the lower modes. This
effectively poses a smoothness condition on the training data.

All of our network models are
trained via minimizing the loss function \eqref{loss} and by using the
open-source Tensorflow library \cite{tensorflow2015}.
%We train the model by minimizing the mean squared loss function \eqref{loss}. The optimization problem is solved with the Adam algorithm \cite{Adam}.
The training data set is divided into mini-batches of size $10$, and the learning rate is taken as $0.001$.  
%And
%we typically train the model for
%$500$ epochs and reshuffle the training data in each epoch.
The network eights are initialized randomly from Gaussian
distributions and biases are initialized to zeros. 
%All the model
%construction, training and prediction are performed with the
%open-source Tensorflow library \cite{tensorflow2015}. 
Upon successful training, the neural network models are then marched
forward in time with the
$\Delta$ time step, as in \eqref{model-predict}. The results, labelled as
``prediction'', are compared against the reference solutions, labelled
as ``exact'', of the true underlying governing equations. 

%In the following figures, 
%the true solutions of the governing equations will be labeled by ``exact'', 
%and the prediction of the neural network models will be labeled by ``prediction''. The expansion coefficients given by 
%the orthogonal projection of the exact solutions will be labeled by ``optimal''. 

%We march the trained network systems
%up to $t\gg \Delta$ to examine their (relatively) long-term
%behaviors. For the two linear examples, we set $t=2$; and for the two
%nonlinear examples, we set $t=20$.

\subsection{Example 1: Advection Equation}

%We first study two linear ODE systems, as textbook examples. In both
%examples, our one-step ResNet method has 3 hidden layers, each of
%which has 30 neurons. For the multi-step RT-ResNet and RS-ResNet
%methods, they both have $3$ ResNet blocks ($K=3$), each of which
%contains 3 hidden layers with 20 neurons in each layer.

We first consider a one-dimensional advection equation with periodic boundary condition: 
\begin{equation}
	\label{eq:example1}
	\begin{cases}
		u_t + u_x = 0,\quad (x,t) \in (0,2\pi) \in \mathbb R^+, \\
		u(0,t)=u(2\pi,t), \quad t\in \mathbb R^+.
	\end{cases}
\end{equation}

\begin{figure}[htbp]
	\centering
{\includegraphics[width=0.6\textwidth]{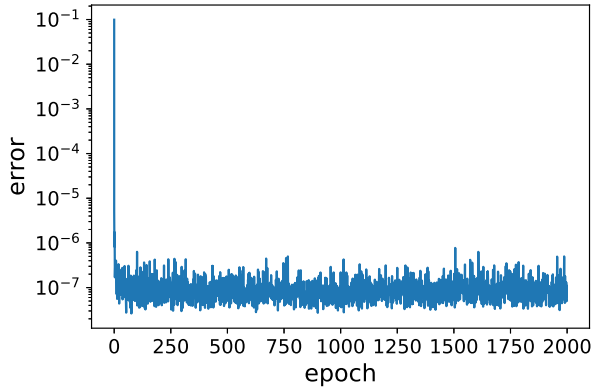}}
	\caption{\small
		Example 1: Training loss history.
	}\label{fig:ex1_loss}
\end{figure}

\begin{figure}[htbp]
	\centering
	{\includegraphics[width=0.48\textwidth]{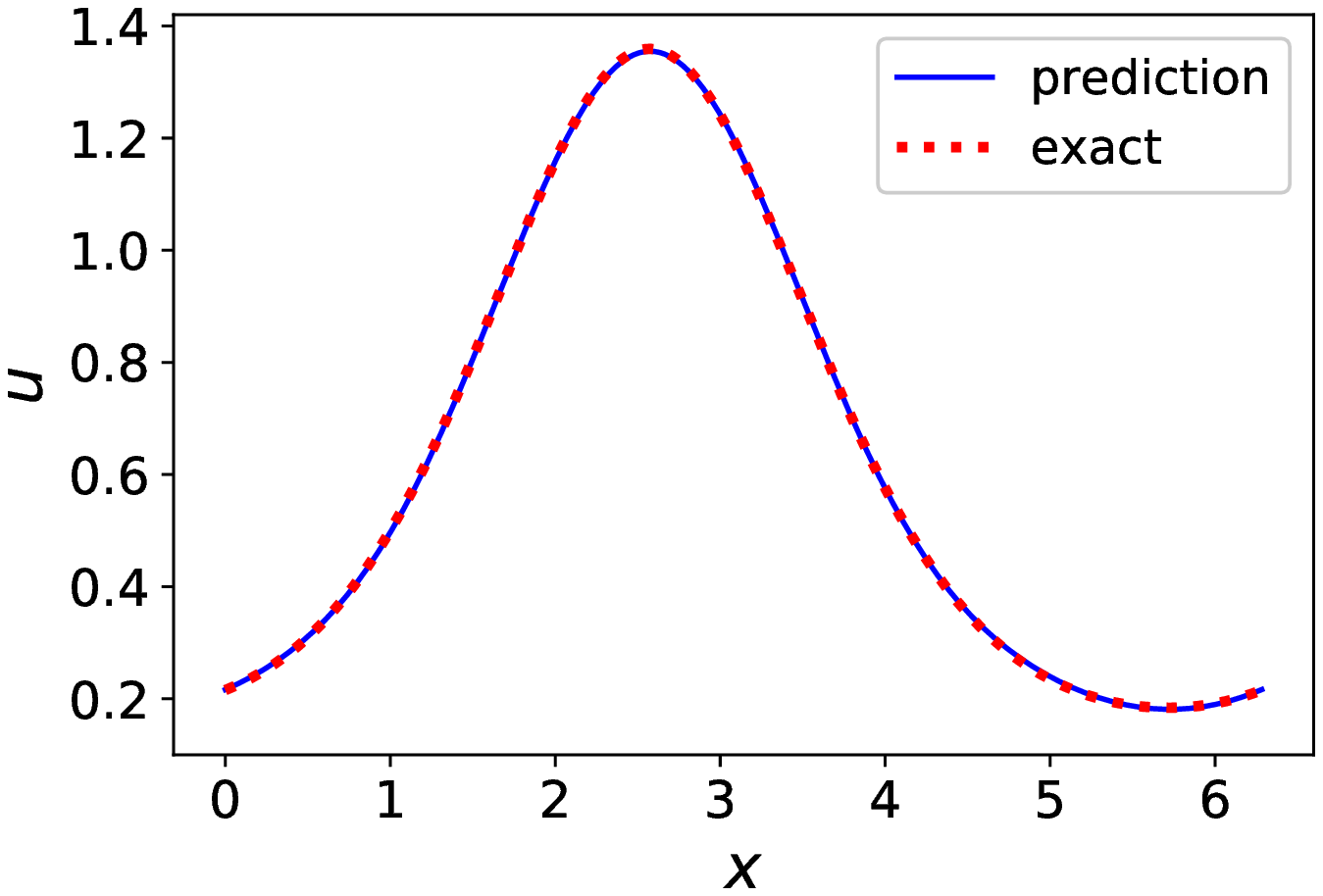}}
	{\includegraphics[width=0.48\textwidth]{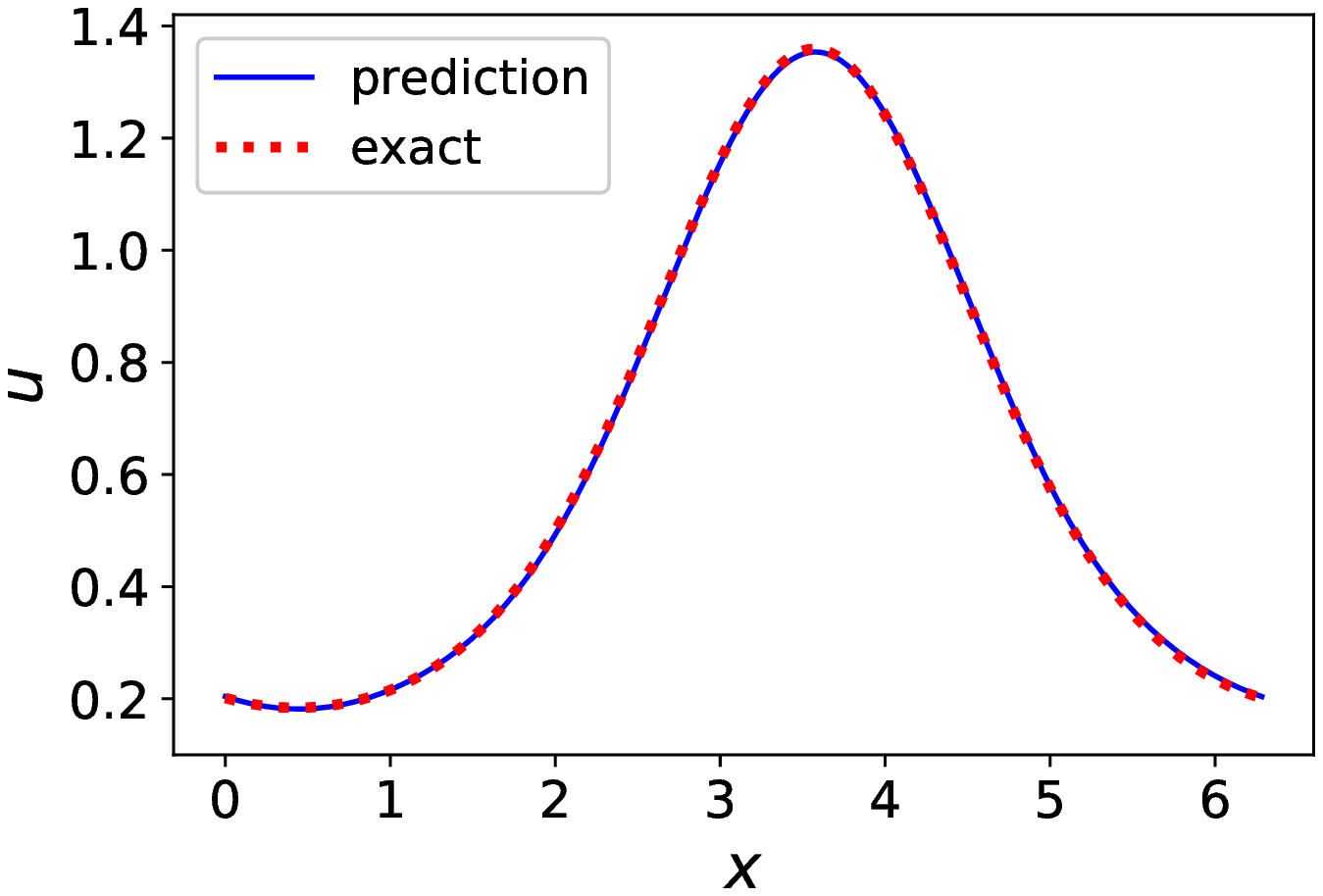}}
	{\includegraphics[width=0.48\textwidth]{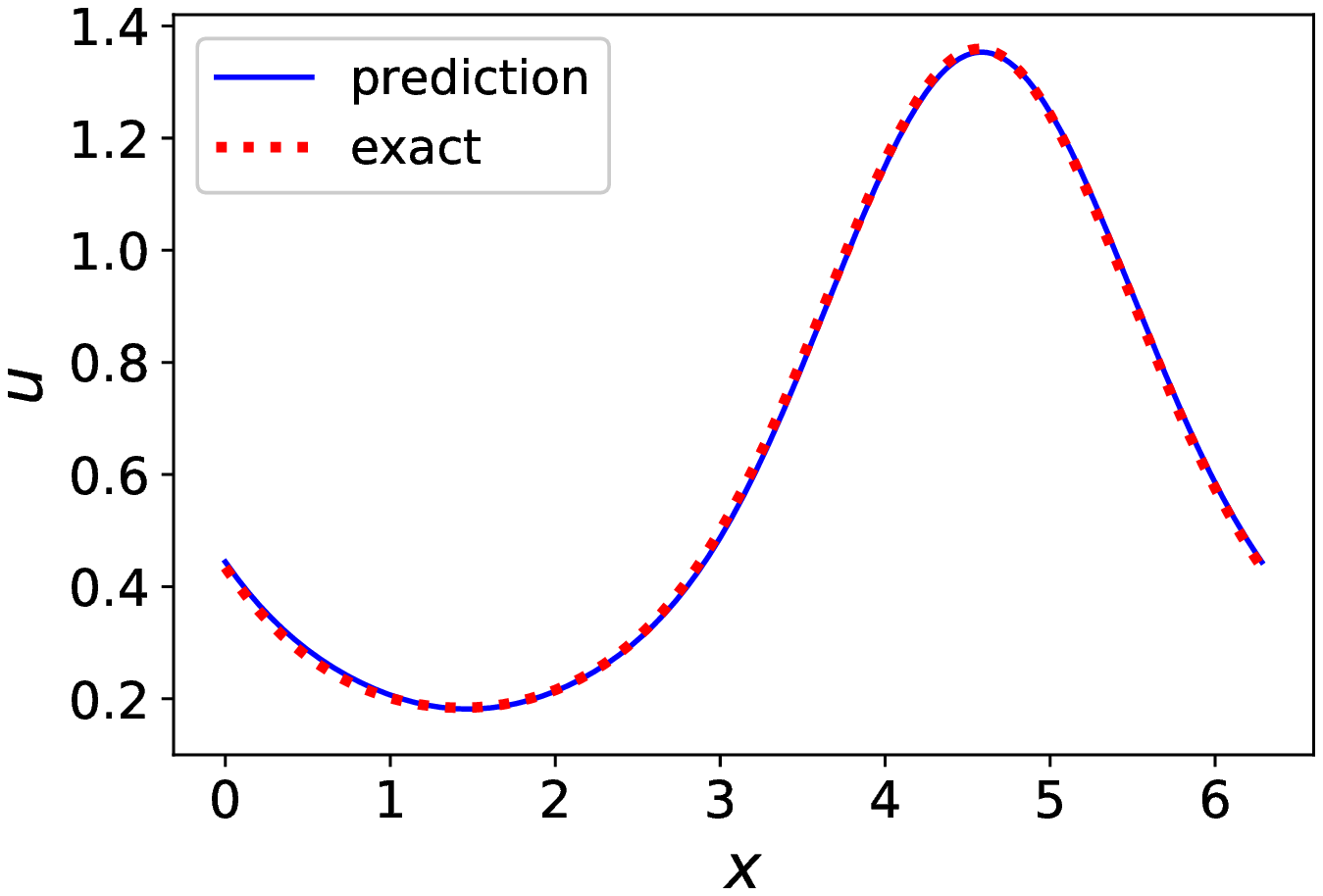}}
	{\includegraphics[width=0.48\textwidth]{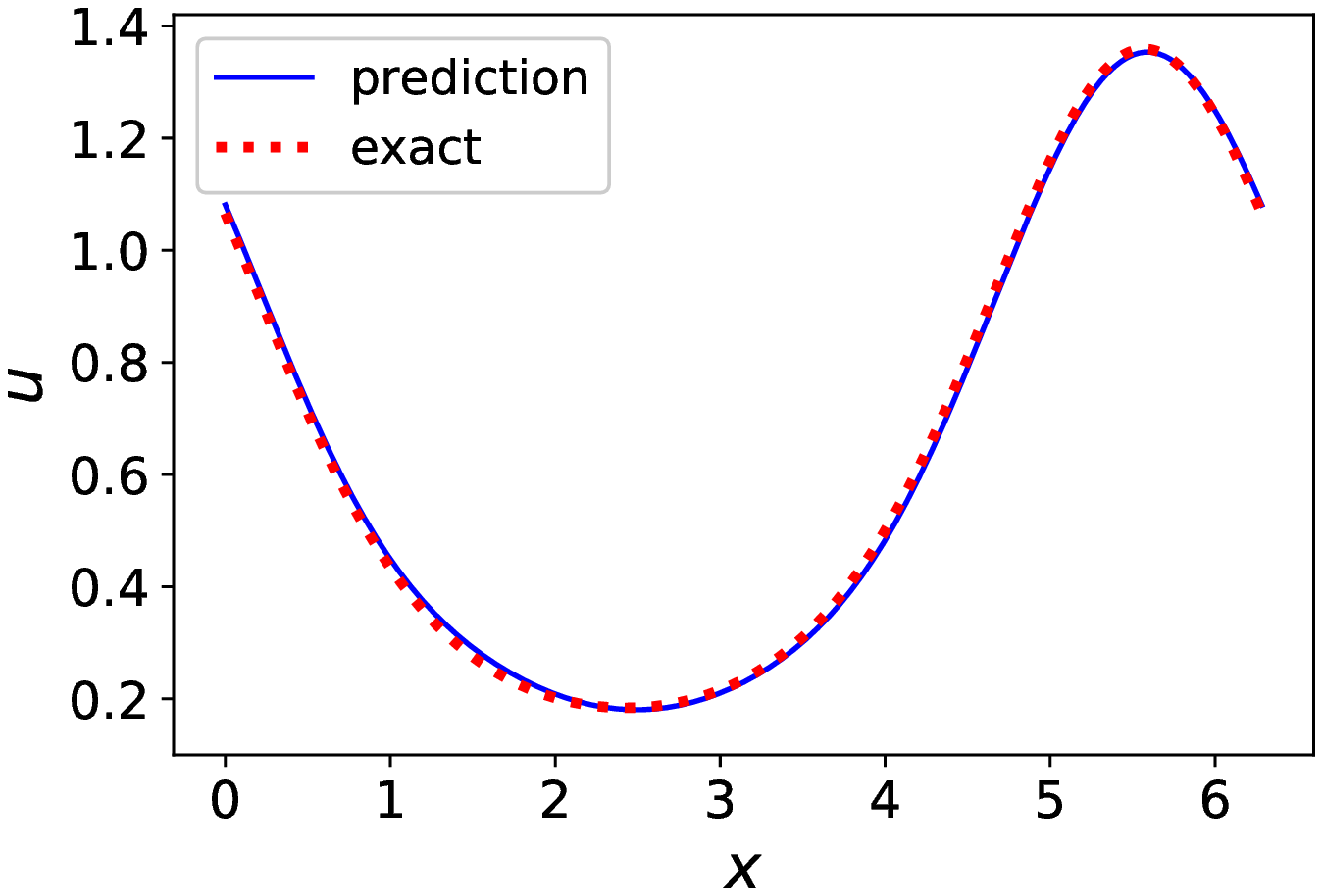}}
	{\includegraphics[width=0.48\textwidth]{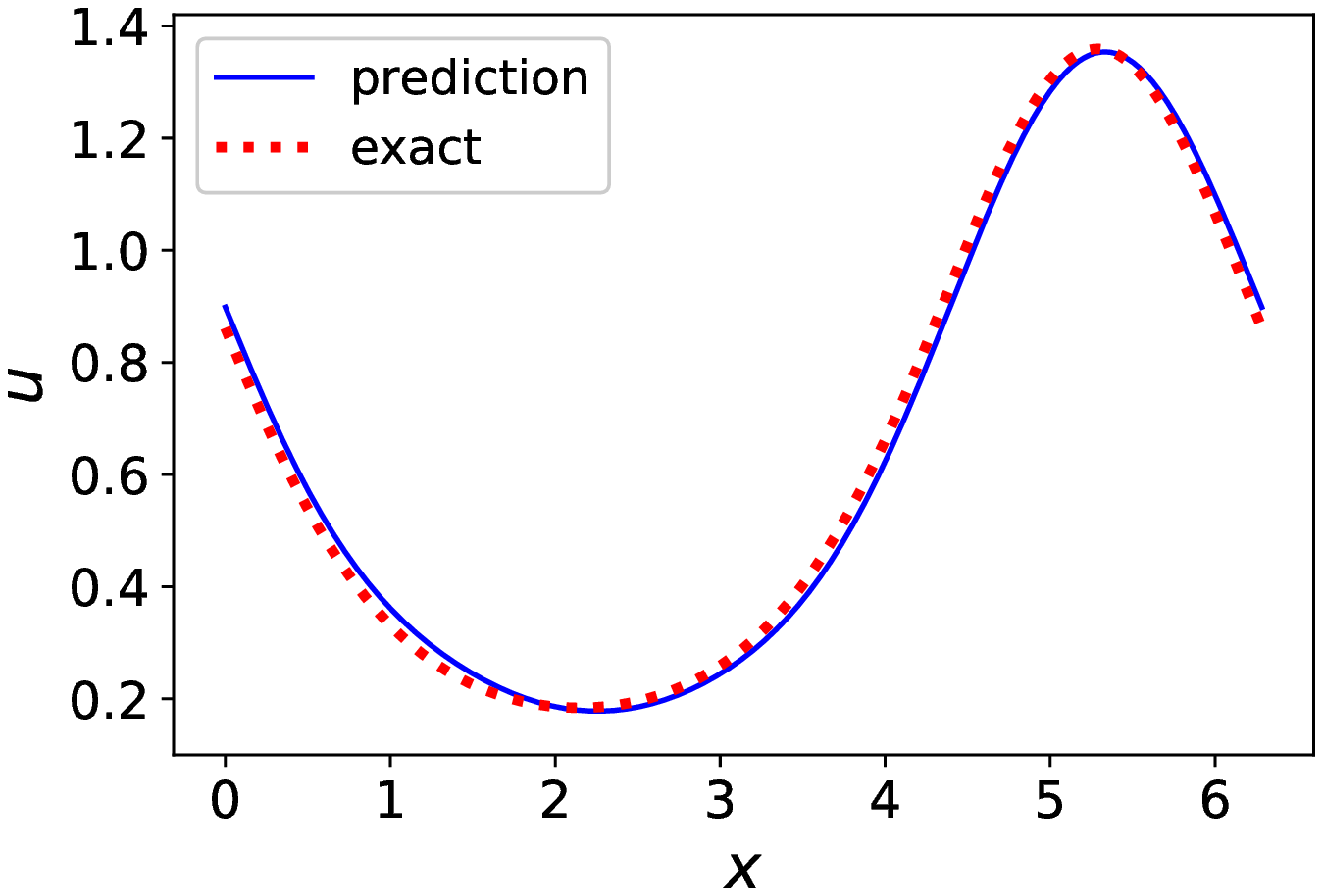}}
	{\includegraphics[width=0.48\textwidth]{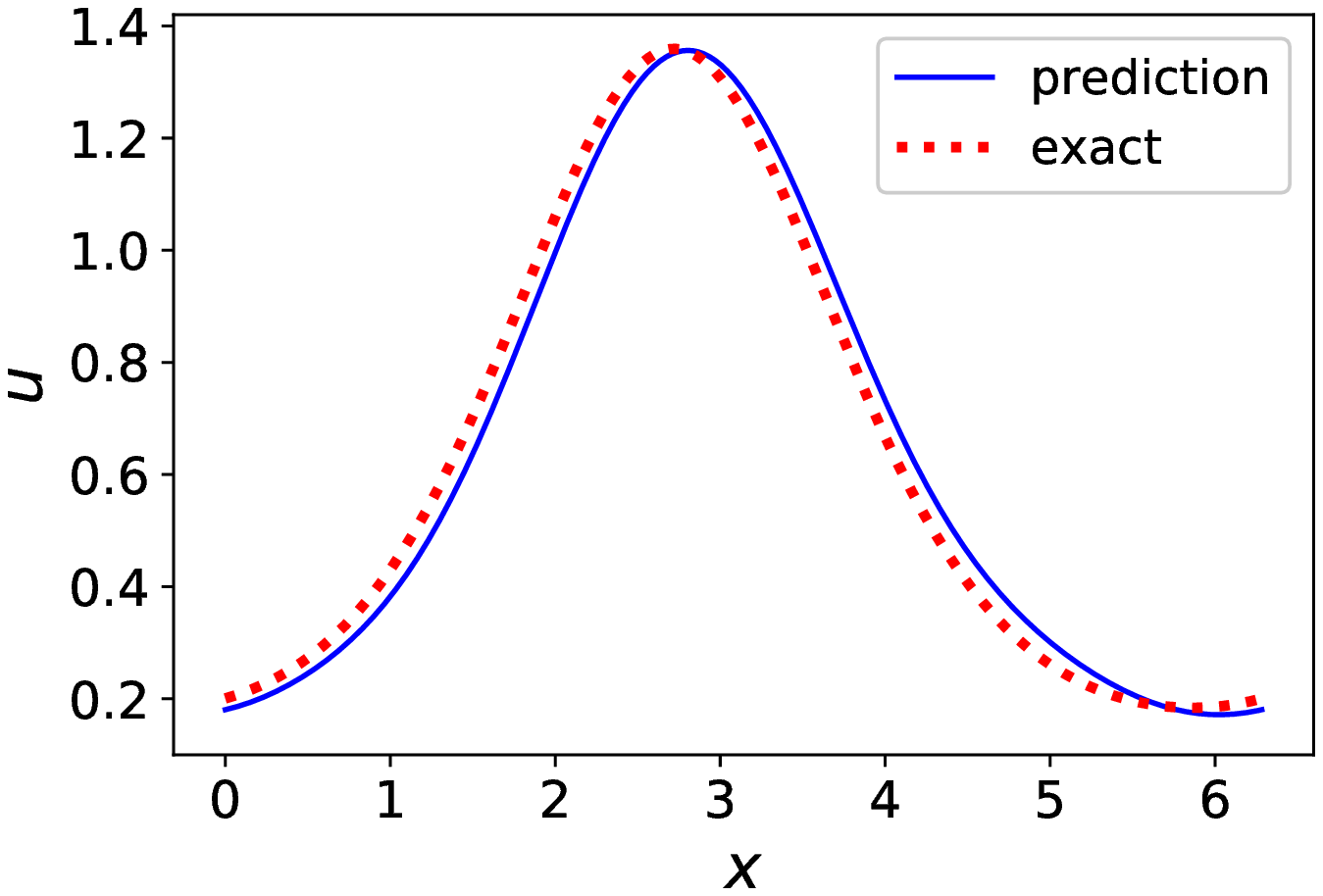}}
	\caption{\small
		Example 1: Comparison of the true solution and the learned model solution at different time. 
		Top-left: $t=1$; top-right: $t=2$; middle-left: $t=3$; middle-right: $t=4$; 
		bottom-left: $t=10$; bottom-right: $t=20$.
	}\label{fig:ex1_solu}
\end{figure}

\begin{figure}[htbp]
	\centering
%	{\includegraphics[width=0.48\textwidth]{Figure/Example1SmallModalDomain/ABSerror.eps}}
	{\includegraphics[width=0.6\textwidth]{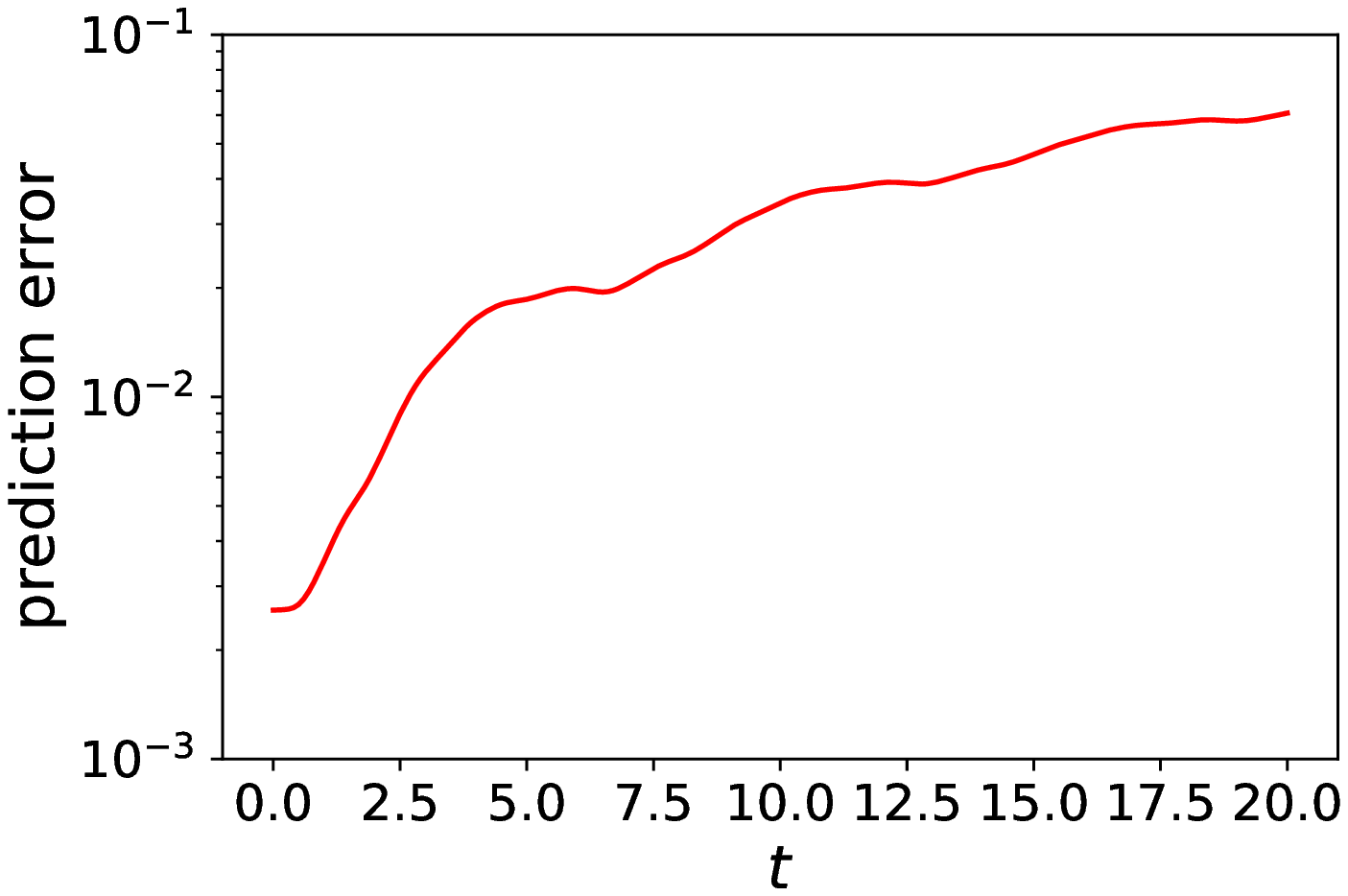}}
	\caption{\small
		Example 1: The evolution of the relative error in the
                prediction in $l^2$-norm.
                %Left: absolute error; right: relative error.  
	}\label{fig:ex1_error}
\end{figure}

\begin{figure}[htbp]
	\centering
	{\includegraphics[width=0.48\textwidth]{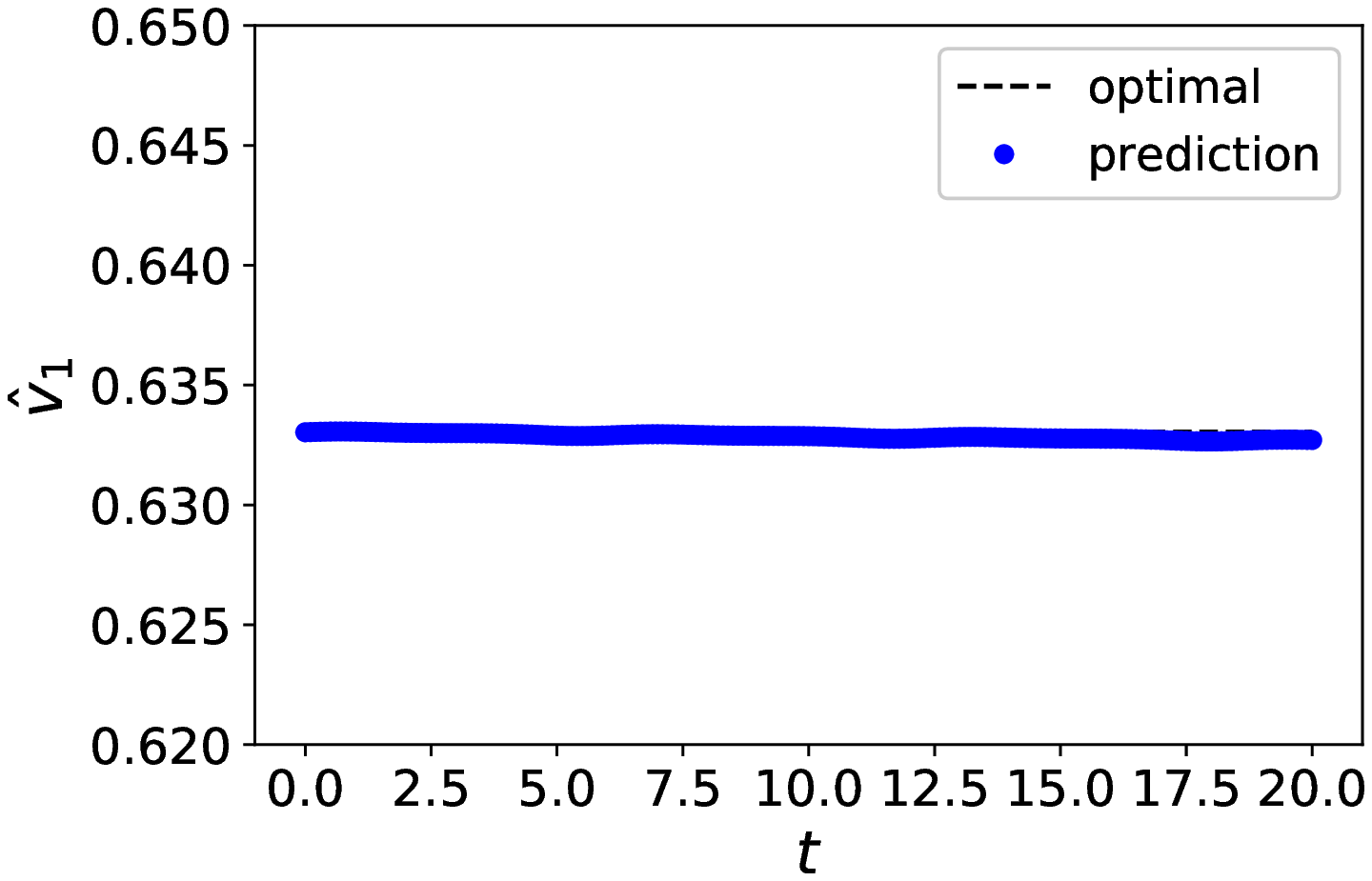}}
	{\includegraphics[width=0.48\textwidth]{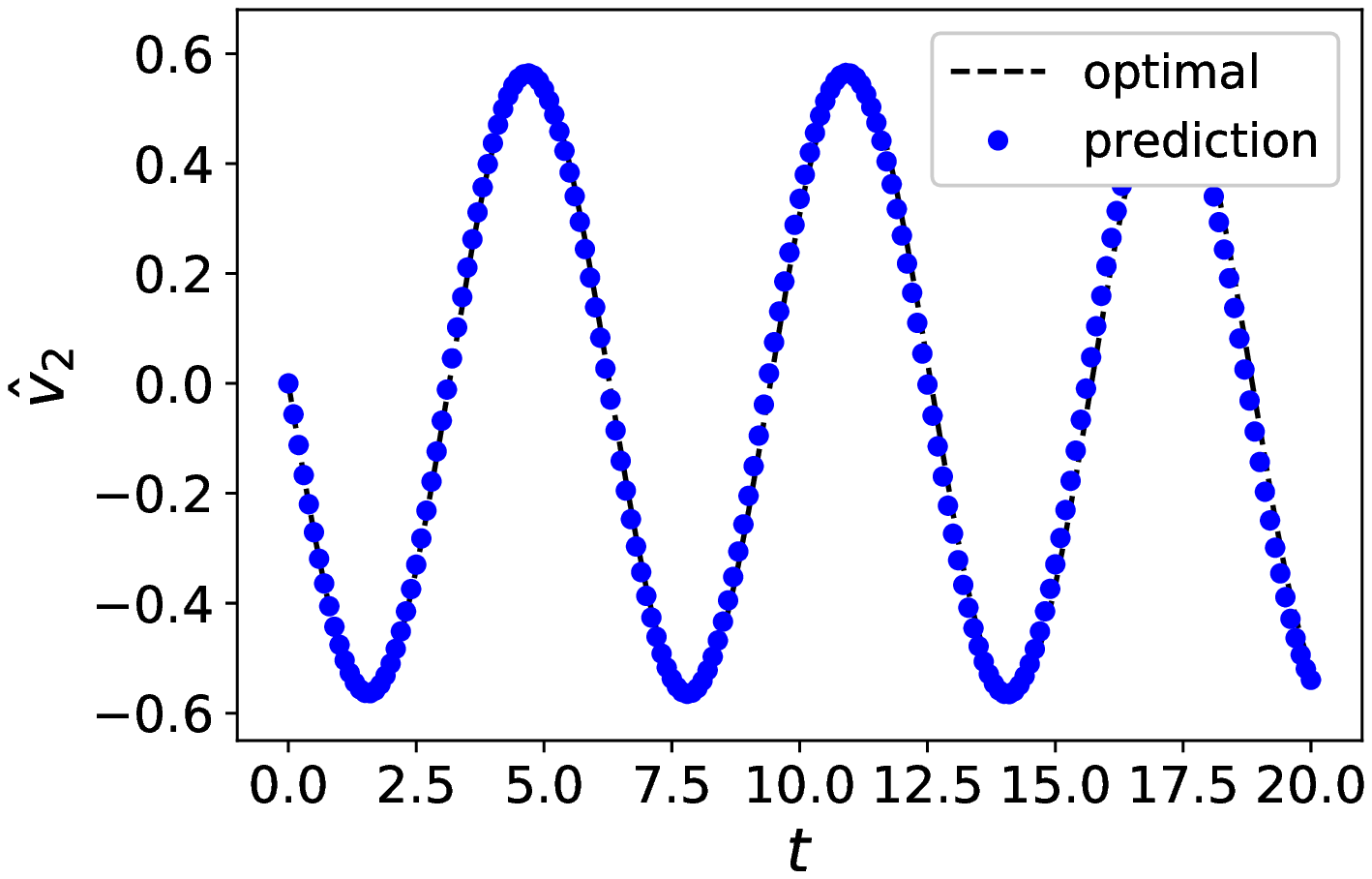}}
	{\includegraphics[width=0.48\textwidth]{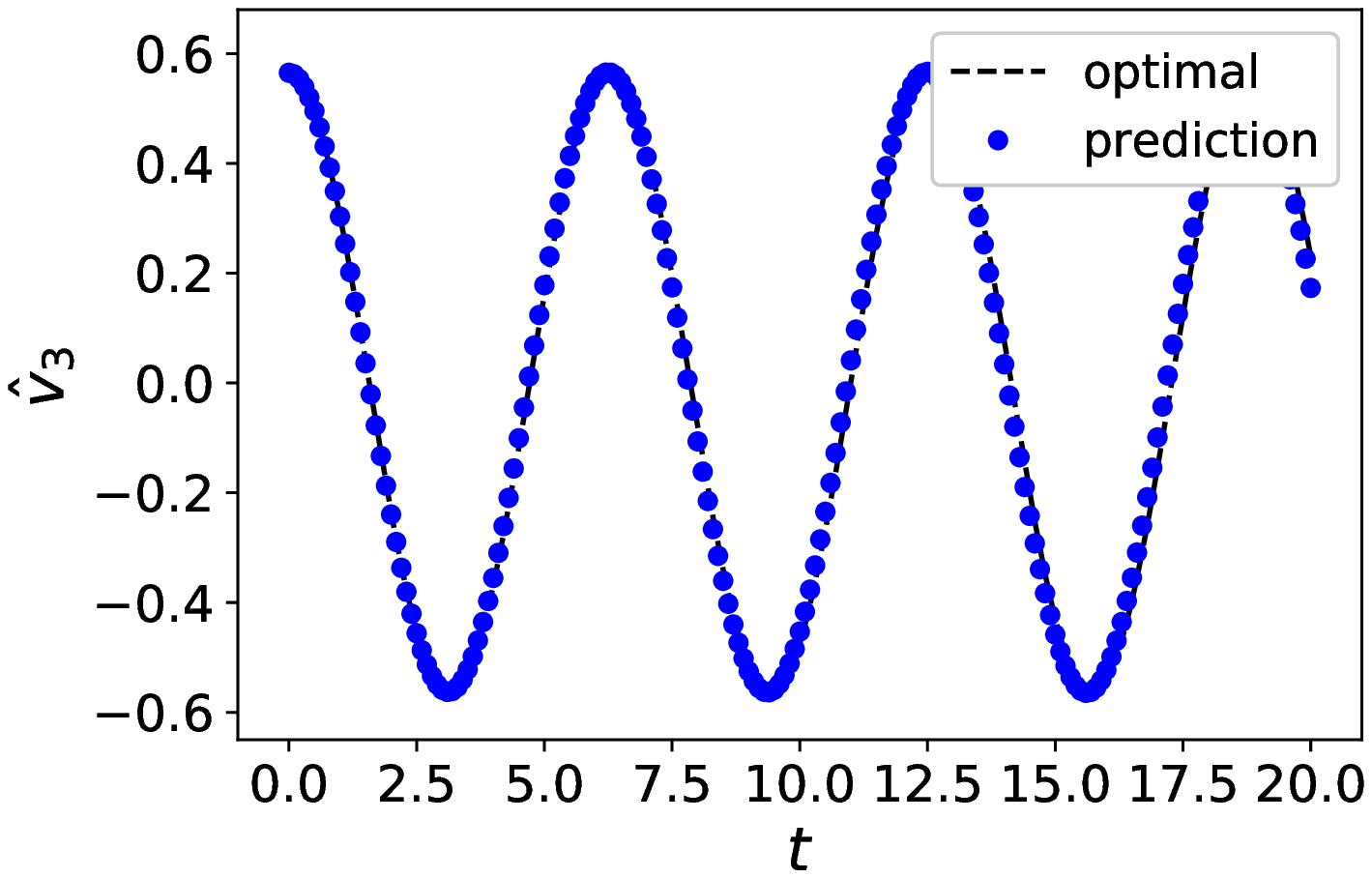}}
	{\includegraphics[width=0.48\textwidth]{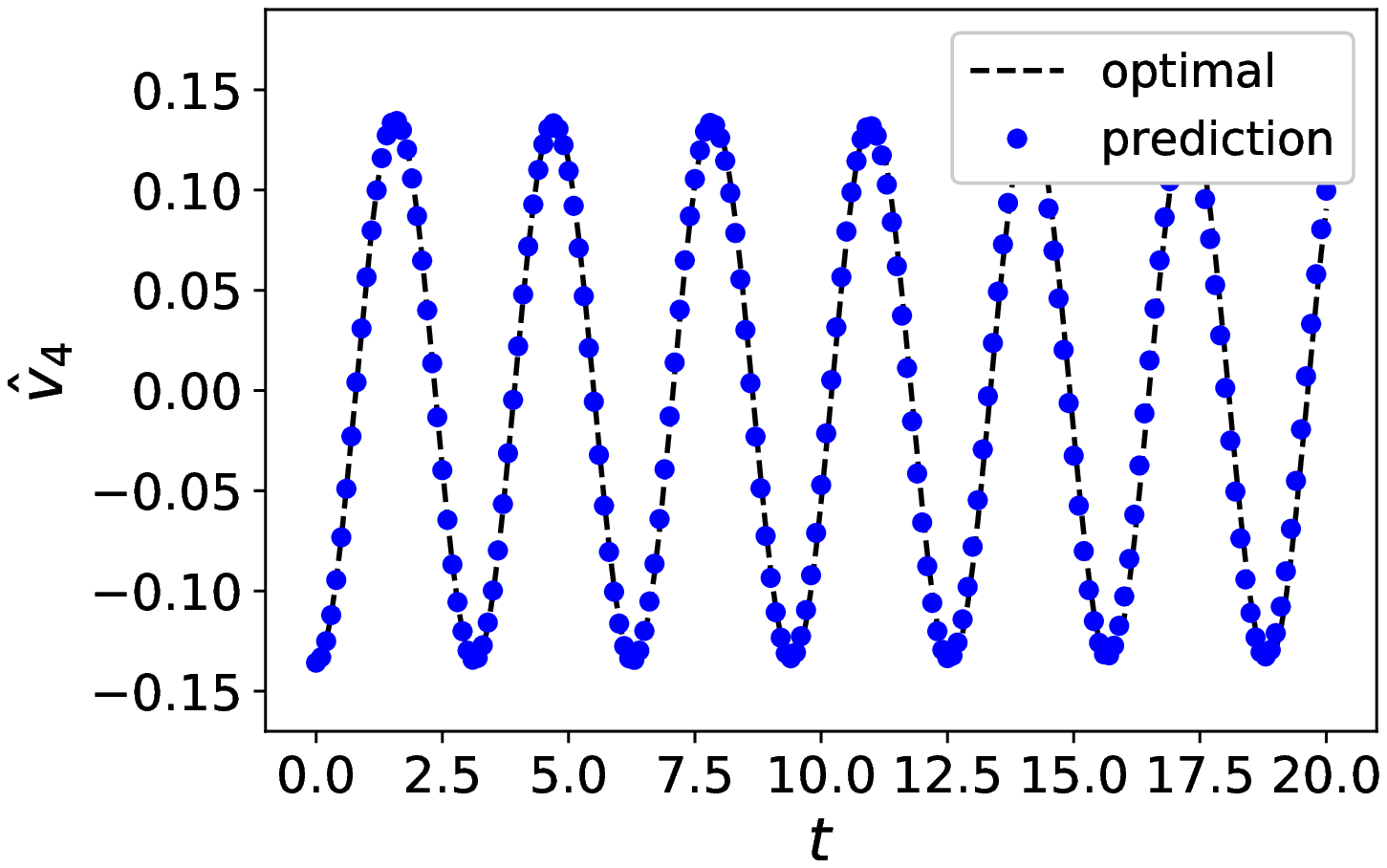}}
	{\includegraphics[width=0.48\textwidth]{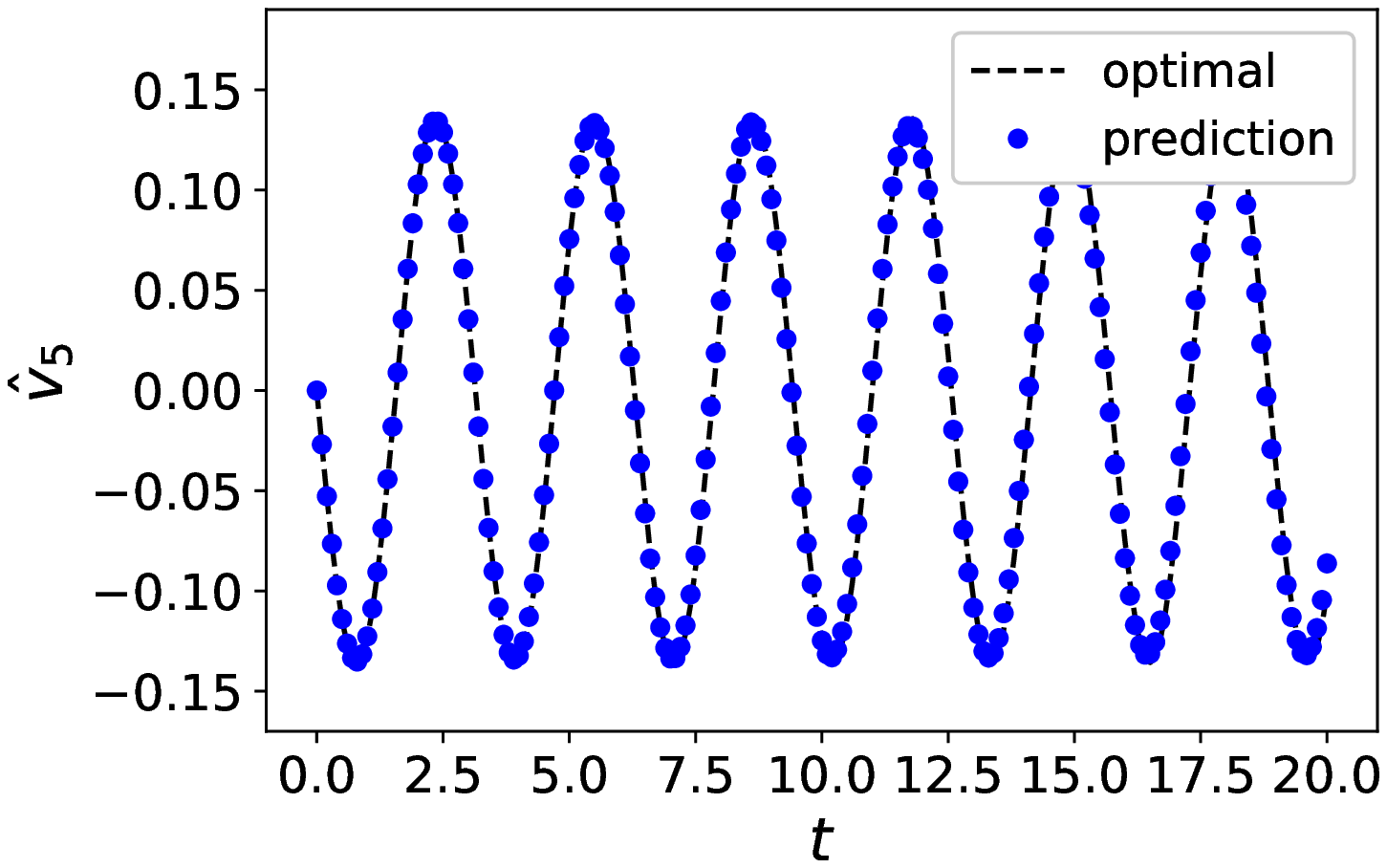}}
	{\includegraphics[width=0.48\textwidth]{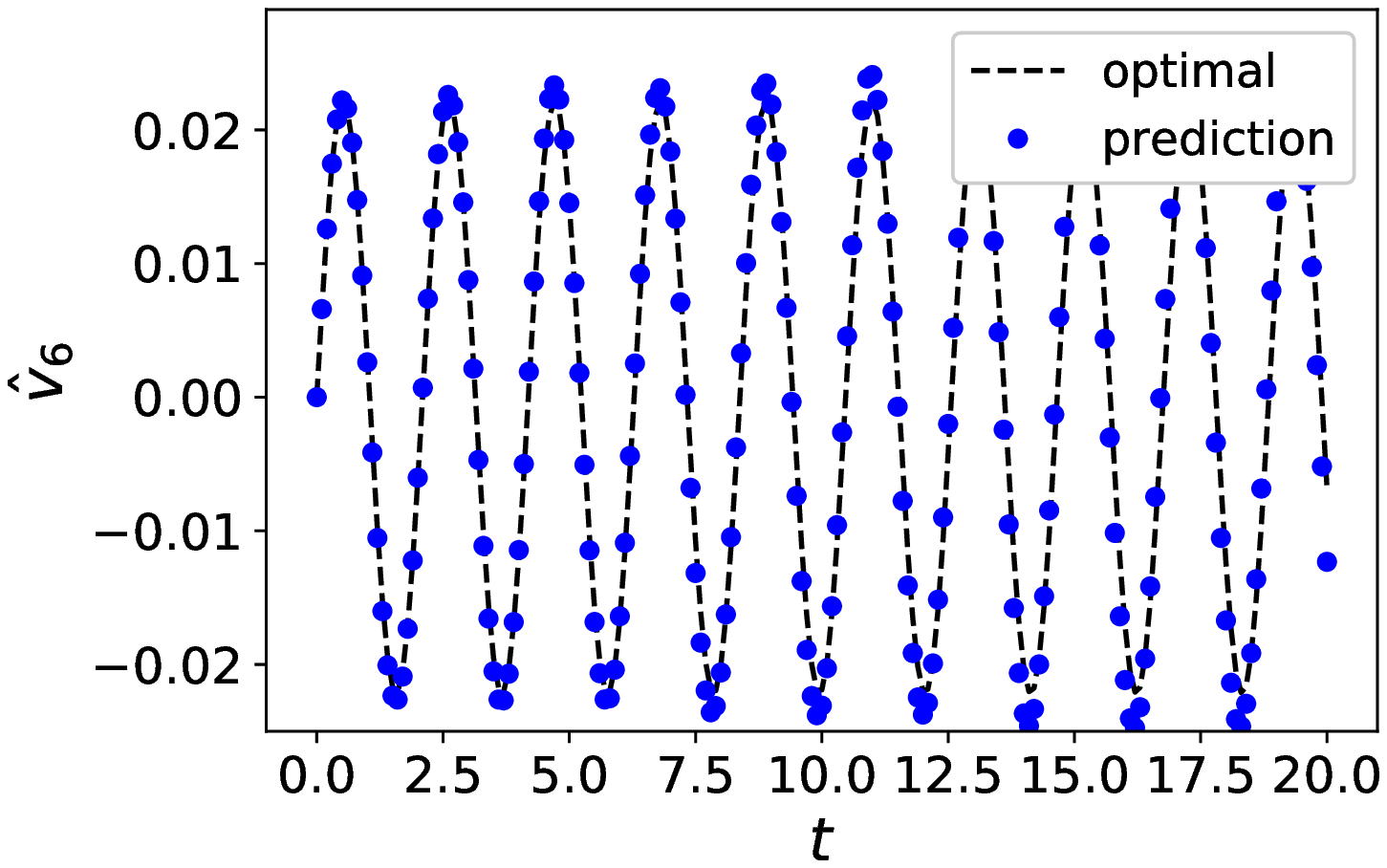}}
	{\includegraphics[width=0.48\textwidth]{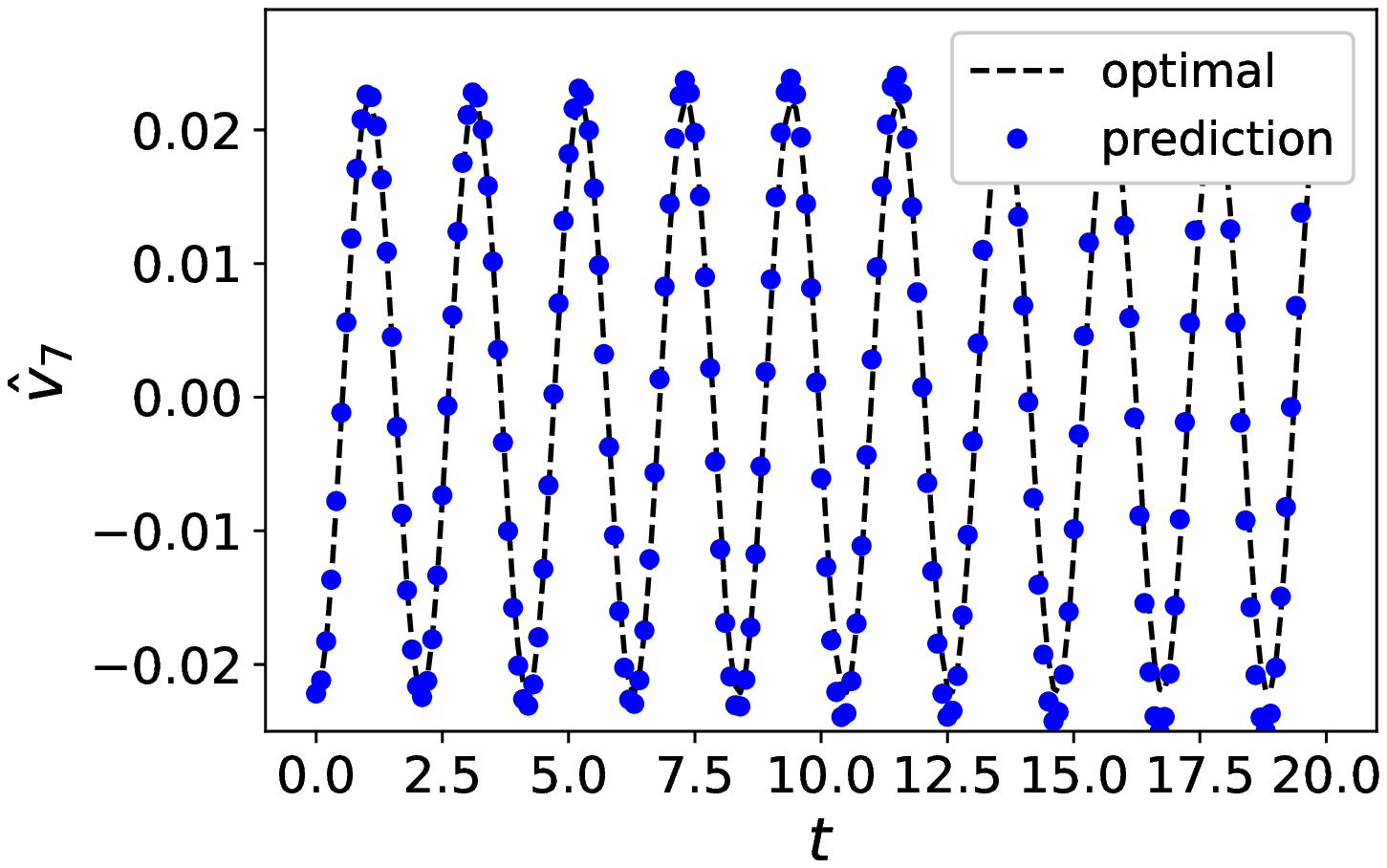}}
	\caption{\small
		Example 1: Evolution of the expansion coefficients for the learned model solution and the projection of the true solution.
	}\label{fig:ex1_coef}
\end{figure}

The finite dimensional  approximation space is set as $\mathbb V_n =
\textrm{span}\{ {e^{ikx}}, k\leq 3\}$, which implies  $n=7$.
%\begin{align*}
%& \phi_1 (x) =1, \quad \phi_2(x) = \cos(x), \quad \phi_3(x) = \sin (x), \quad \phi_4 (x) = \cos(2x),\\
%& \phi_5(x) = \sin(2x), \quad \phi_6(x) = \cos(3 x), \quad \phi_7 (x) = \sin(3x). 
%\end{align*} 
The time lag $\Delta$ is taken as 0.1. 
%The domain $D$ on the modal space is taken as $[-1,1]^7$, from which we sample $50,000$ training data.  
The domain $D$ in the modal space is fixed as $[-0.8,0.8]$ for $k\leq
1$ and $[-0.2,0.2]$ for $k=2$, and $[-0.03,0.03]$ for $k=3$.
By using uniform {distribution}, we generate $80,000$ training data in the
modal space.
For neural network modeling, we employ the block ResNet structure with
two blocks ($K=2$), each of which contains 3 hidden layers of equal
width of 30 neurons.
The loss function training history is shown in
Fig.~\ref{fig:ex1_loss}, where the network is trained for up to
$2,000$ epochs. Convergence can be achieved after about $200$ epochs.
To validate the model, we employ an initial condition $u_0(x)=\frac12
\exp ( \sin (x) )$ and conduct simulations using the trained model, in
the form of \eqref{model-initial}--\eqref{model-final}, for
up to $t=20$. Note that this particular initial condition, albeit
smooth, is fairly representative, as it is not in the approximation
space $\mathbb V_n$.
In Fig.~\ref{fig:ex1_solu}, the solution prediction of the trained
model is plotted at different time instances, along with the true
solution for comparison.
The relative error in the predicted solution is shown in Fig.~\ref{fig:ex1_error}. 
We observe that the network model produces
accurate prediction results for time up to $t=20$. The error grows
over time. This is consistent with the error estimate from Theorem
\ref{thm:errorbound} and is expected from any numerical time integrator.
To further examine the solution property, we also plot in Fig.~\ref{fig:ex1_coef} 
the evolution of the learned expansion coefficients $\hat v_j$, $1\le
j \le 7$. For reference purpose, the optimal coefficients obtained by
the $L_{\mathcal B}^2(\Omega)$ orthogonal projection of 
the true solution onto ${\mathbb V}_n$ are also plotted. We observe
good agreement between the two solutions.
%%%%%%%%%%%%%%%%%

\subsection{Example 2: Diffusion Equation}

We now consider the following diffusion equation with Dirichlet boundary condition: 
\begin{equation}
	\label{eq:example2}
	\begin{cases}
		u_t = \sigma u_{xx},\quad (x,t) \in (0,\pi) \times \mathbb R^+, \\
		u(0,t)=u(\pi,t)=0,\quad t \in \mathbb R^+.
	\end{cases}
\end{equation}

\begin{figure}[htbp]
	\centering
	{\includegraphics[width=0.6\textwidth]{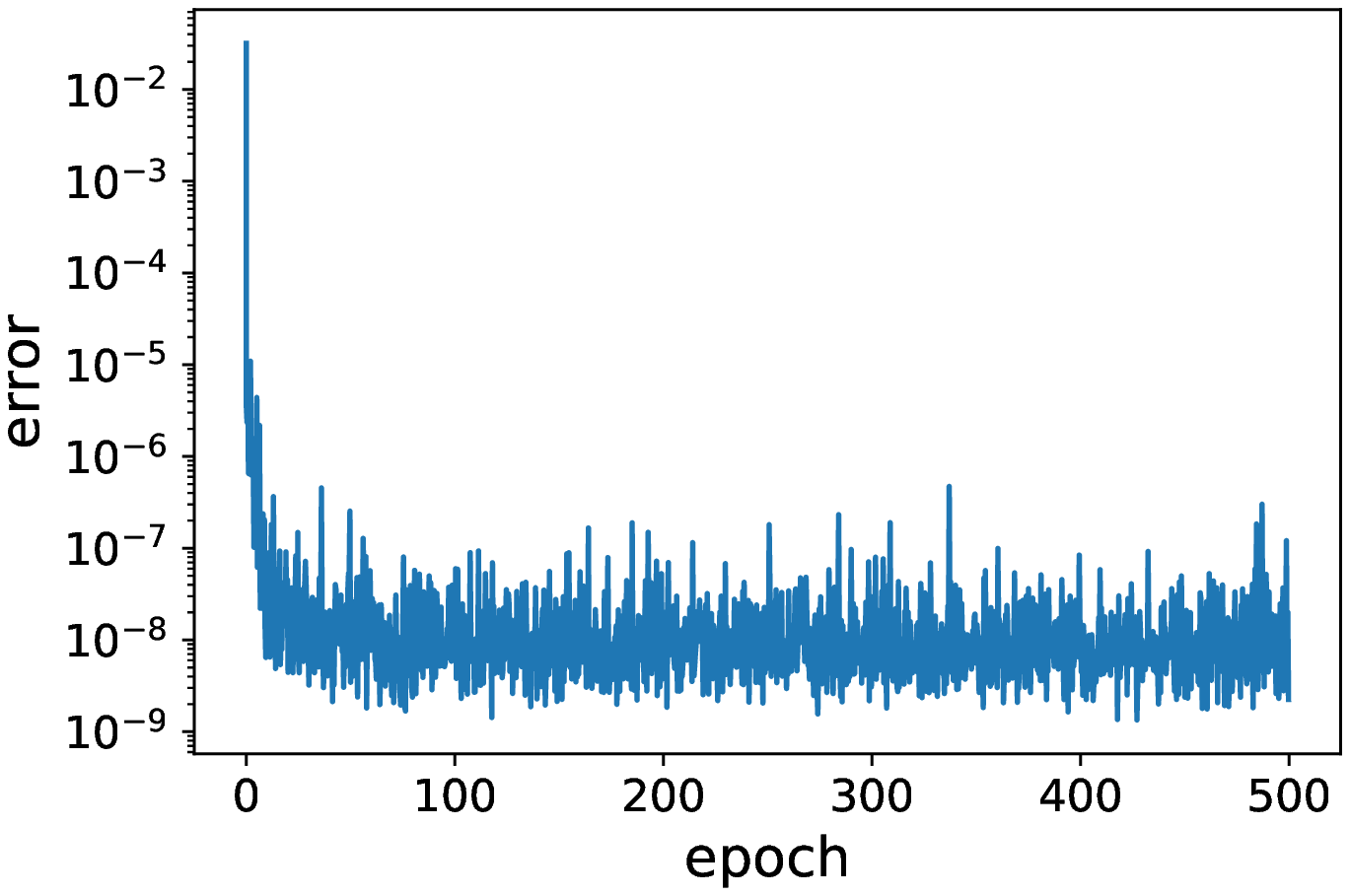}}
	\caption{\small
		Example 2: Training loss history.
	}\label{fig:ex2_loss}
\end{figure}

\begin{figure}[htbp]
	\centering
	{\includegraphics[width=0.48\textwidth]{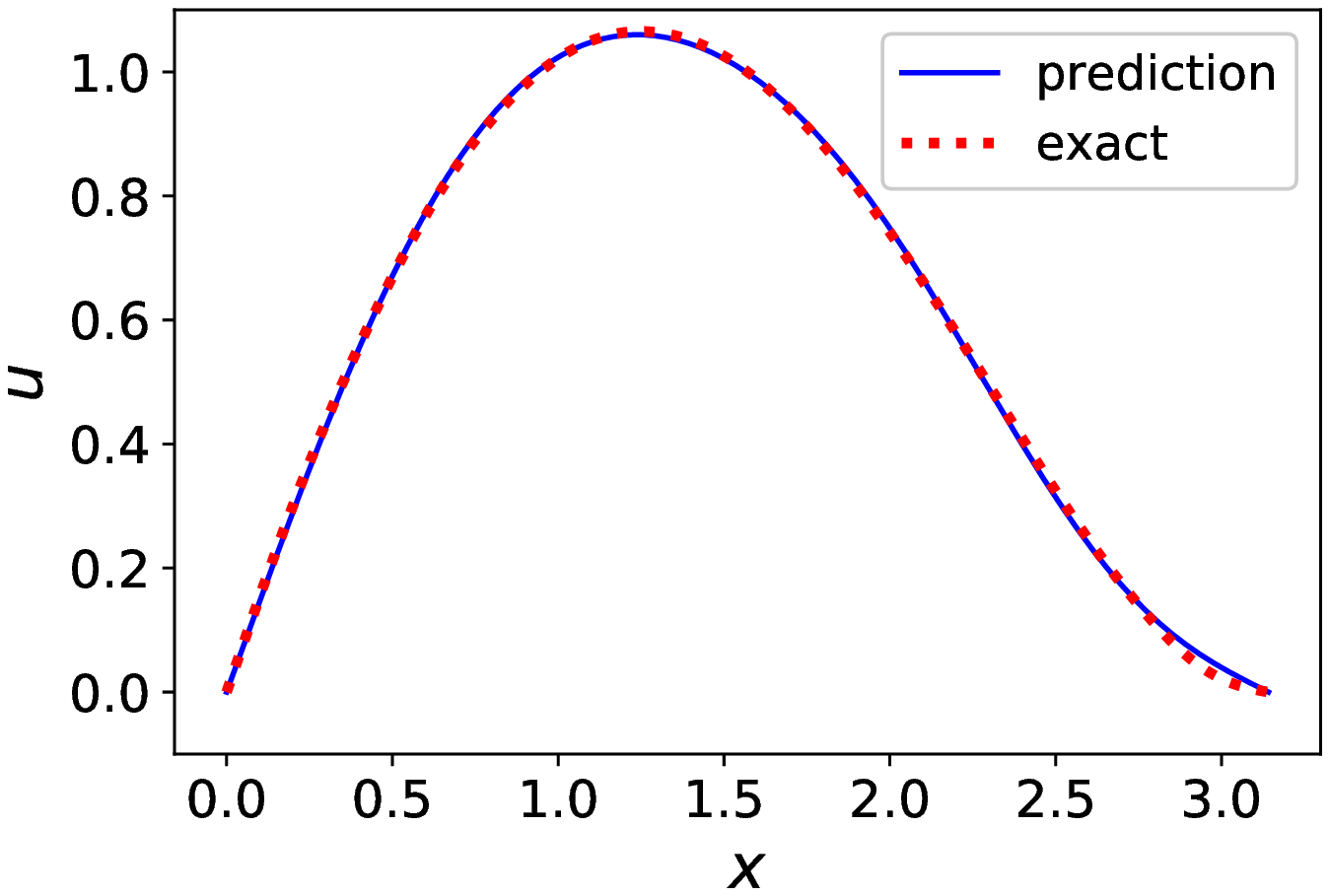}}
	{\includegraphics[width=0.48\textwidth]{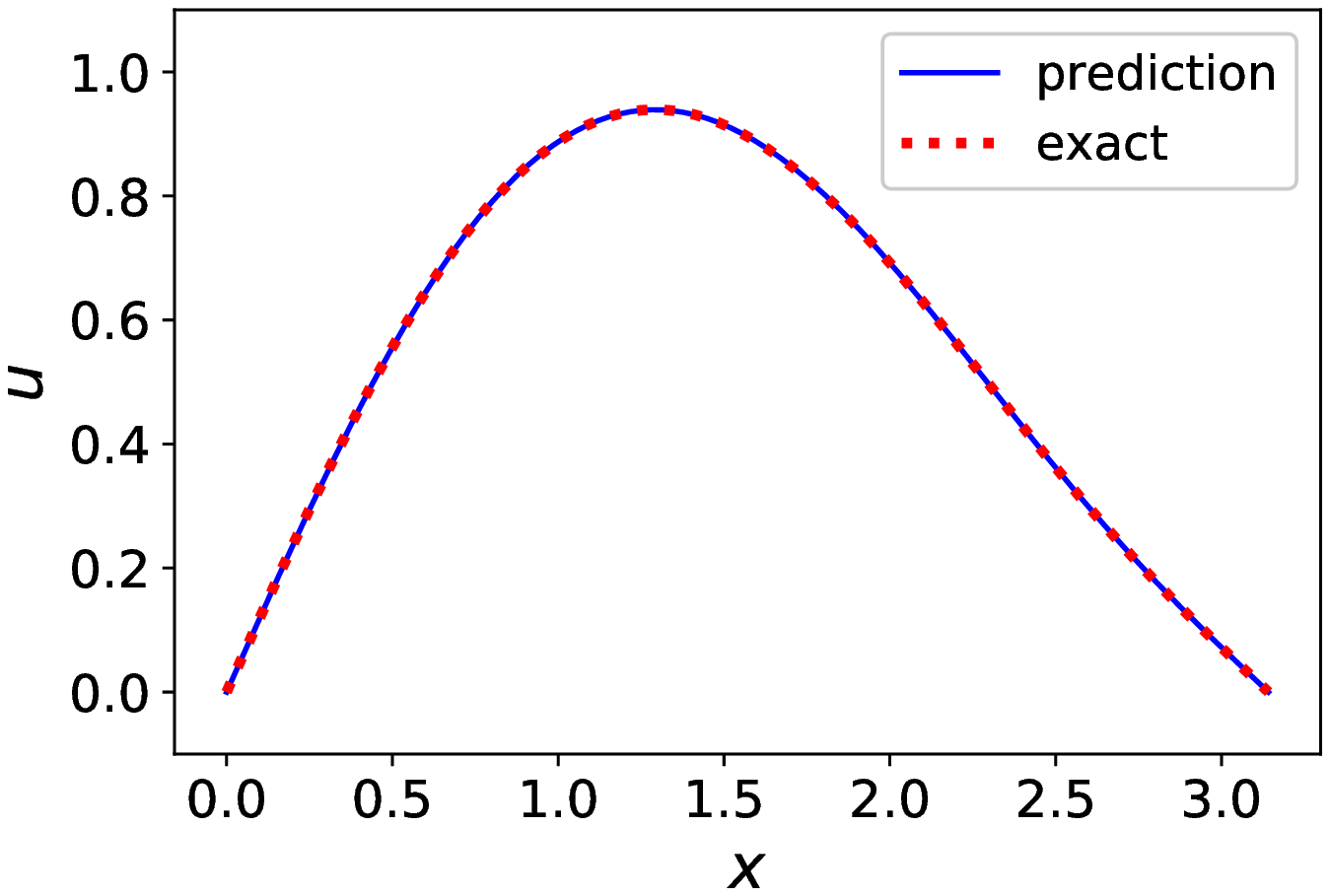}}
	{\includegraphics[width=0.48\textwidth]{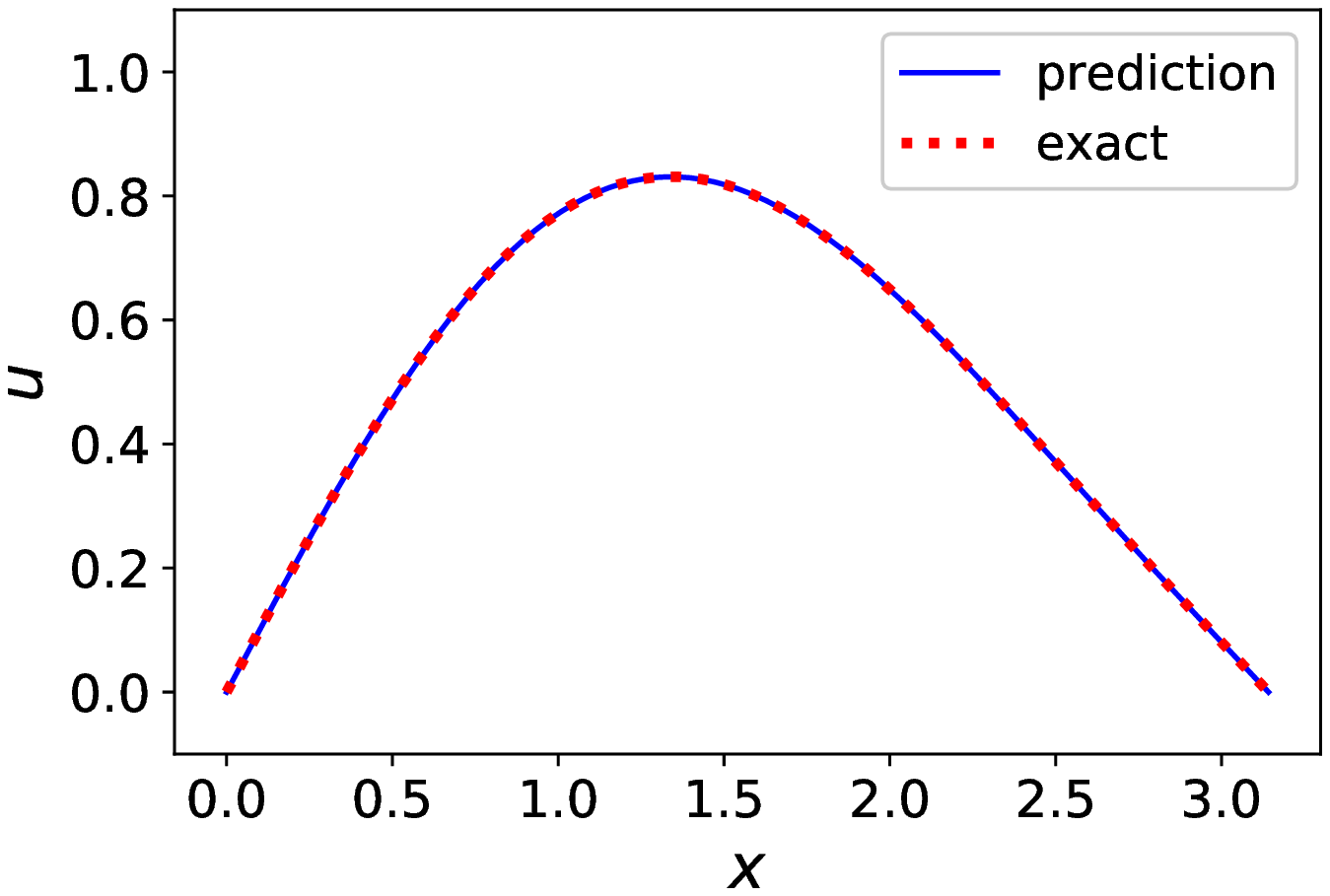}}
	{\includegraphics[width=0.48\textwidth]{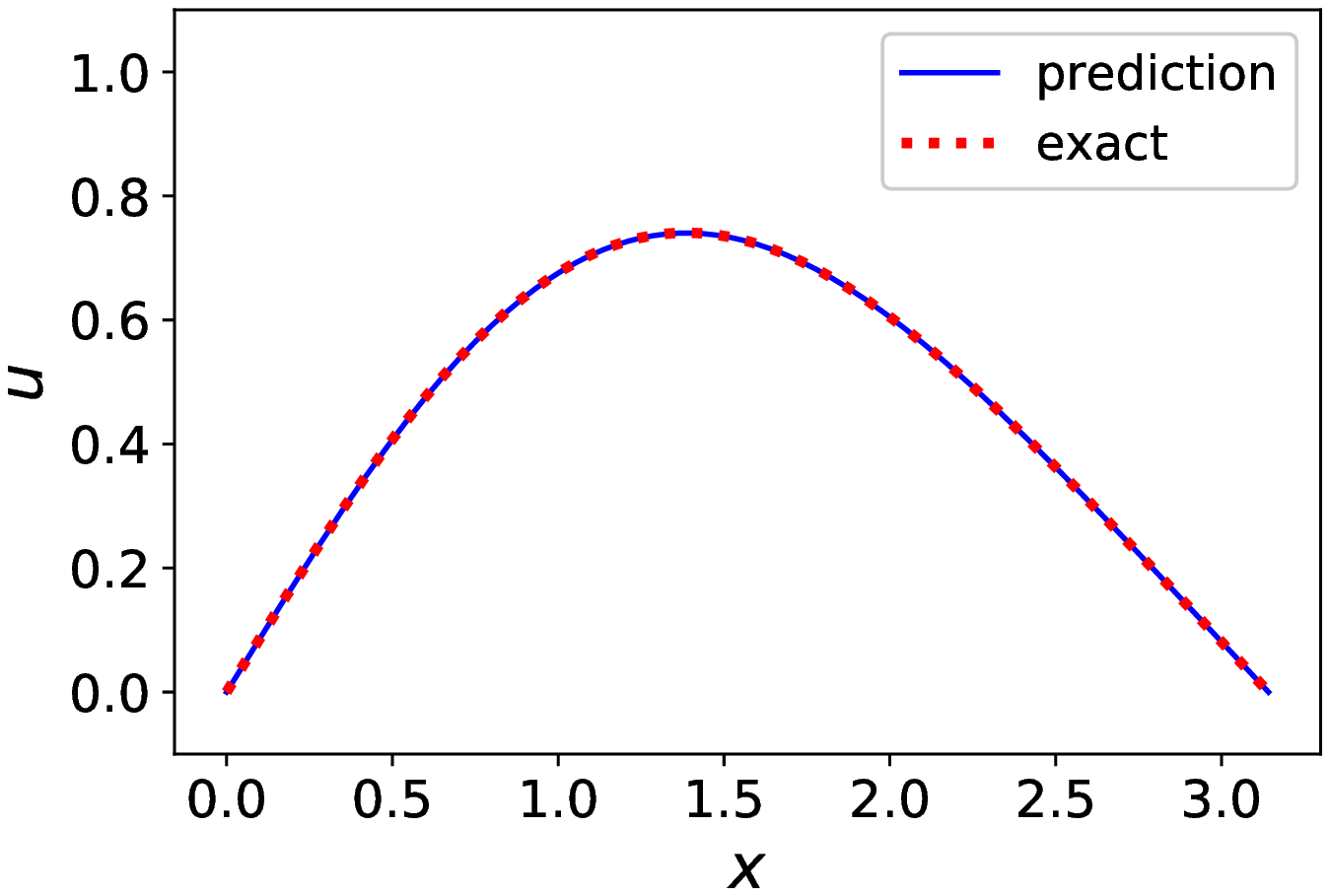}}
	\caption{\small
		Example 2: Comparison of the true solution and the learned model solution at different time. Top-left: $t=0$; top-right: $t=1$; bottom-left: $t=2$; bottom-right: $t=3$.   
	}\label{fig:ex2_solu}
\end{figure}

\begin{figure}[htbp]
	\centering
	%{\includegraphics[width=0.48\textwidth]{Figure/Example2/ABSerror.eps}}
	{\includegraphics[width=0.6\textwidth]{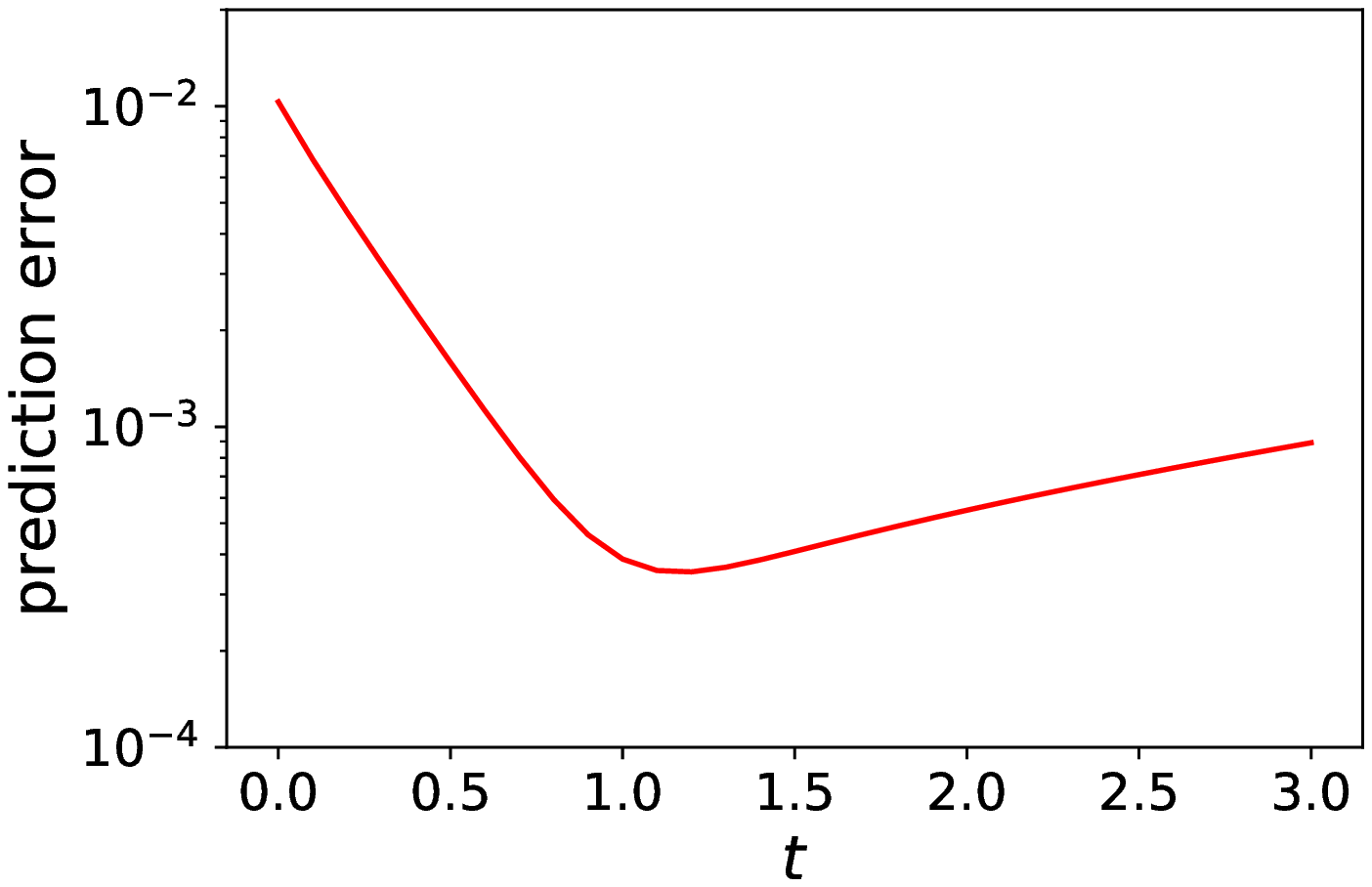}}
	\caption{\small
		Example 2: The evolution of the relative error in the
                prediction in $l^2$-norm.
                %Left: absolute error; right: relative error.  
	}\label{fig:ex2_error}
\end{figure}

\begin{figure}[htbp]
	\centering
	{\includegraphics[width=0.325\textwidth]{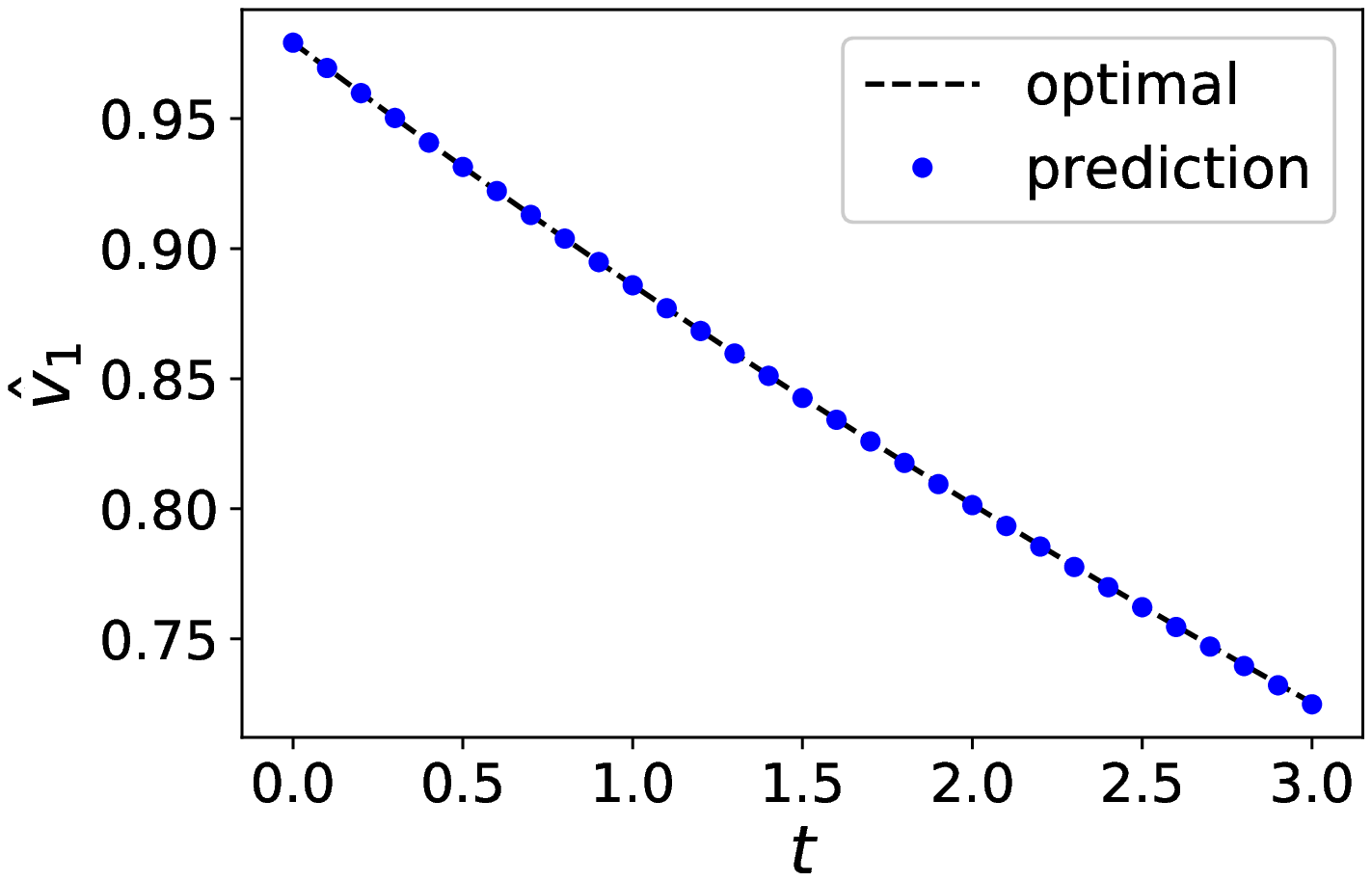}}
	{\includegraphics[width=0.325\textwidth]{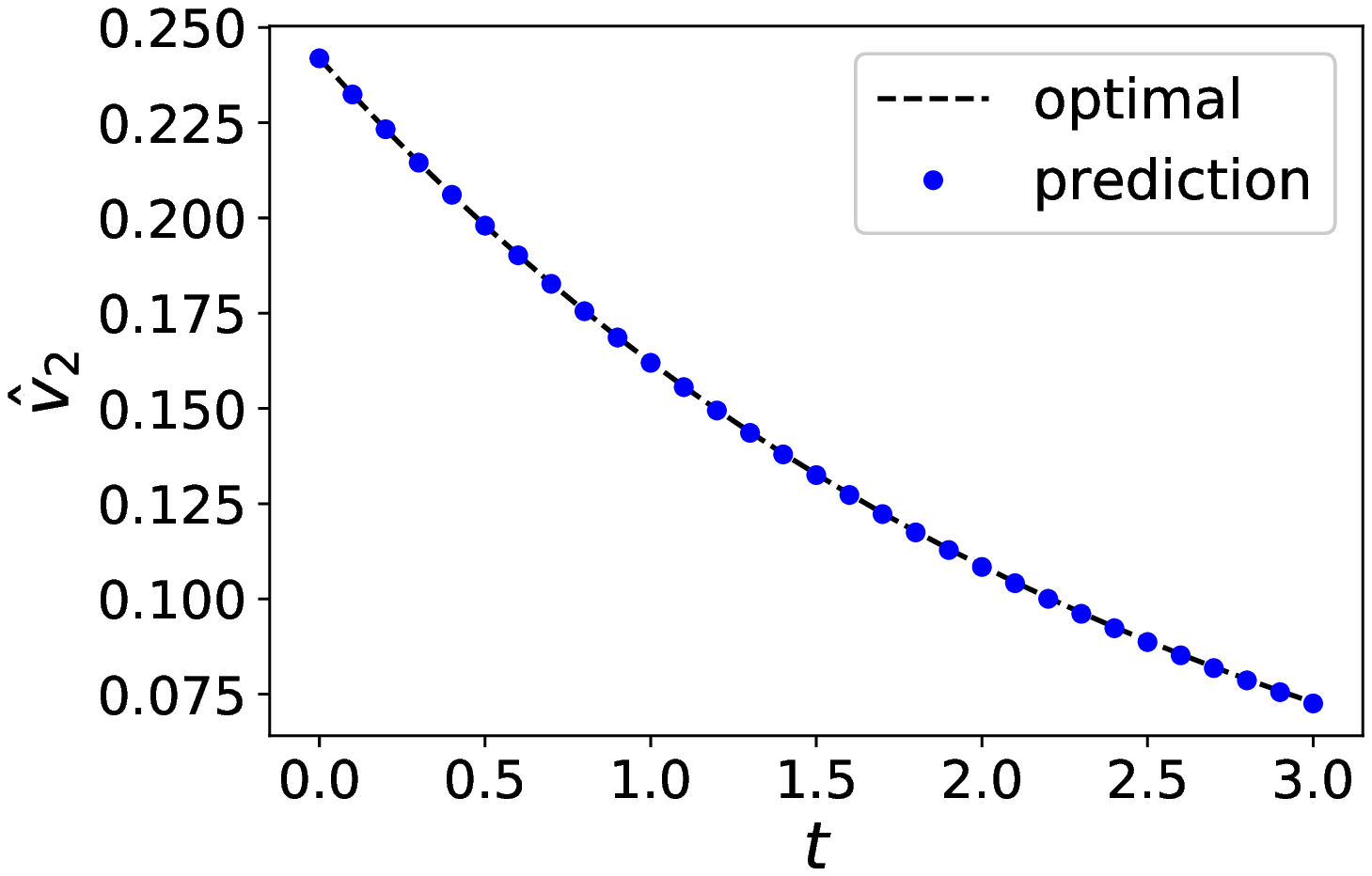}}
	{\includegraphics[width=0.325\textwidth]{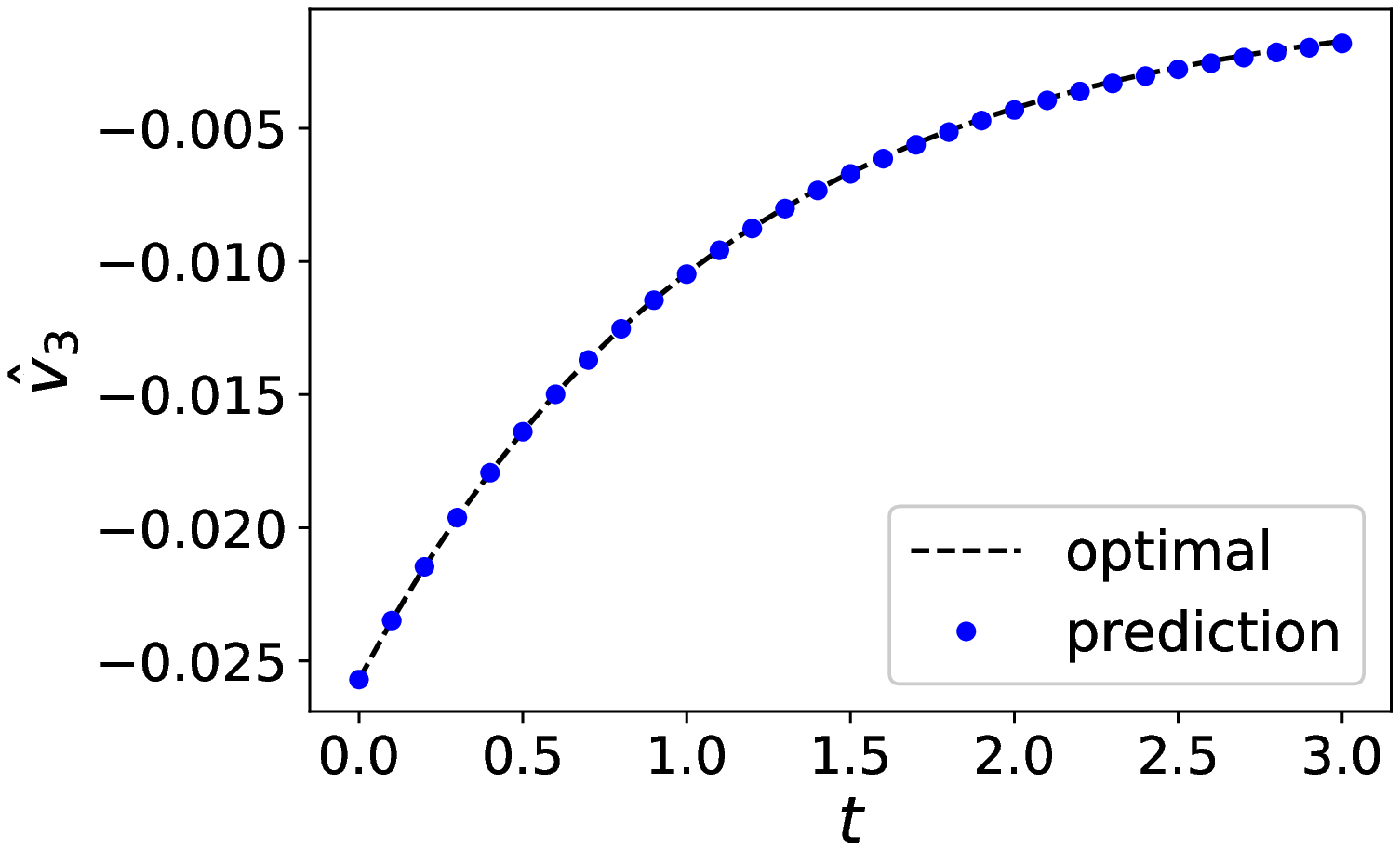}}
	{\includegraphics[width=0.325\textwidth]{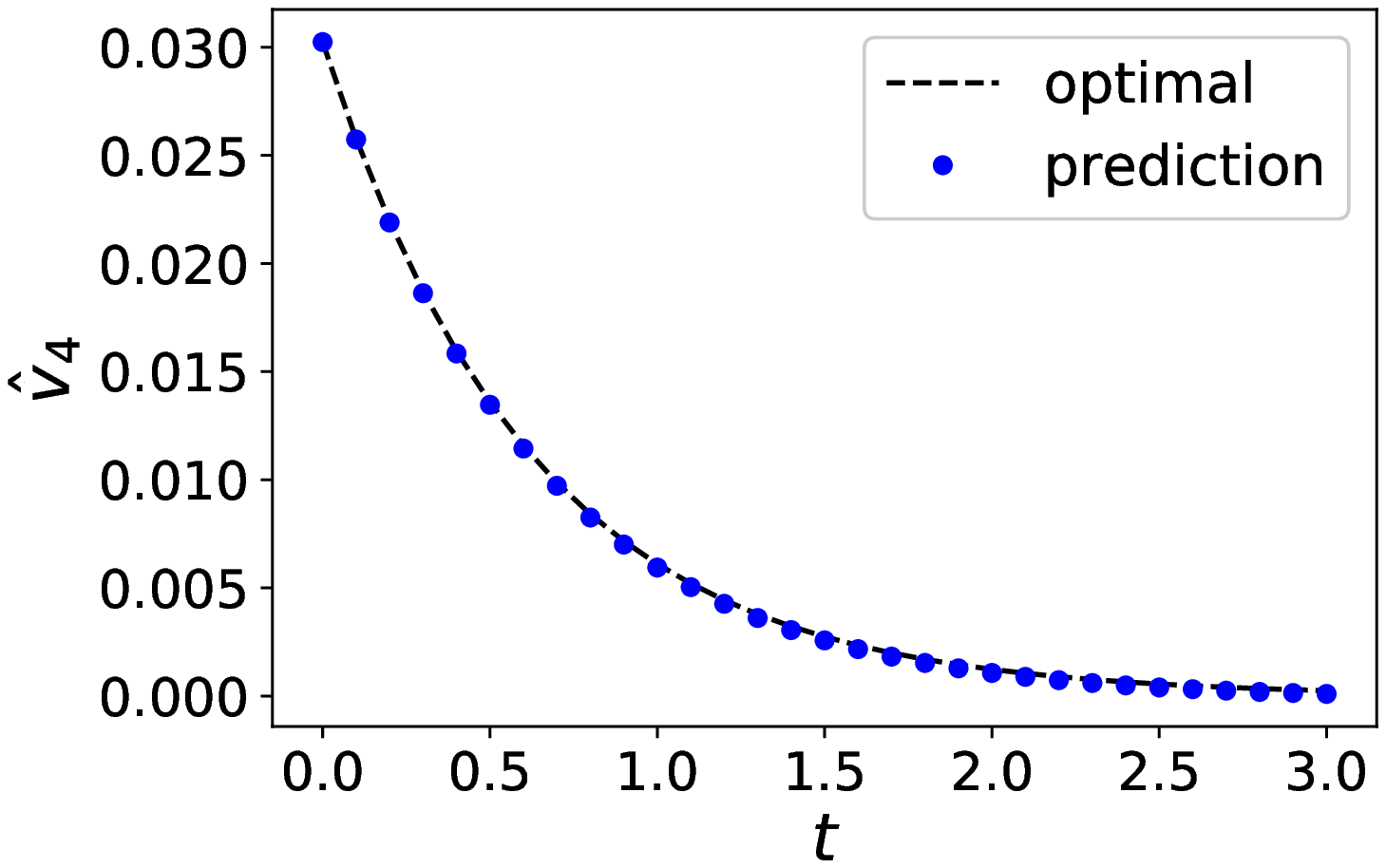}}
	{\includegraphics[width=0.325\textwidth]{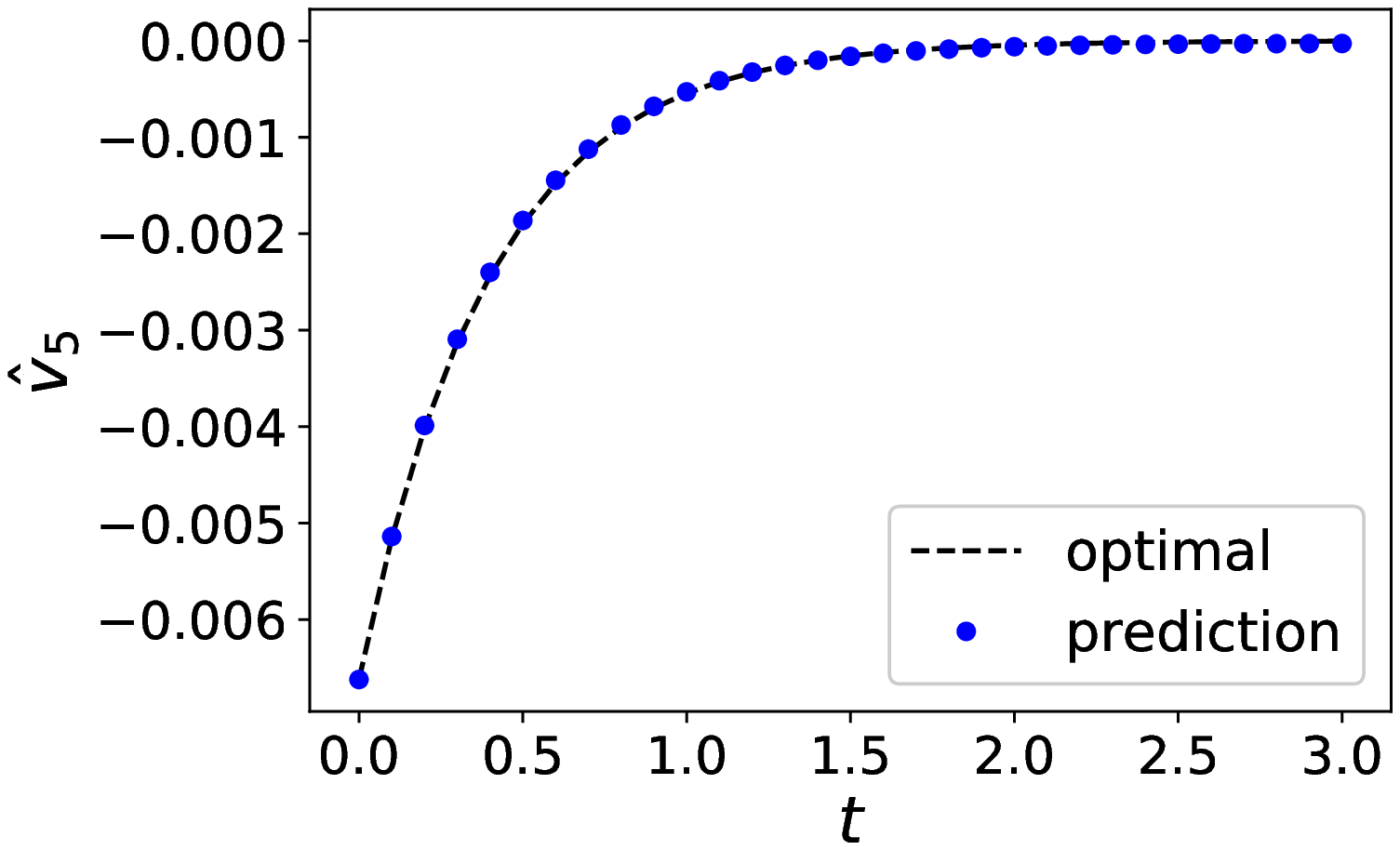}}
	\caption{\small
		Example 2: Evolution of the expansion coefficients for the learn model and the projection of the true solution.
	}\label{fig:ex2_coef}
\end{figure}

The finite dimensional approximation space is chosen as $\mathbb V_n =
{\rm span} \{  \sin(jx), 1\le j \le 5 \}$, where $n=5$. (The symmetry
of the problem allows this simplified choice to facilitate the
computation and comparison.)
The time lag $\Delta$ is taken as 0.1. 
The domain $D$ in the modal space is taken as $[-1,  1]\times[-0.5,  0.5]\times[-0.2,  0.2]\times[-0.05, 0.05]\times[-0.01,  0.01]
$, from which we sample $30,000$ training data.  
In this example, we employ a single-block ResNet method ($K=1$)
containing 3 hidden layers of equal width of 30 neurons. 
The training of the network model is conducted for up to 500 epochs,
where satisfactory convergence is established, as shown in
Fig.~\ref{fig:ex2_loss}.
For validation and accuracy test, we conduct numerical predictions of
the trained network model using initial condition
$$u_0(x)= \frac{x}{\pi} \left[  5 - \frac{4 x}{\pi} - 7 \left( \frac{x}{\pi} \right)^2 + 6 \left( \frac{ x } {\pi} \right)^3  \right],$$ 
for up to time $t=3$. % in the form of \eqref{model-initial}--\eqref{model-final}. 
Fig.~\ref{fig:ex2_solu} shows the solution predicted by the trained network model, along with the exact solution of the true equation \eqref{eq:example2}. 
%We see that the magnitude of the solution decays as time increases, due to the effect of diffusion. 
It can be seen that the predicted solution agrees well with the exact
solution. The relative error in the numerical prediction is shown in
Fig.~\ref{fig:ex2_error},
in term of $l^2$-norm. Finally,
the evolution of the expansion coefficients is also shown given in
Fig.~\ref{fig:ex2_coef}.
We observe that they agree well with the optimal coefficients obtained
by projecting the exact solution onto the linear space  $\mathbb V_n$.
%%%%%%%%%%%%%

\begin{figure}[htbp]
	\centering
	{\includegraphics[width=0.48\textwidth]{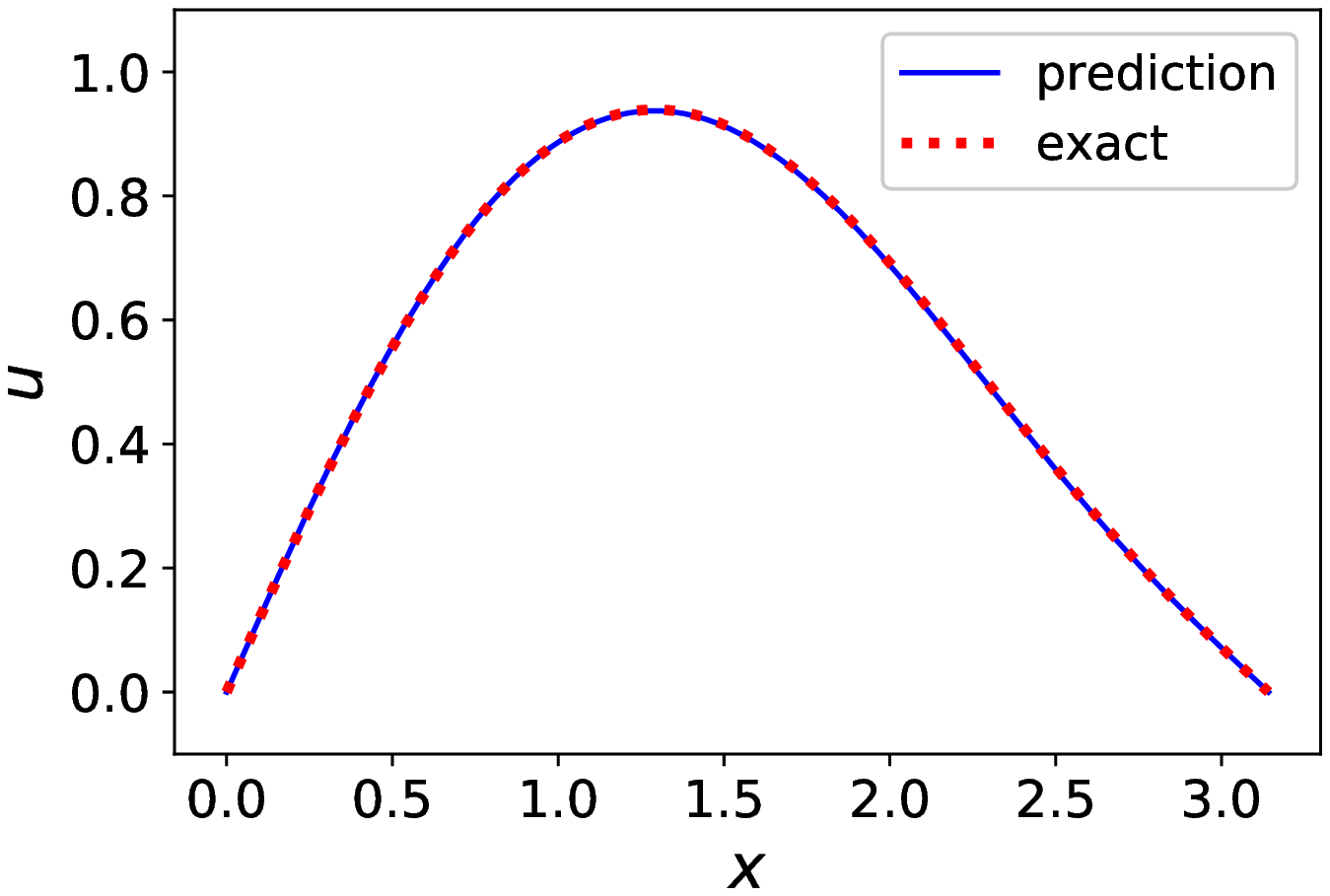}}
	{\includegraphics[width=0.48\textwidth]{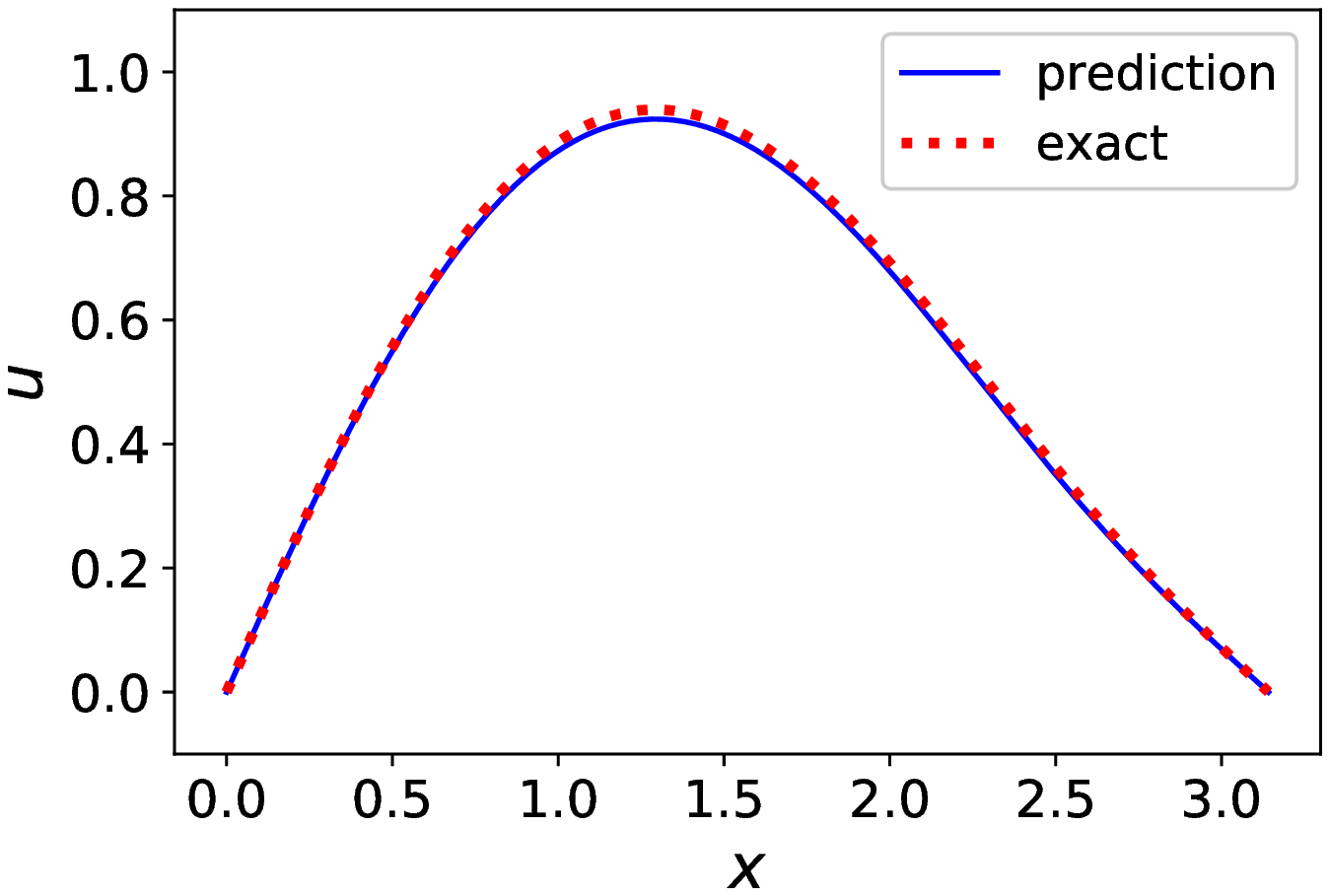}}
	{\includegraphics[width=0.48\textwidth]{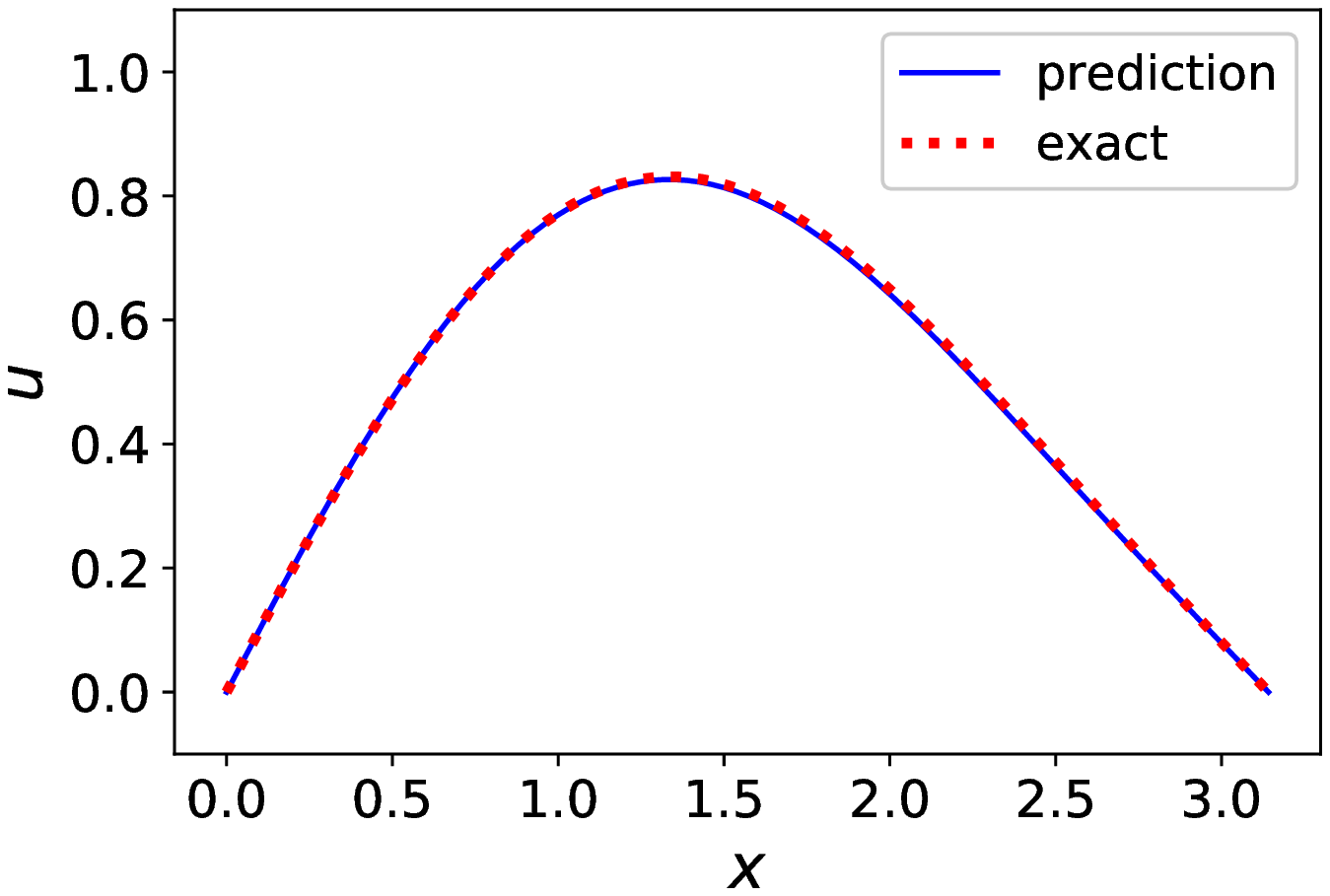}}
	{\includegraphics[width=0.48\textwidth]{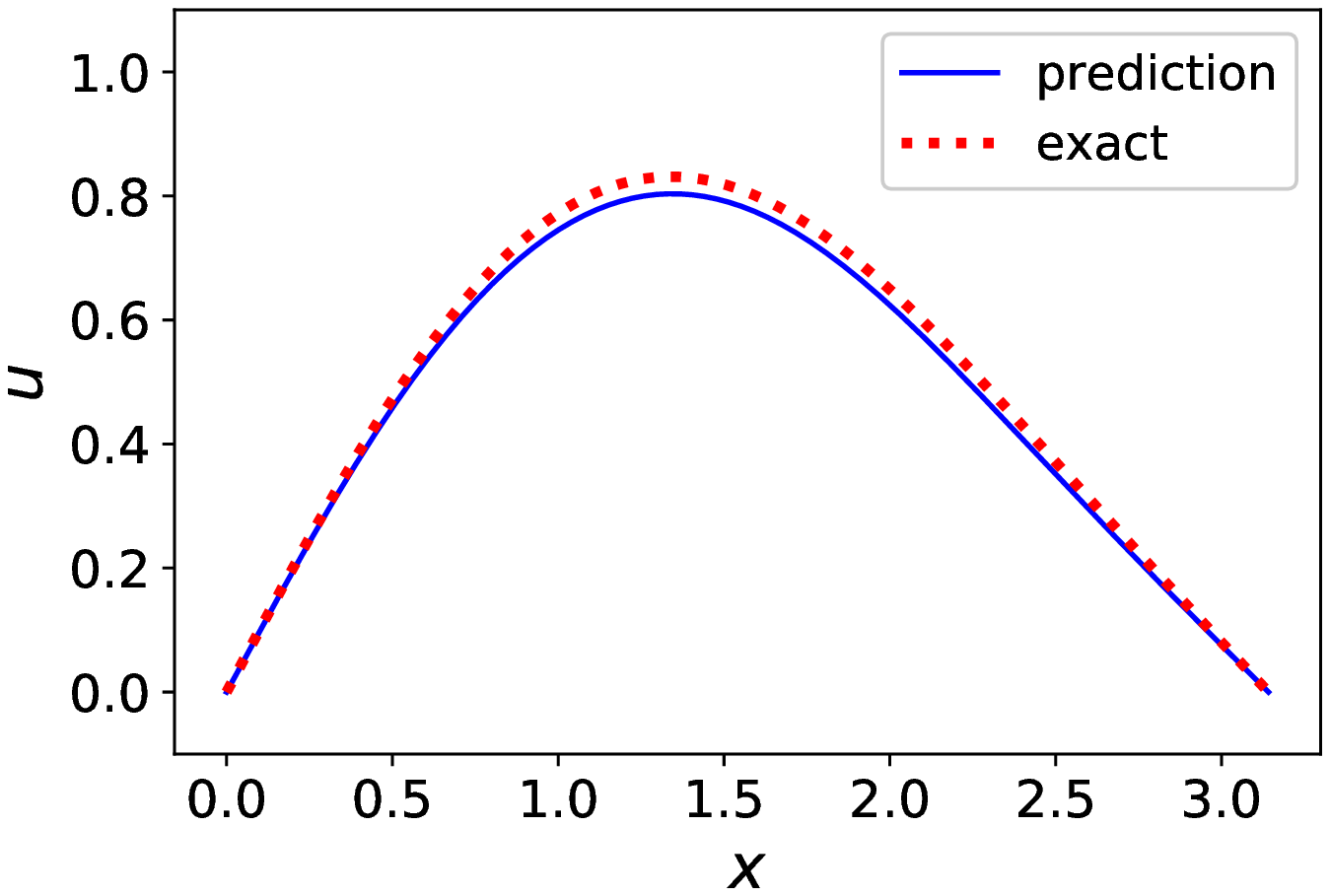}}
	{\includegraphics[width=0.48\textwidth]{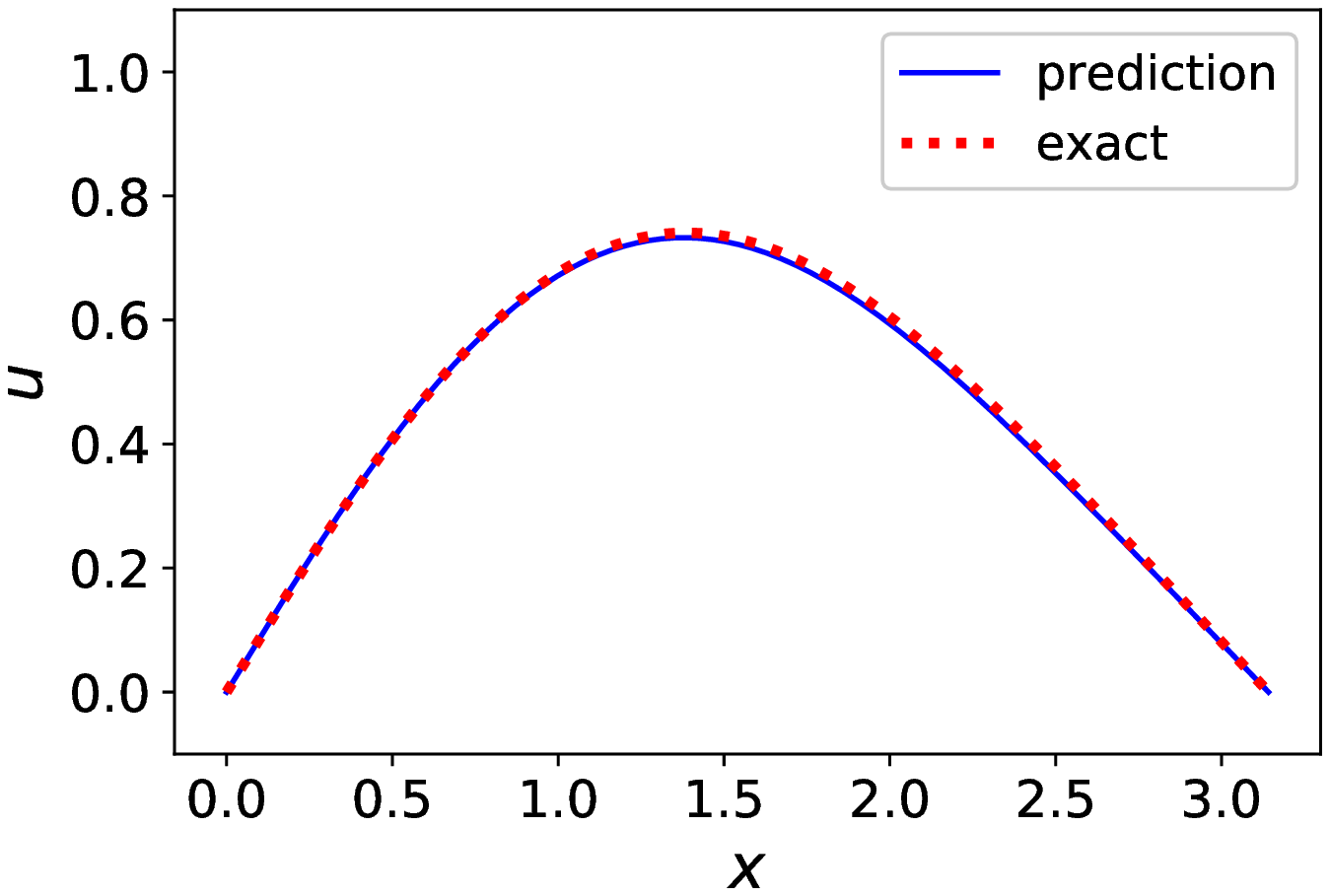}}
	{\includegraphics[width=0.48\textwidth]{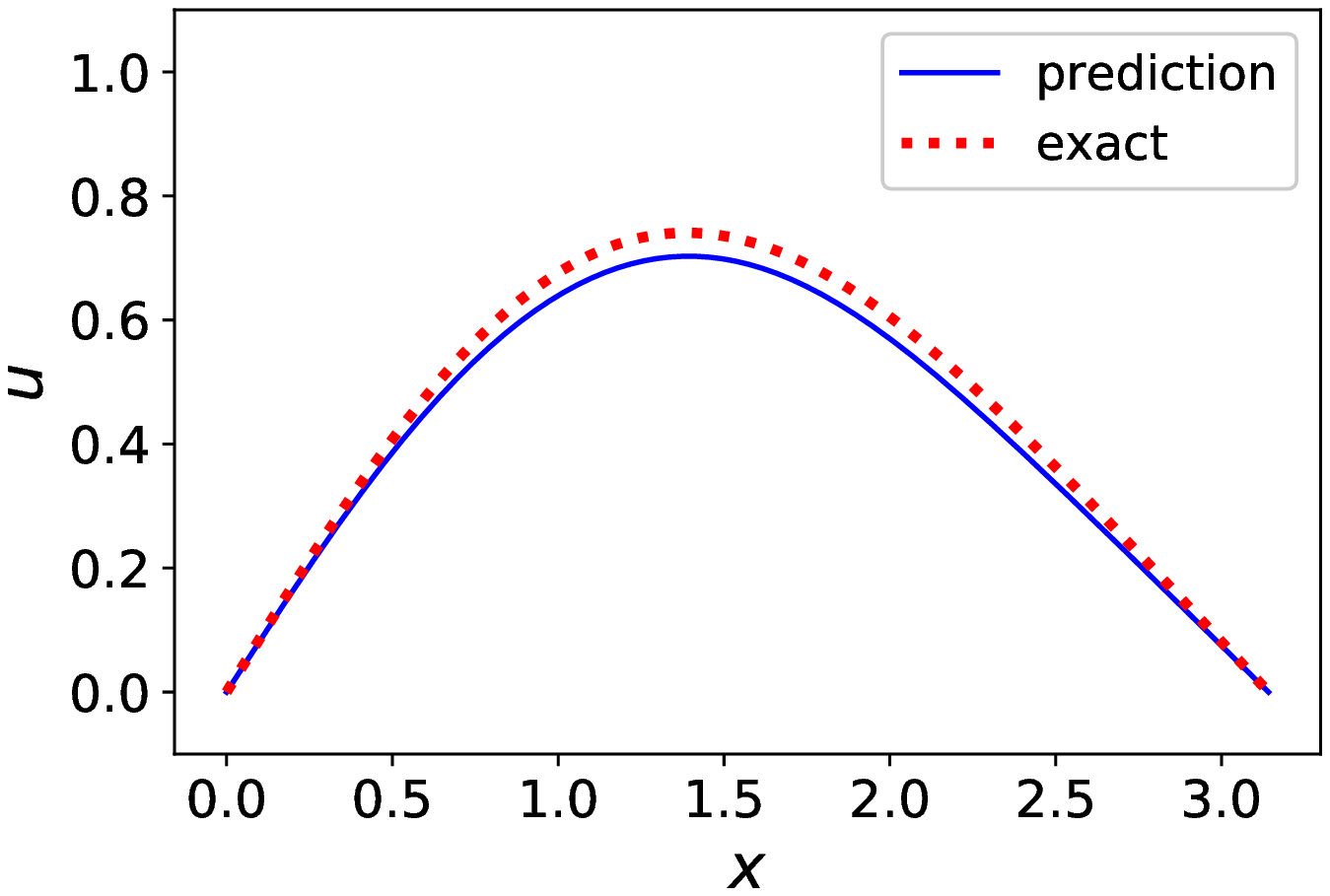}}
	\caption{\small
		Example 2: The solution at different time predicted by the neural network model trained with noisy data. 
		Left: $2 \%$ noise; right: $5 \%$ noise. 
		From top to bottom: solutions at $t=1$, $t=2$ and $t=3$, respectively.  
	}\label{fig:ex2_solu_noisy}
\end{figure}

We now consider the case of noisy data. All training data are then perturbed by 
a multiplicative factor $(1+\epsilon)$, where $\epsilon\sim
[-\eta,\eta]$ follows uniform distribution. We consider two cases of
 $\eta = 0.02$ and $\eta = 0.05$, which respectively correspond to
 $\pm 2 \%$ and $\pm 5 \%$ relative noises in all data.
 In Fig.~\ref{fig:ex2_solu_noisy}, 
 the numerical solutions produced by the neural netwrok models are
 presented, after network training of  $500$ epochs, along with the
 exact solution.
 We observe that the predictions of the network model are fairly
 robust against data noise. At higher level noises in data, the
 predictive results contain relatively larger numerical errors, as expected.
% It can be seen that the proposed neural network model delivers accurate prediction and is robust against noises. As the noise level increases, 
 %the neural network prediction deviates more from the exact solutions as expected, while the main structure of the solution is still well resolved.  Since the time lag $\Delta = 0.1$ is relatively coarse, such noisy data cases are not amenable to some standard equation recovery approaches that require the estimate of temporal derivative of the solution.   

\subsection{Example 3: Viscous Burgers' Equation}

We now consider the viscous Burgers' equation with Dirichlet boundary condition: 
\begin{equation}
\label{eq:example3}
\begin{cases}
u_t + \left( \frac{u^2}2 \right)_x = \sigma u_{xx},\quad (x,t) \in (-\pi,\pi) \times \mathbb R^+, \\
u(-\pi,t)=u(\pi,t)=0,\quad t \in \mathbb R^+.
\end{cases}
\end{equation}

\begin{figure}[htbp]
	\centering
	{\includegraphics[width=0.6\textwidth]{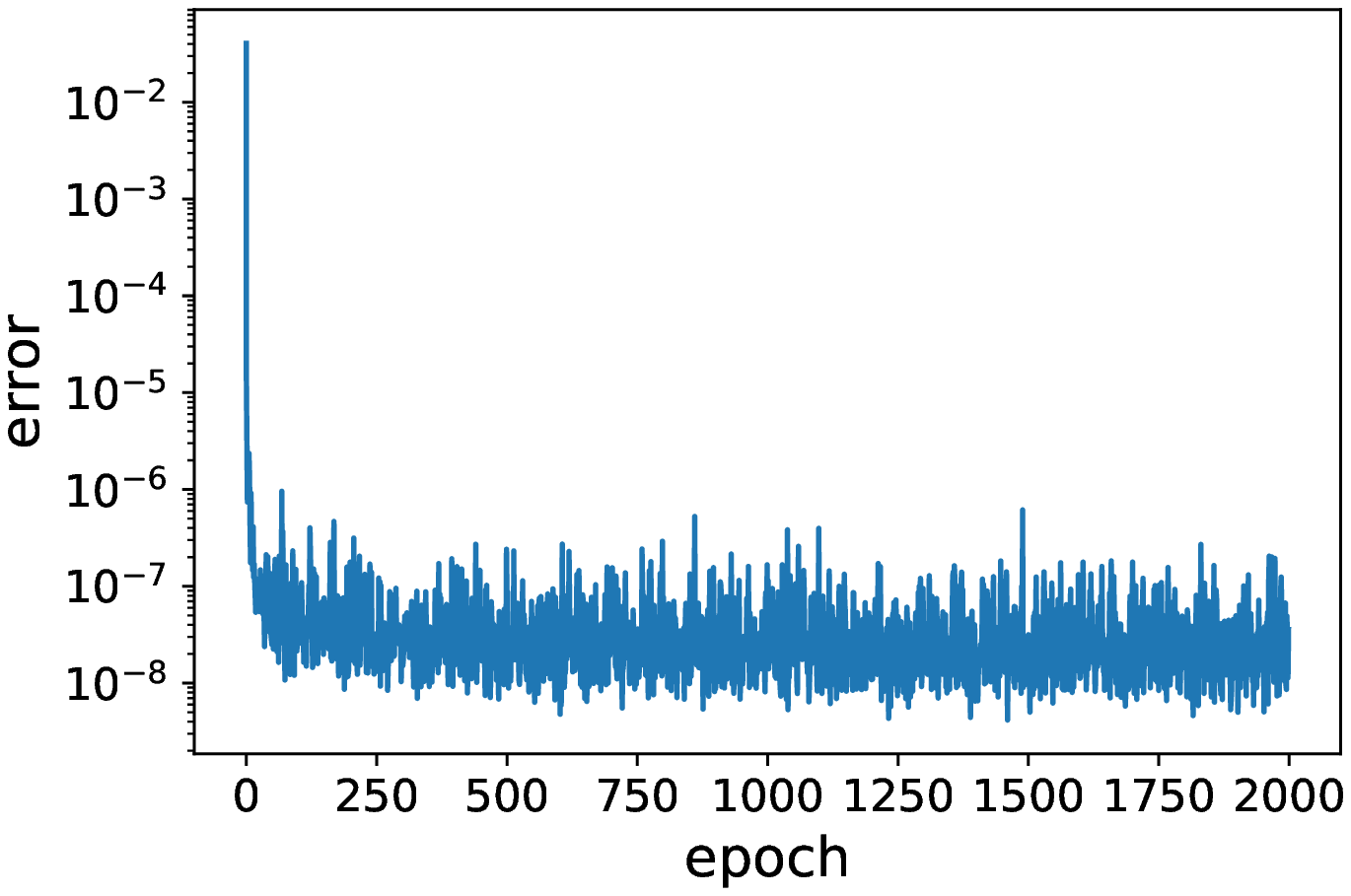}}
	\caption{\small
		Example 3: Training loss history.
	}\label{fig:ex3_loss}
\end{figure}

\begin{figure}[htbp]
	\centering
	{\includegraphics[width=0.48\textwidth]{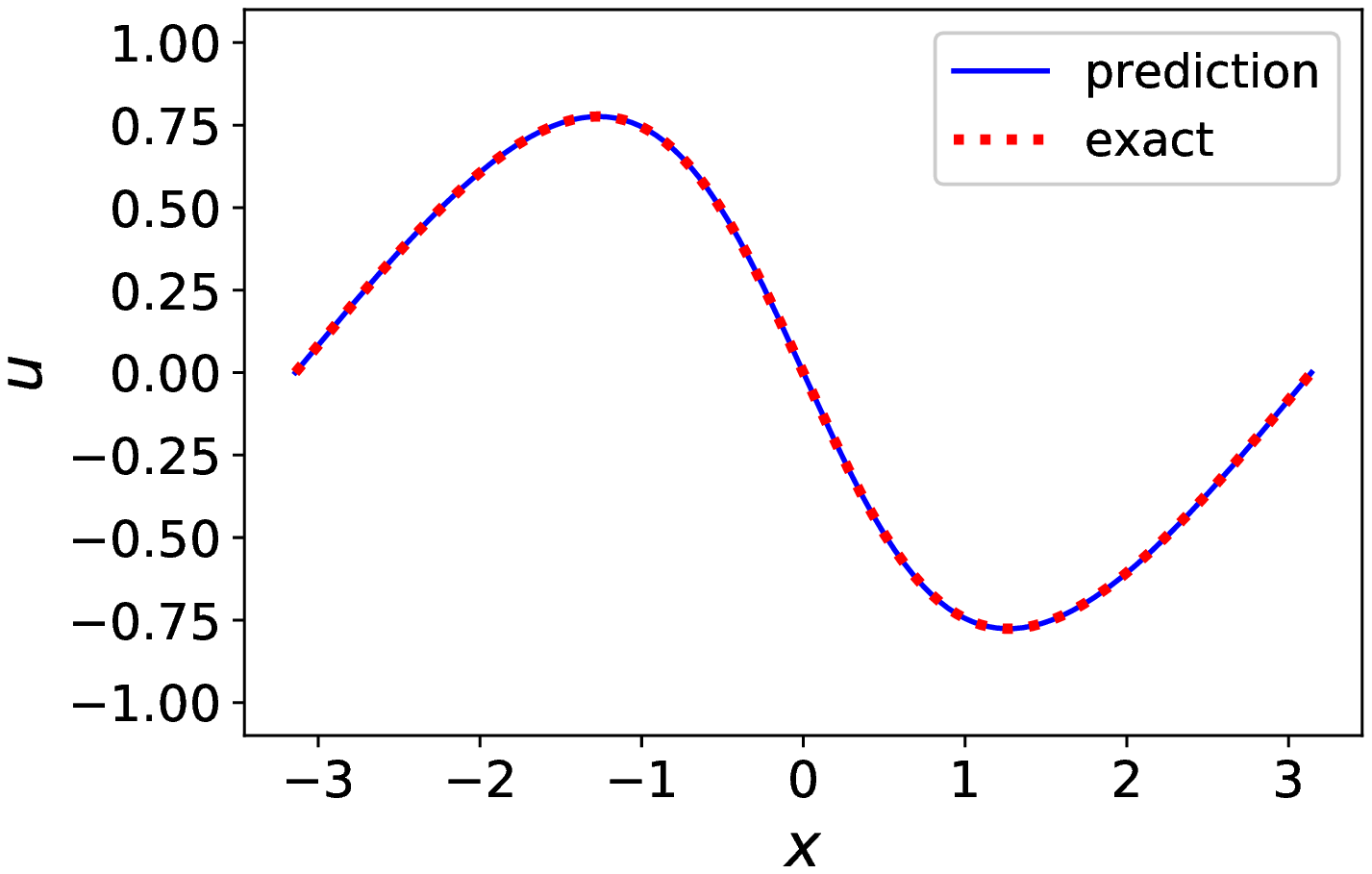}}
	{\includegraphics[width=0.48\textwidth]{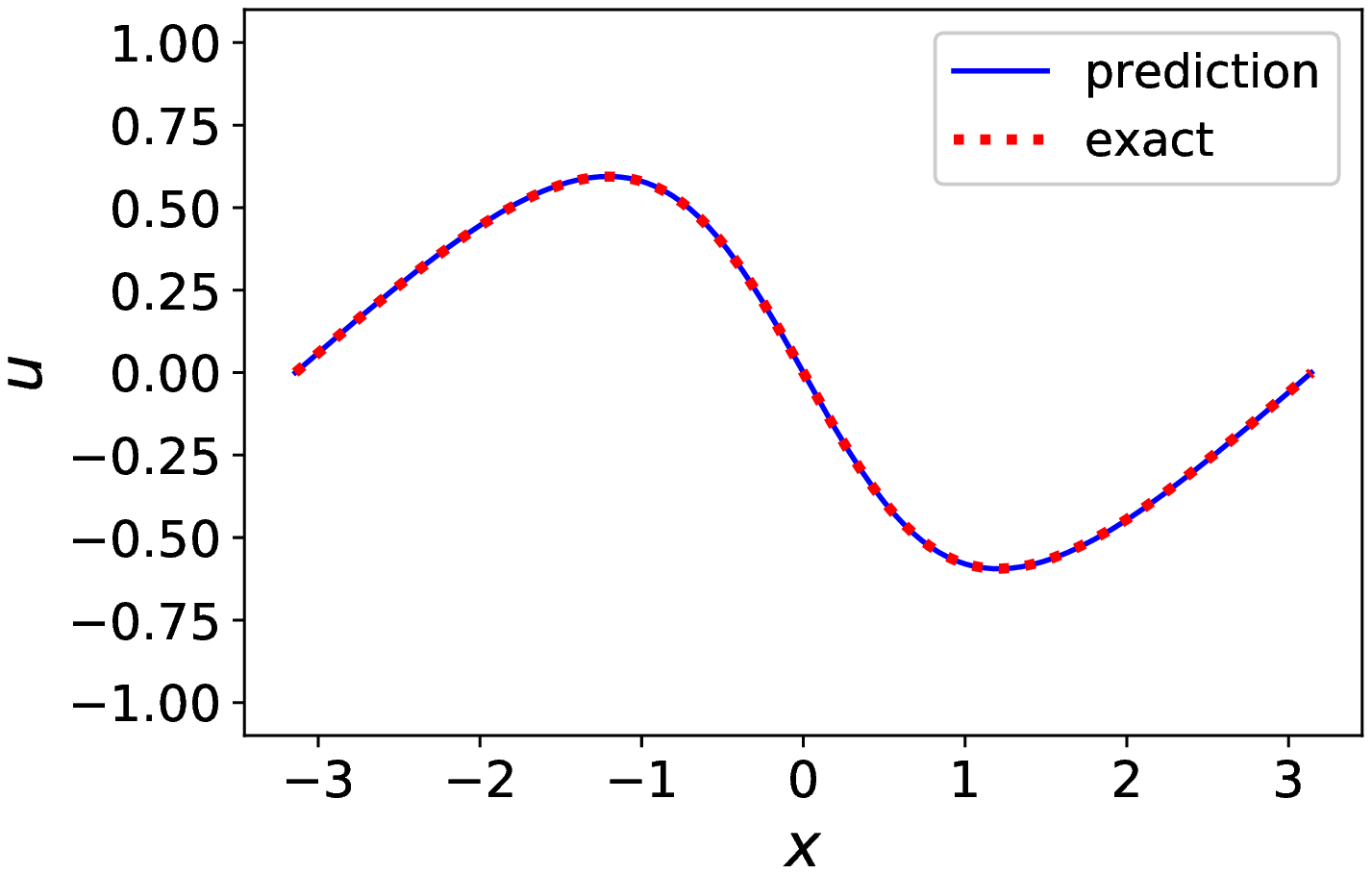}}
	{\includegraphics[width=0.48\textwidth]{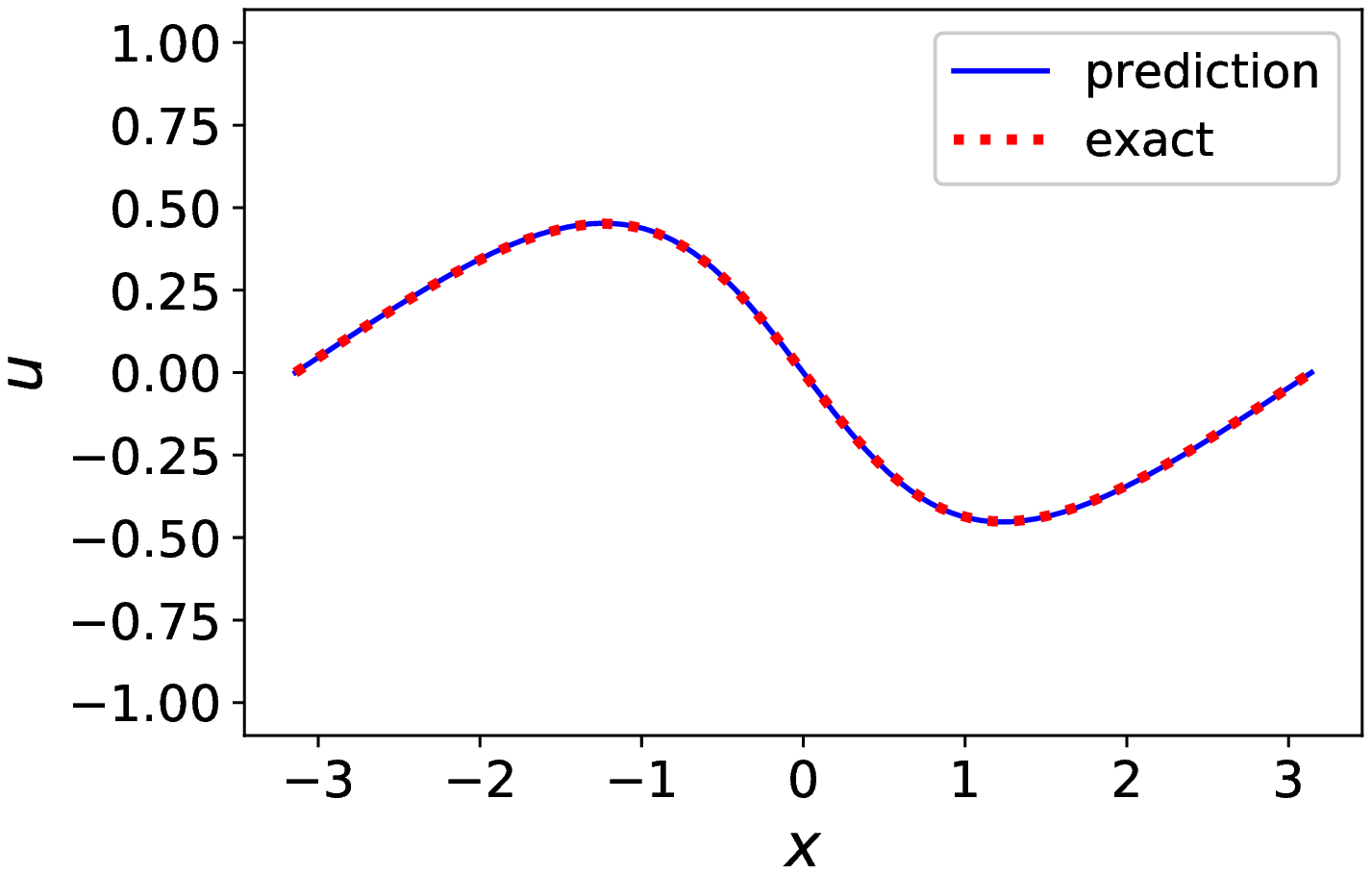}}
	{\includegraphics[width=0.48\textwidth]{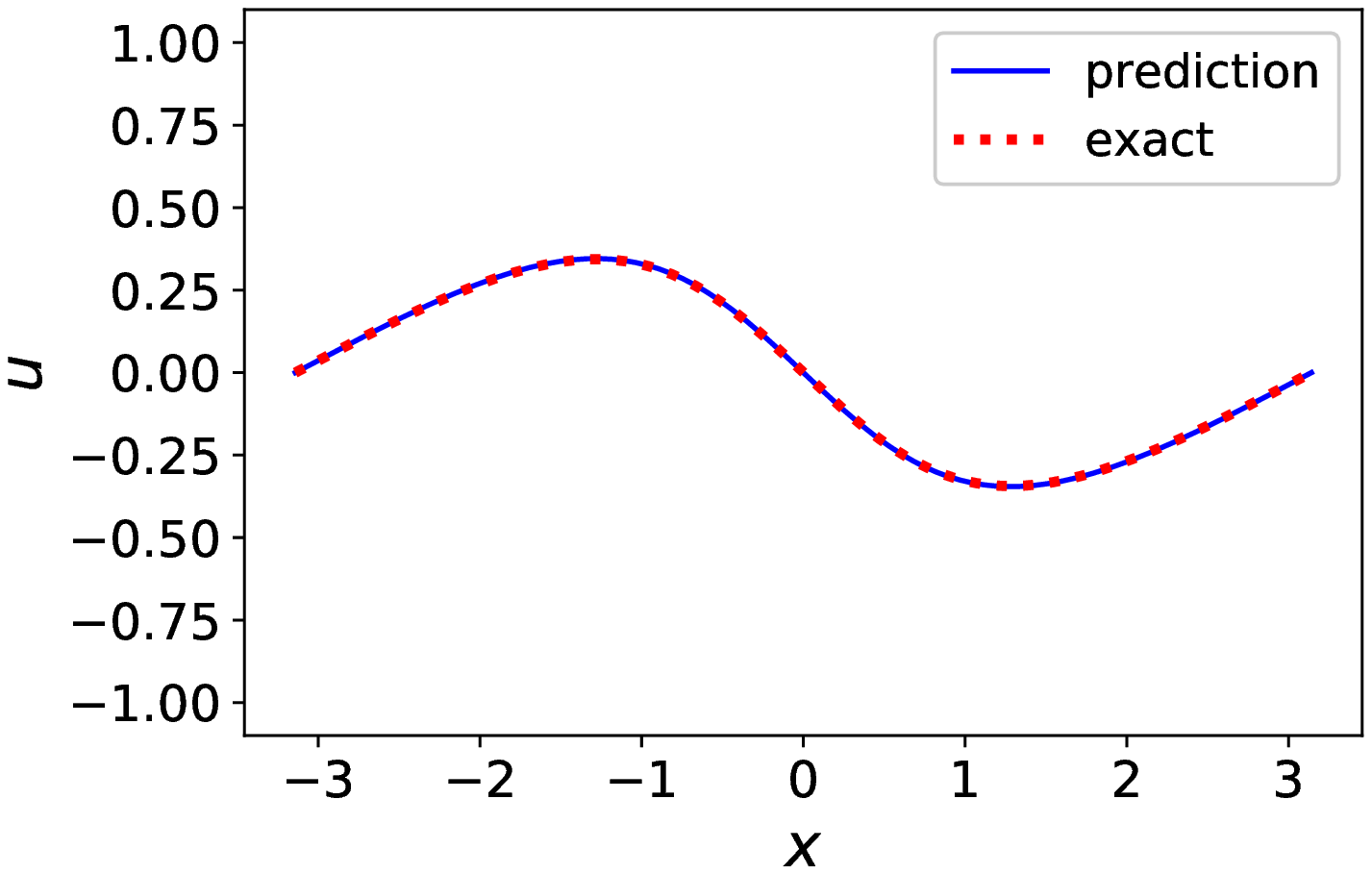}}	
	\caption{\small
	Example 3: Comparison of the true solution, the learned model solution and the solution by Galerkin method at different time. Top-left: $t=0.5$; top-right: $t=1$; bottom-left: $t=1.5$; bottom-right: $t=2$.  
	}\label{fig:ex3_solu}
\end{figure}

\begin{figure}[htbp]
	\centering
%	{\includegraphics[width=0.48\textwidth]{Figure/Example3a/ABSerror.eps}}
	{\includegraphics[width=0.6\textwidth]{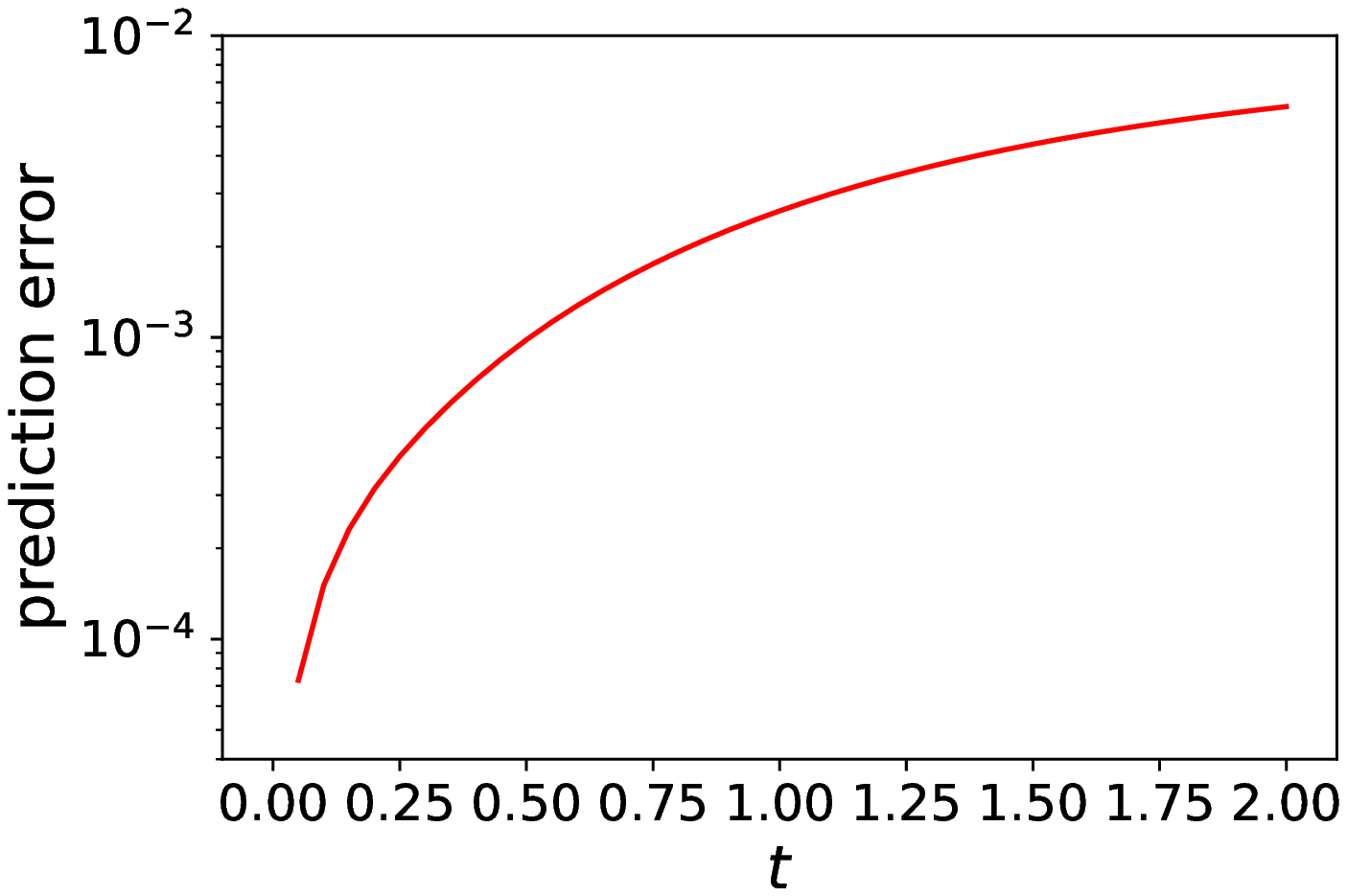}}
	\caption{\small
		Example 3: The evolution of the relative  errors in
                the prediction in $l^2$-norm.
                %Left: absolute error; right: relative error.  
	}\label{fig:ex3_error}
\end{figure}

\begin{figure}[htbp]
	\centering
	{\includegraphics[width=0.325\textwidth]{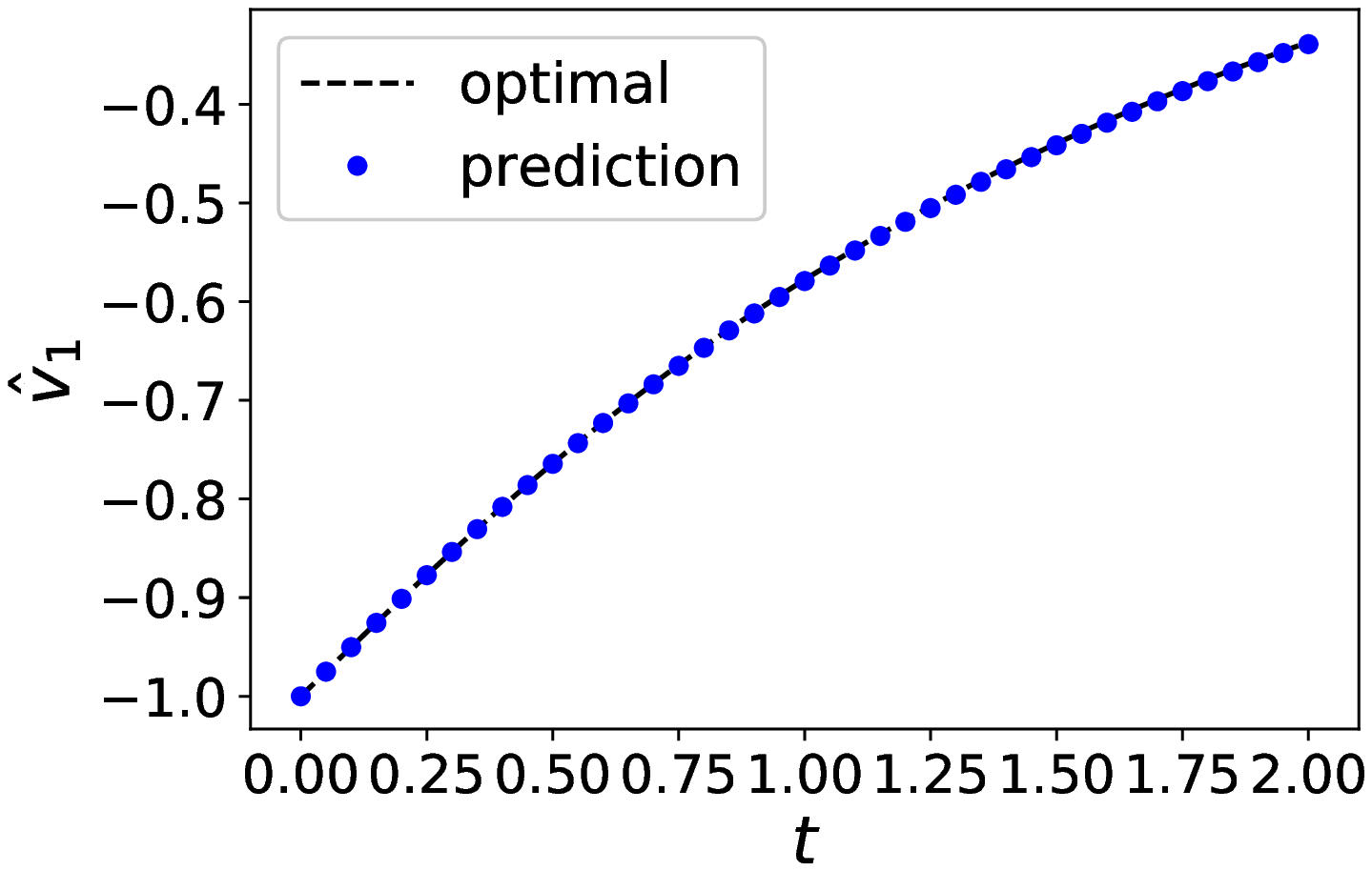}}
	{\includegraphics[width=0.325\textwidth]{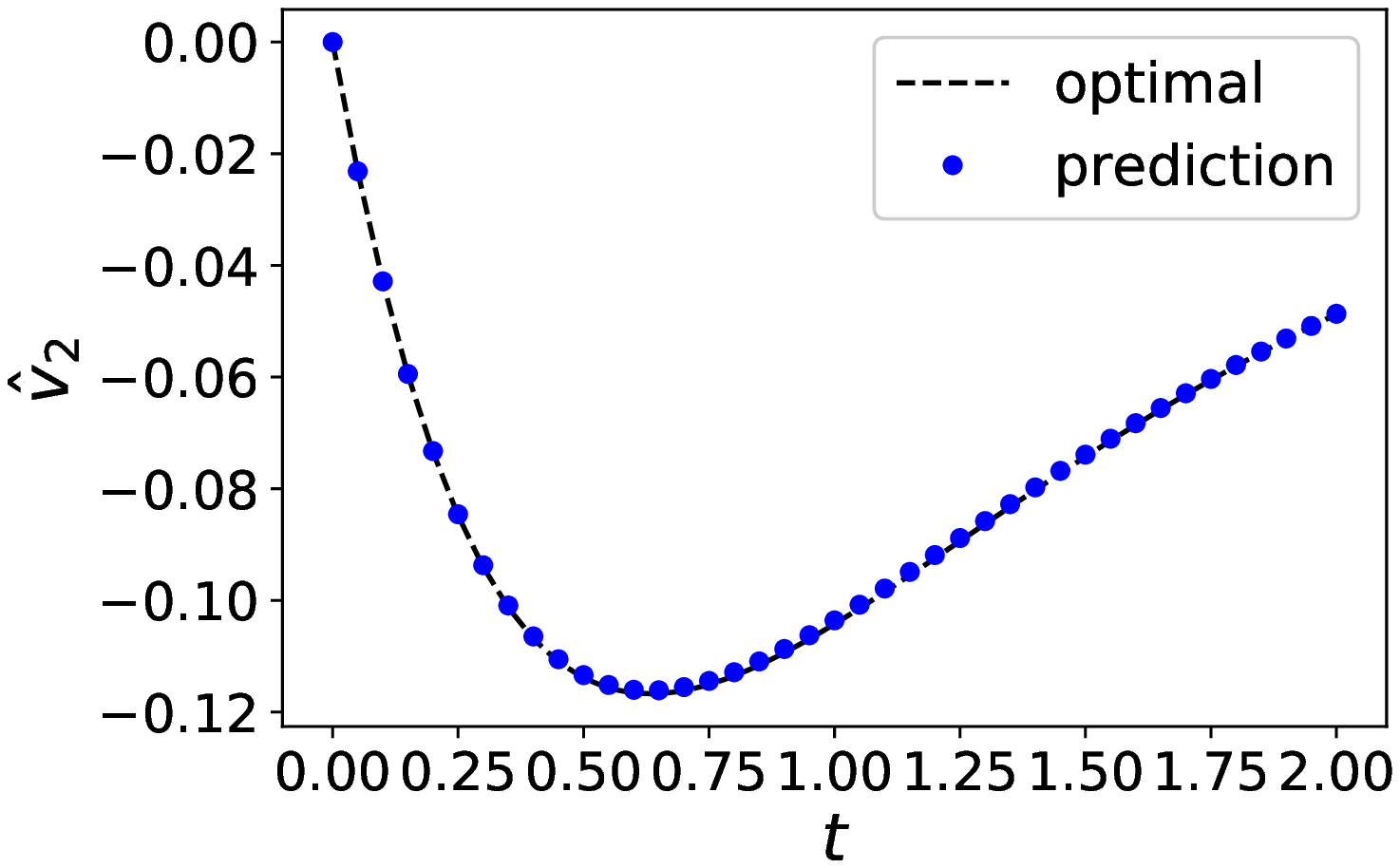}}
	{\includegraphics[width=0.325\textwidth]{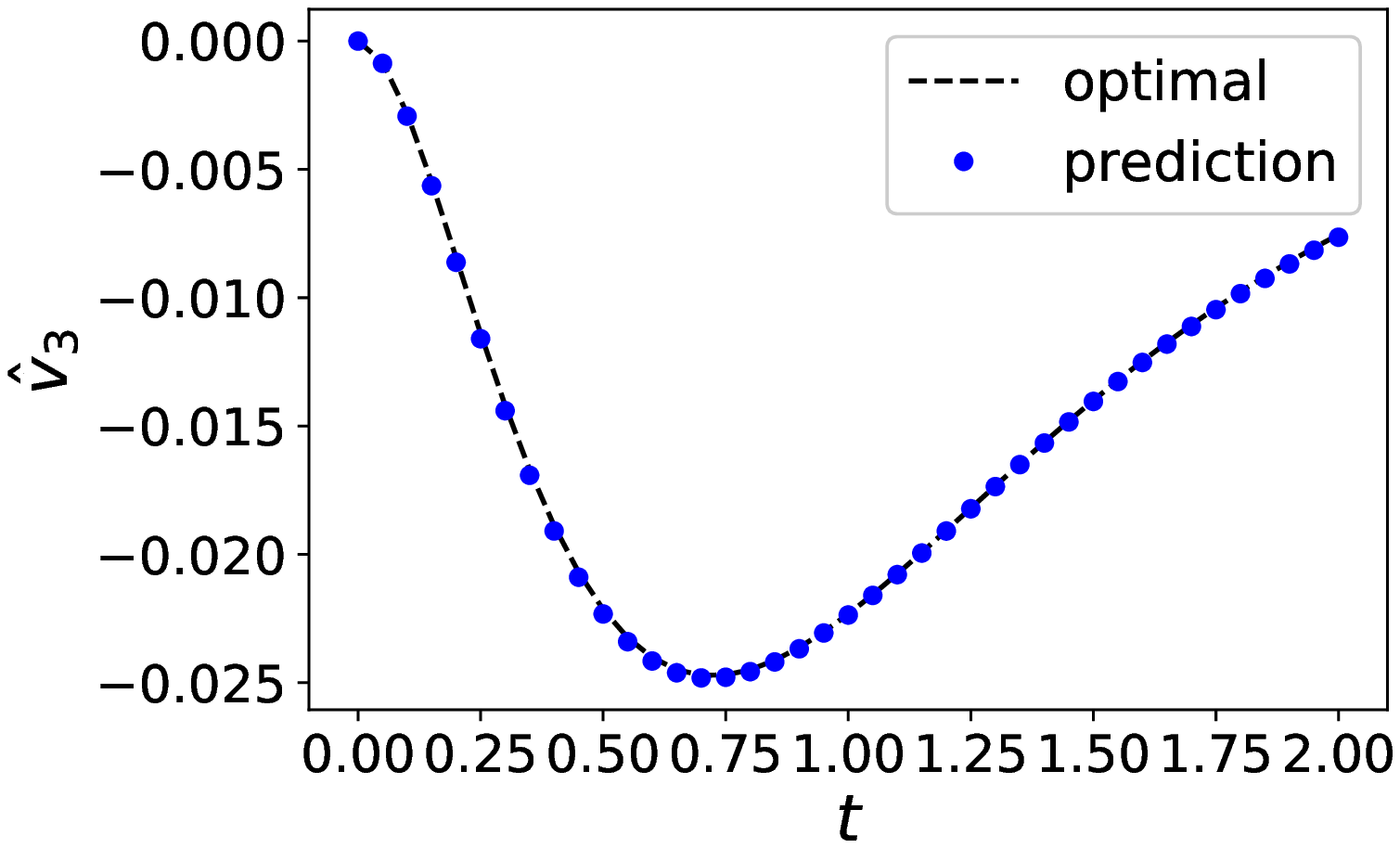}}
	{\includegraphics[width=0.325\textwidth]{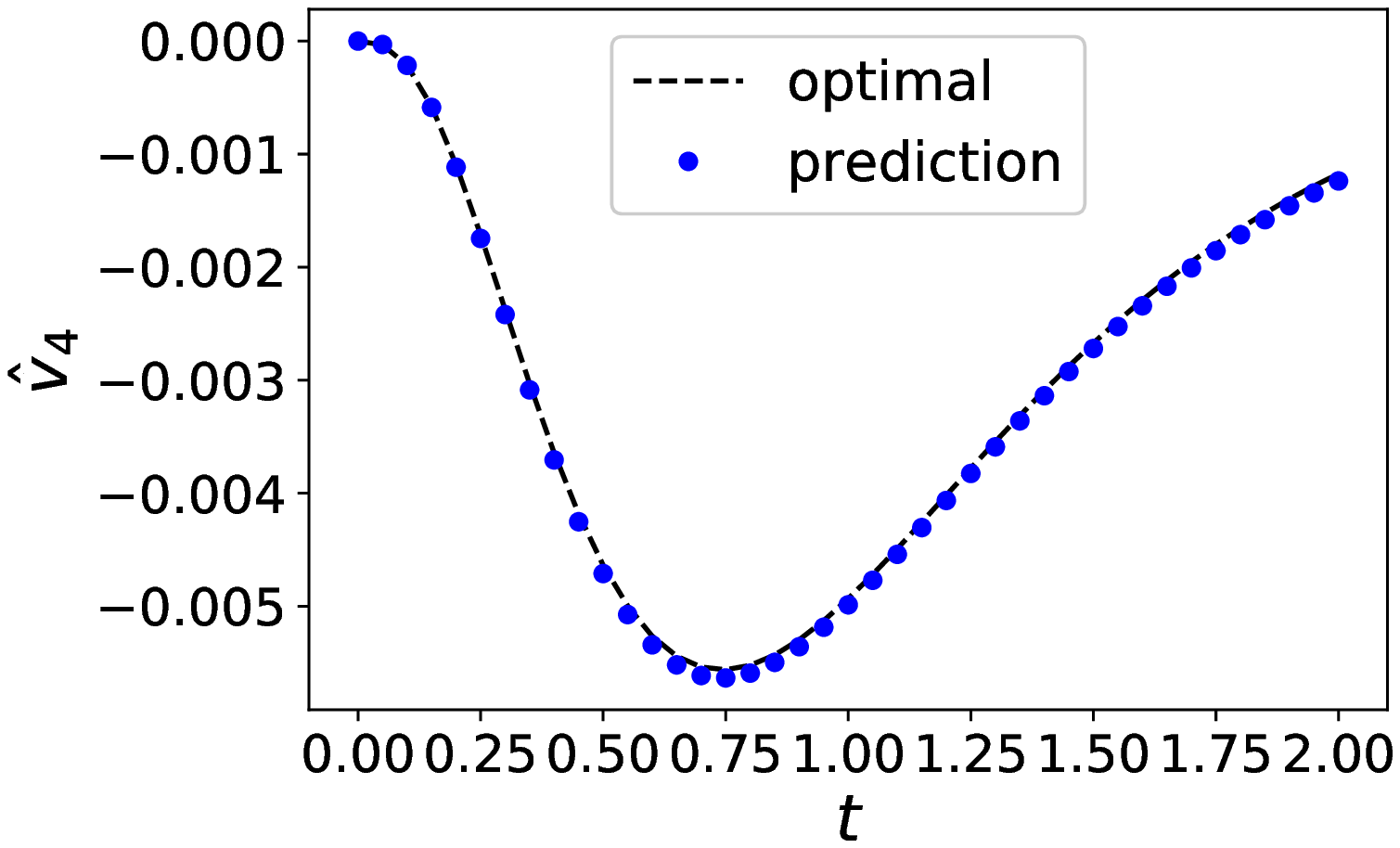}}
	{\includegraphics[width=0.325\textwidth]{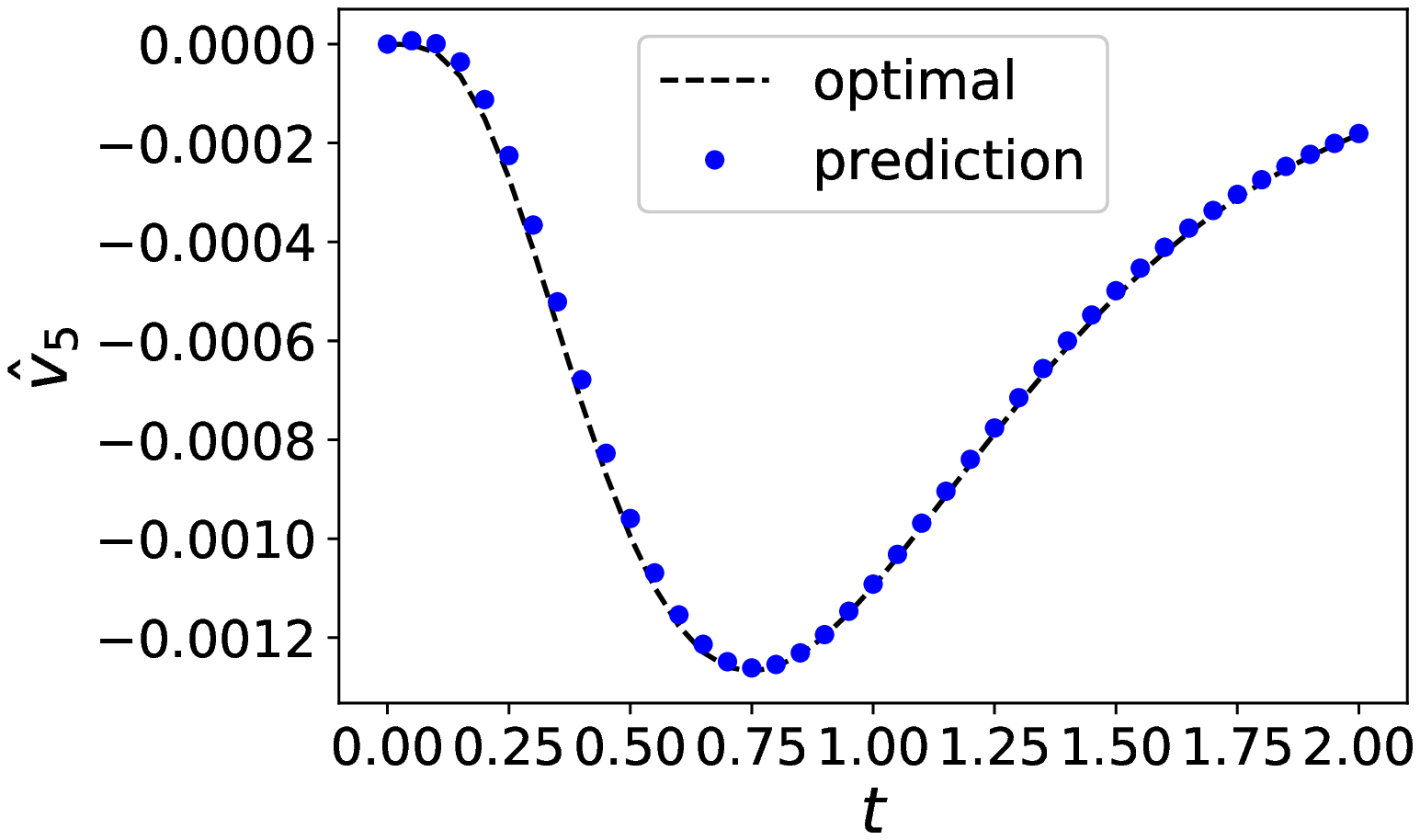}}
	\caption{\small
		Example 3: Evolution of the expansion coefficients for the learned model and the projection of the true solution.
	}\label{fig:ex3_coef}
\end{figure}

We first consider a modestly large viscosity $\sigma = 0.5$. The
approximation space is chosen as $\mathbb V_n = {\rm span} \{
\sin(jx), 1\le j \le 5 \}$ with $n=5$. 
The time lag $\Delta$ is fixed at $\Delta= 0.05$. 
The domain $D$ in the modal space is chosen as 
%$[-1.5,  1.5]\times[-0.5,  0.5]\times[-0.2,  0.2]\times[-0.2, 0.2]\times[-0.1,  0.1]$
$[-1.5,  1.5]\times[-0.2,  0.2]\times[-0.05,  0.05]\times[-0.01,
0.01]\times[-0.002,  0.002]$, from which $100,000$ training data are generated.
The block ResNet method with two blocks ($K=2$) is used, where each
block contains 3 hidden layers of equal width of 30 neurons. 
Upon training the network model satisfactorily (see Fig.~\ref{fig:ex3_loss} for the training loss history), we validate the
trained model for the initial condition 
$$u_0(x)= -\sin (x),$$ 
for time up to $t=2$. %by using \eqref{model-initial}--\eqref{model-final}. 
In Fig.~\ref{fig:ex3_solu}, we compare the predicted solution against the exact solution at different time. The error of the prediction is computed and displayed in Fig.~\ref{fig:ex3_error}. 
We observe that the network model produces 
accurate prediction results. 
The learned expansion coefficients are  shown in Fig.~\ref{fig:ex3_coef} and agree well with the optimal coefficients given by orthogonal projection of the exact solution.

\begin{figure}[htbp]
	\centering
	{\includegraphics[width=0.6\textwidth]{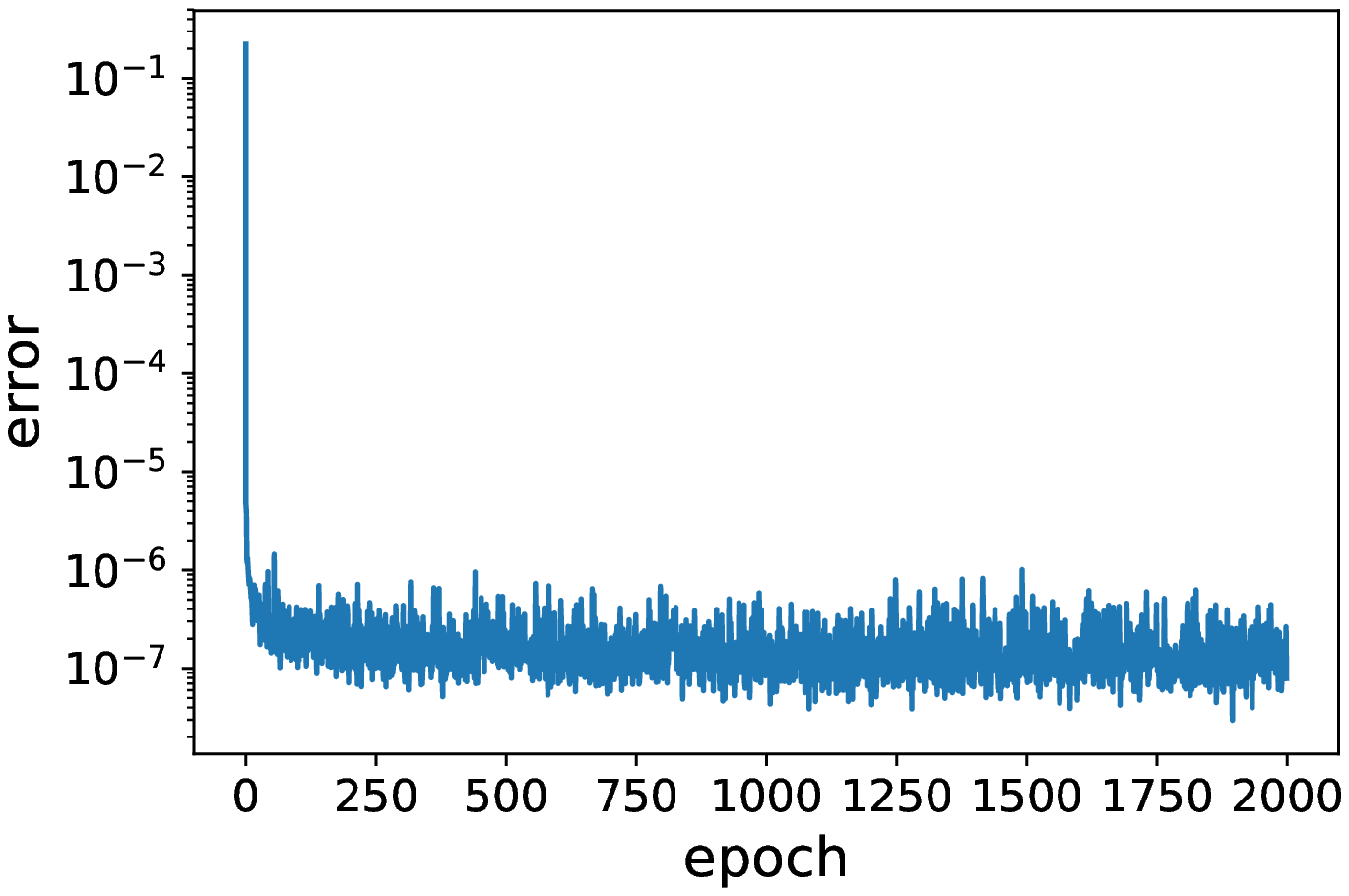}}
	\caption{\small
		Example 3: Training loss history.
	}\label{fig:ex3b_loss}
\end{figure}

\begin{figure}[htbp]
	\centering
	{\includegraphics[width=0.48\textwidth]{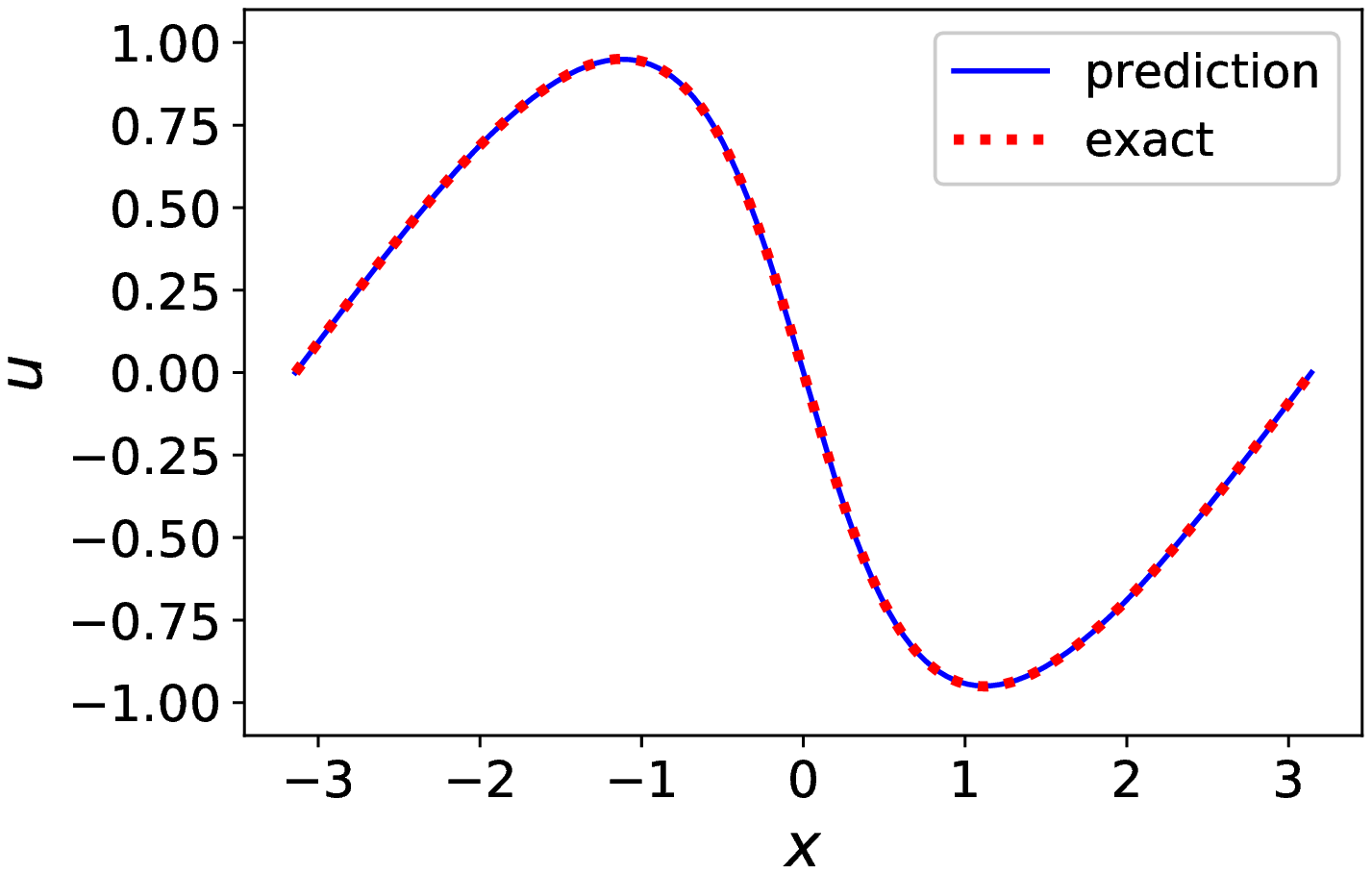}}
	{\includegraphics[width=0.48\textwidth]{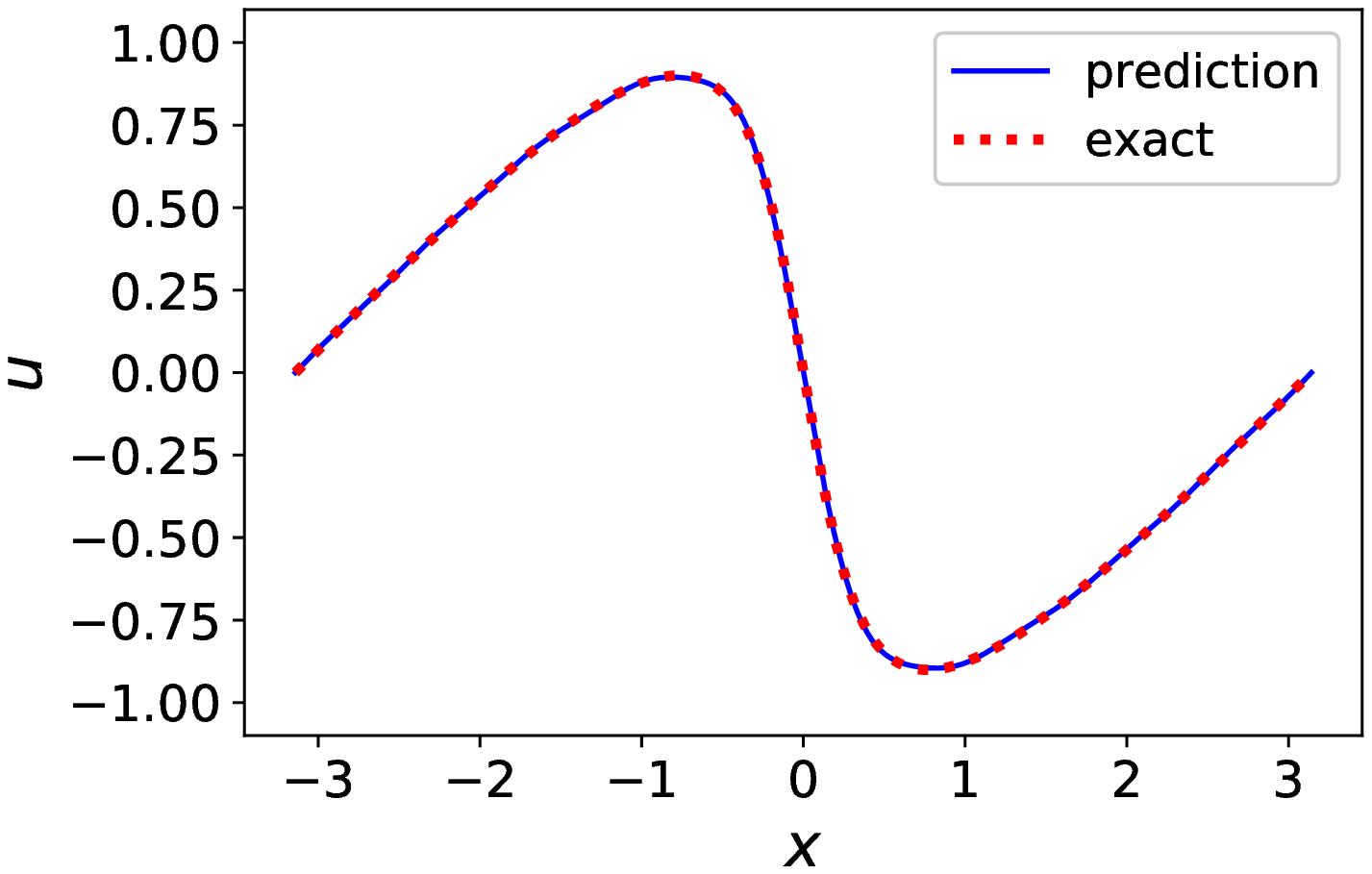}}
	{\includegraphics[width=0.48\textwidth]{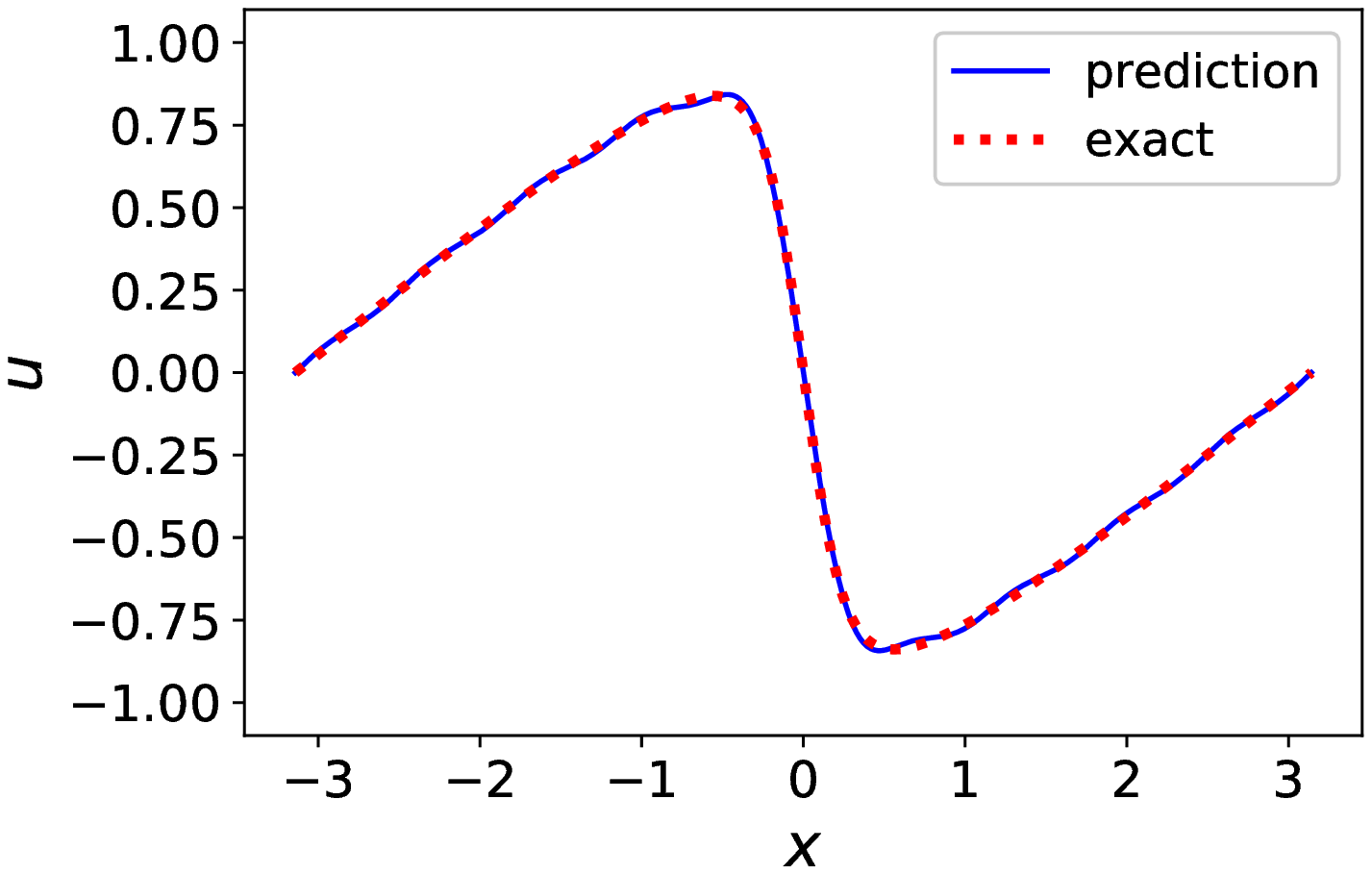}}
	{\includegraphics[width=0.48\textwidth]{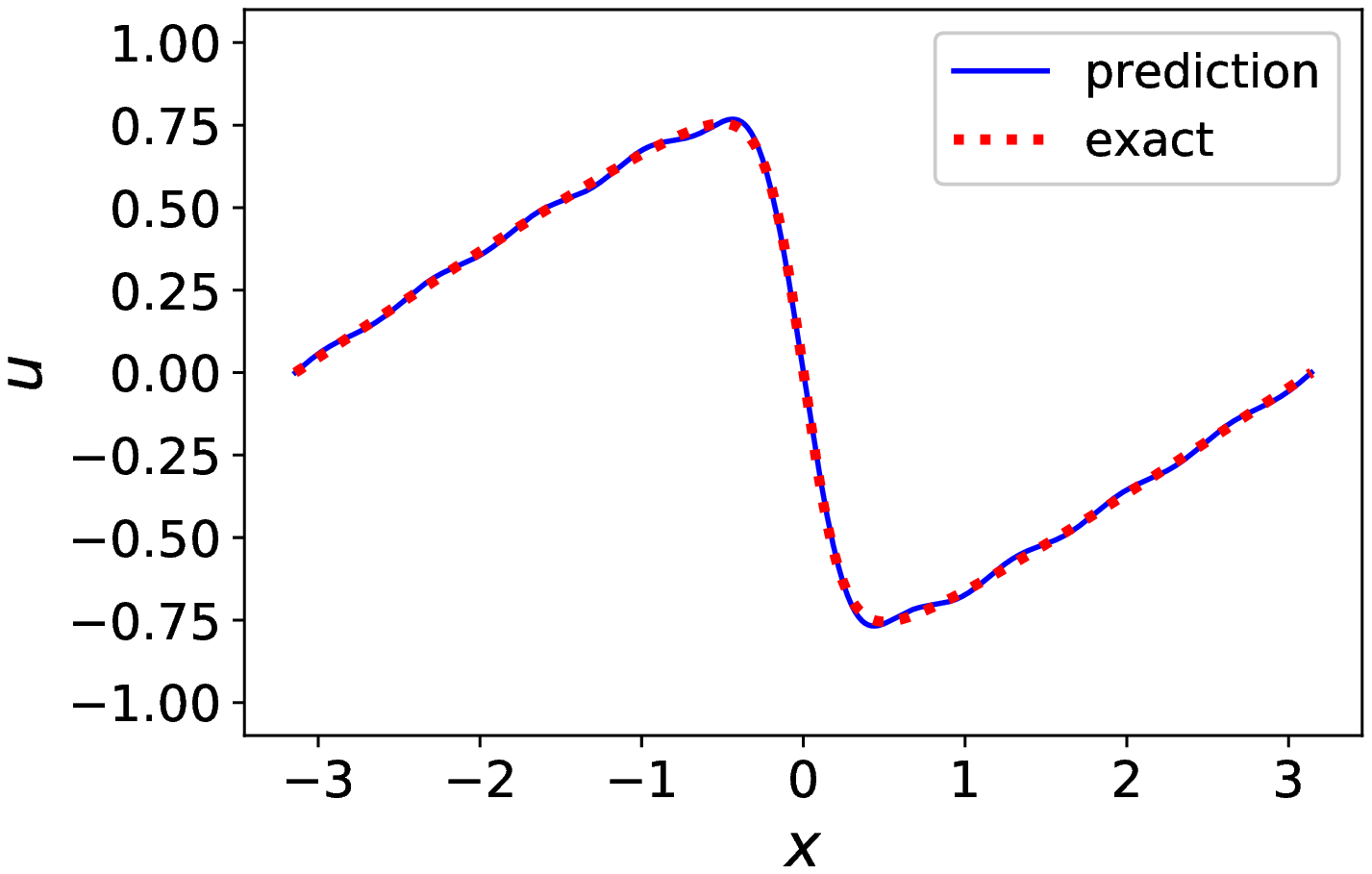}}	
	\caption{\small
	Example 3: Comparison of the true solution, the learned model solution and the solution by Galerkin method at different time. Top-left: $t=0.5$; top-right: $t=1$; bottom-left: $t=1.5$; bottom-right: $t=2$.  
	}\label{fig:ex3b_solu}
\end{figure}

\begin{figure}[htbp]
	\centering
%	{\includegraphics[width=0.48\textwidth]{Figure/Example3b/ABSerror.eps}}
	{\includegraphics[width=0.6\textwidth]{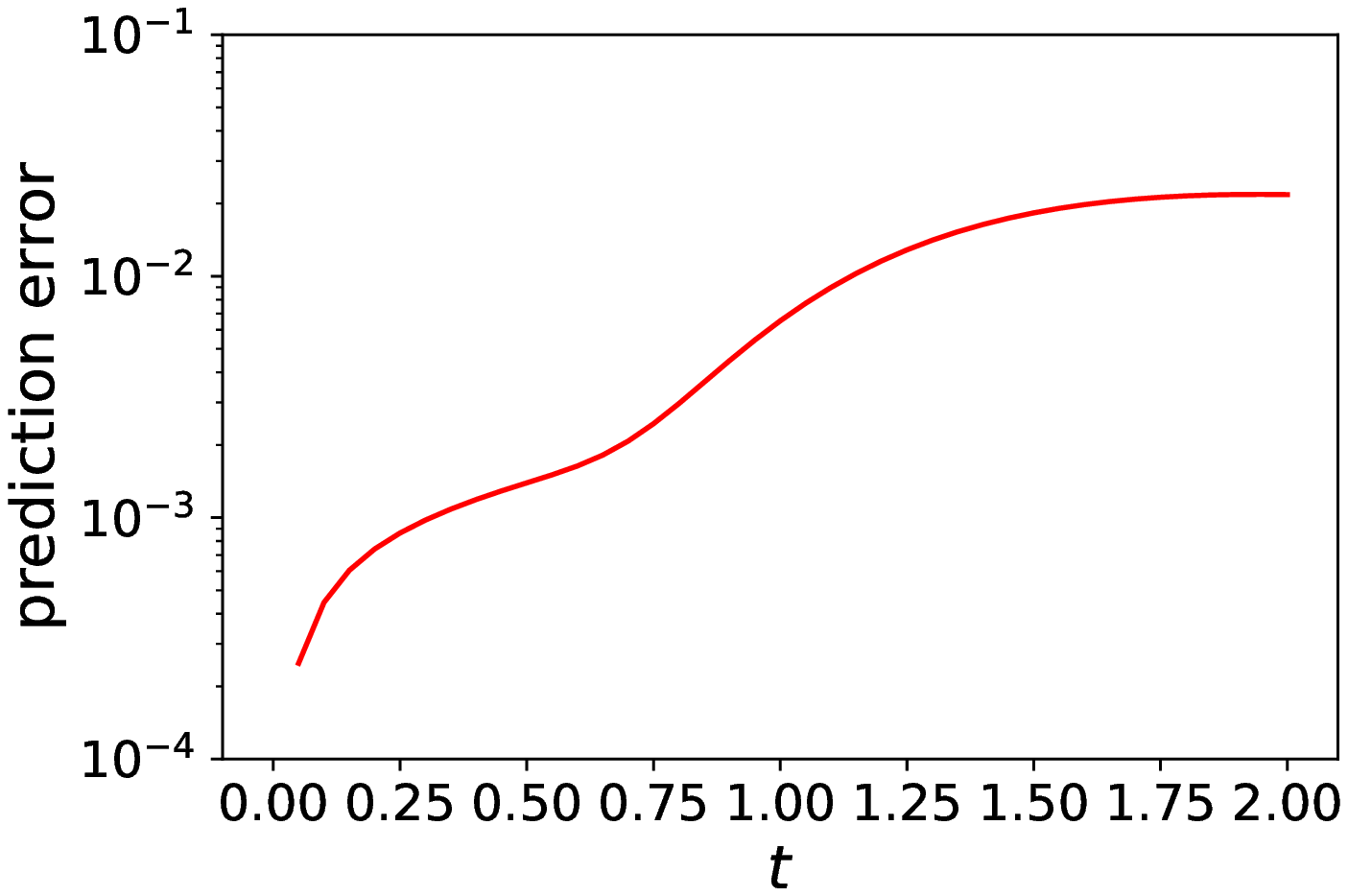}}
	\caption{\small
		Example 3: The evolution of the relative error in
                prediction in $l^2$-norm.
                %Left: absolute error; right: relative error.  
	}\label{fig:ex3b_error}
\end{figure}

\begin{figure}[htbp]
	\centering
	{\includegraphics[width=0.325\textwidth]{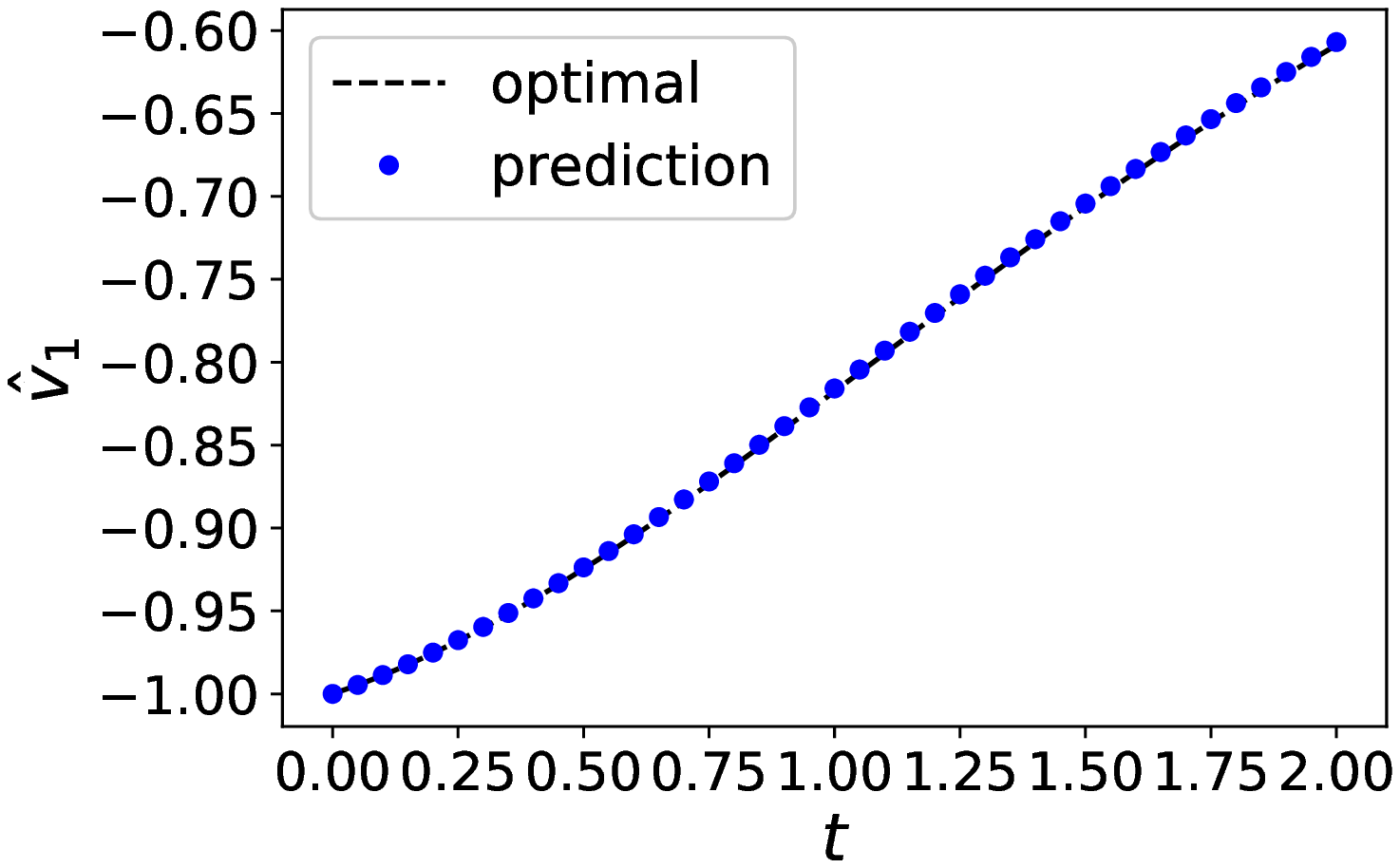}}
	{\includegraphics[width=0.325\textwidth]{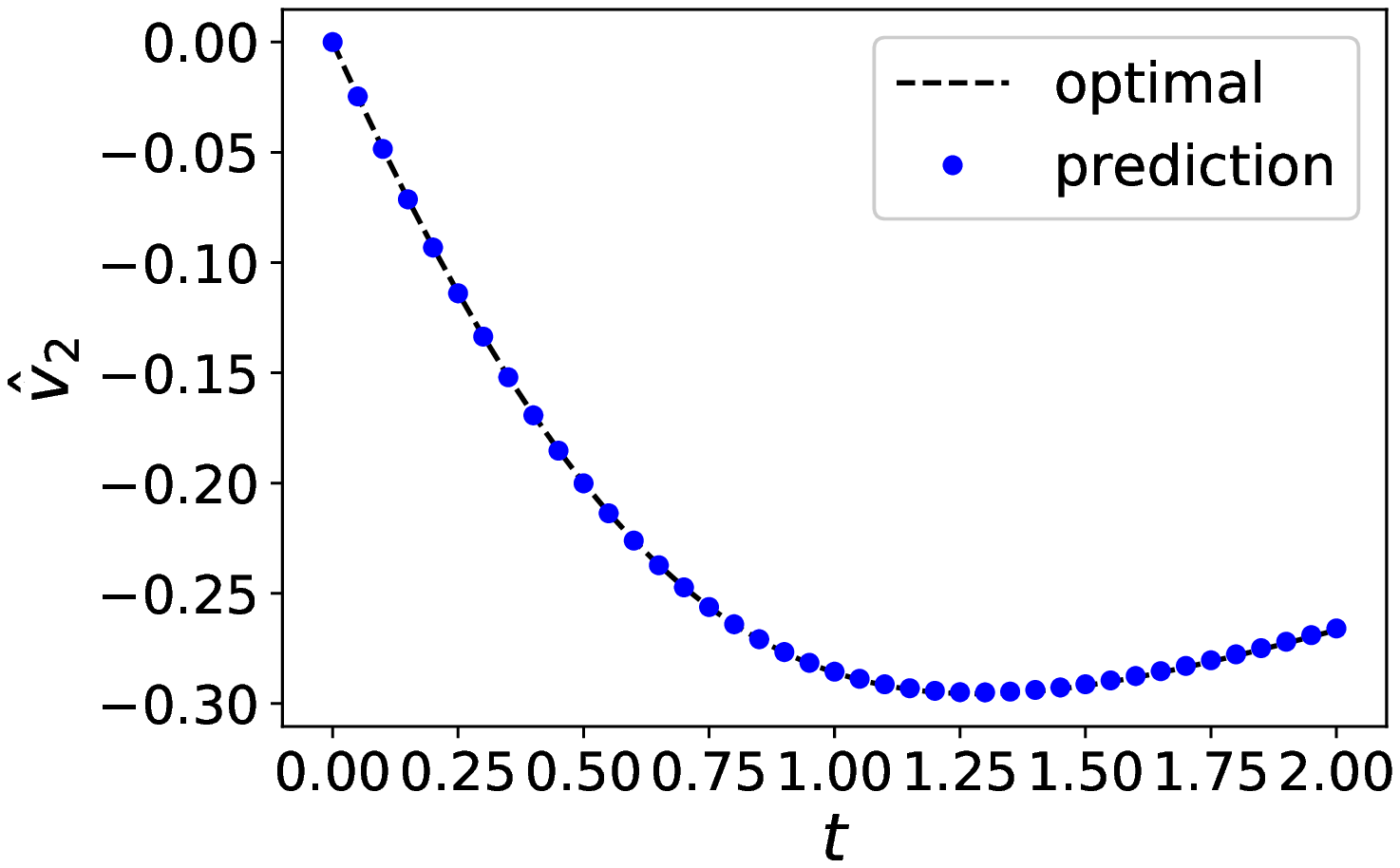}}
	{\includegraphics[width=0.325\textwidth]{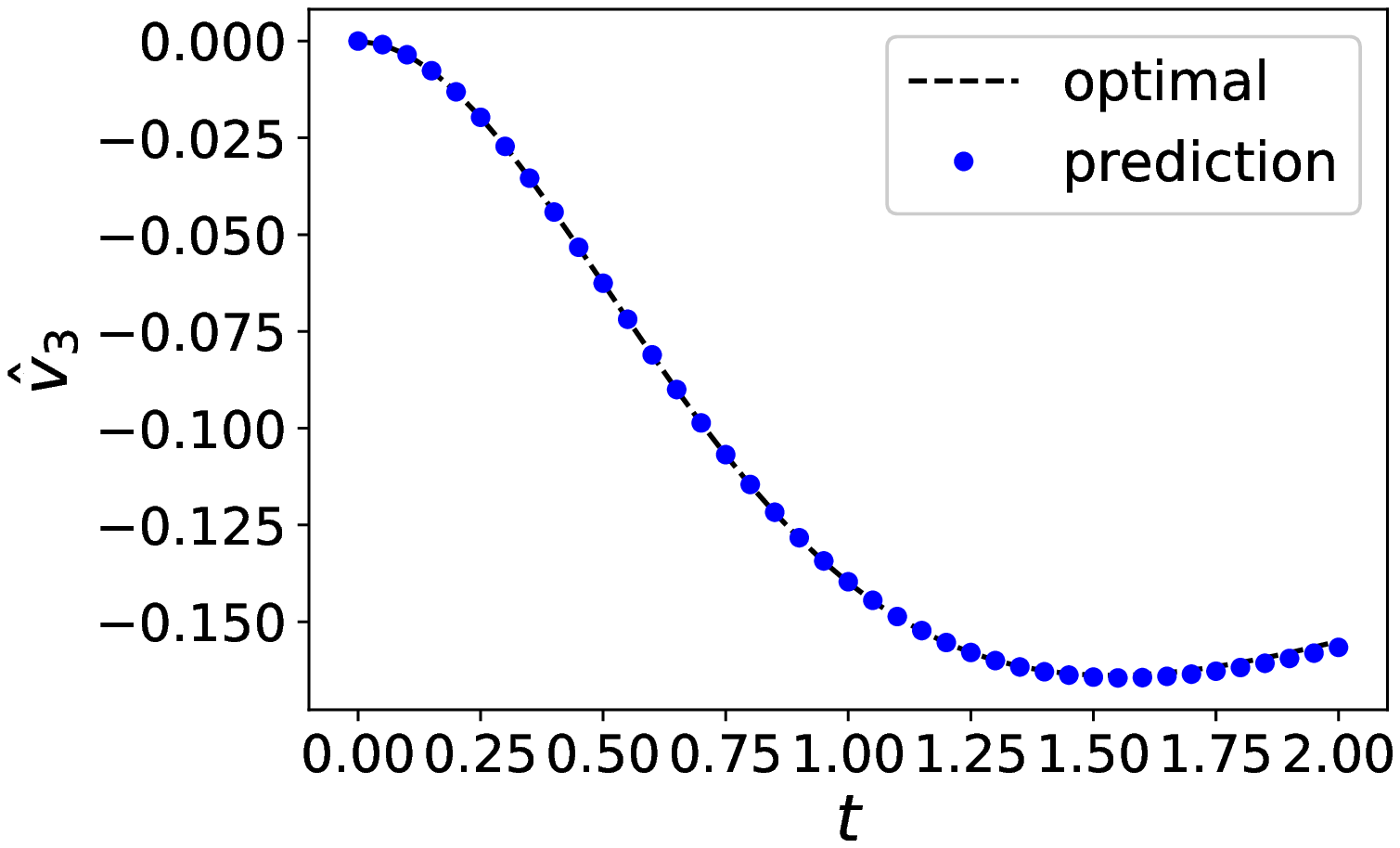}}
	{\includegraphics[width=0.325\textwidth]{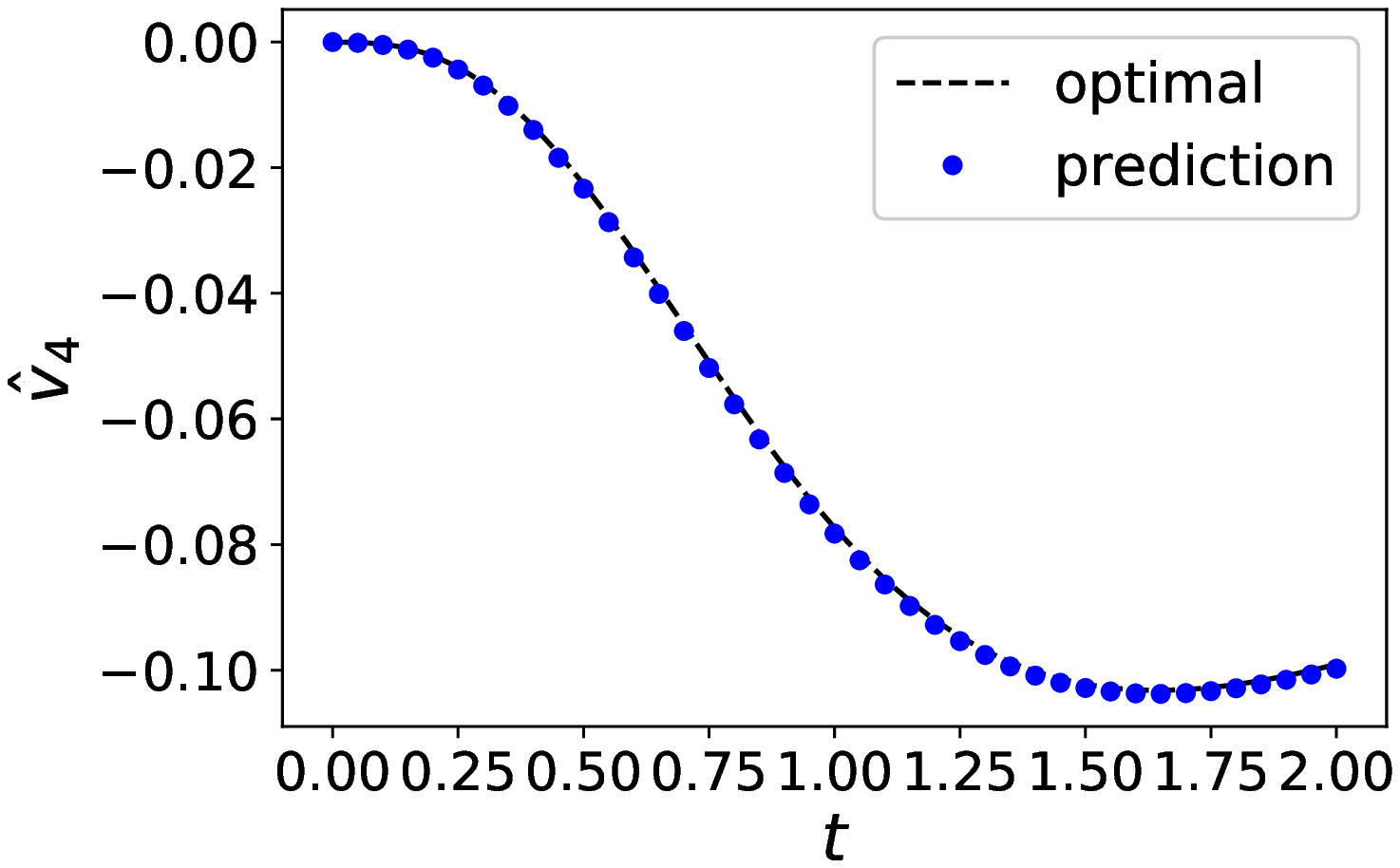}}
	{\includegraphics[width=0.325\textwidth]{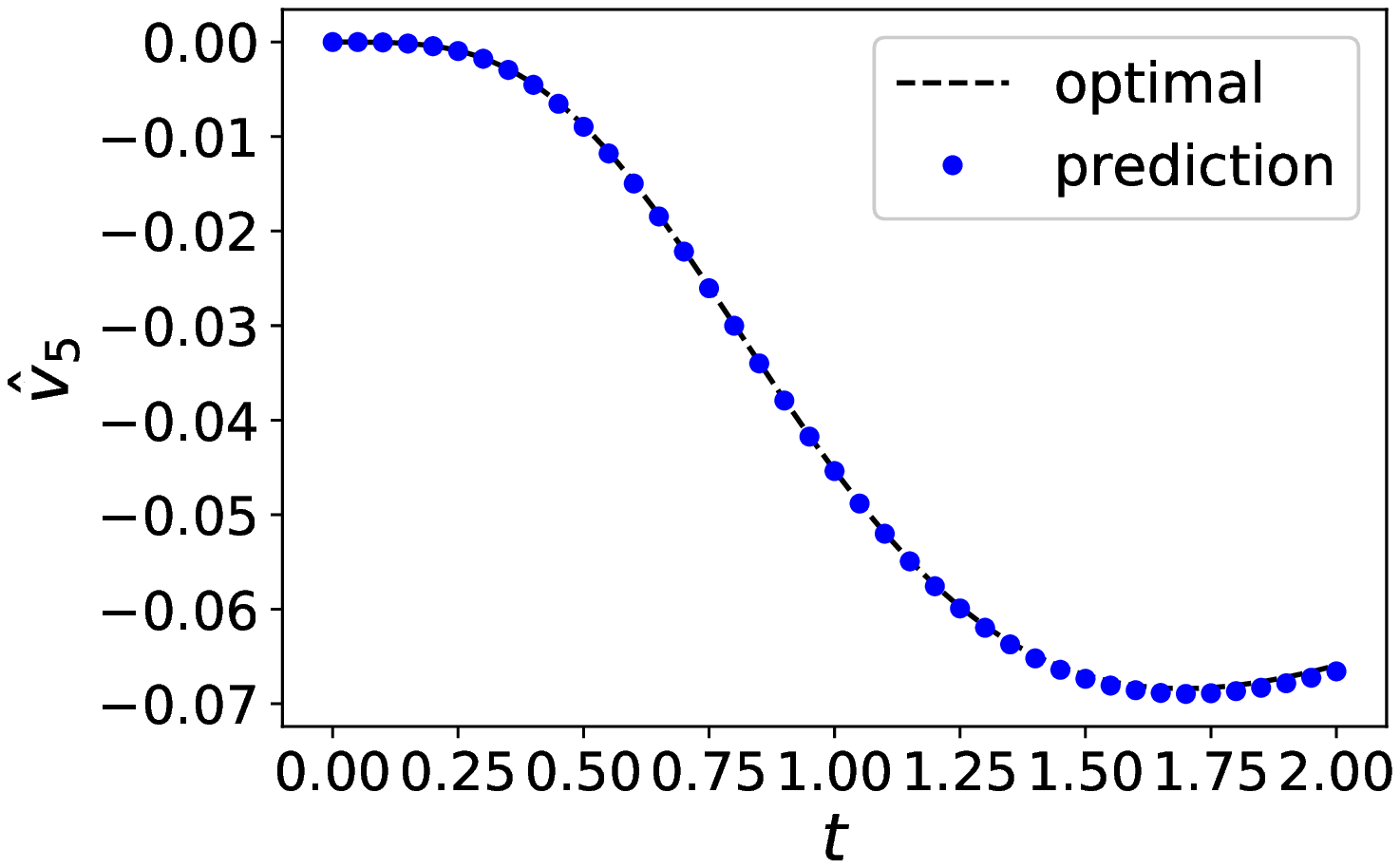}}
	{\includegraphics[width=0.325\textwidth]{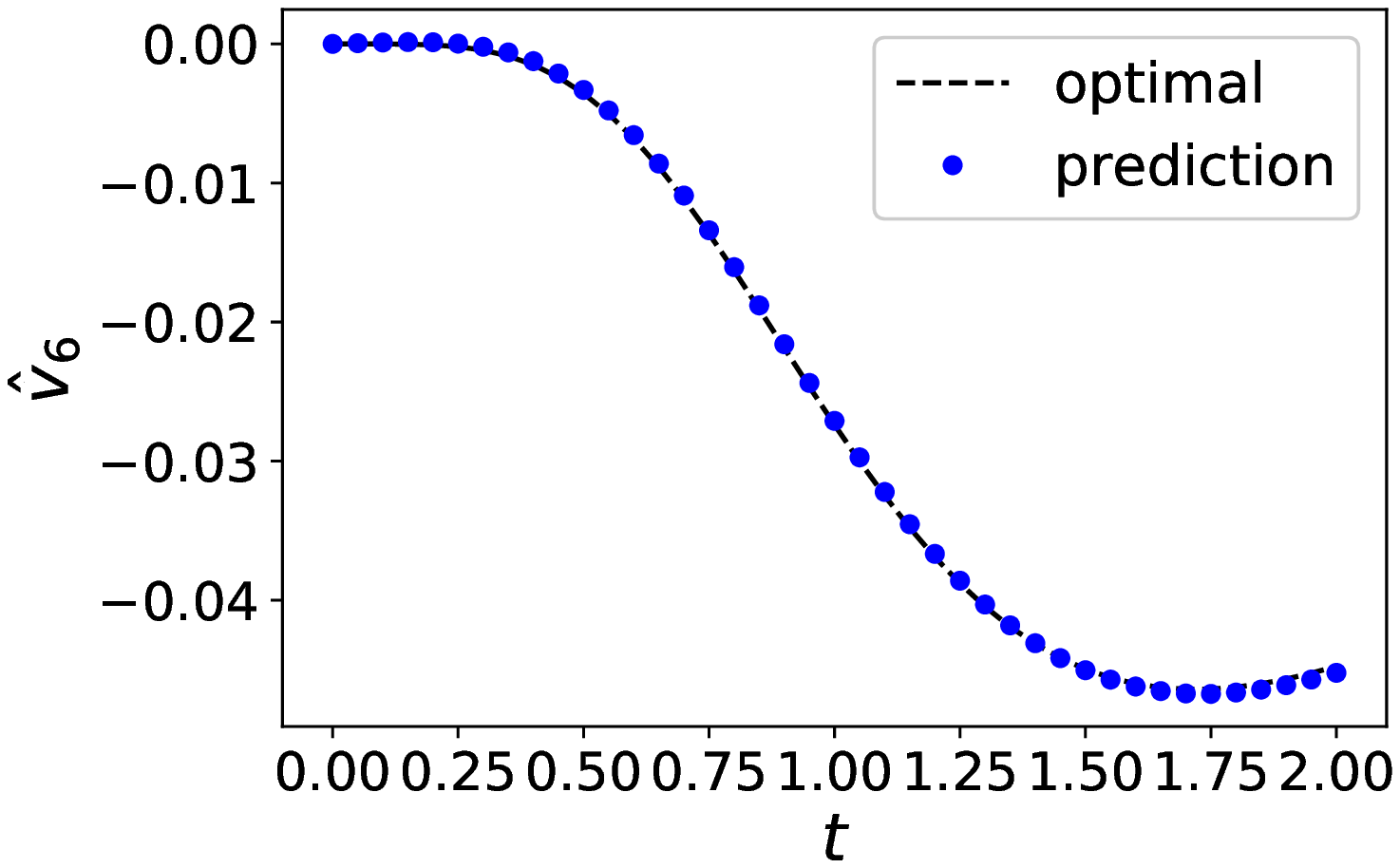}}
	{\includegraphics[width=0.325\textwidth]{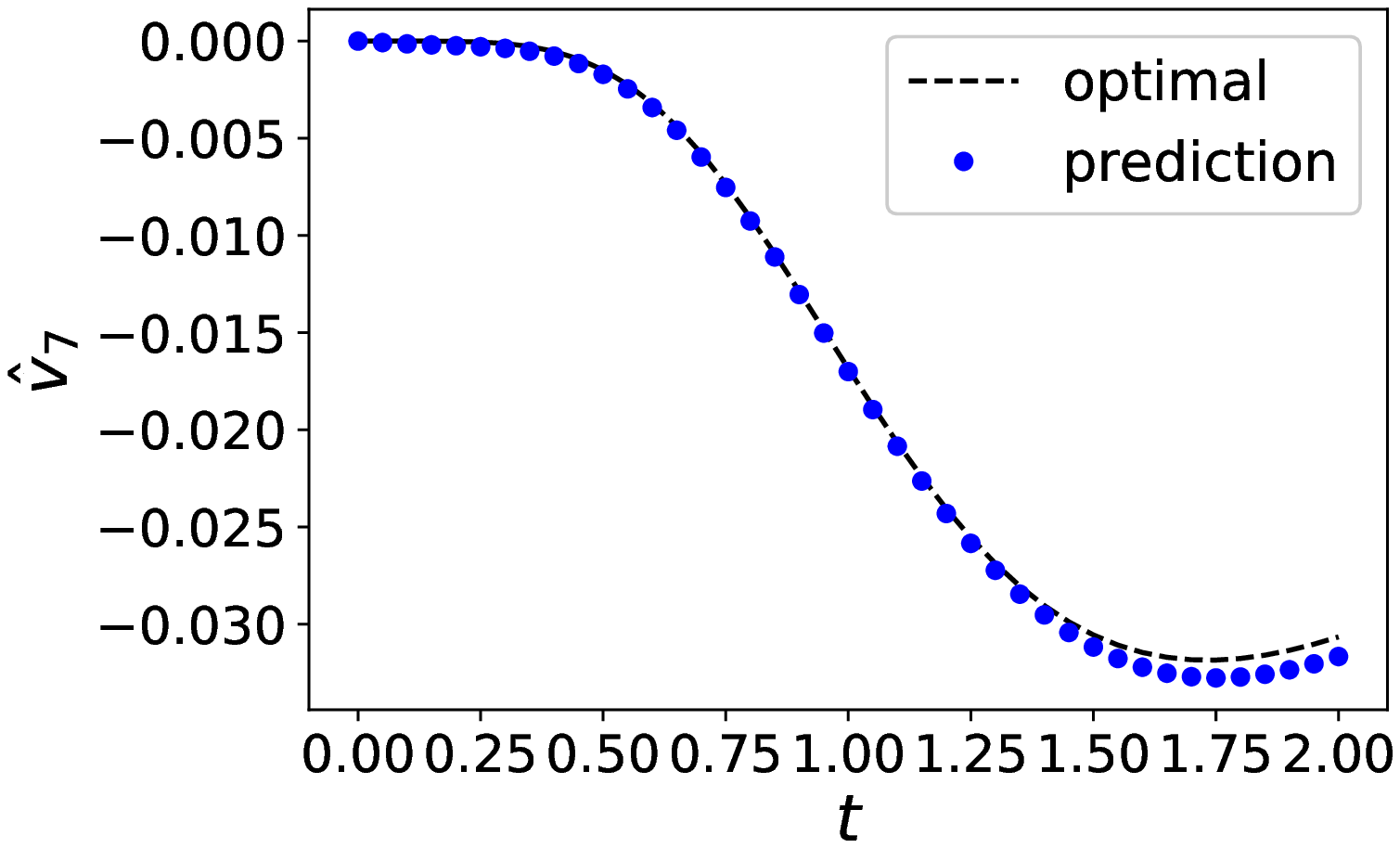}}
	{\includegraphics[width=0.325\textwidth]{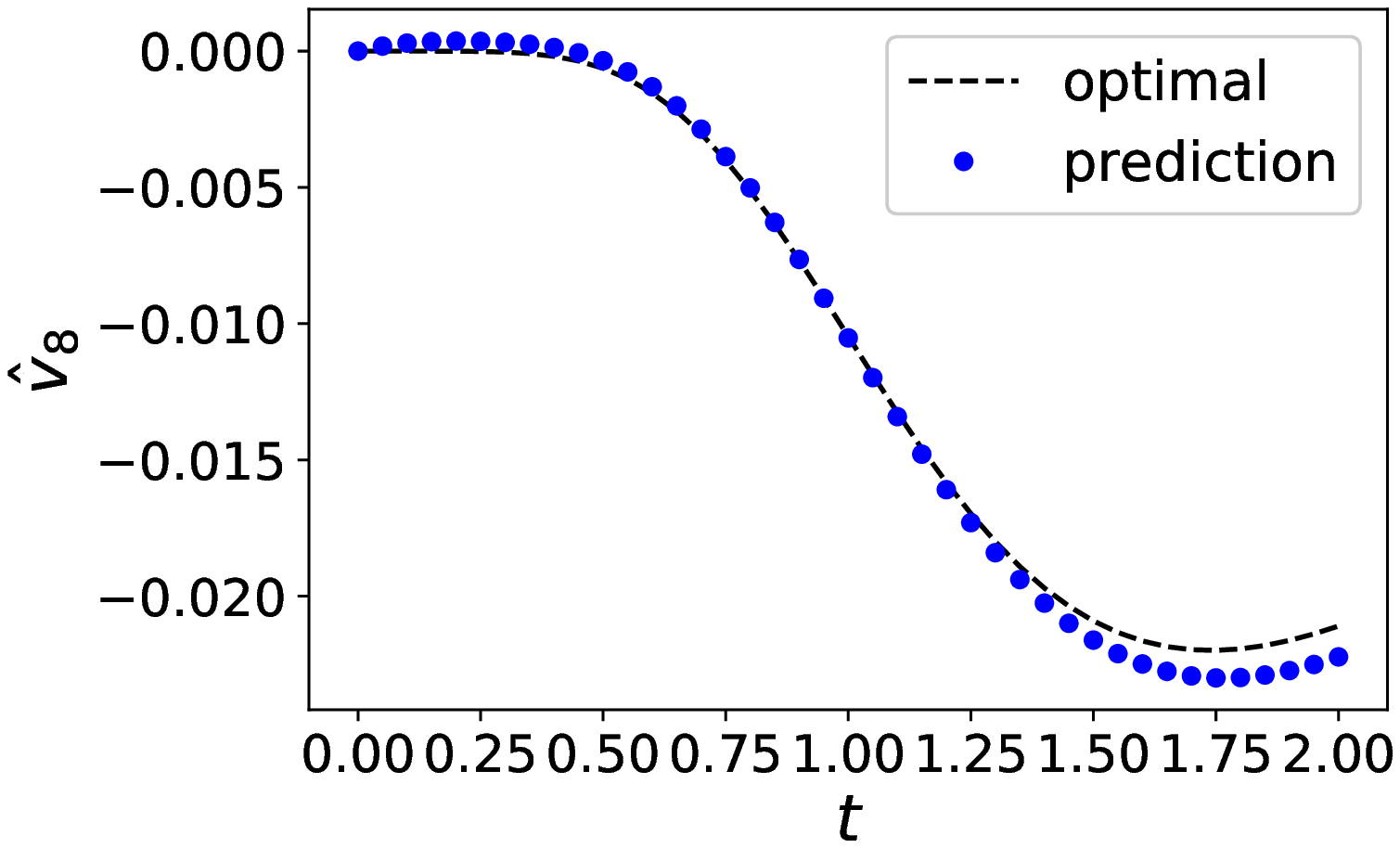}}
	{\includegraphics[width=0.325\textwidth]{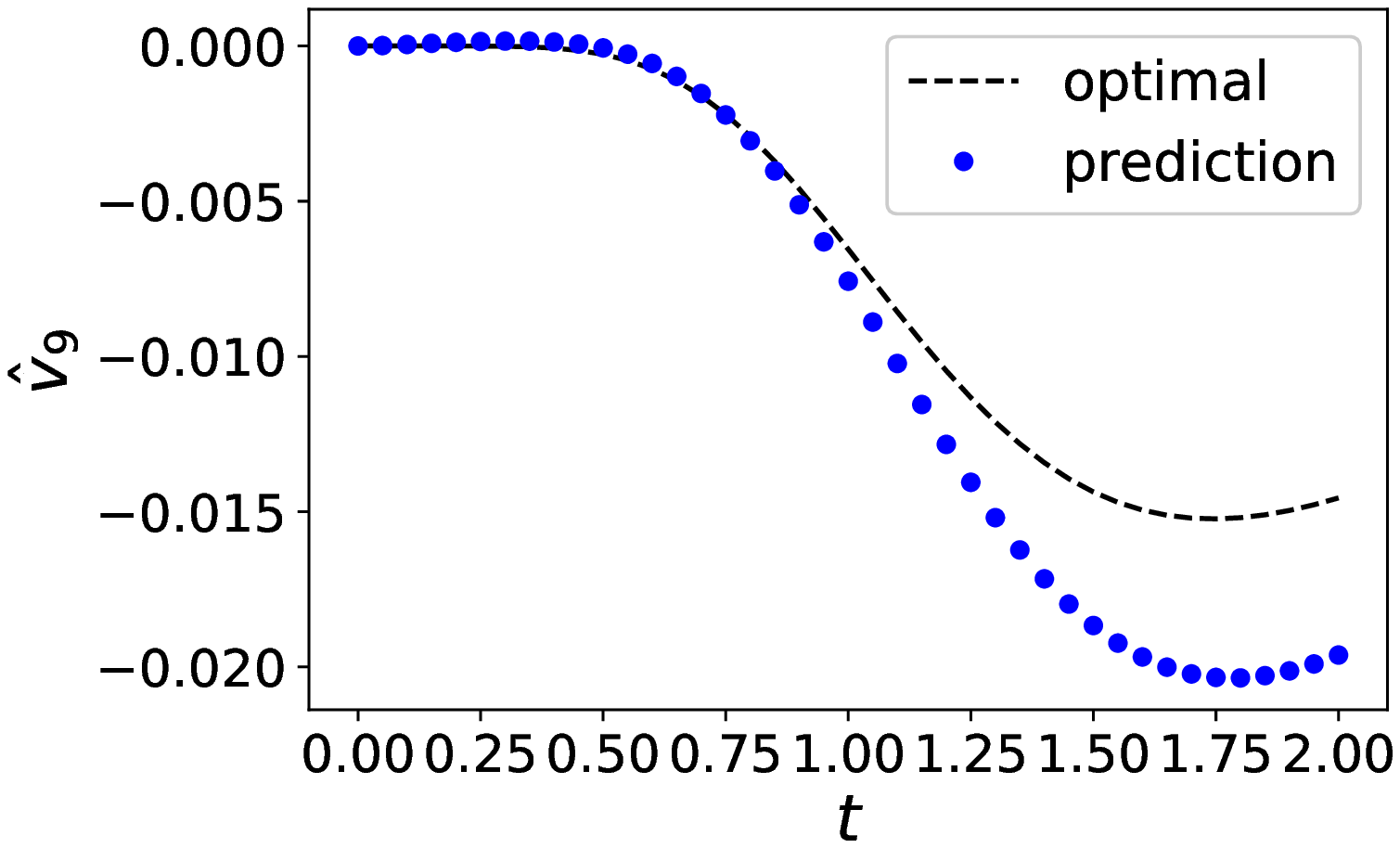}}
	\caption{\small
		Example 3: Evolution of the expansion coefficients for the learned model solution and the projection of the true solution.
	}\label{fig:ex3b_coef}
\end{figure}

We then consider a smaller viscosity $\sigma = 0.1$. The approximation
space is chosen to be relatively larger as $\mathbb V_n = {\rm span}\{  \sin(jx), 1\le j \le 9 \}$ with $n=9$. 
The time lag $\Delta$ is taken as 0.05. 
The domain $D$ in the modal space is taken as $[-1.5,  1.5]\times[-0.5,  0.5]\times[-0.2,  0.2]^2 \times[-0.1,  0.1]^2 \times [-0.05,0.05]^2 \times [-0.02,0.02] 
$, from which we sample $500,000$ training data.  
In this example, we use the four-block ResNet method ($K=4$) with each
block containing 3 hidden layers of equal width of 30 neurons. 
The network model is trained for up to  2,000 epochs, and training loss history is shown in Fig.~\ref{fig:ex3b_loss}. 
Then we validate the  
trained model for the initial condition 
$$u_0(x)= -\sin (x).$$ 
%by using \eqref{model-initial}--\eqref{model-final}. 
In Fig.~\ref{fig:ex3b_solu}, we present the prediction results
generated by the trained network, for time up to $t=2$. One can see
that the predicted solutions are very close to the exact ones,
This can  also be seen in the relative error of the  prediction from
Fig.~\ref{fig:ex3b_error}.
We note that at time $t=2$ the exact solution develops a relatively
sharp (albeit still smooth) transition layer at this relatively low
visocity. The numerical solution starts to exhibit Gibbs' type small
oscillations. This is a rather common feature for global type
approximation and is not unexpected.
Comparison between the learned and optimal expansion coefficients
%against the optimal orthogonal projection coefficients 
is also shown
in Fig.~\ref{fig:ex3b_coef}.

\subsection{Example 4: Inviscid Burgers' Equation}

We now consider the inviscid Burgers' equation with Dirichlet boundary condition: 
\begin{equation}
\label{eq:example4}
\begin{cases}
u_t + \left( \frac{u^2}2 \right)_x = 0,\quad (x,t) \in (-\pi,\pi) \times \mathbb R^+, \\
u(-\pi,t)=u(\pi,t)=0,\quad t \in \mathbb R^+.
\end{cases}
\end{equation}
This represents a challenging problem, as the nonlinear
hyperbolic nature of this equation can produce shocks over time even if the initial solution is smooth. 
%This feature makes the prediction and data-driven modeling for this problem challenging. 

\begin{figure}[htbp]
	\centering
	{\includegraphics[width=0.6\textwidth]{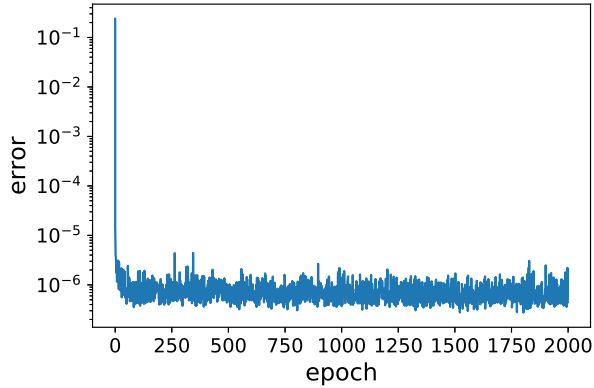}}
	\caption{\small
		Example 4: Training loss history.
	}\label{fig:ex4_loss}
\end{figure}

 The approximation space is chosen as $\mathbb V_n = {\rm span} \{
 \sin(jx), 1\le j \le 9\}$ with $n=9$. 
The time lag $\Delta$ is taken as $0.05$. 
The domain $D$ in the modal space is taken as $[-1.1,  1.1]\times[-0.5,  0.5] \times [-0.3 \times 0.3]^7
$, from which we sample $1,000,000$ training data.
We remark that by sampling in this manner, all of our training data
are smooth.
In this example, we use the block ResNet method with $K=4$ blocks,
each of which contains 3 hidden layers of equal width of 30.
The network training is conducted for up to $2,000$ epochs, when it was
deemed satisfactory, as shown in the loss history in Fig.~\ref{fig:ex4_loss}. 
For validation, we conduct system prediction using
 initial condition 
$$u_0(x)= -\sin (x),$$ 
for time up to $t=2$.
Although the initial condition is smooth,
the exact solution will start to develop shock 
at $t=1$.

\begin{figure}[htbp]
	\centering
	{\includegraphics[width=0.48\textwidth]{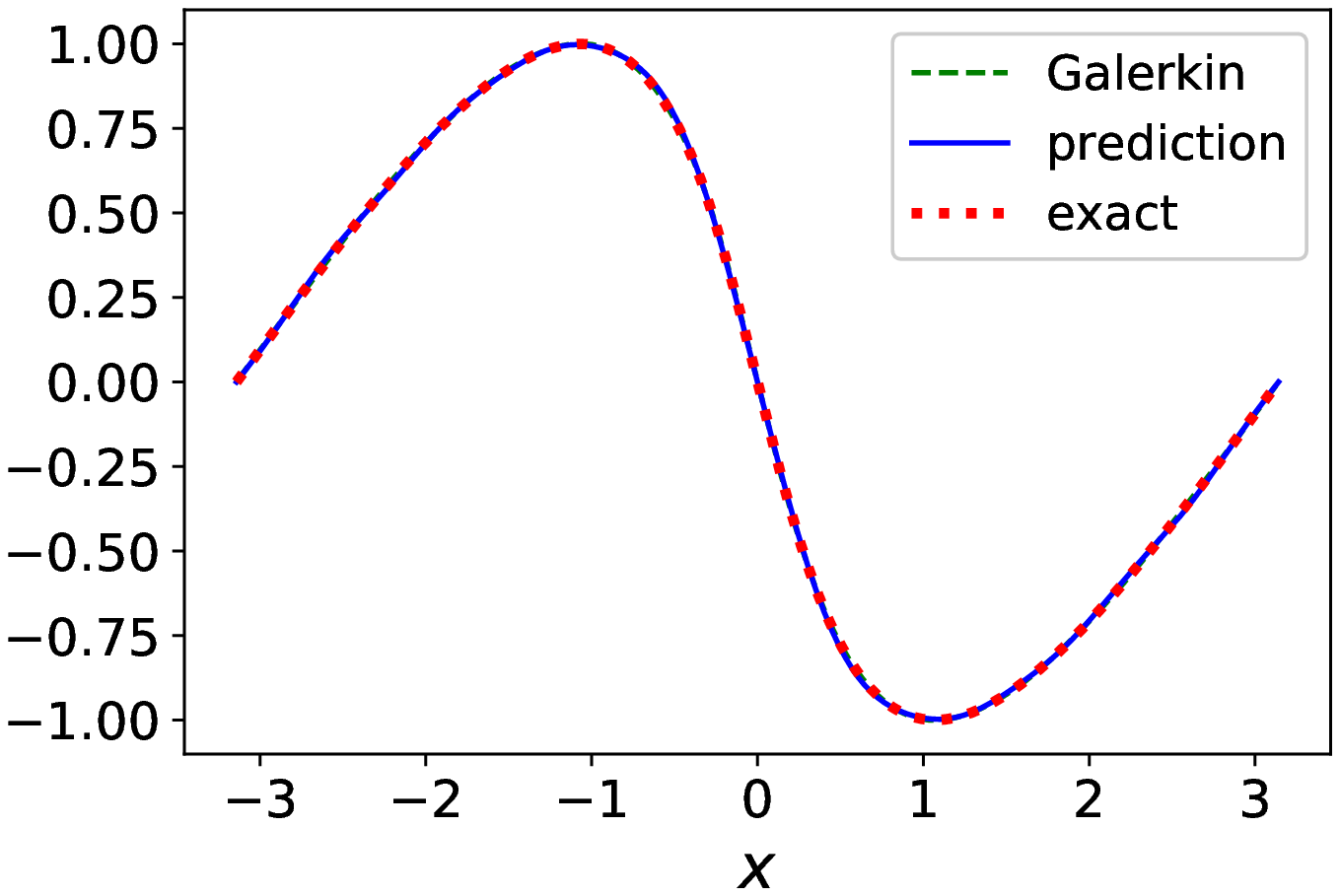}}
	{\includegraphics[width=0.48\textwidth]{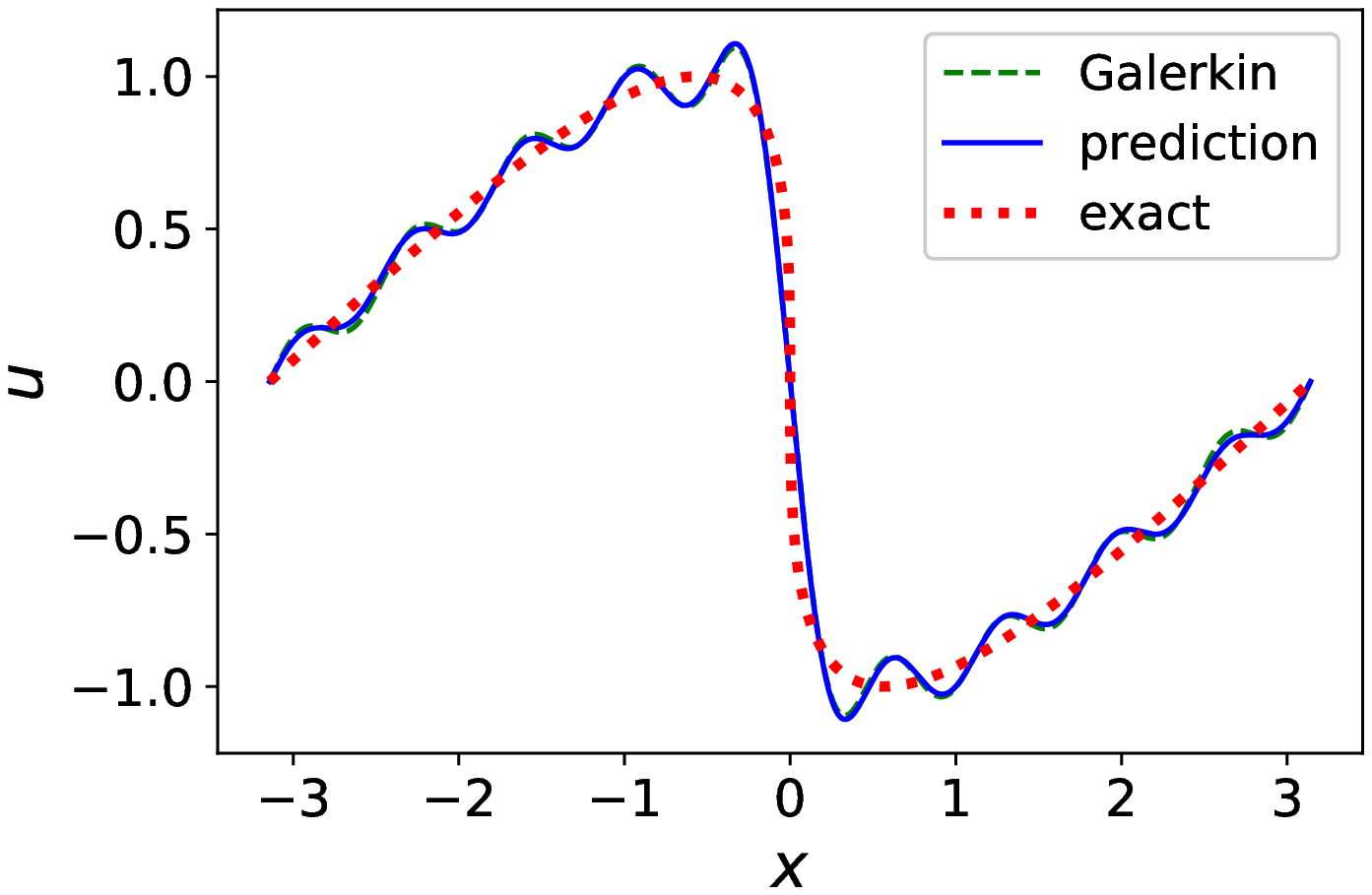}}
	{\includegraphics[width=0.48\textwidth]{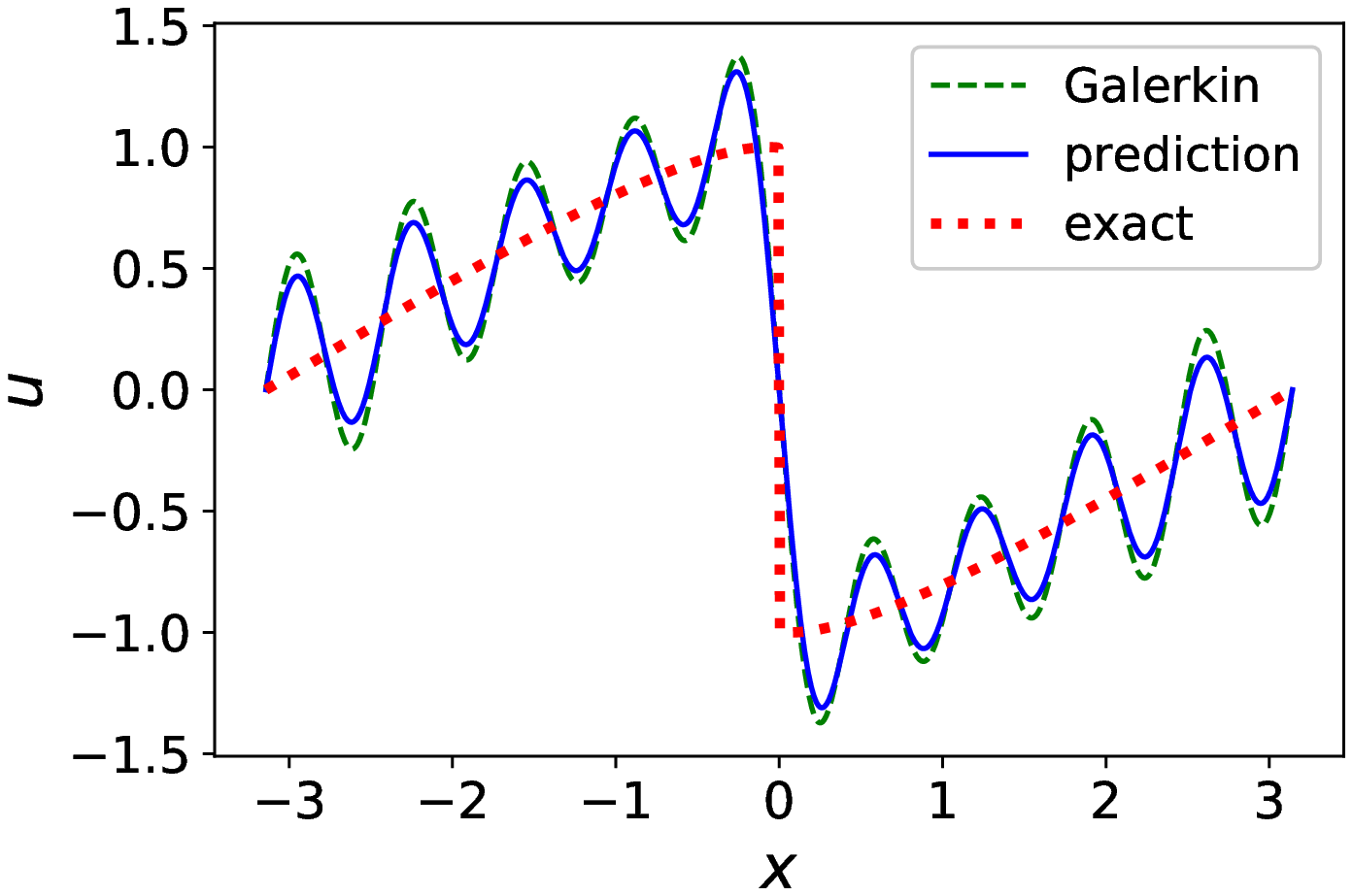}}
	{\includegraphics[width=0.48\textwidth]{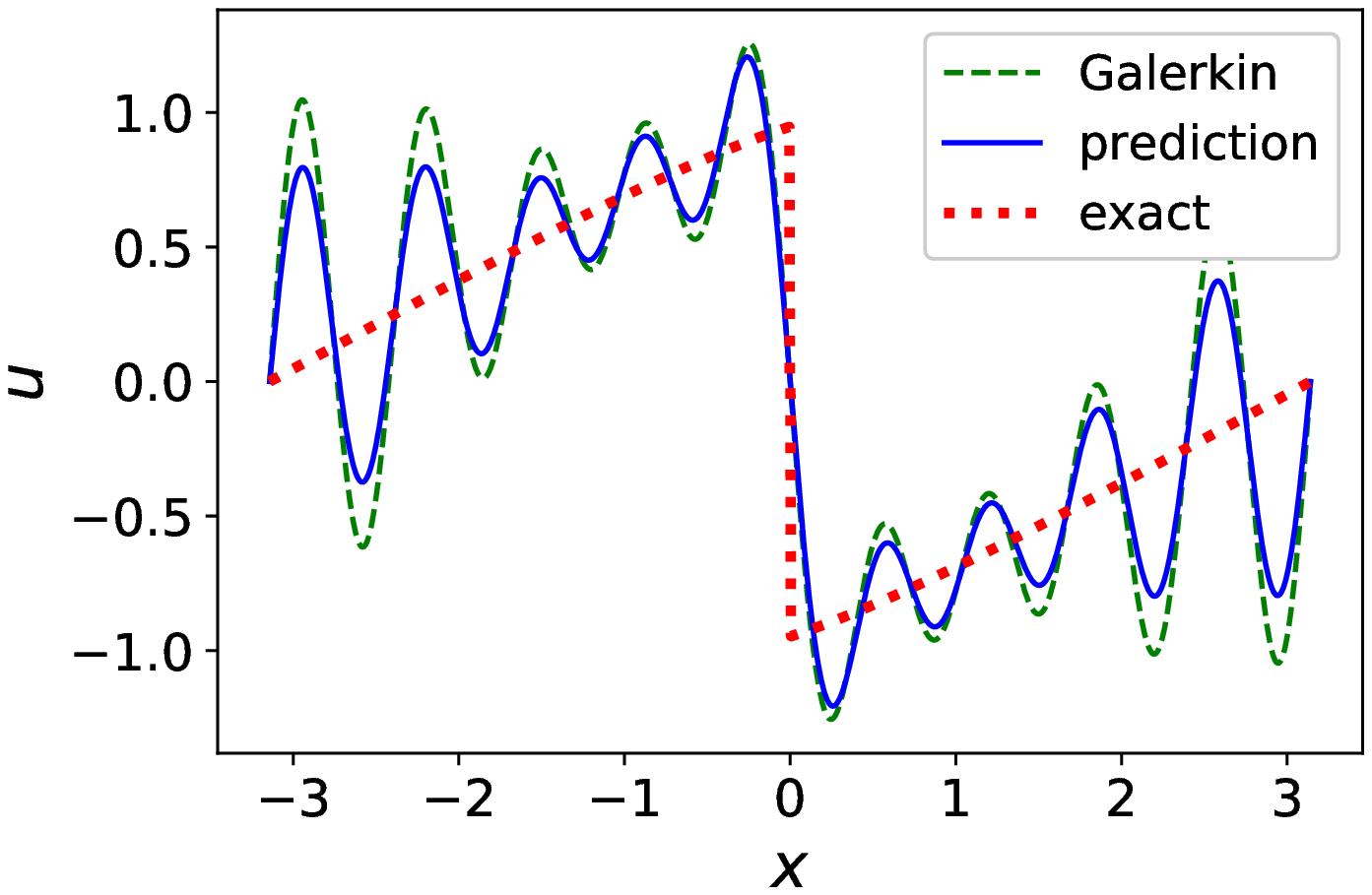}}	
	\caption{\small
		Example 4: Comparison of the true solution, the learned model solution and the solution by Galerkin method at different time. Top-left: $t=0.5$; top-right: $t=1$; bottom-left: $t=1.5$; bottom-right: $t=2$.  
	}\label{fig:ex4_solu}
\end{figure}

\begin{figure}[htbp]
	\centering
	{\includegraphics[width=0.325\textwidth]{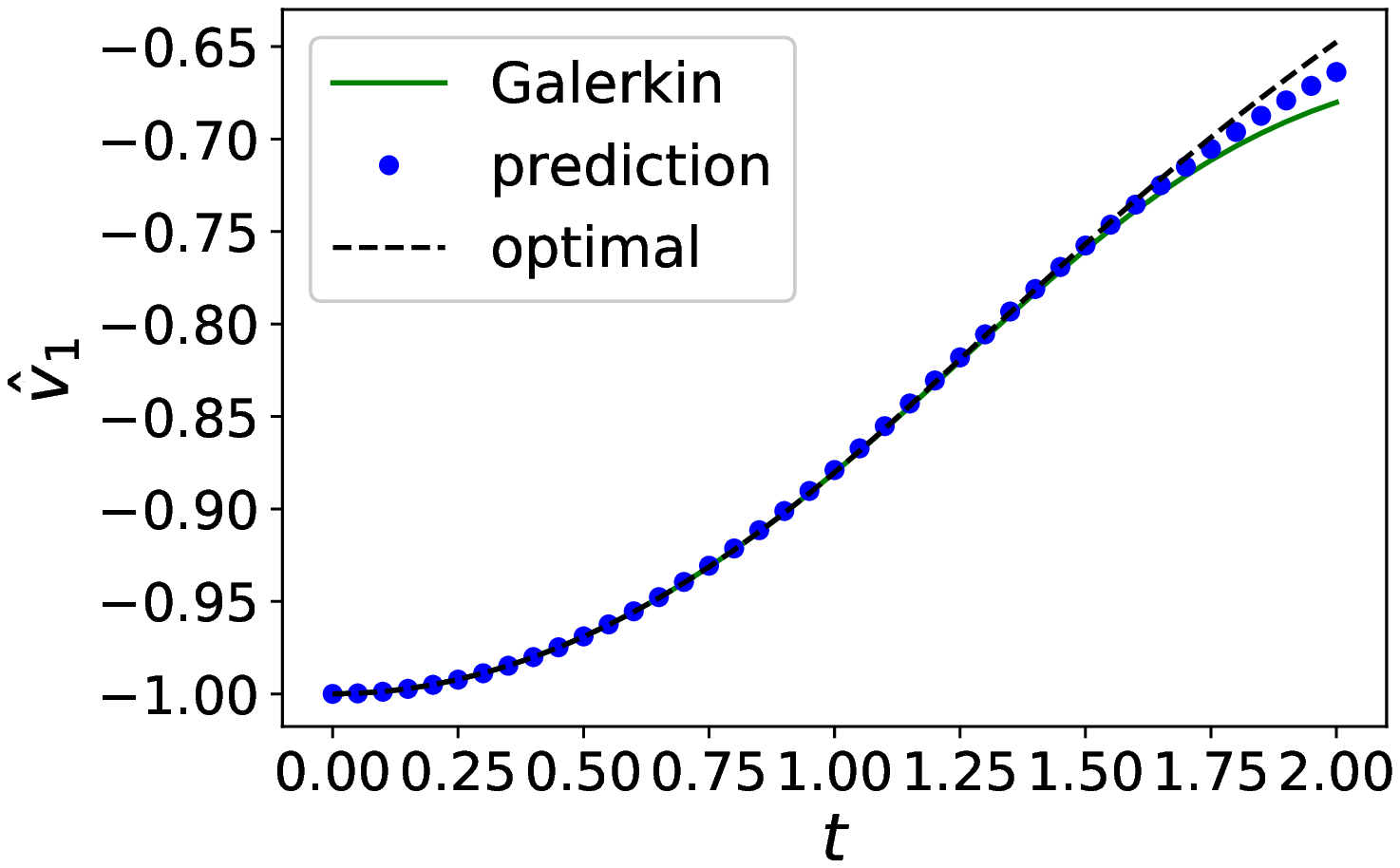}}
	{\includegraphics[width=0.325\textwidth]{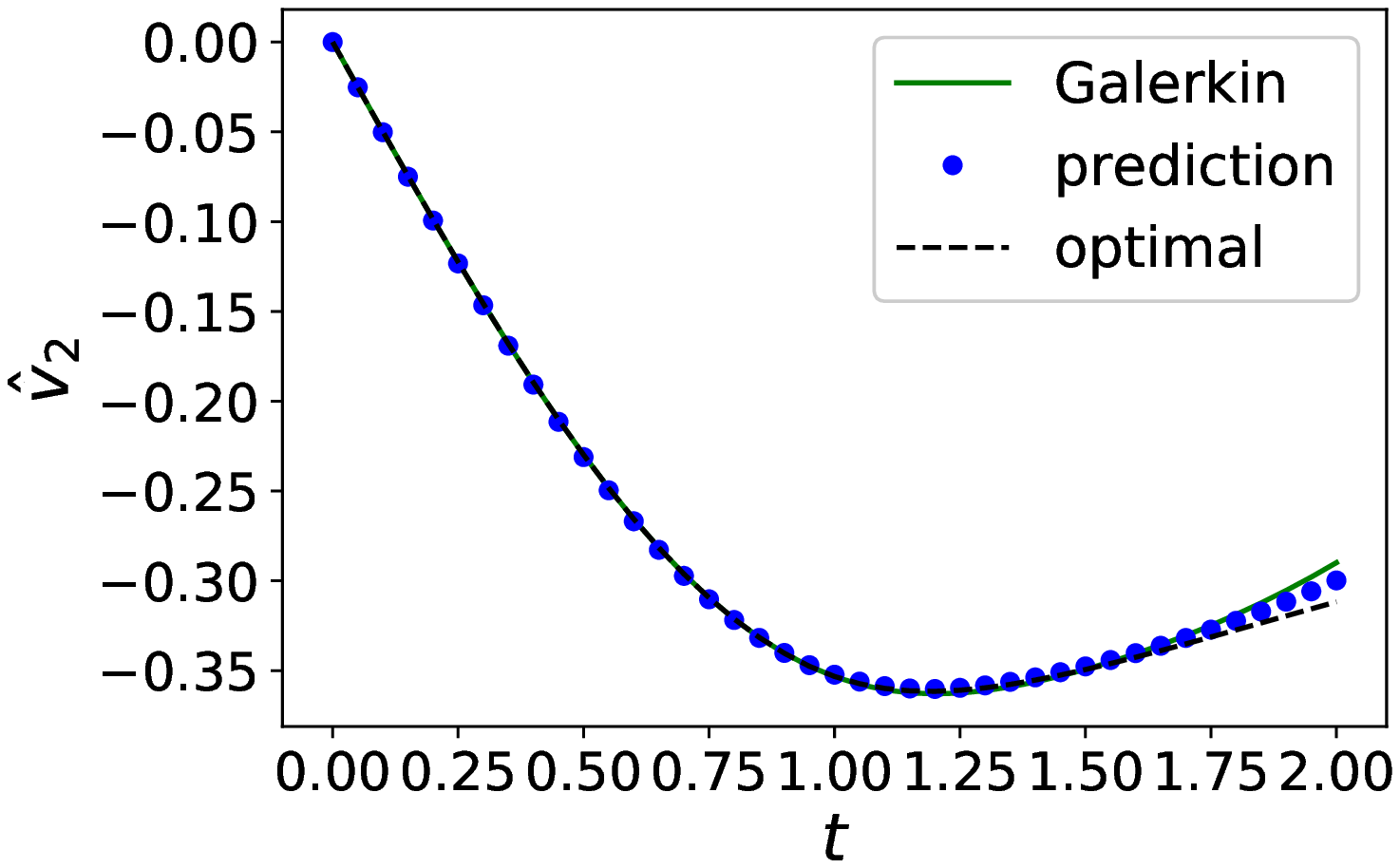}}
	{\includegraphics[width=0.325\textwidth]{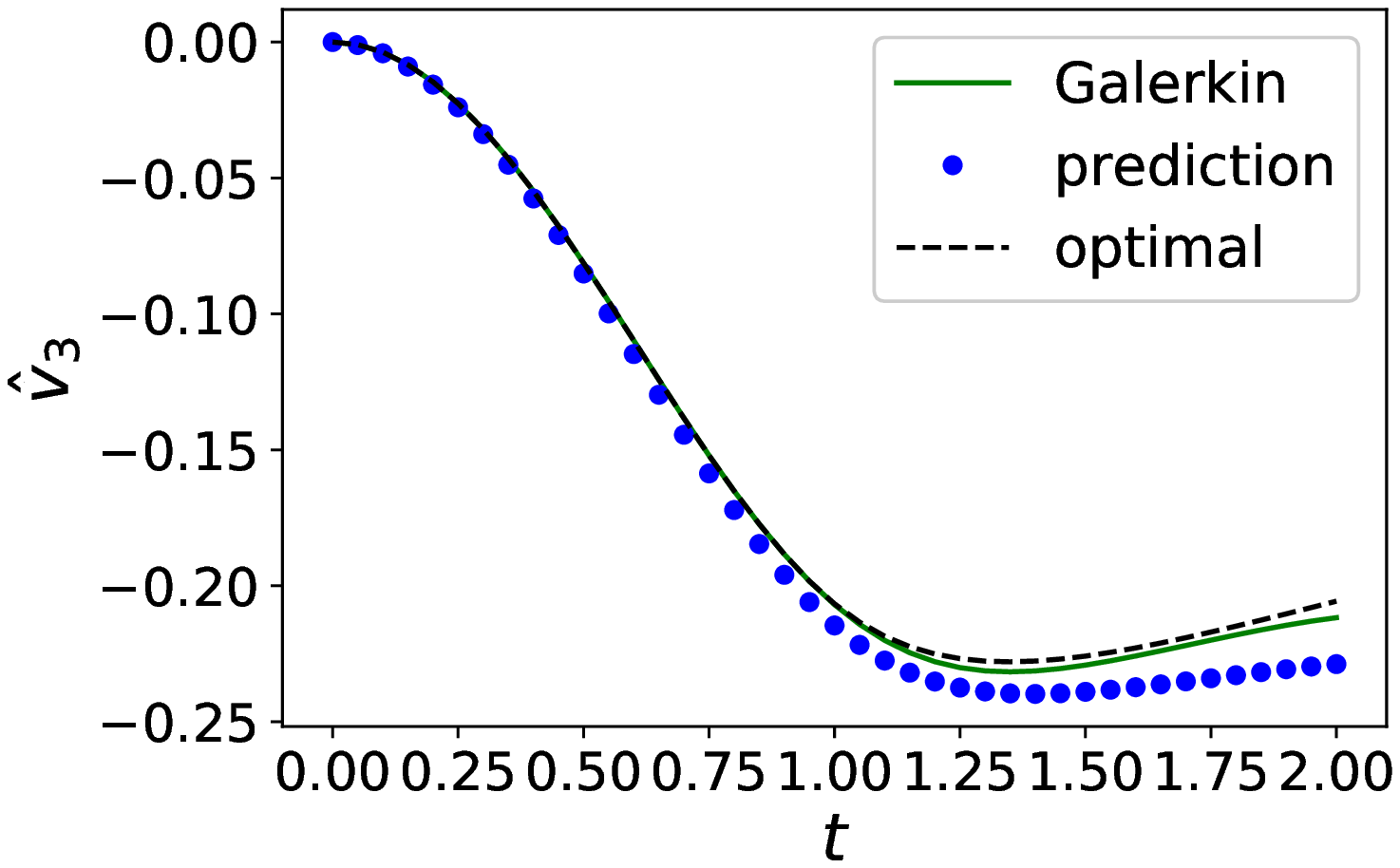}}
	{\includegraphics[width=0.325\textwidth]{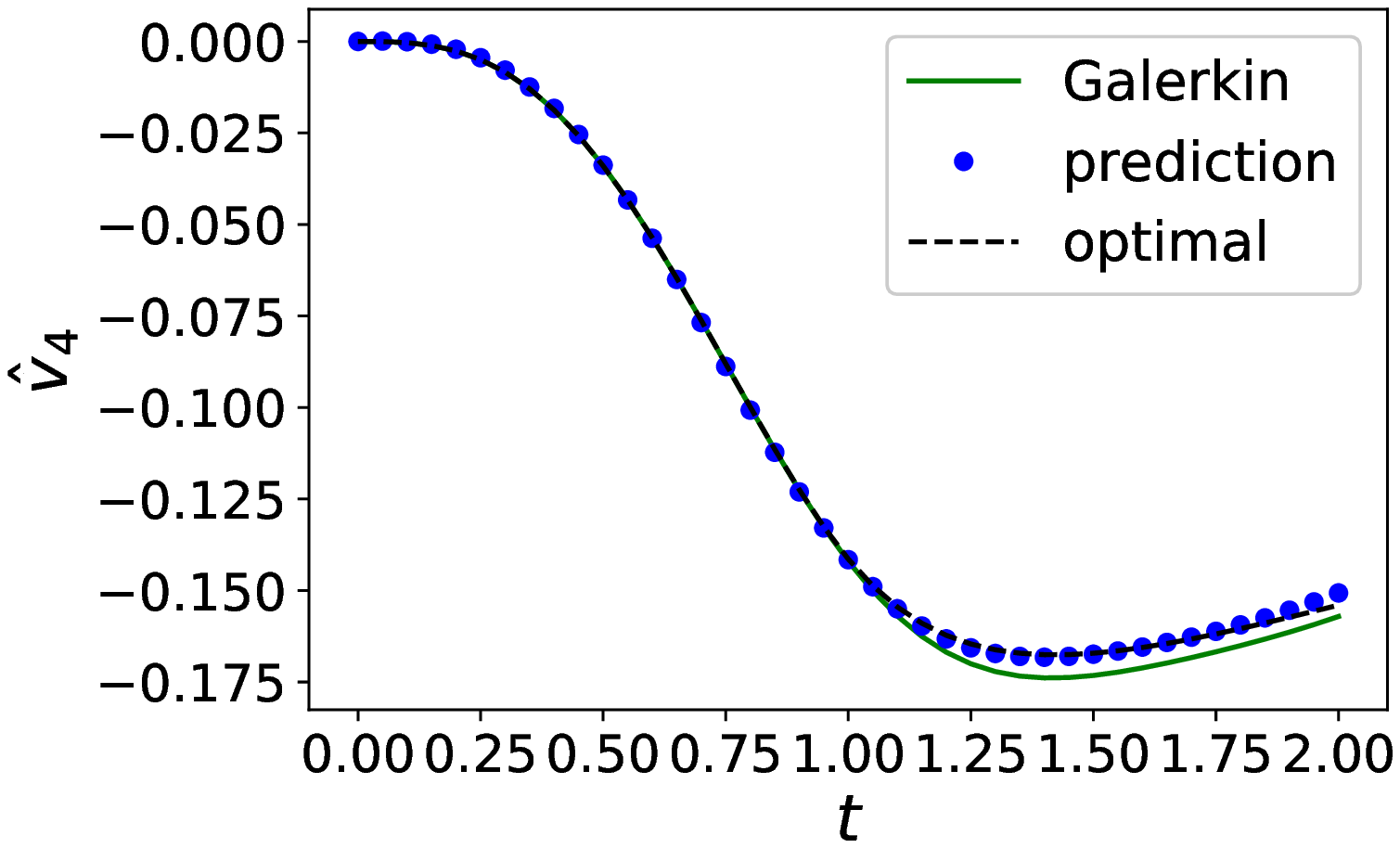}}
	{\includegraphics[width=0.325\textwidth]{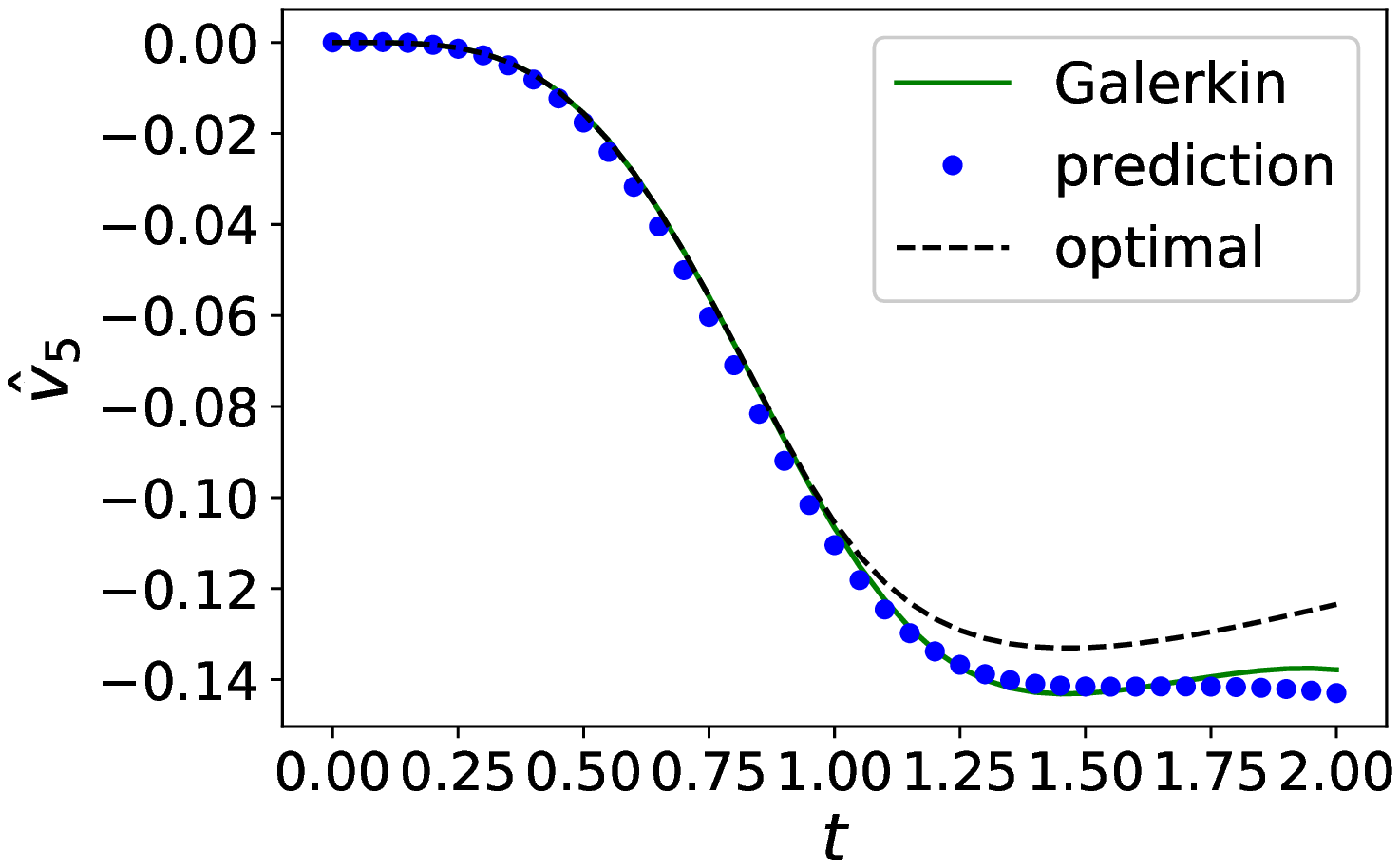}}
	{\includegraphics[width=0.325\textwidth]{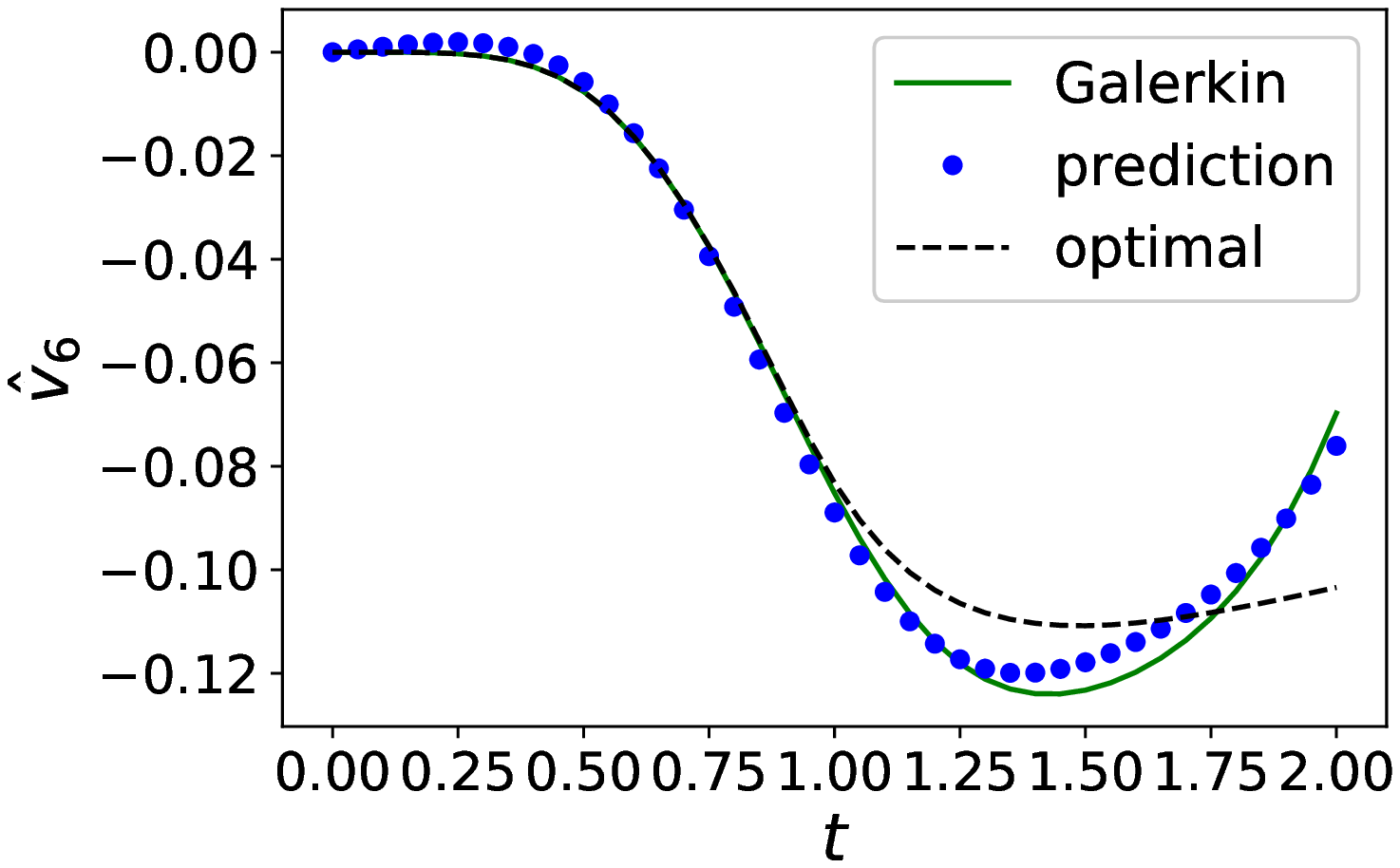}}
	{\includegraphics[width=0.325\textwidth]{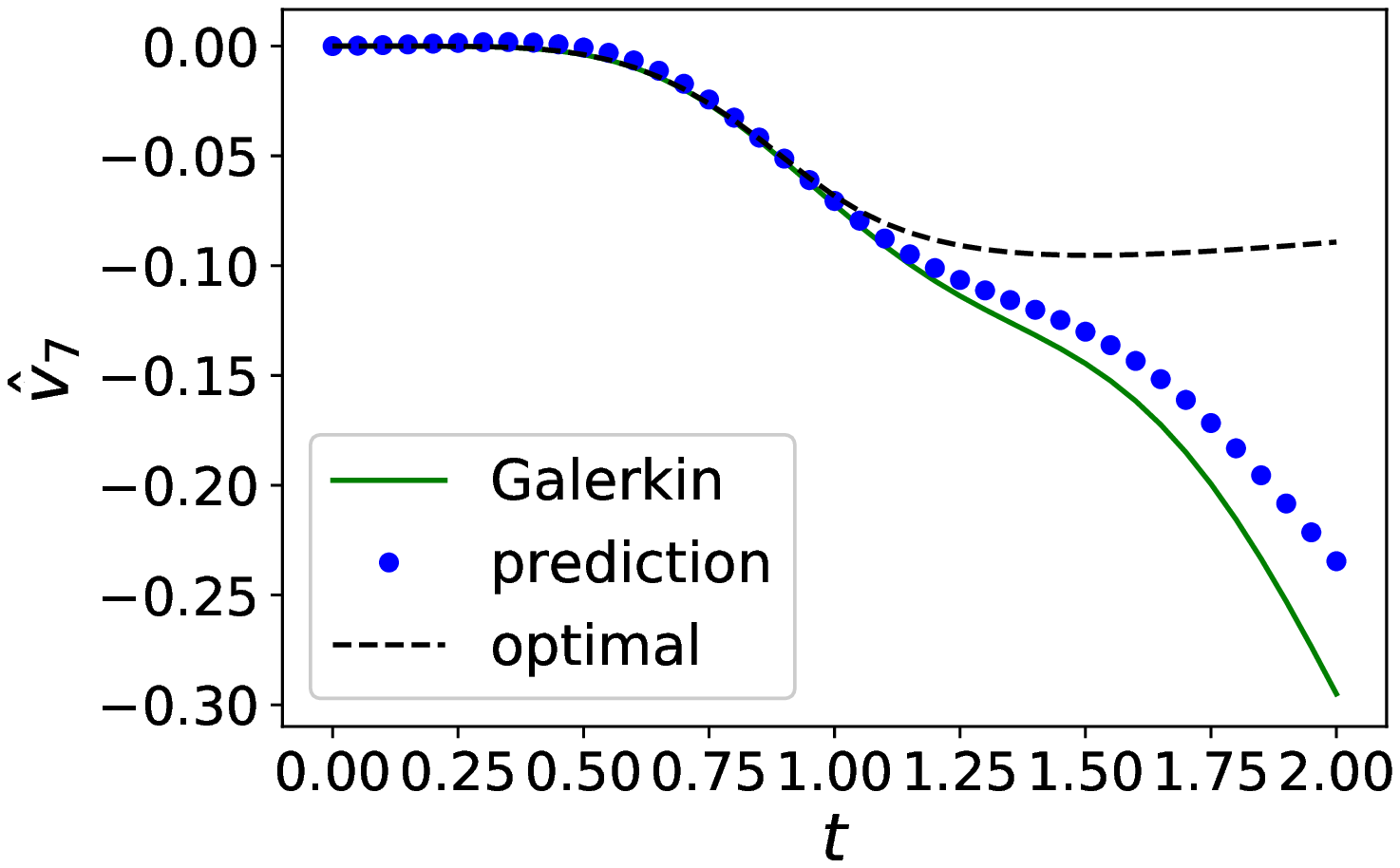}}
	{\includegraphics[width=0.325\textwidth]{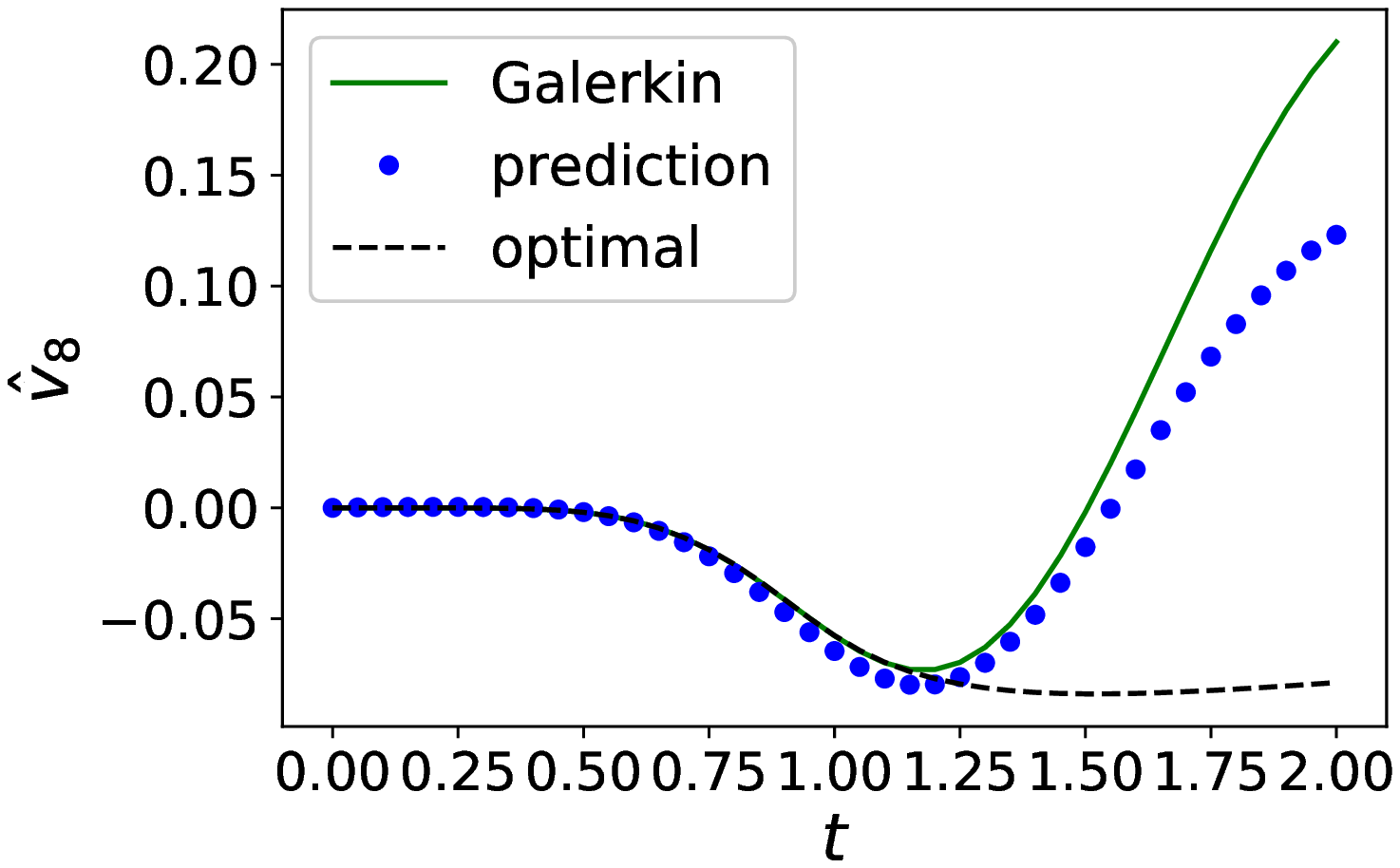}}
	{\includegraphics[width=0.325\textwidth]{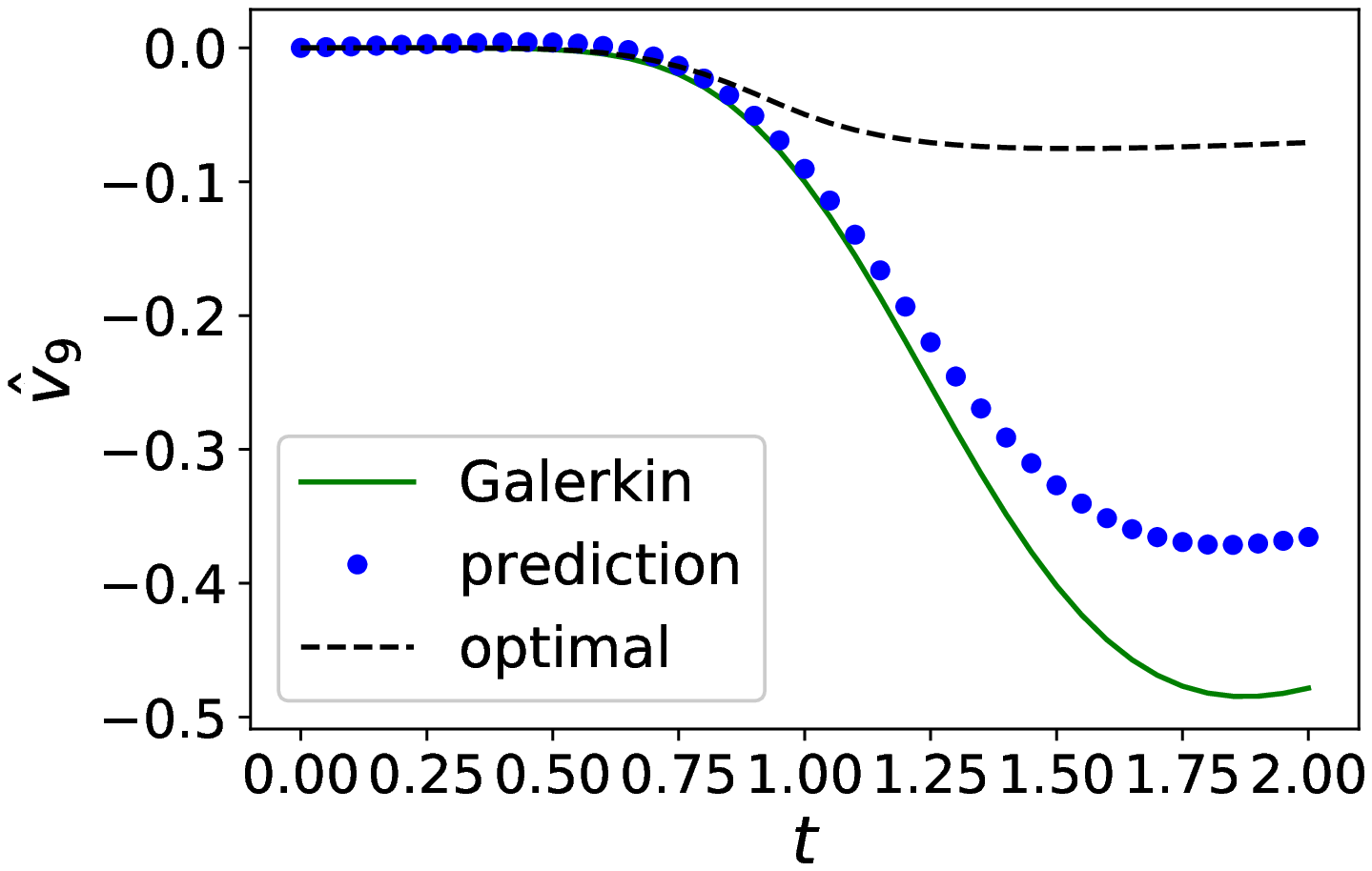}}
	\caption{\small
		Example 4: Evolution of the expansion coefficients for the learned model solution and the projection of the true solution.
	}\label{fig:ex4_coef}
\end{figure}

The prediction results generated by our network model are shown in Fig.~\ref{fig:ex4_solu}, 
along with the exact solution, and the Galerkin solution of
the Burgers' equation using the same linear space $\mathbb V_n$.
Due to the discontinuity in the solution, Gibbs type oscillations  
appear in the predicted solutions. This is not unexpected, as the best
representation of a discontinuous function using the linear subspace
$\mathbb V_n$ will naturally produce oscillations, unless special
treatment such as filtering is utilized (which is not pursued in this
work).
The Galerkin solution of the equation \eqref{eq:example4} exhibits the
similar Gibbs' oscillations for precisely the same reason.
We remark that it can be seen that our neural network prediction is visibly
better than the Galerkin solution. While the network prediction is
entirely data driven and does not require knowledge of the equation,
the Galerkin solution is attainable only after knowing the precise
form of the governing equation.
%However, the location of the shock is approximately captured. 
%We also see that our network model produces more accurate approximation than the Galerkin method. 
%This is impressive, because the network model is purely data-driven
%without using any knowledge of the governing equation, while the
%Galerkin method has taken the form of the equation into account.  
In Fig.~\ref{fig:ex4_coef}, we plot the evolution of the expansion
coefficients in the modal space obtained by the neural network model
prediction, Galerkin solution of the Burgers' equation, and
orthogonal projection of the exact solution (denoted as ``optimal''), the last of which serves
as the reference solution. It is clearly seen that the neural network
model produces more accurate results than the Galerkin solver. The
accuracy improvement is especially visible at higher modes such as
$\hat v_7$, $\hat v_8$, and $\hat v_9$. The cause of the accuracy
improvement over Galerkin method will be pursued in a separate work.

\begin{figure}[htbp]
	\centering
%	{\includegraphics[width=0.48\textwidth]{Figure/Example4/ABSerror.eps}}
	{\includegraphics[width=0.6\textwidth]{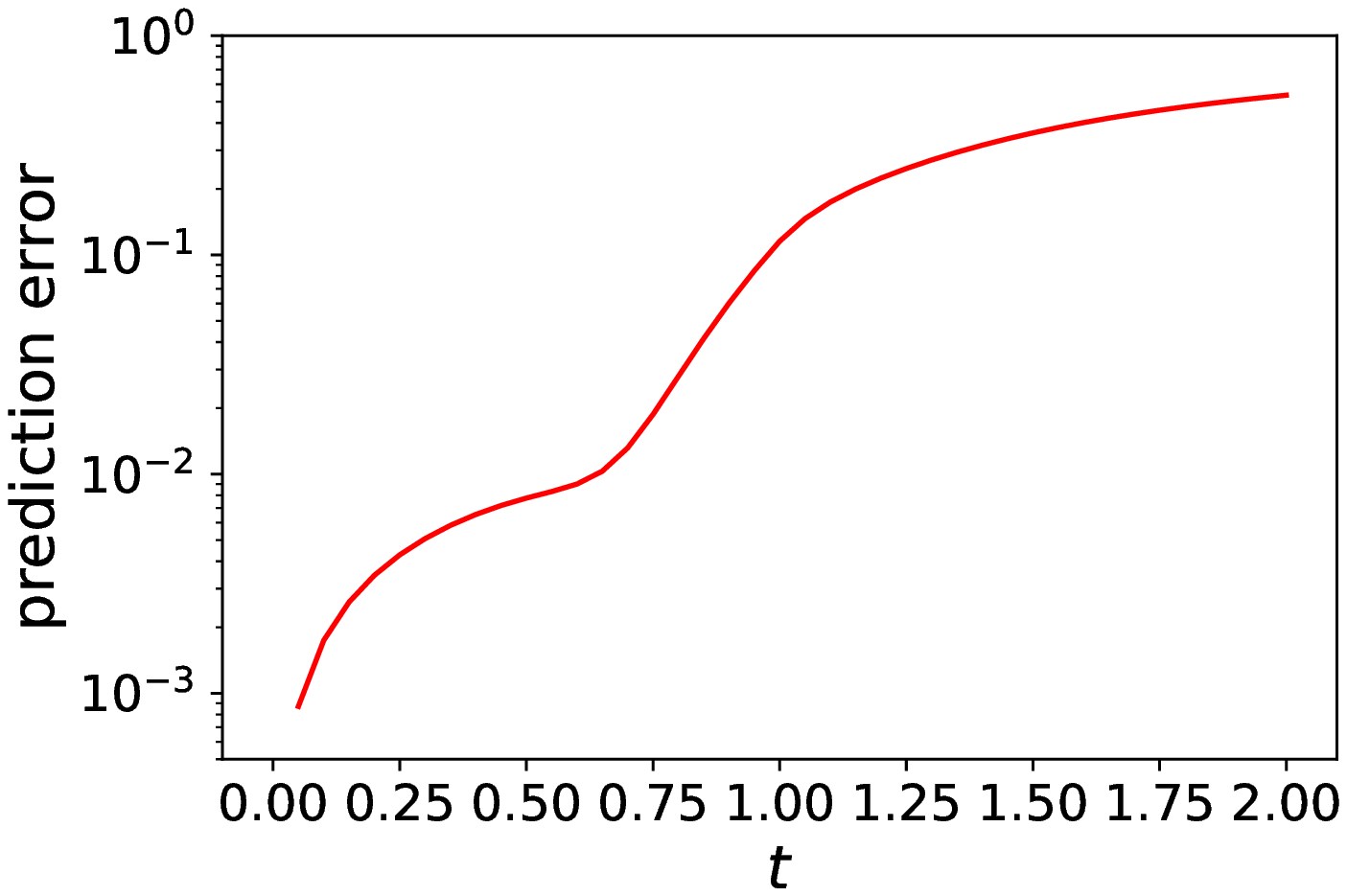}}
	\caption{\small
		Example 4: The evolution of the relative error in the
                prediction in $l^2$-norm.
                %Left: absolute error; right: relative error.  
	}\label{fig:ex4_error}
\end{figure}

\subsection{Example 5: Two-Dimensional Convection-Diffusion Equation}

In our last example, we consider a two-dimensional
convection-diffusion equation to demonstrate the applicability of the
proposed algorithm for multiple dimensions. The equation set up is as follows.
\begin{equation}
\label{eq:example5}
\begin{cases}
u_t + \alpha_1 u_x + \alpha_2 u_y = \sigma_1 u_{xx} + \sigma_2 u_{yy},\quad (x,y,t) \in (-\pi,\pi)^2 \in \mathbb R^+ \\
u(-\pi,y,t)=u(\pi,y,t),~u_x(-\pi,y,t)=u_x(\pi,y,t),~~(y,t) \in (-\pi,\pi) \times \mathbb R^+,
\\
u(x,-\pi,t)=u(x,\pi,t),~u_y(x,-\pi,t)=u_y(x,\pi,t),~~(x,t) \in (-\pi,\pi) \times \mathbb R^+,
\end{cases}
\end{equation}
where the parameters are set as $\alpha_1 = 1 $, $\alpha_2=0.7$,
$\sigma_1=0.1$, and $\sigma_2 = 0.16$ in the test.

\begin{figure}[htbp]
	\centering
	{\includegraphics[width=0.6\textwidth]{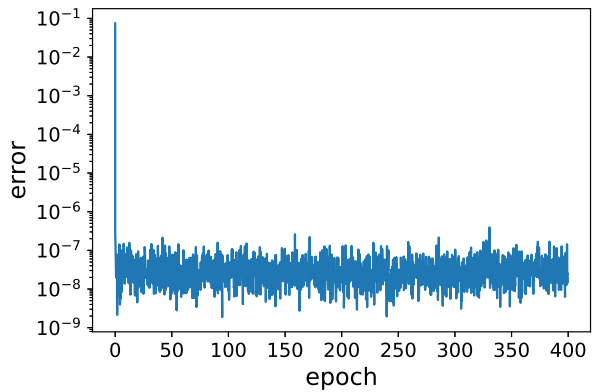}}
	\caption{\small
		Example 5: Training loss history.
	}\label{fig:ex5_loss}
\end{figure}

We chose the finite dimensional approximation space $\mathbb V_n$ as a
span by $n=25$ basis functions:
\begin{align*}
& \phi_1(x,y)=1, \quad \phi_2(x,y)=\cos(x), \quad \phi_3(x,y)=\sin(x), 
\\
& \phi_4(x,y)=\cos(2 x), \quad \phi_5(x,y)=\sin(2 x), \quad \phi_6(x,y)=\cos(3 x), \quad \phi_7(x,y)=\sin(3 x),
\\
& \phi_8(x,y)=\cos(y), \quad \phi_9(x,y)=\sin(y), \quad \phi_{10}(x,y)=\cos(2 y), \quad \phi_{11}(x,y)=\sin(2 y),
\\
& \phi_{12}(x,y)=\cos(3 y), \quad \phi_{13}(x,y)=\sin(3 y),  \quad \phi_{14}(x,y)= \cos(x)  \cos(y),  
\\
& \phi_{15}(x,y)= \cos(x)  \sin(y), \quad   \phi_{16}(x,y)= \sin(x)  \cos(y), \quad \phi_{17}(x,y)= \sin(x)  \sin(y),
\\
& \phi_{18}(x,y)= \cos(x)  \cos(2y),  \quad \phi_{19}(x,y)= \cos(x)  \sin(2y), \quad   \phi_{20}(x,y)= \sin(x)  \cos(2y),
\\
& \phi_{21}(x,y)= \sin(x)  \sin(2y), \quad    \phi_{22}(x,y)= \cos(2x)  \cos(y), \quad   \phi_{23}(x,y)= \cos(2x)  \sin(y), 
\\
&   \phi_{24}(x,y)= \sin(2x)  \cos(y), \quad \phi_{25}(x,y)= \sin(2x)  \sin(y).
\end{align*}
The time lag $\Delta$ is taken as 0.1 and
$1,000,000$ training data are generated.
The neural network model is a
five-block ResNet method ($K=5$) with each block containing 3 hidden
layers of equal width of 40 neurons.
Upon training the network model satisfactorily for 400 epochs (see
Fig.~\ref{fig:ex5_loss} for the training loss history), we validate the 
trained models by using the initial condition 
$$u_0(x)= \frac{2}{5} \exp \left( \frac{ \sin (x) - \cos(y) } 2 \right),$$ 
for time up to $t=3$.

\begin{figure}[htbp]
	\centering
	{\includegraphics[width=0.48\textwidth]{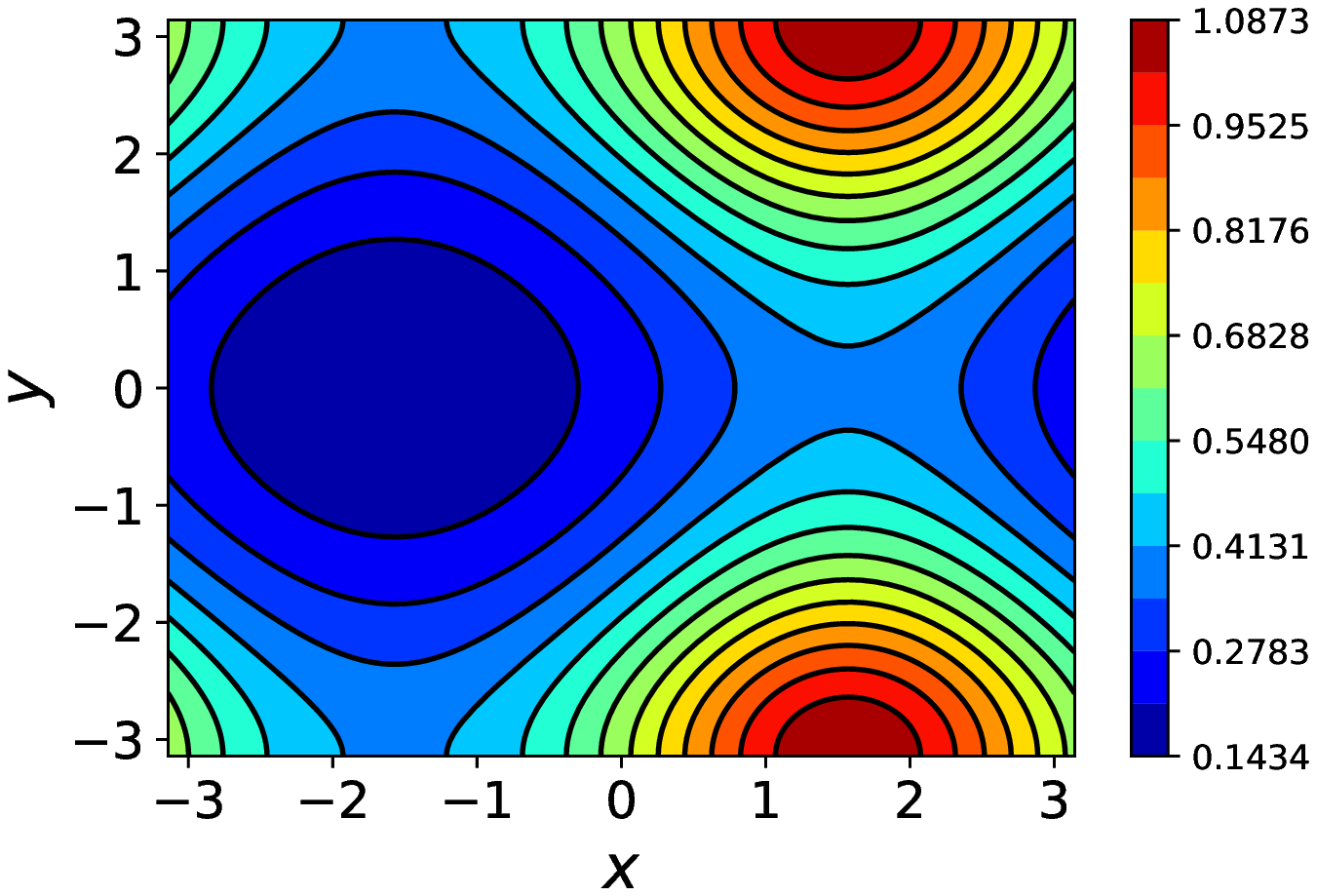}}
	{\includegraphics[width=0.48\textwidth]{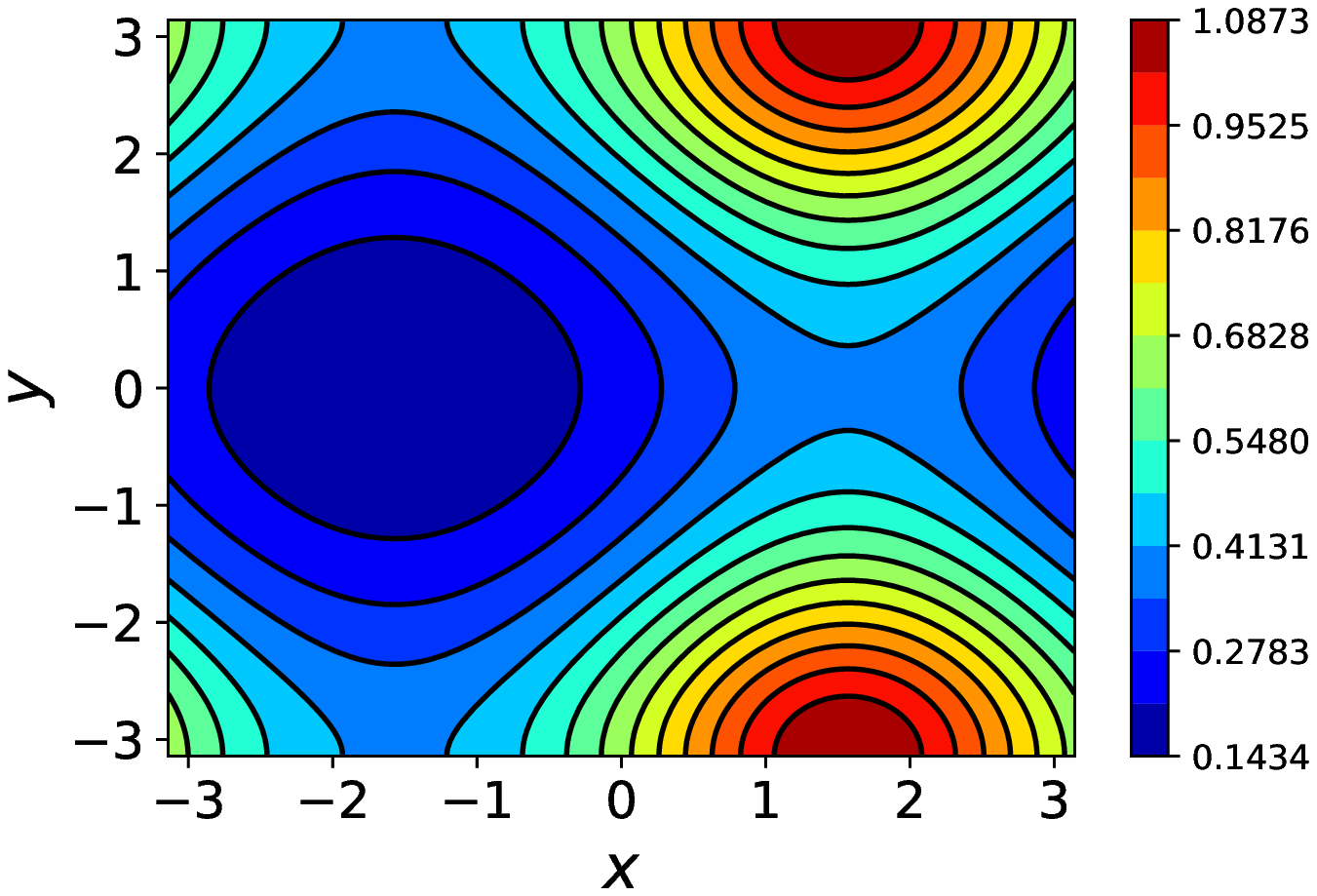}}
	{\includegraphics[width=0.48\textwidth]{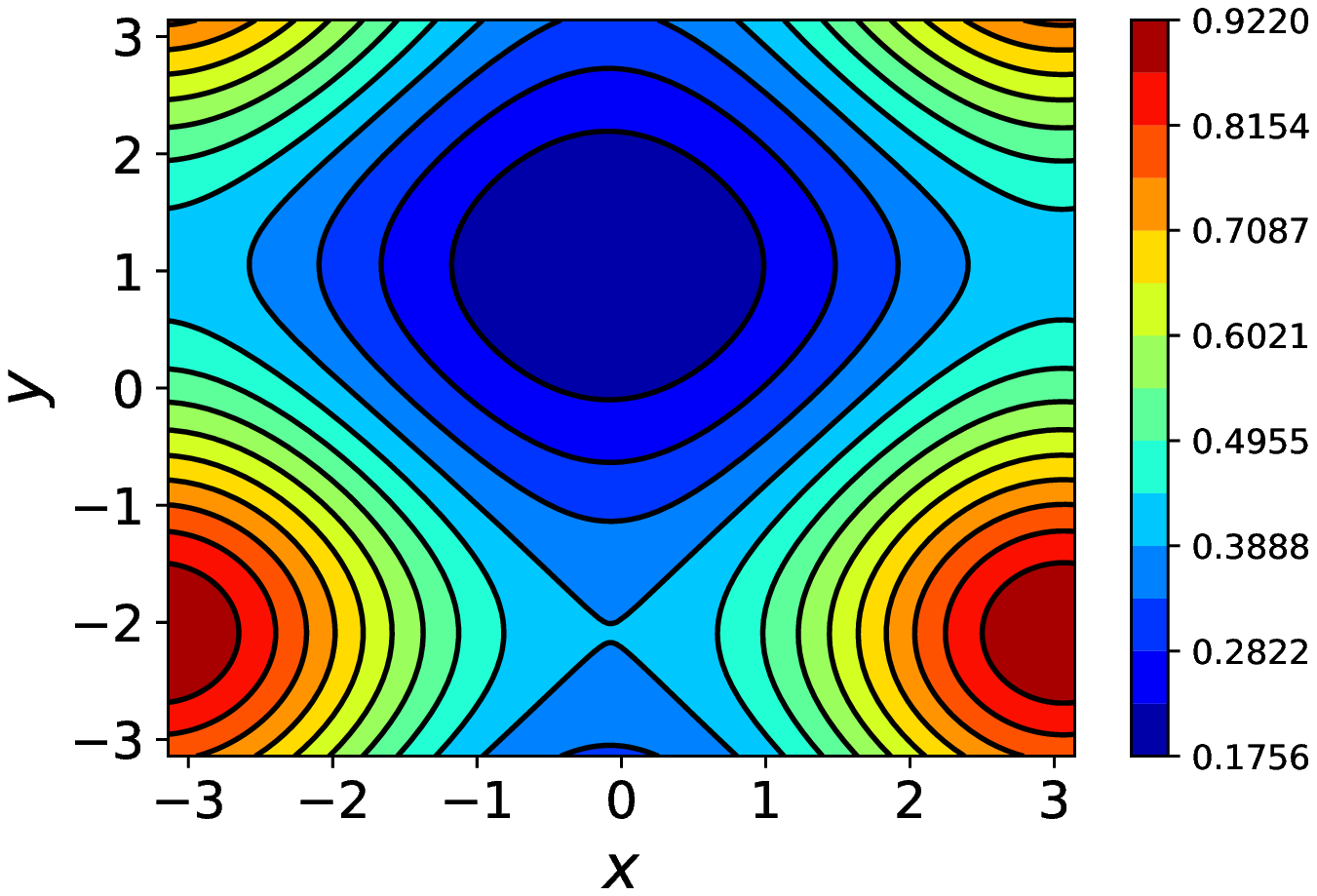}}
	{\includegraphics[width=0.48\textwidth]{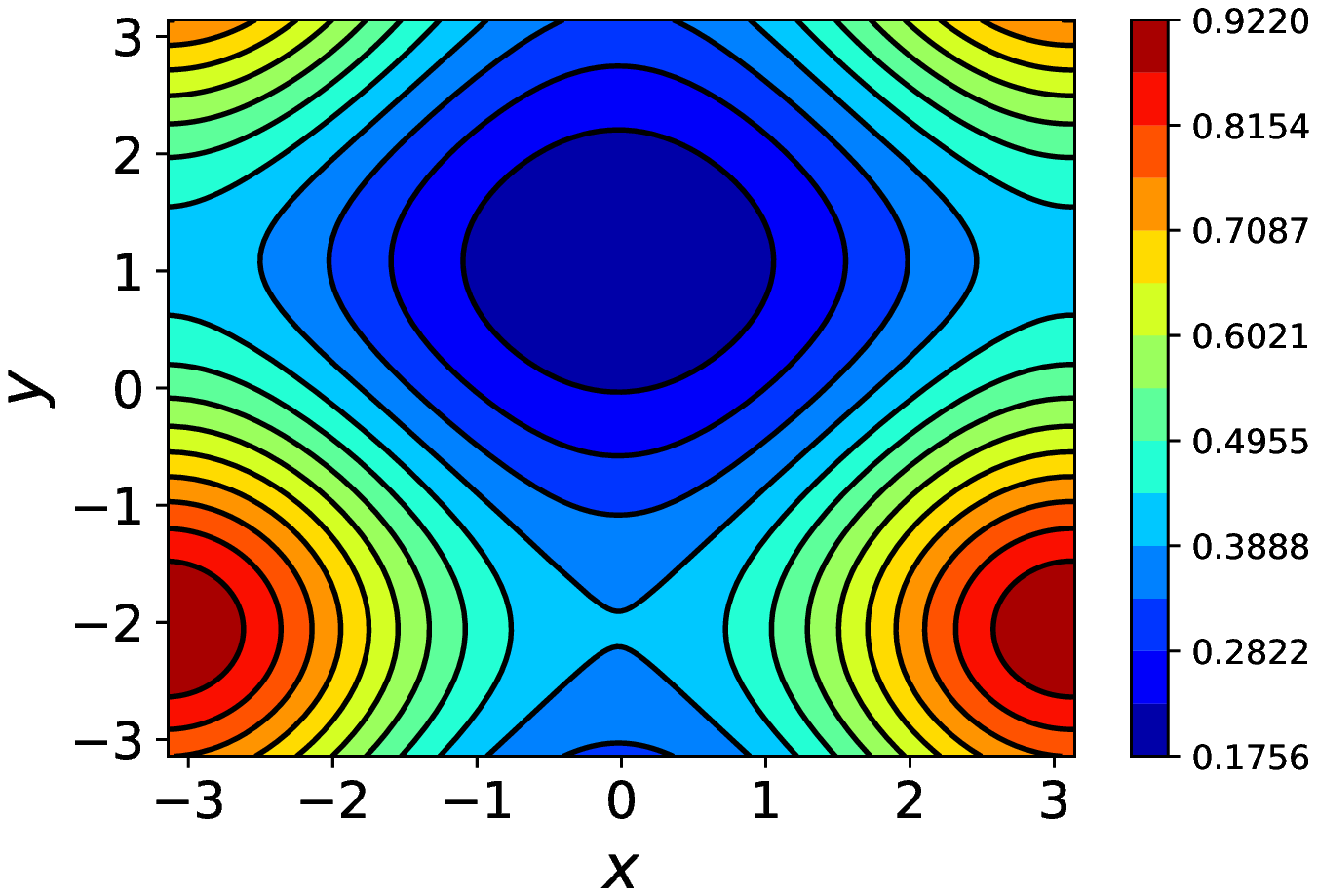}}
	{\includegraphics[width=0.48\textwidth]{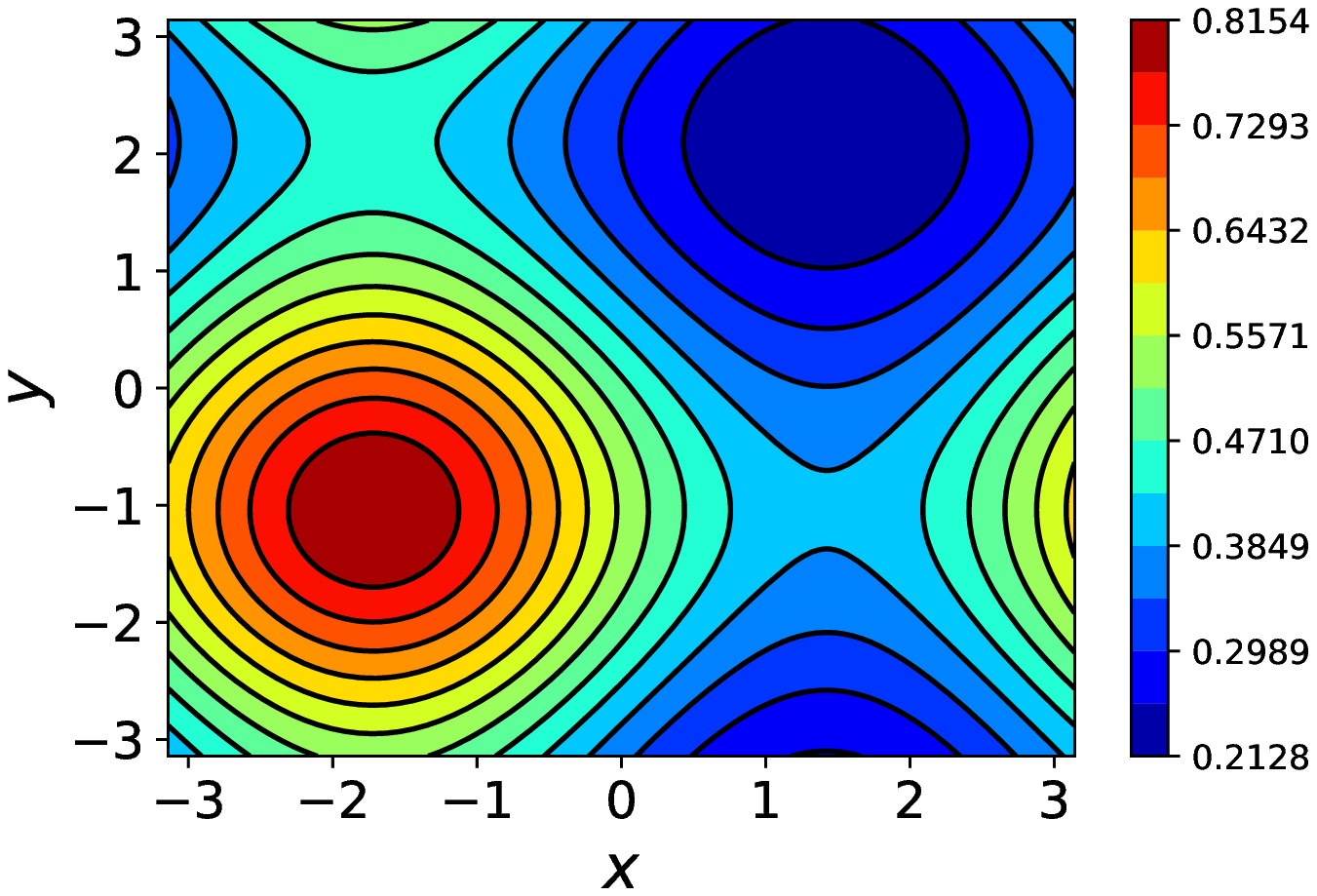}}
	{\includegraphics[width=0.48\textwidth]{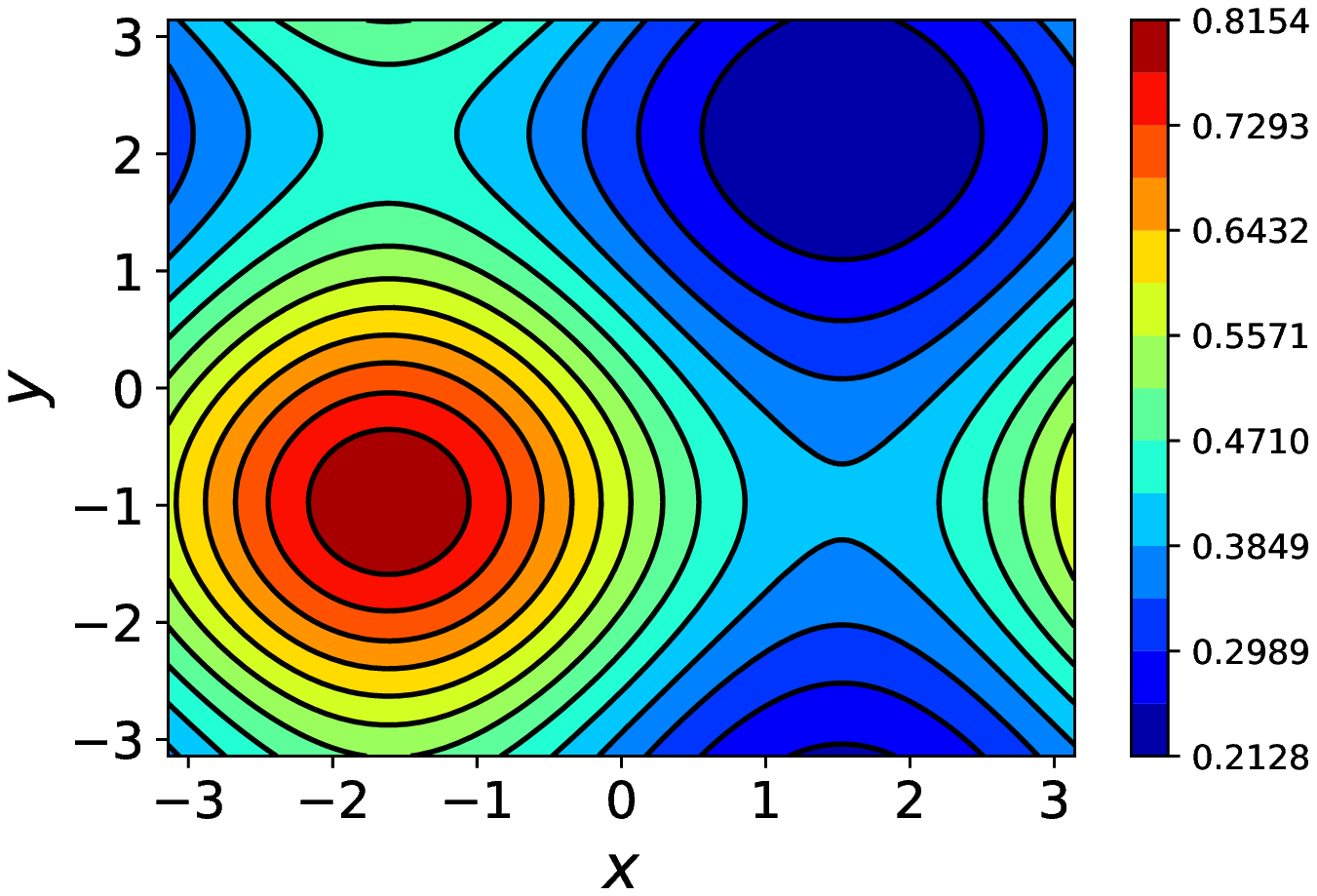}}
	\caption{\small
		Example 5: Contour plots of the solutions at different time. From top to bottom $t=0,$ $t=1.5$ and $t=3$. 
		Left: learned model solution; right: true solution.  15 equally
		spaced contour lines are shown at the same levels for the learned model solution and true solution. 
	}\label{fig:ex5_soluC}
\end{figure}

\begin{figure}[htbp]
	\centering
	{\includegraphics[width=0.48\textwidth]{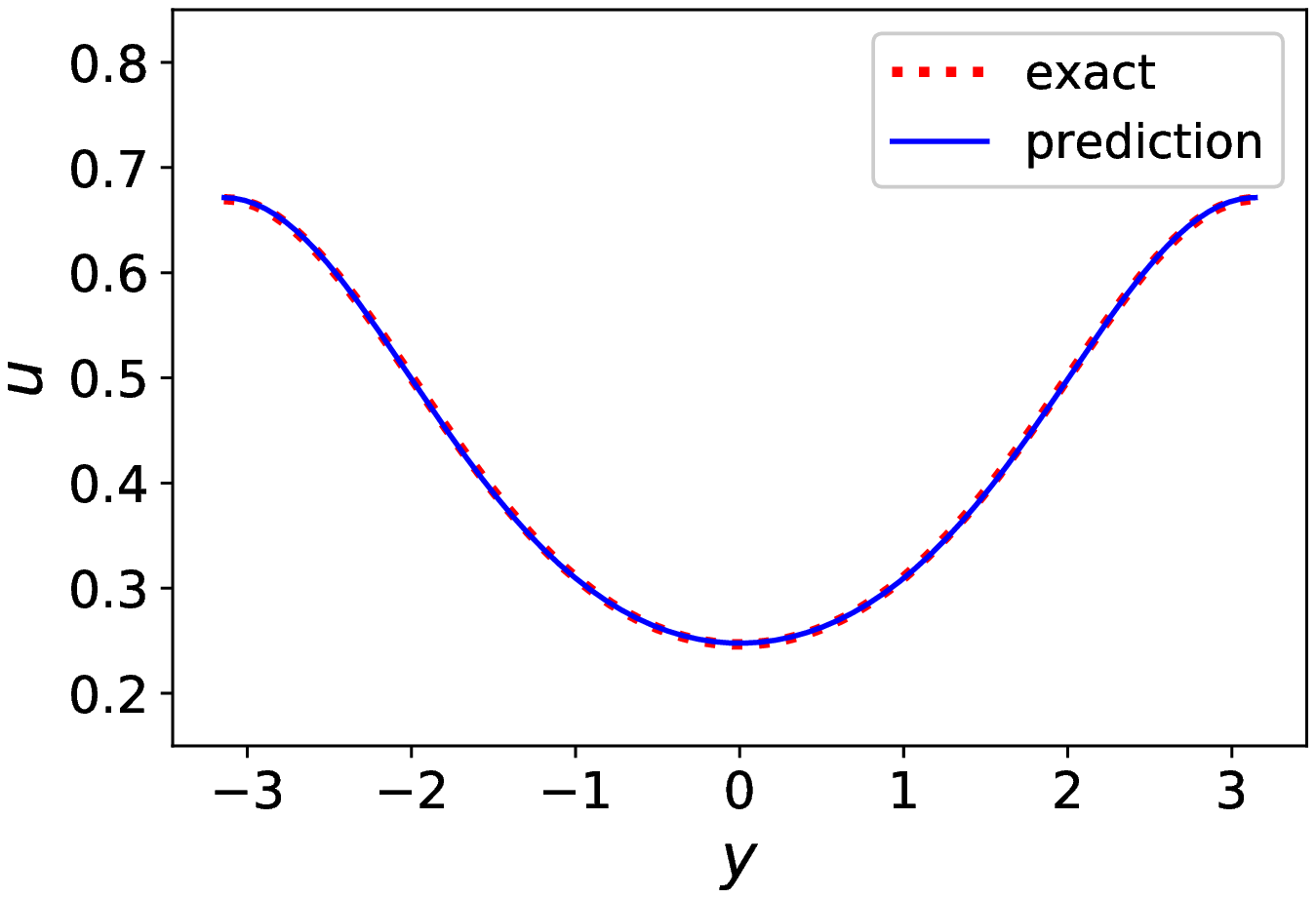}}
	{\includegraphics[width=0.48\textwidth]{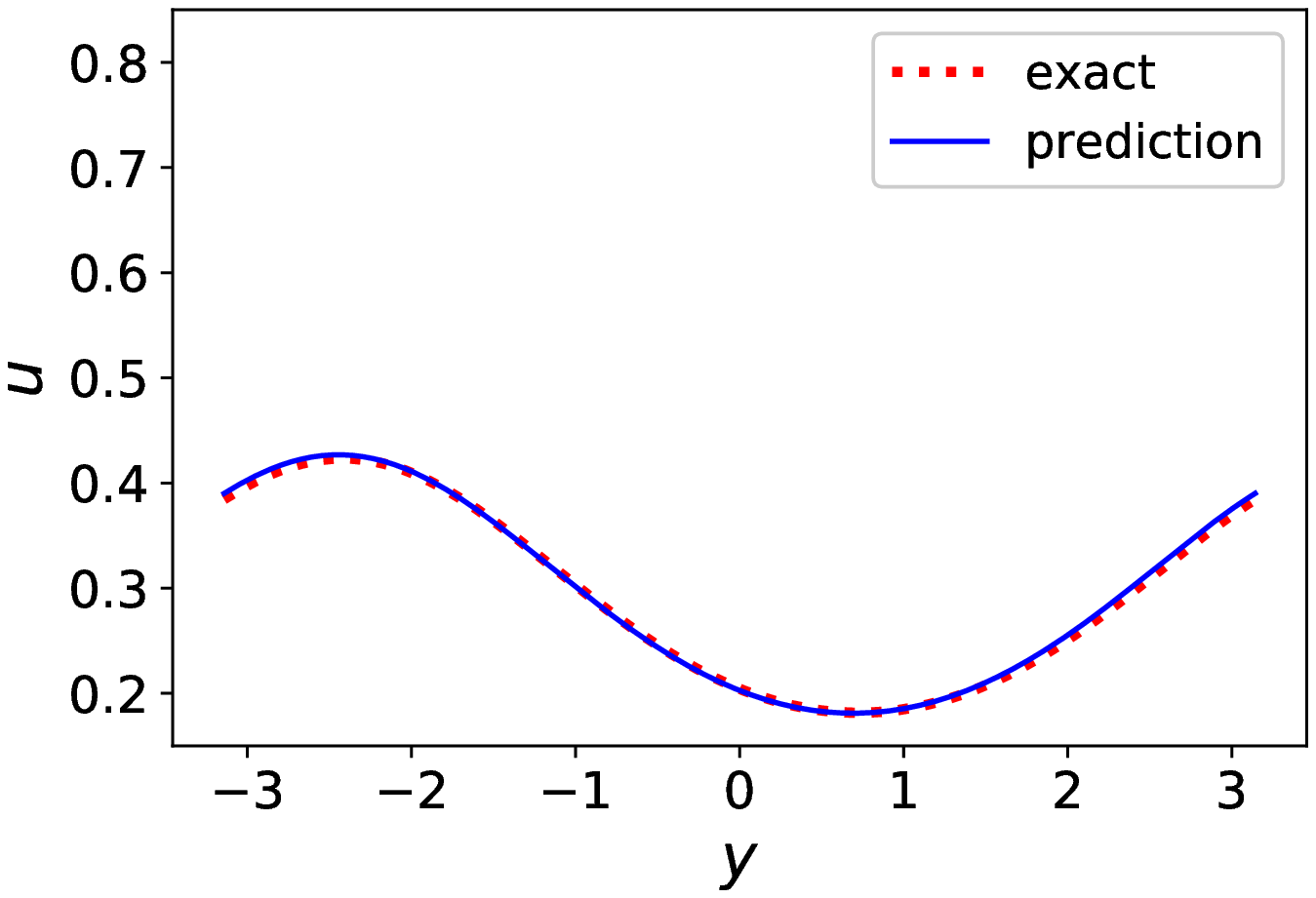}}
	{\includegraphics[width=0.48\textwidth]{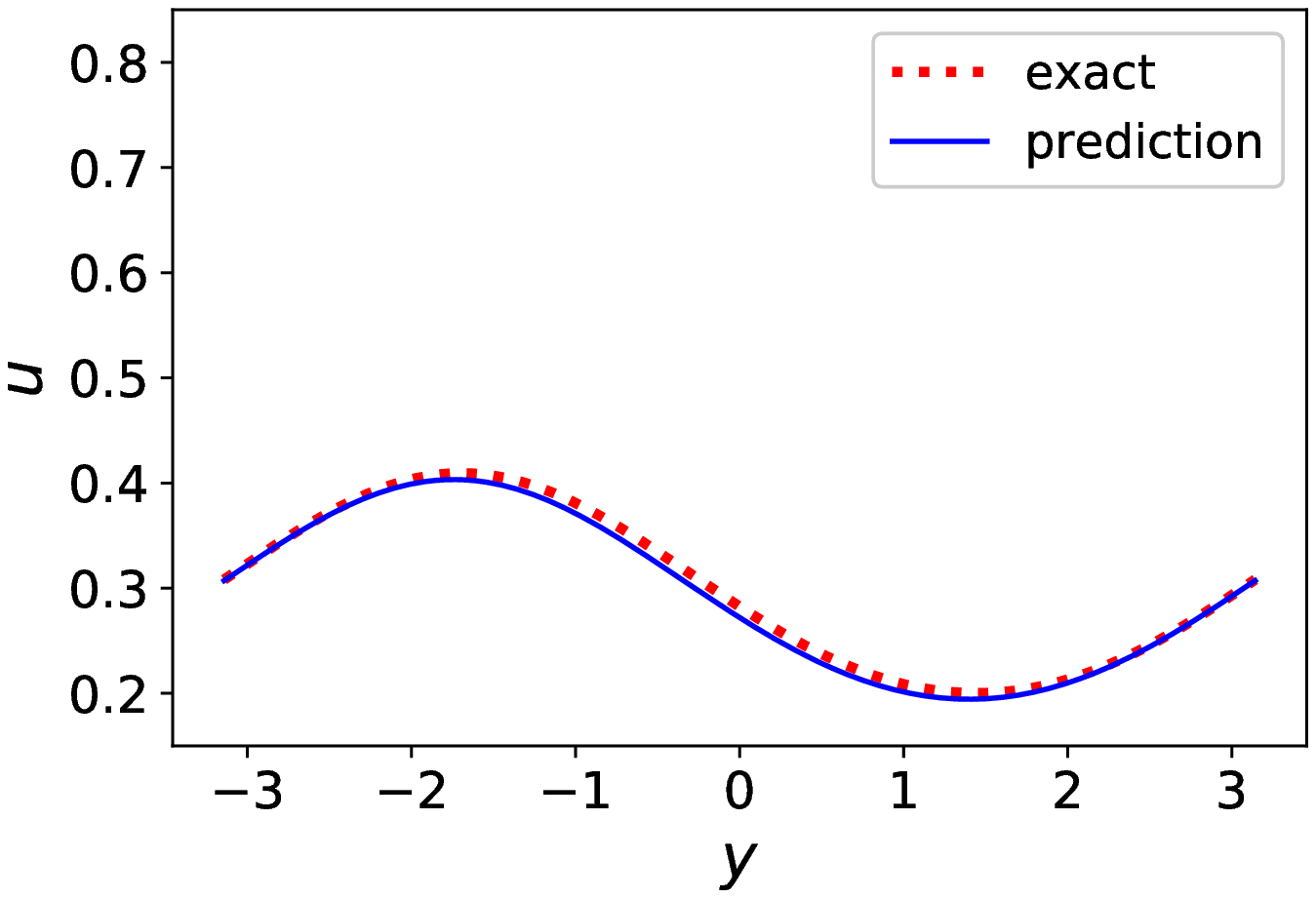}}
	{\includegraphics[width=0.48\textwidth]{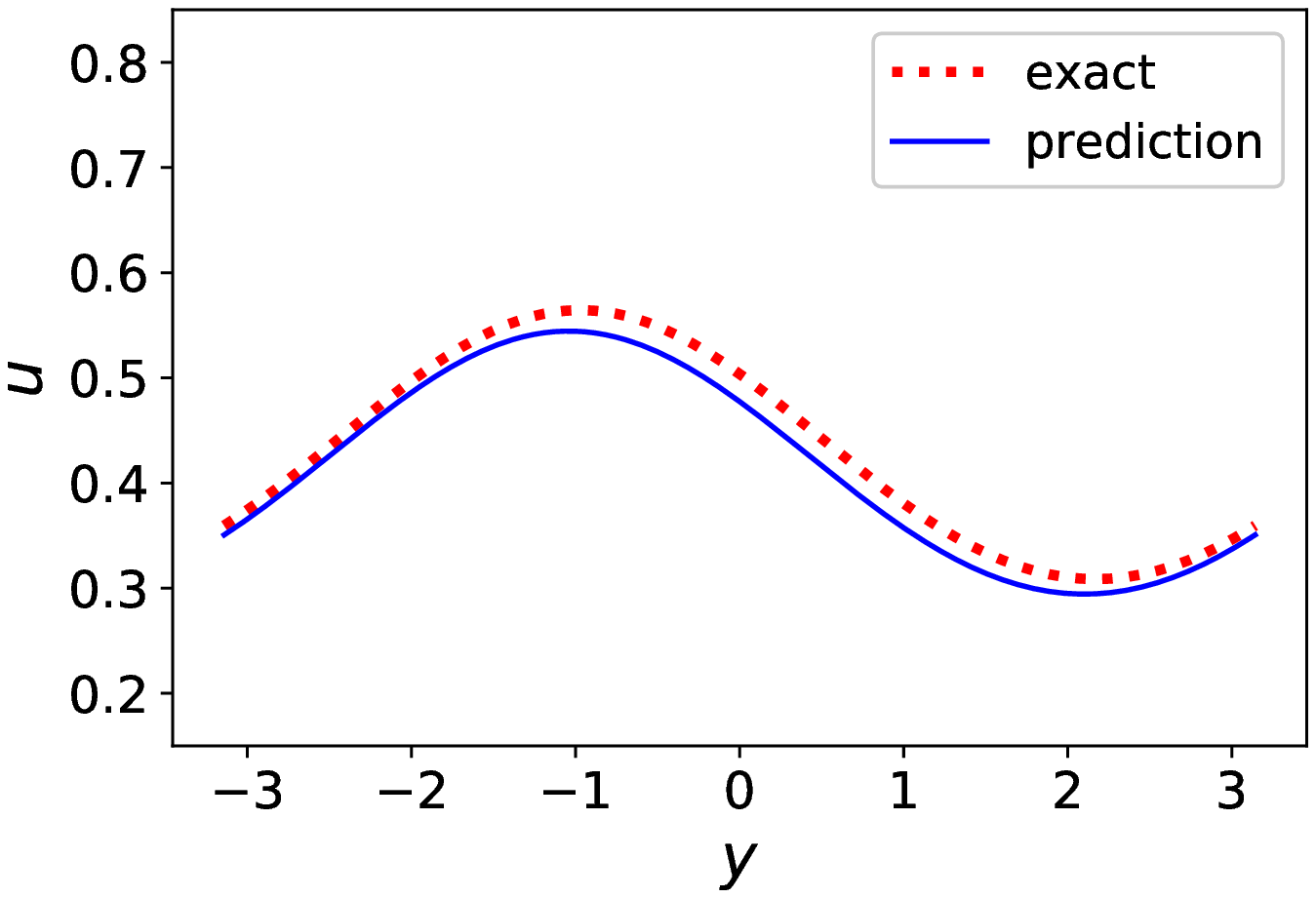}}
	\caption{\small
		Example 5: Comparison of the true solution and the learned model solution at different time along the line $x=0$. Top-left: $t=0$; top-right: $t=1$; bottom-left: $t=2$; bottom-right: 
		$t=3$.  
	}\label{fig:ex5_soluA}
\end{figure}

\begin{figure}[htbp]
	\centering
%	{\includegraphics[width=0.48\textwidth]{Figure/Example5/ABSerror.eps}}
	{\includegraphics[width=0.6\textwidth]{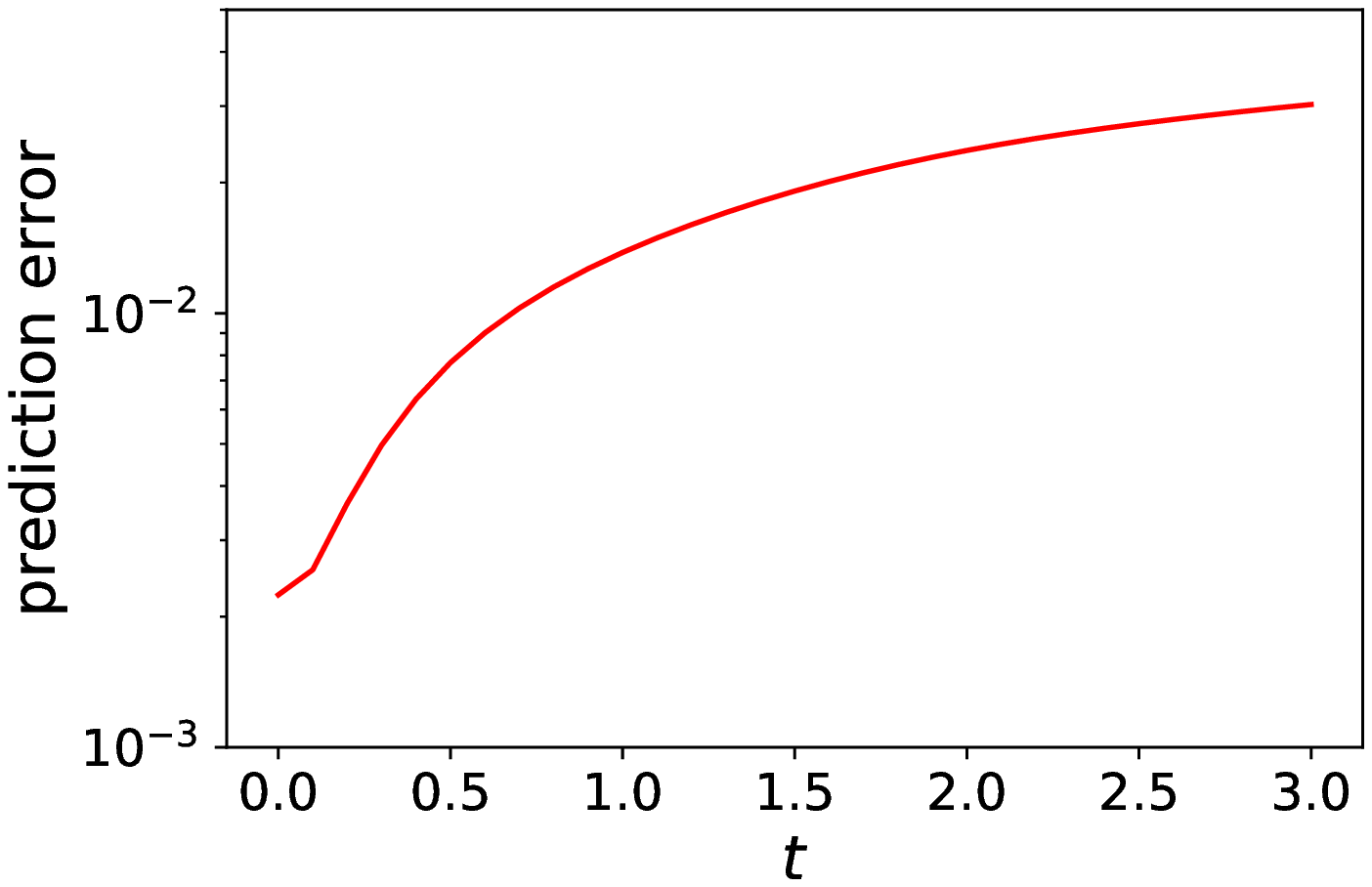}}
	\caption{\small
		Example 5: The evolution of the relative error in the
                prediction  in $l^2$-norm.
                %Left: absolute error; right: relative error.  
	}\label{fig:ex5_error}
\end{figure}

\begin{figure}[htbp]
	\centering
	{\includegraphics[width=0.24\textwidth]{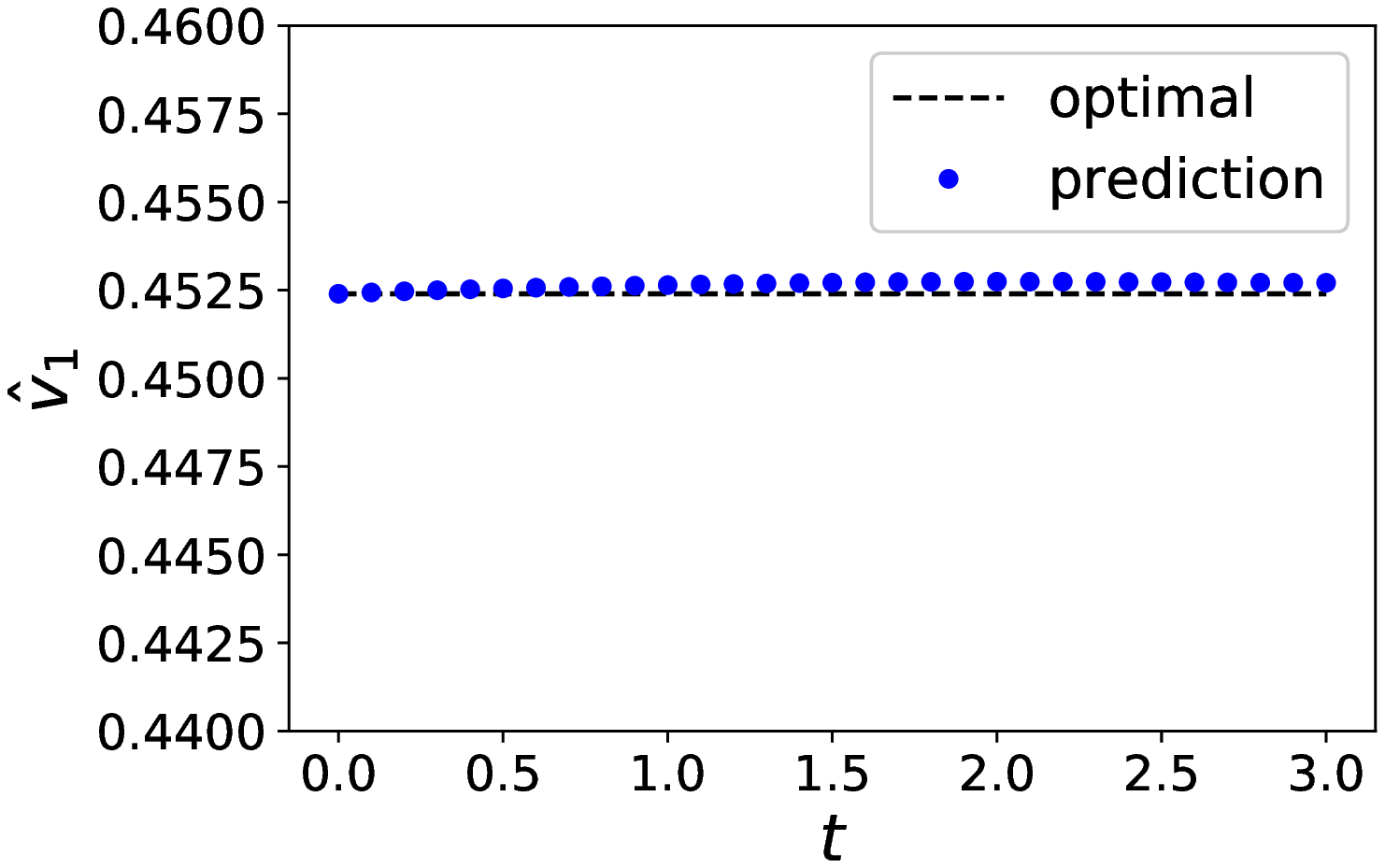}}
	{\includegraphics[width=0.24\textwidth]{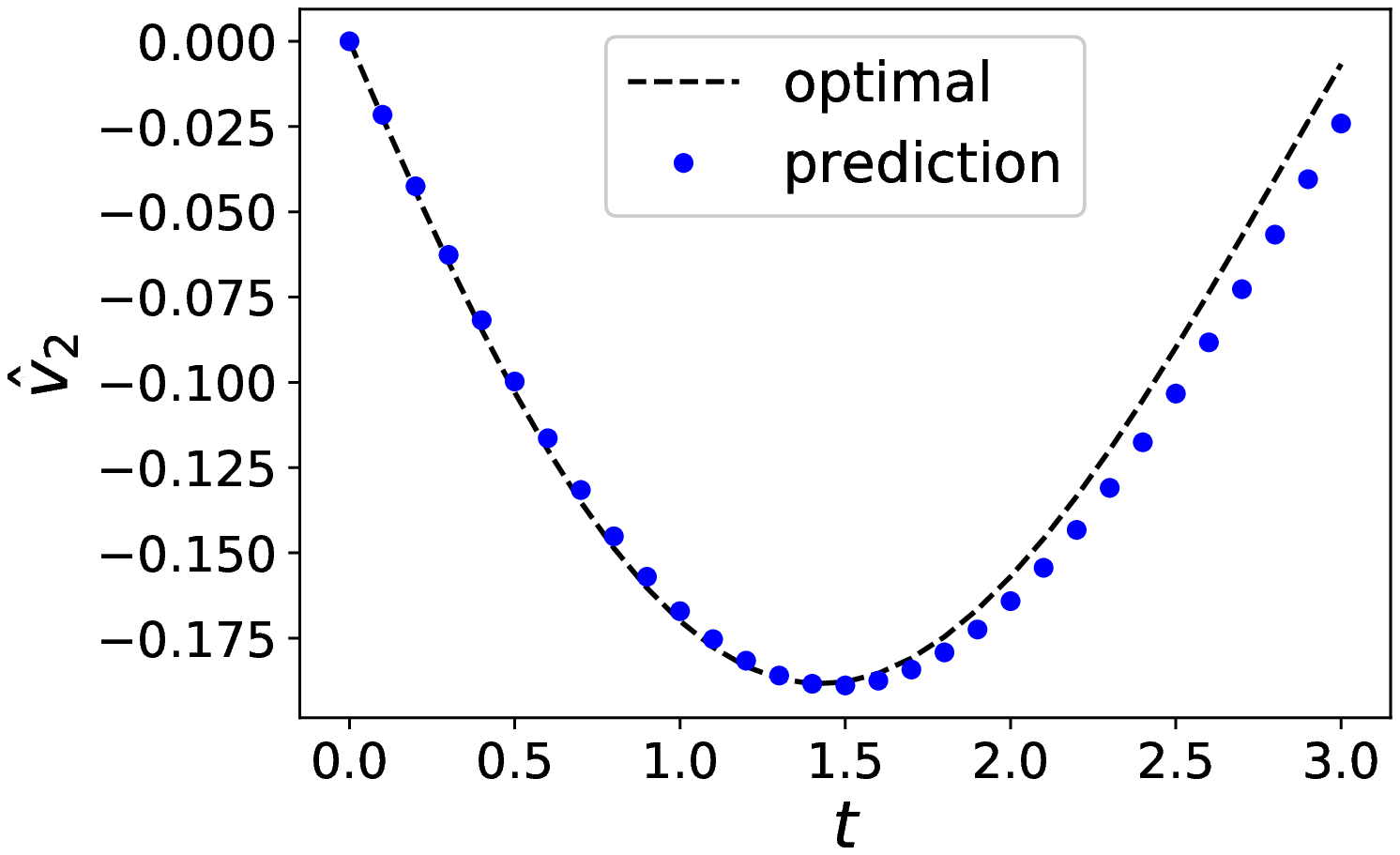}}
	{\includegraphics[width=0.24\textwidth]{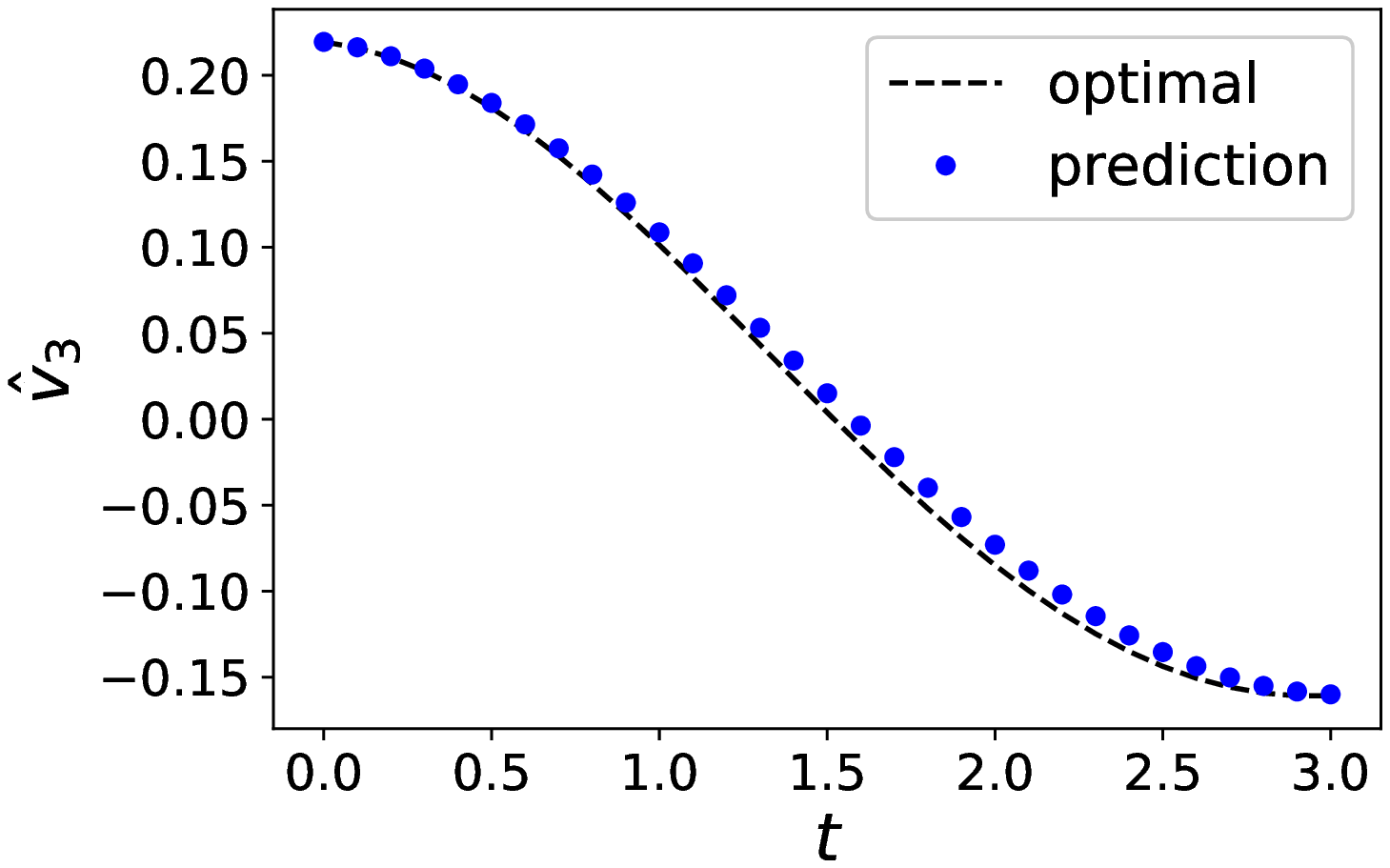}}
	{\includegraphics[width=0.24\textwidth]{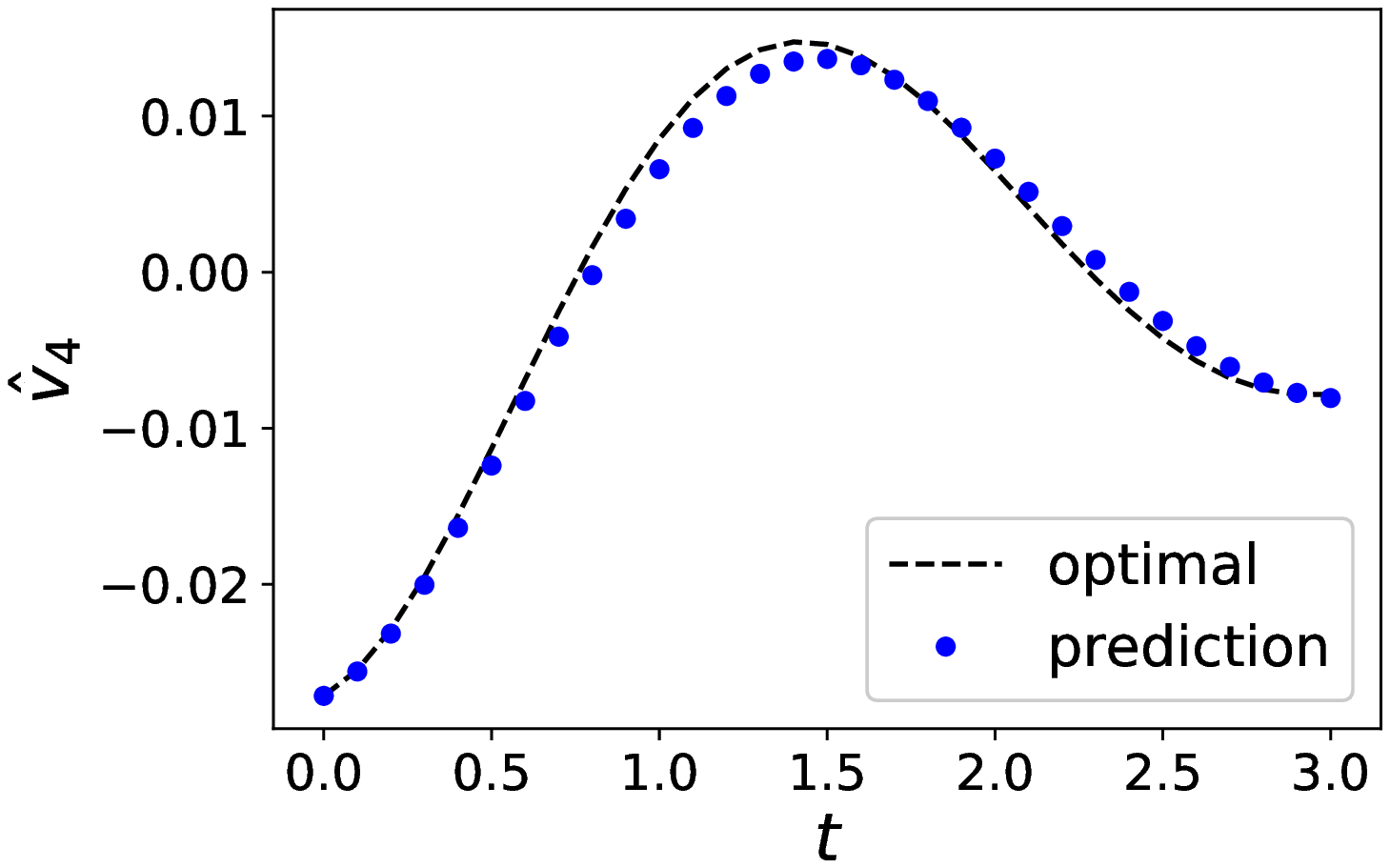}}
	{\includegraphics[width=0.24\textwidth]{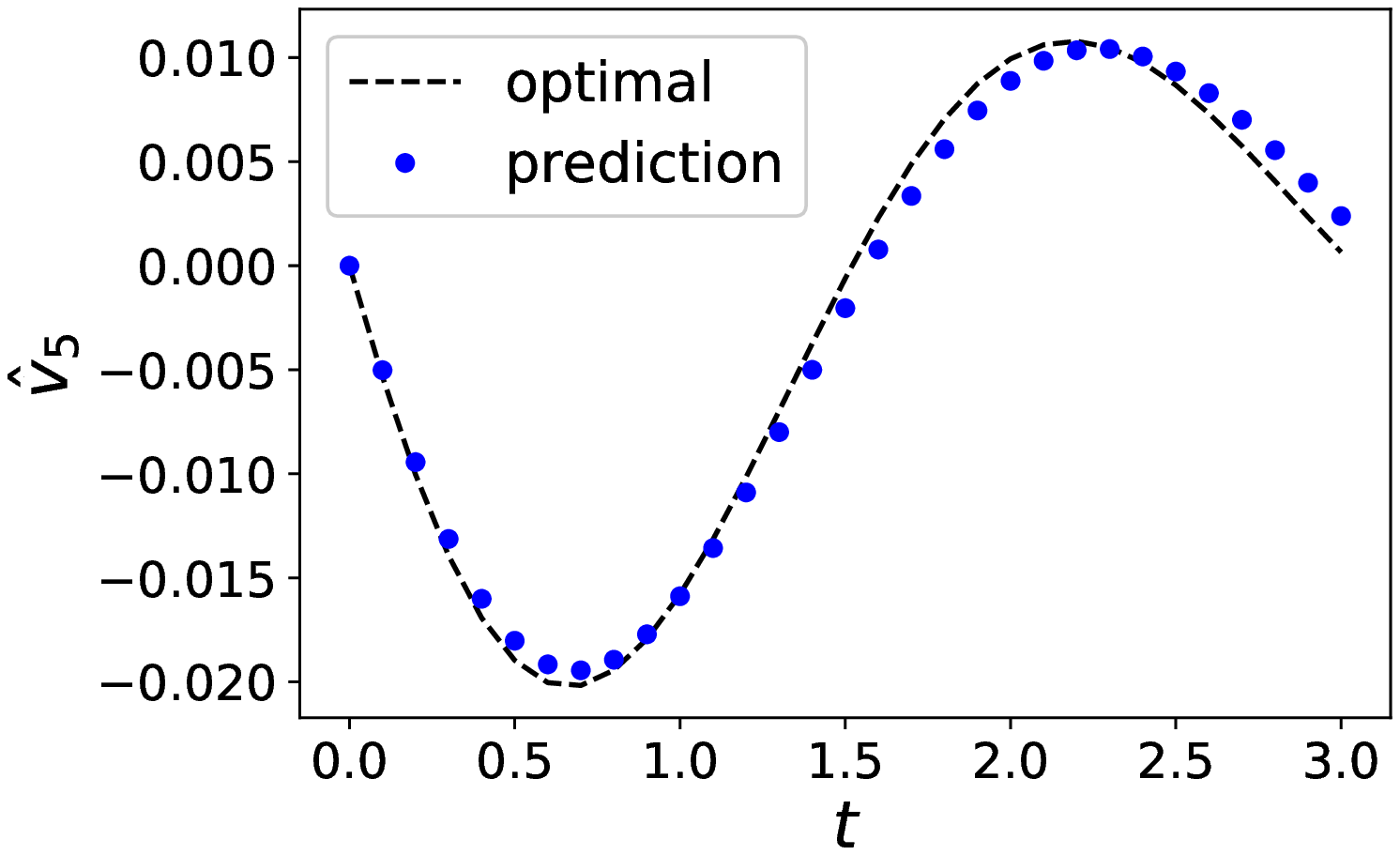}}
	{\includegraphics[width=0.24\textwidth]{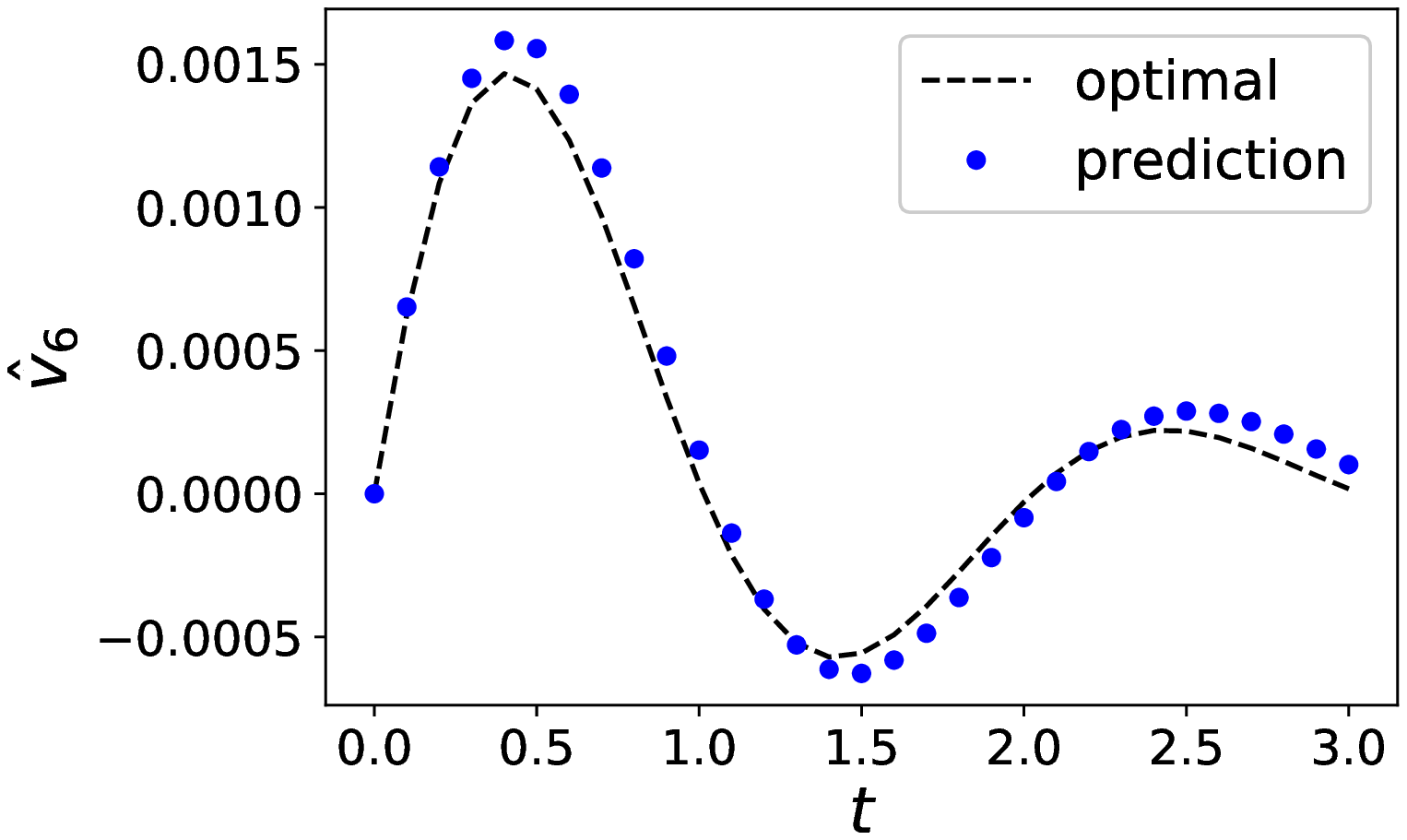}}
	{\includegraphics[width=0.24\textwidth]{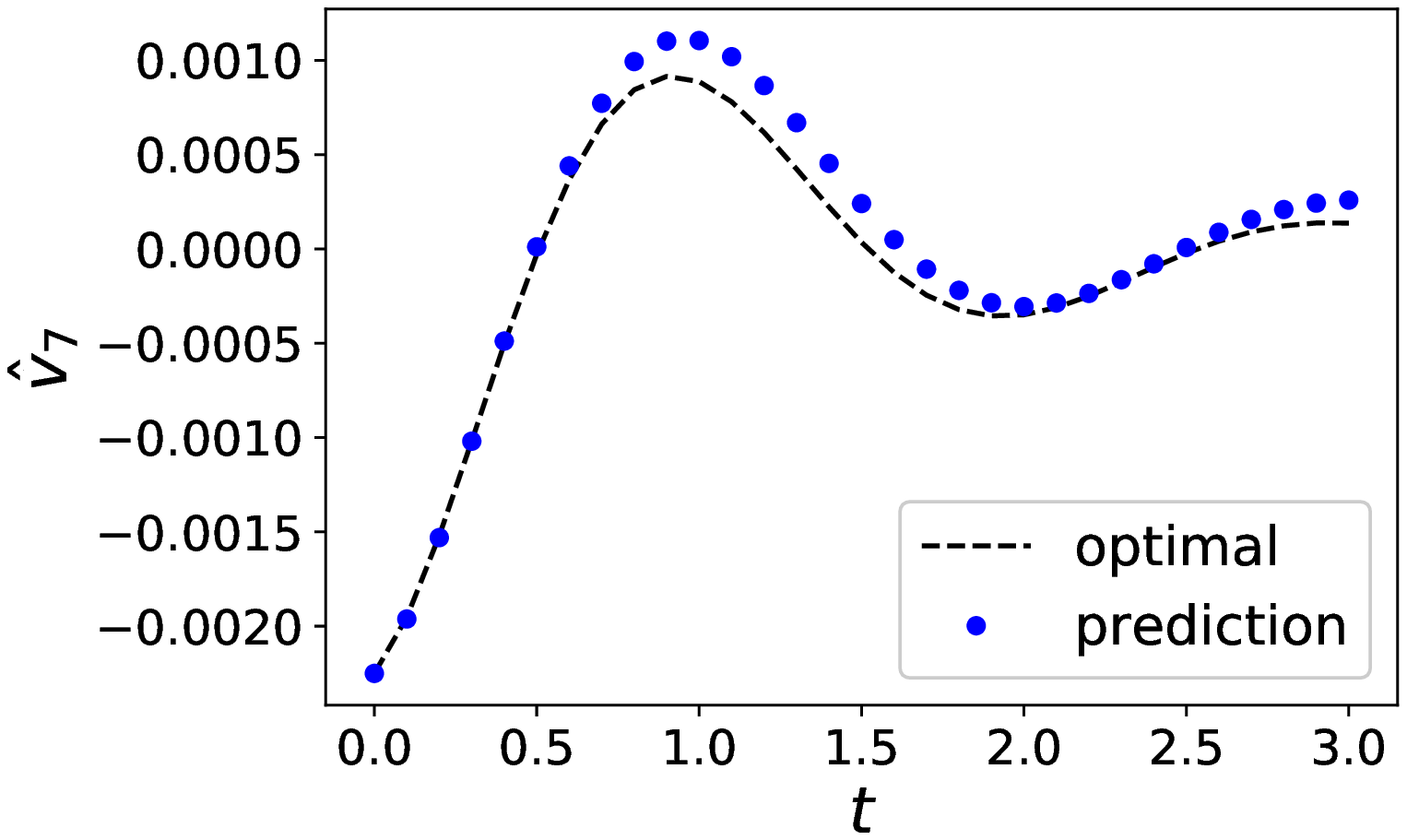}}
	{\includegraphics[width=0.24\textwidth]{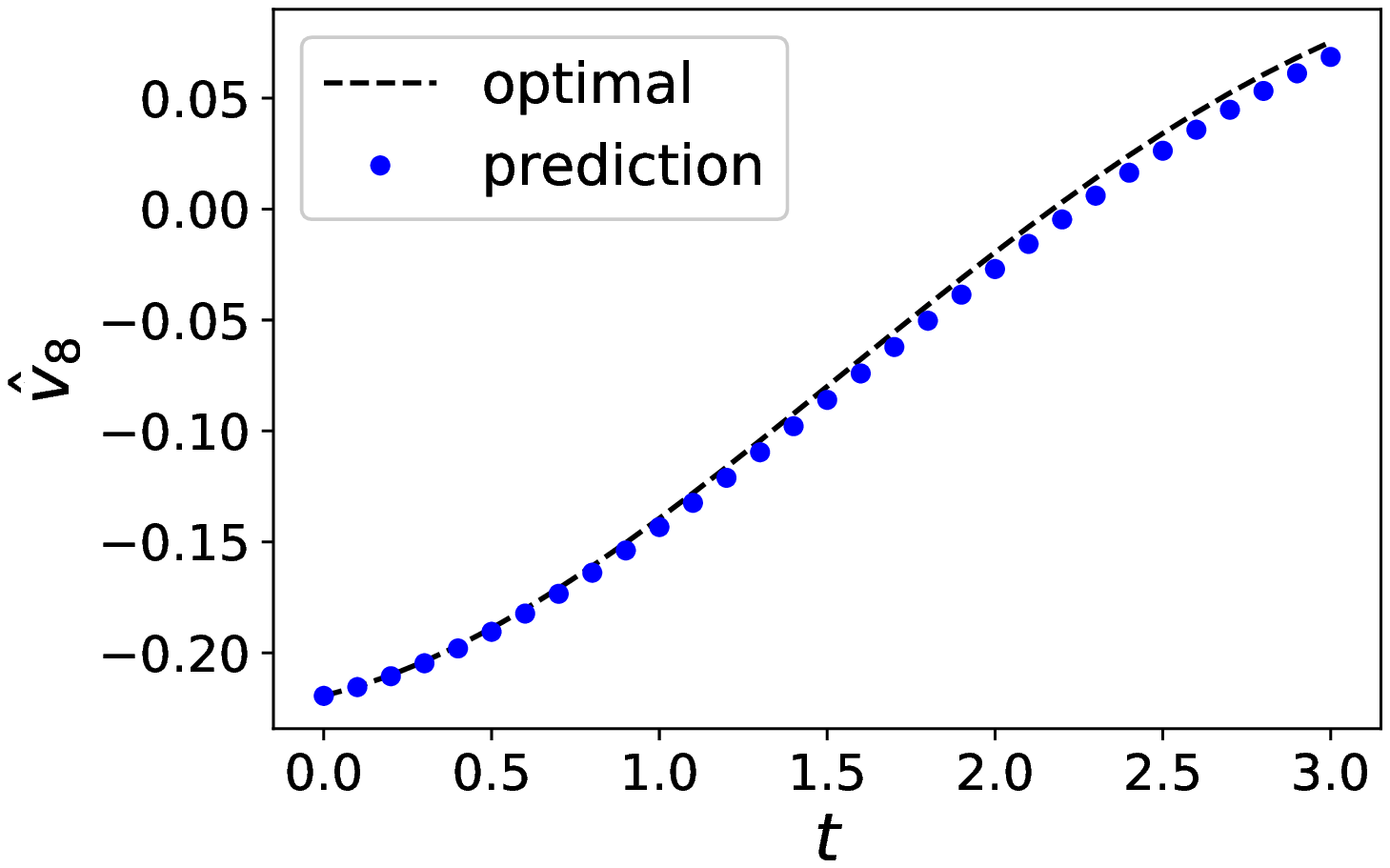}}
	{\includegraphics[width=0.24\textwidth]{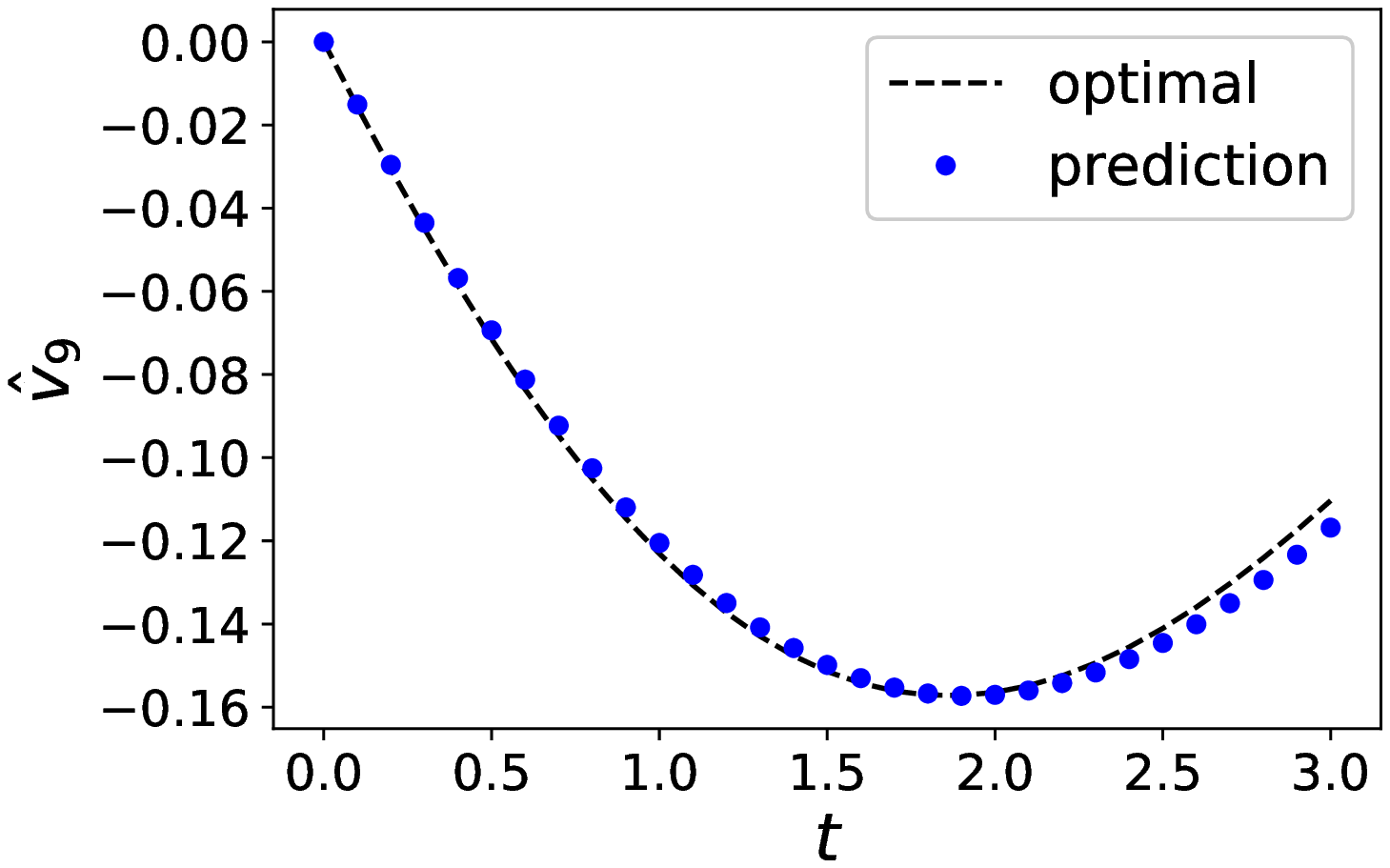}}
	{\includegraphics[width=0.24\textwidth]{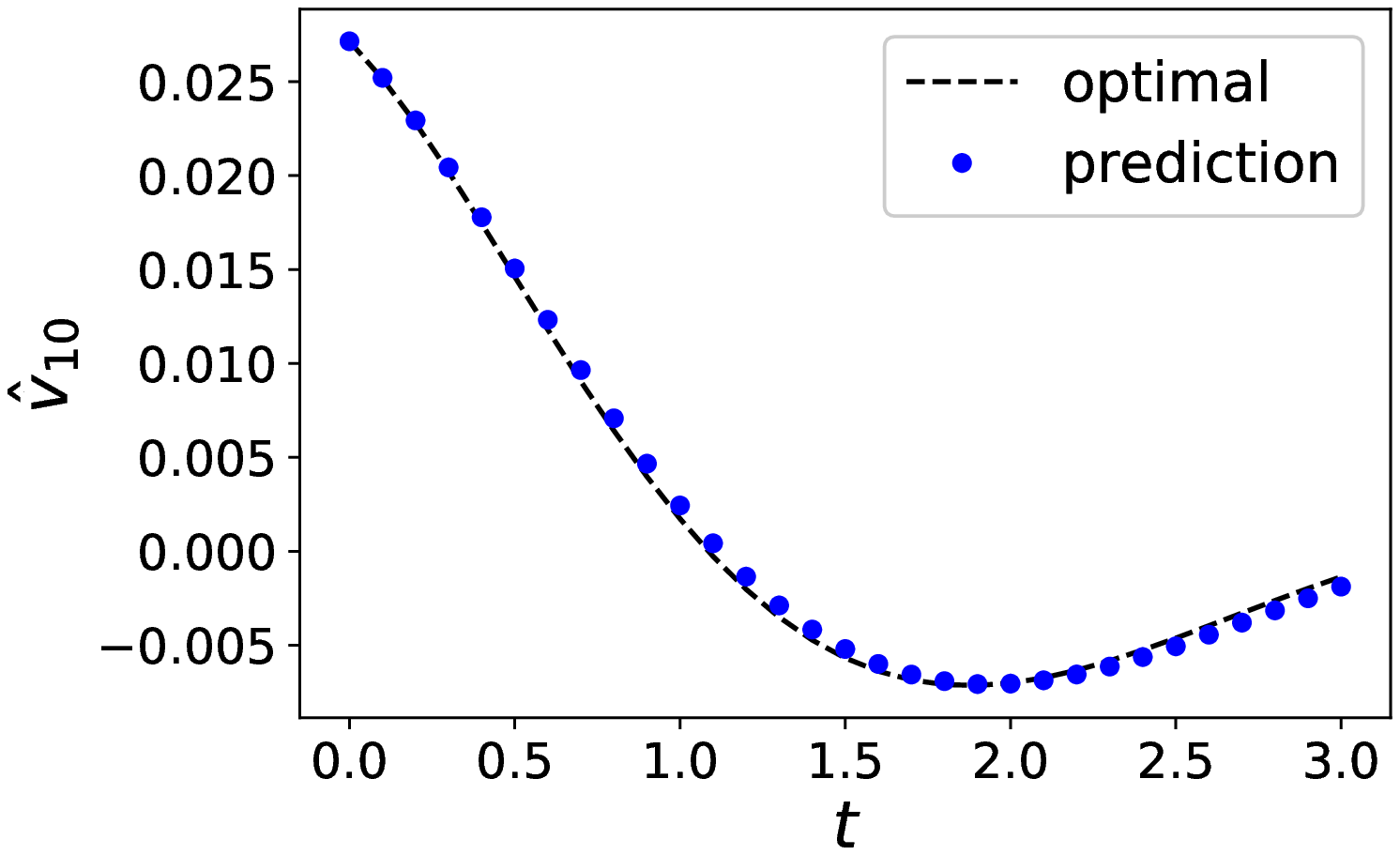}}
	{\includegraphics[width=0.24\textwidth]{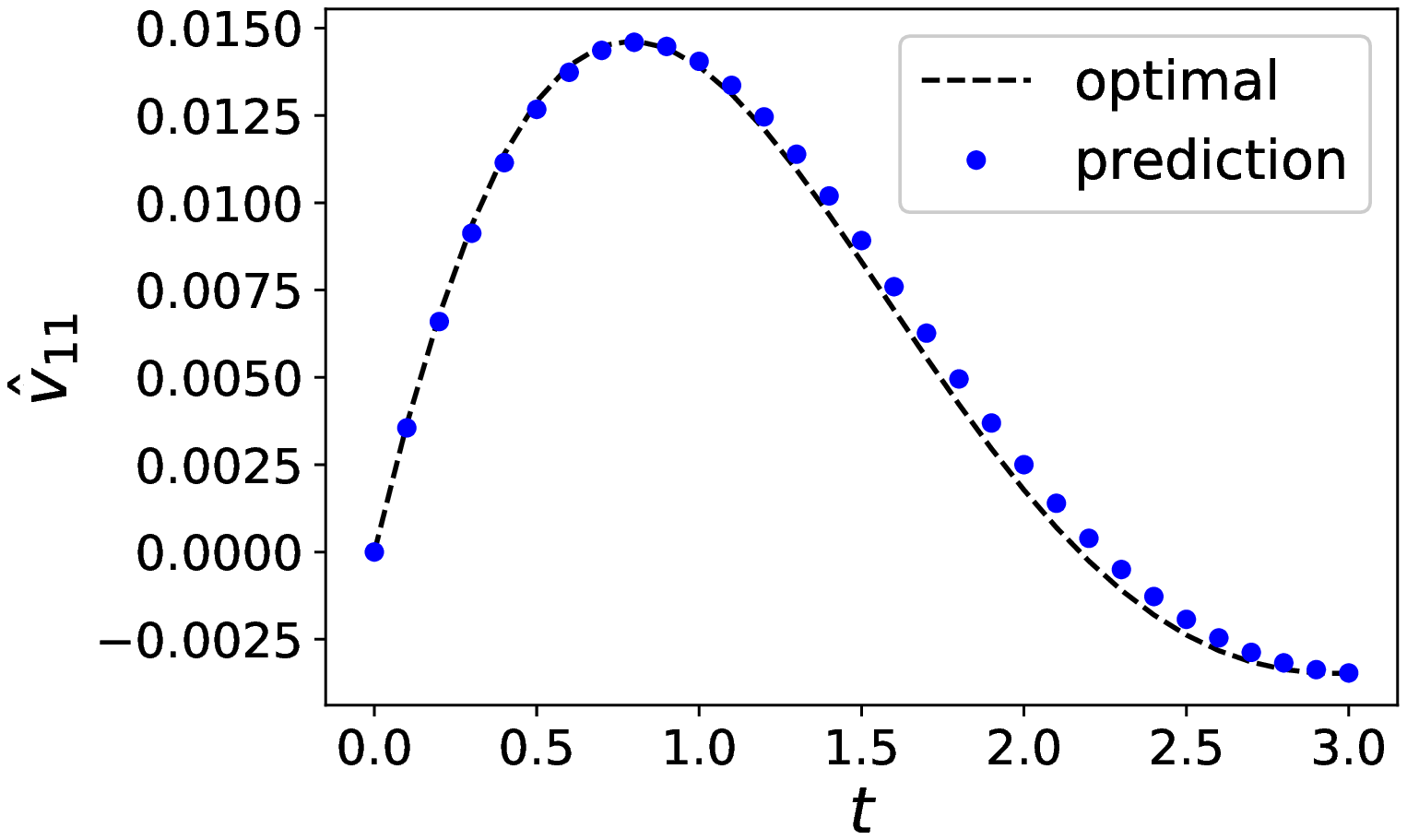}}
	{\includegraphics[width=0.24\textwidth]{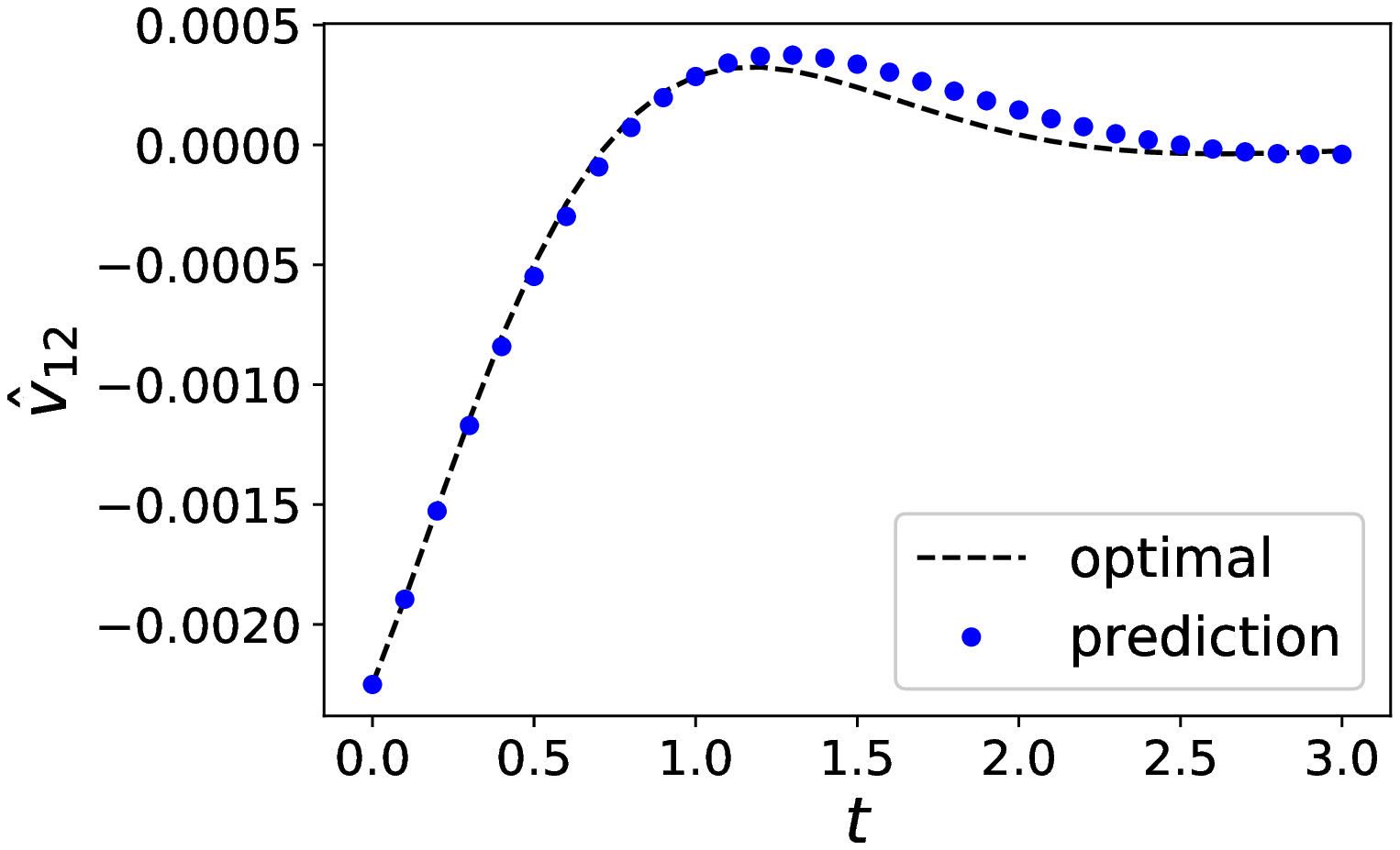}}
	{\includegraphics[width=0.24\textwidth]{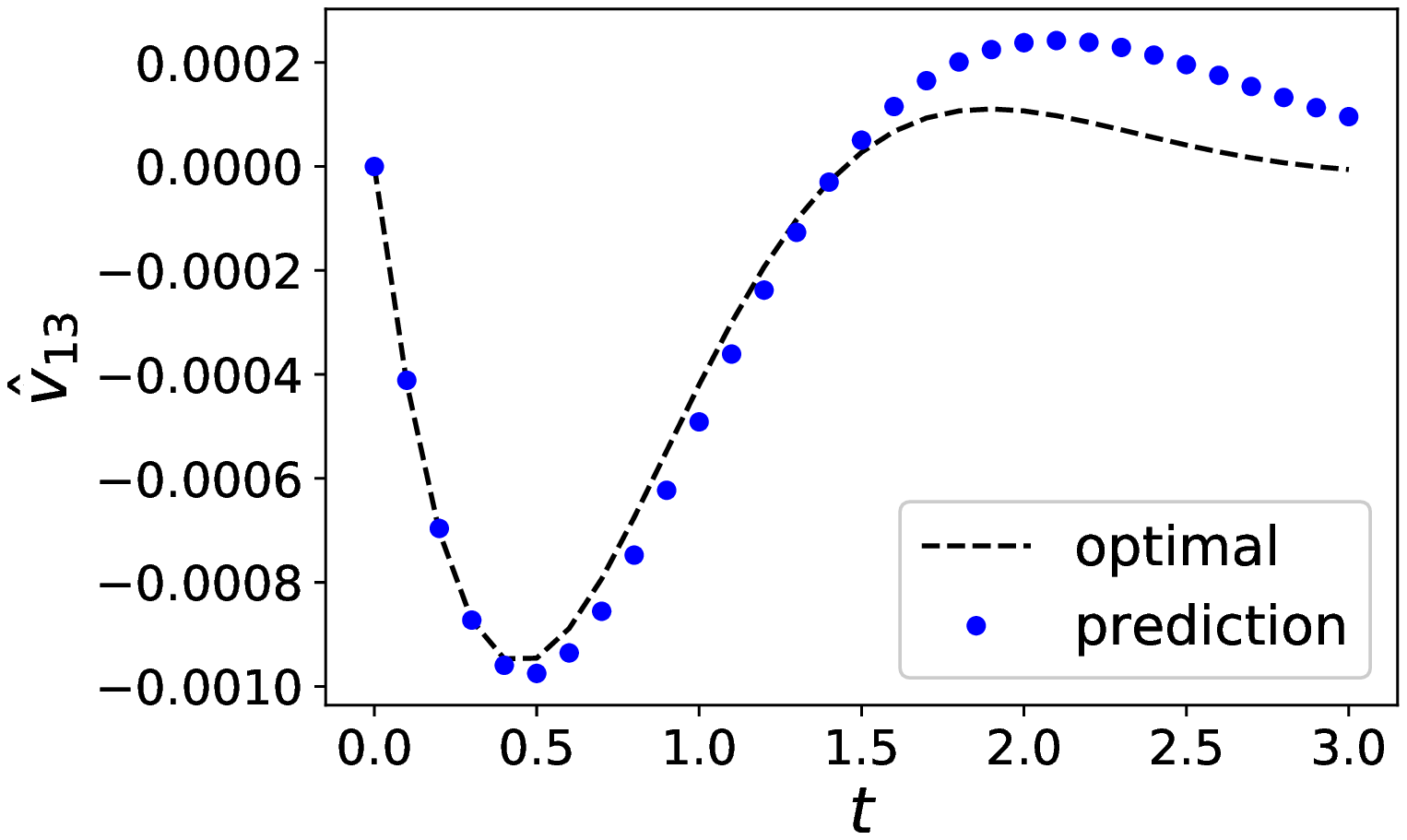}}
	{\includegraphics[width=0.24\textwidth]{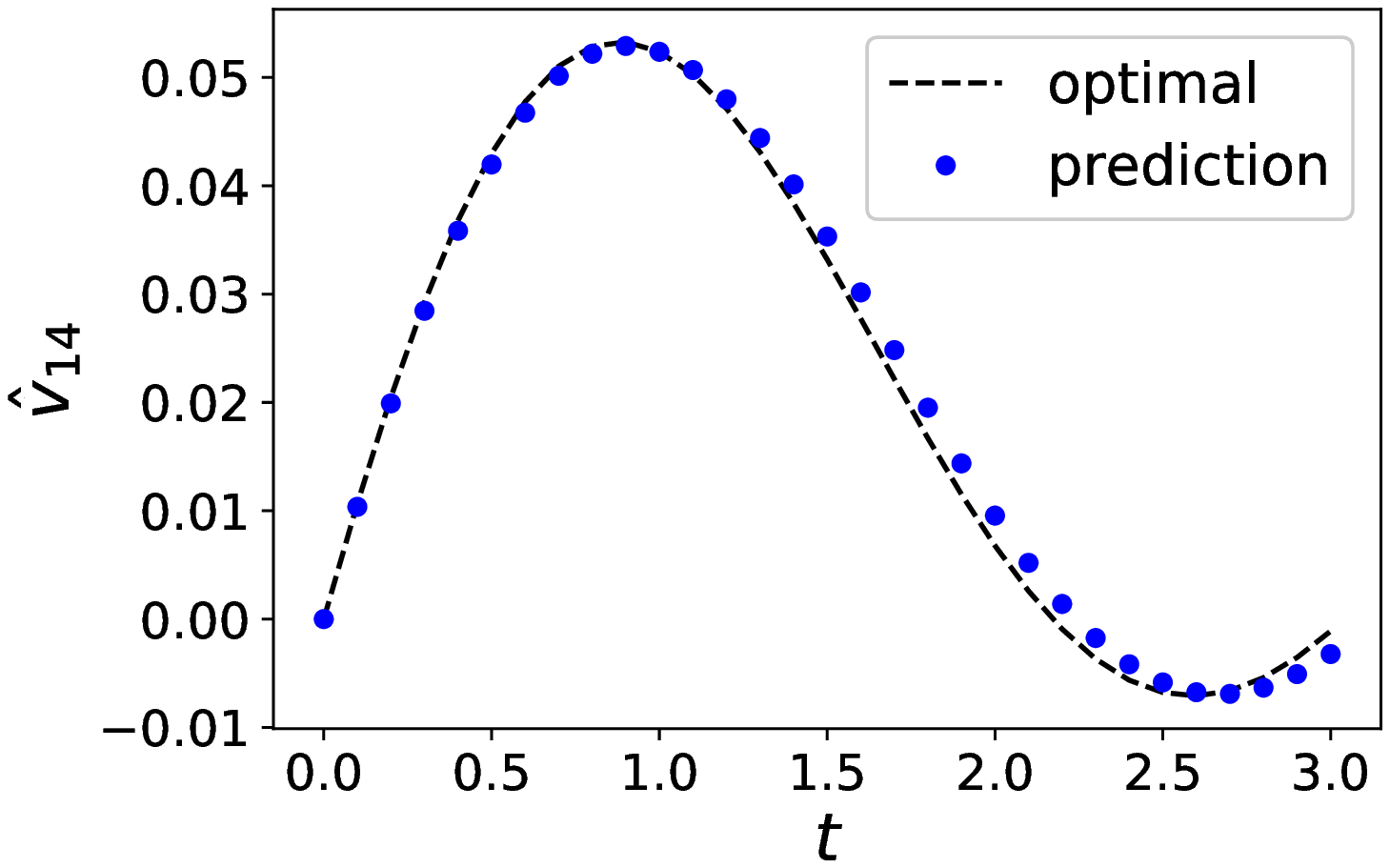}}
	{\includegraphics[width=0.24\textwidth]{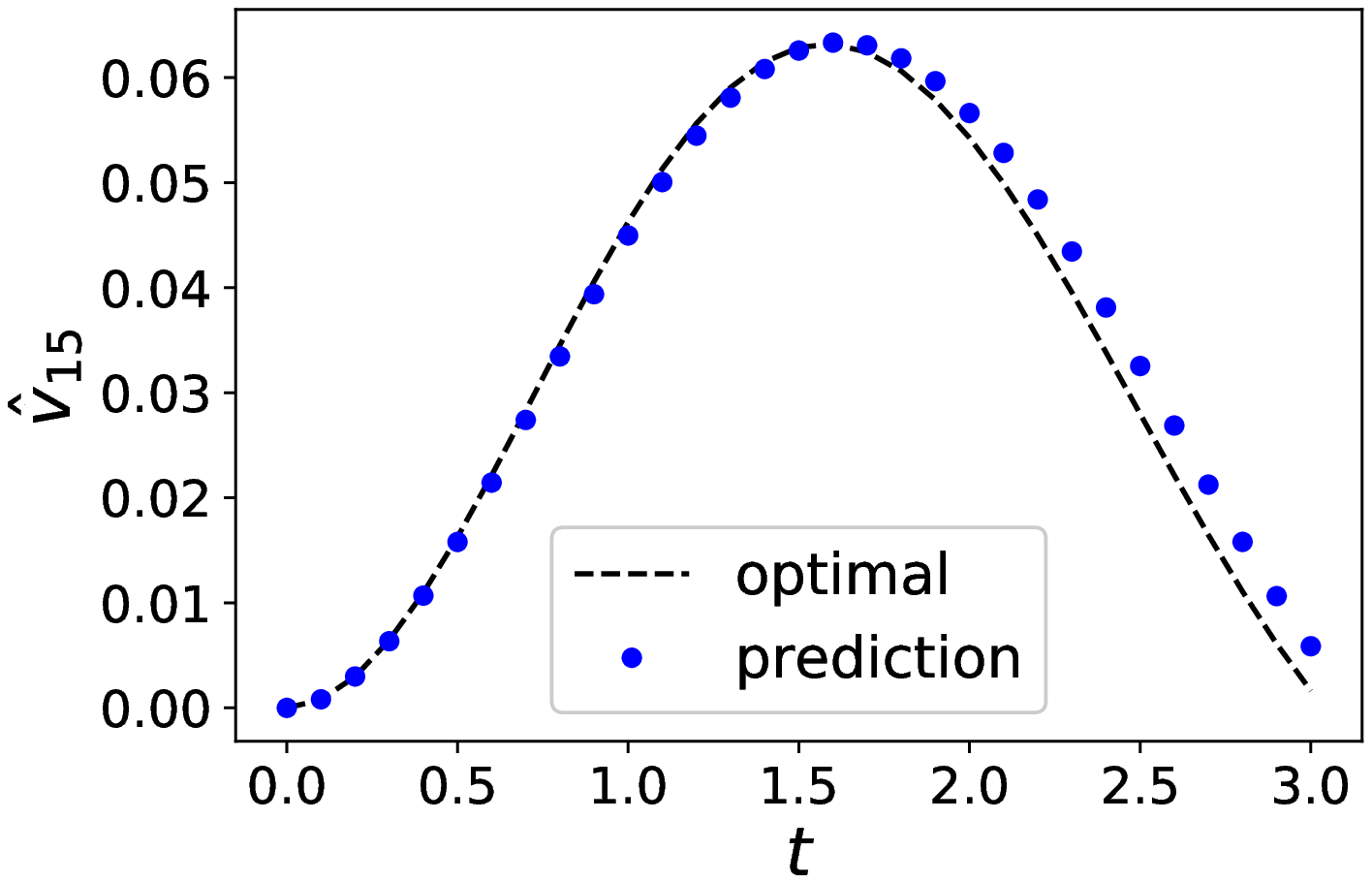}}
	{\includegraphics[width=0.24\textwidth]{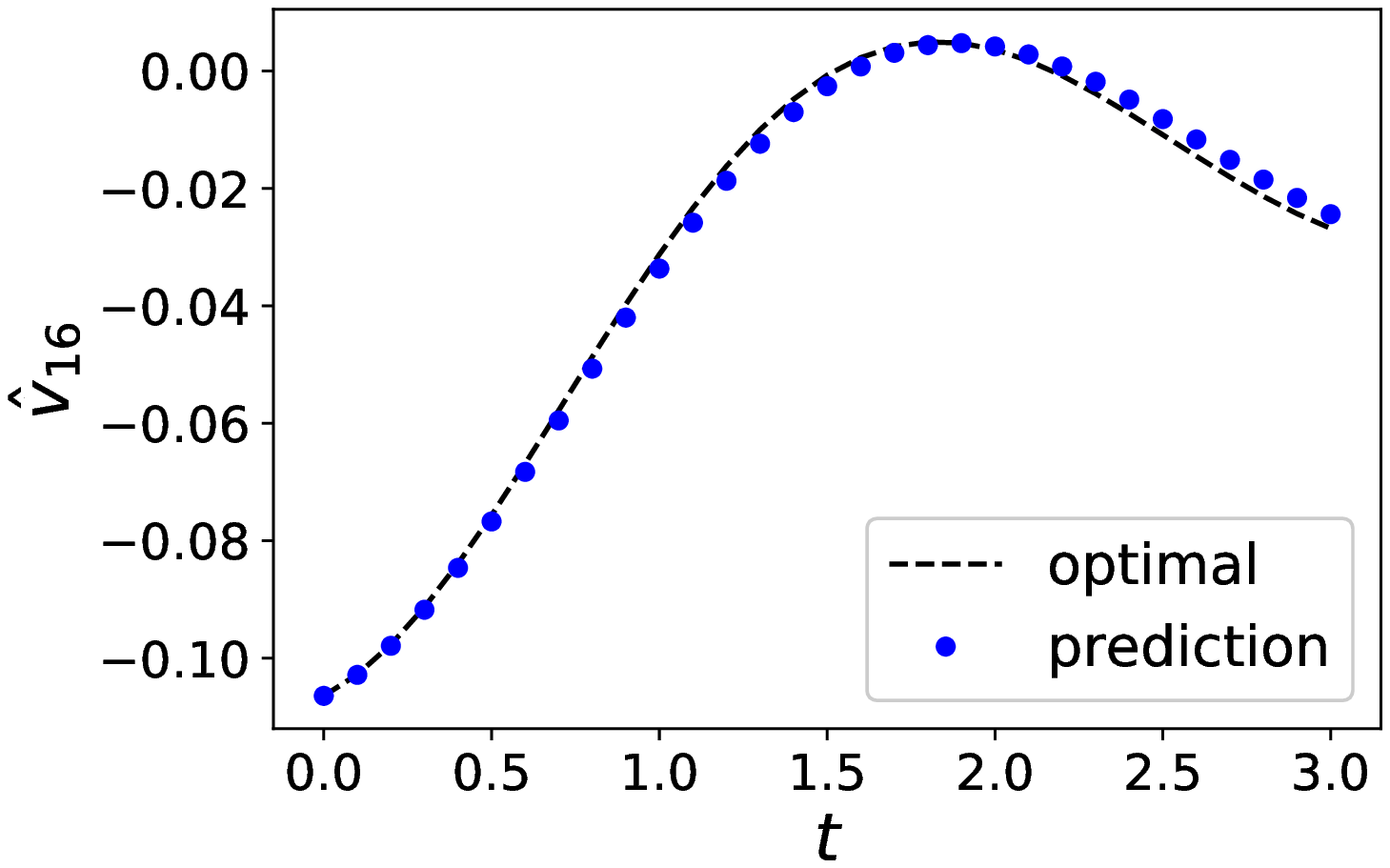}}
	{\includegraphics[width=0.24\textwidth]{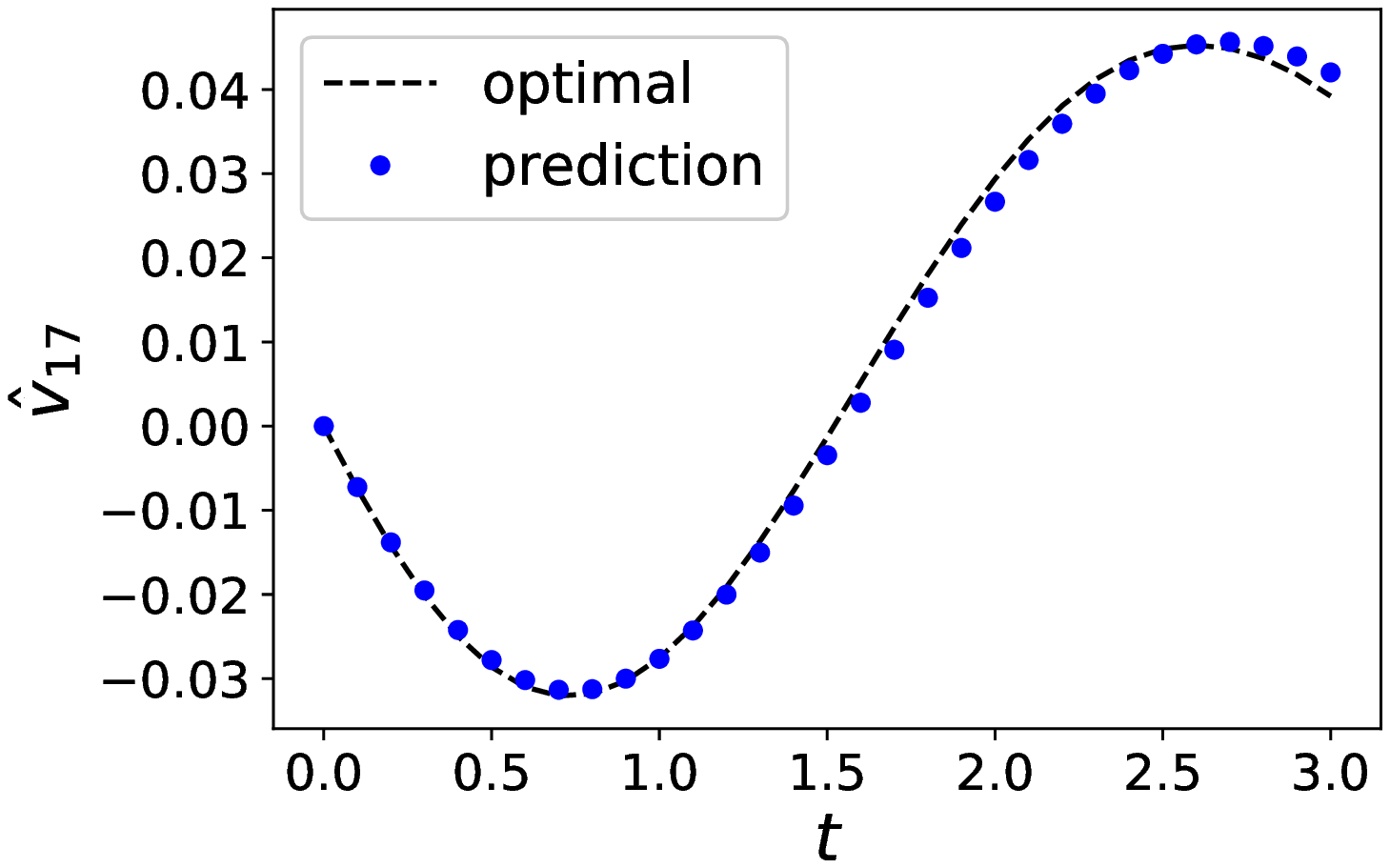}}
	{\includegraphics[width=0.24\textwidth]{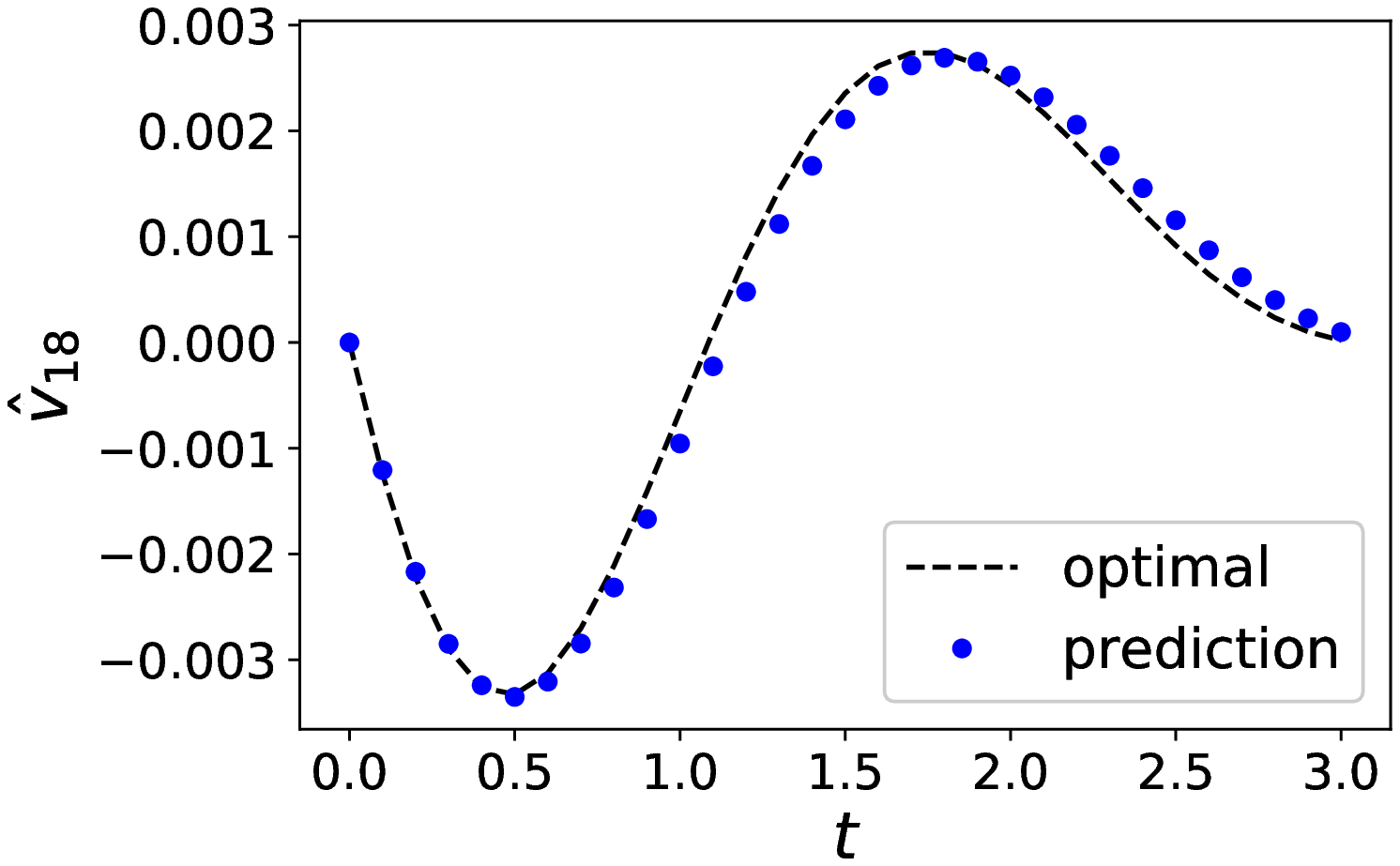}}
	{\includegraphics[width=0.24\textwidth]{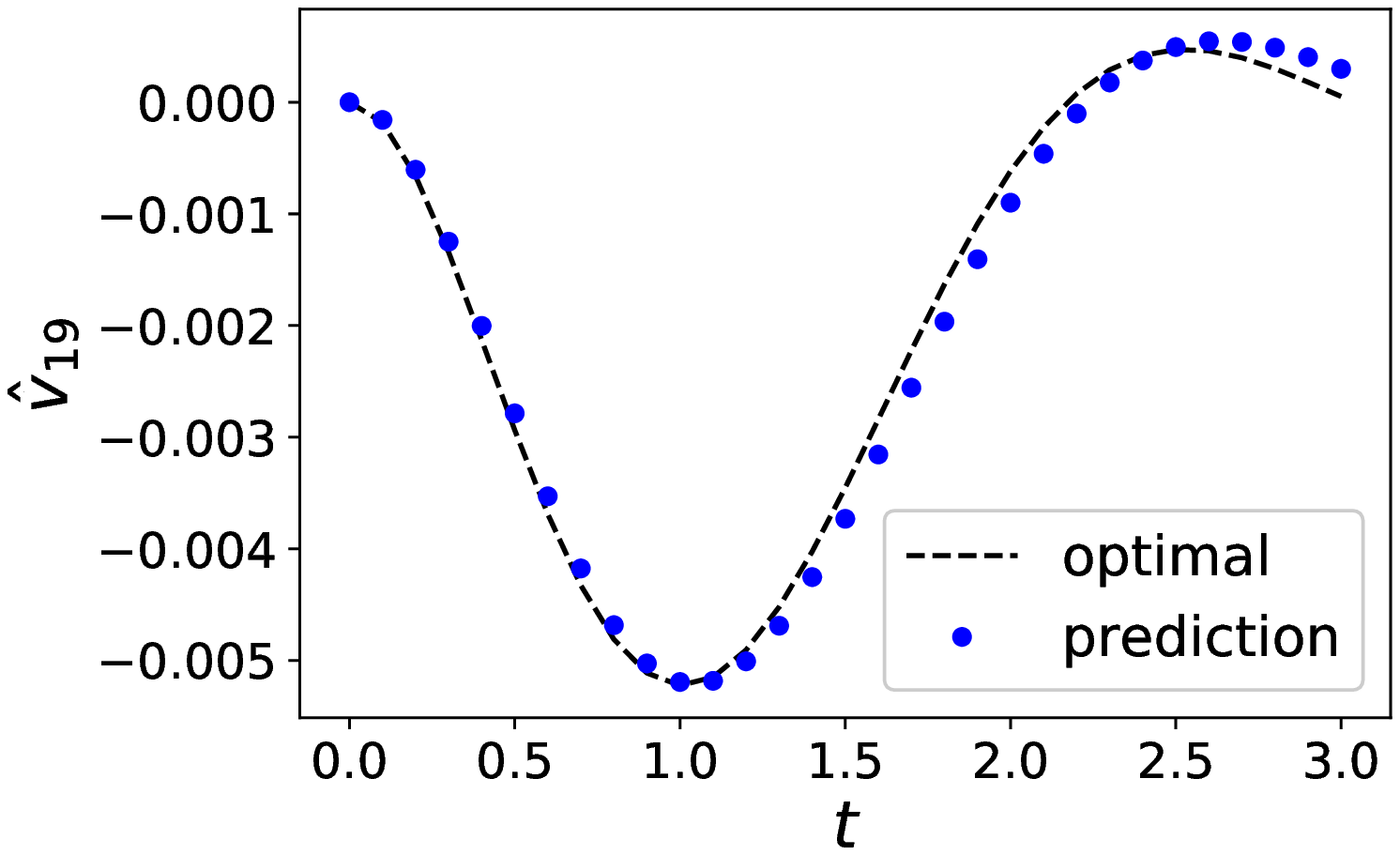}}
	{\includegraphics[width=0.24\textwidth]{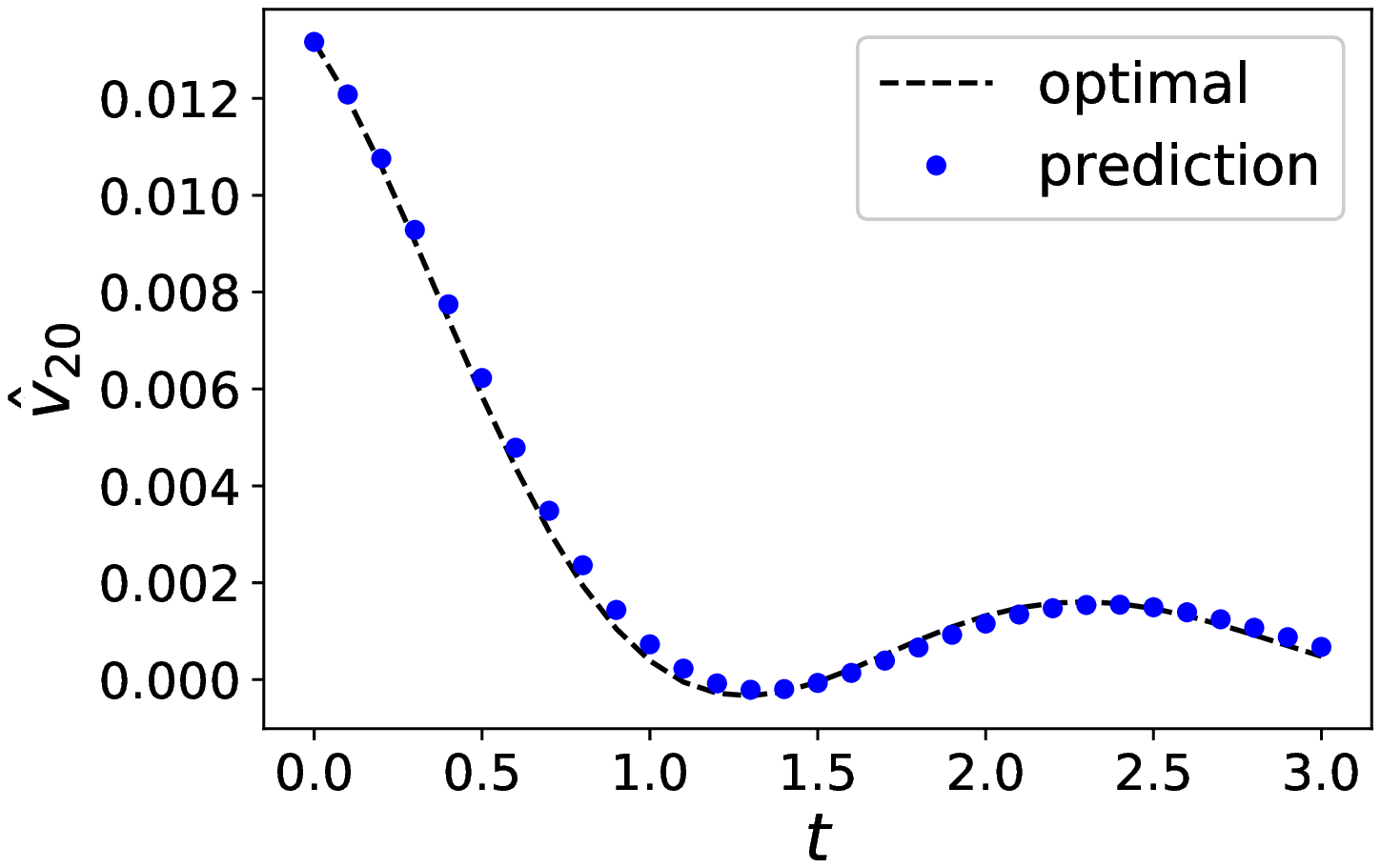}}
	{\includegraphics[width=0.24\textwidth]{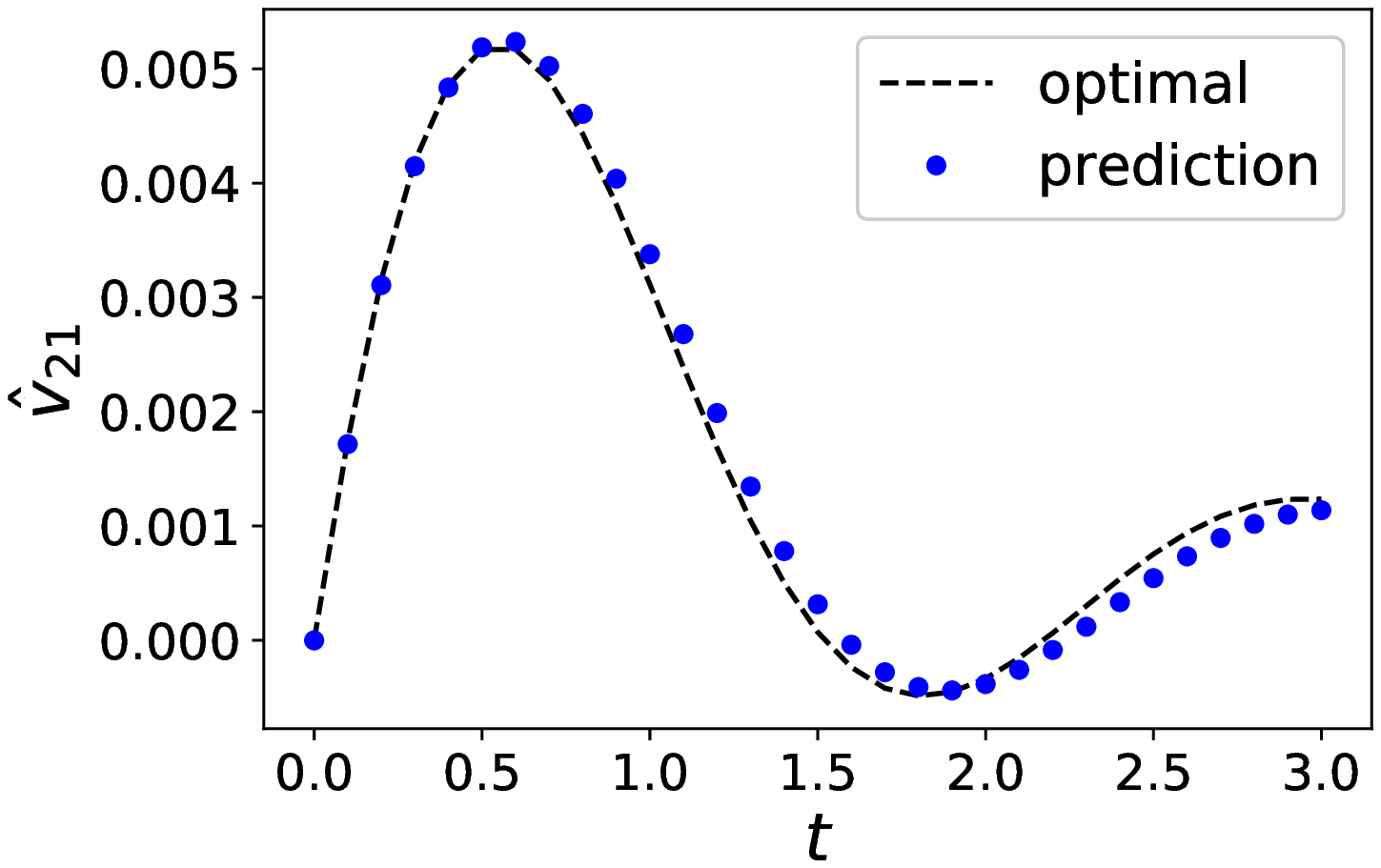}}
	{\includegraphics[width=0.24\textwidth]{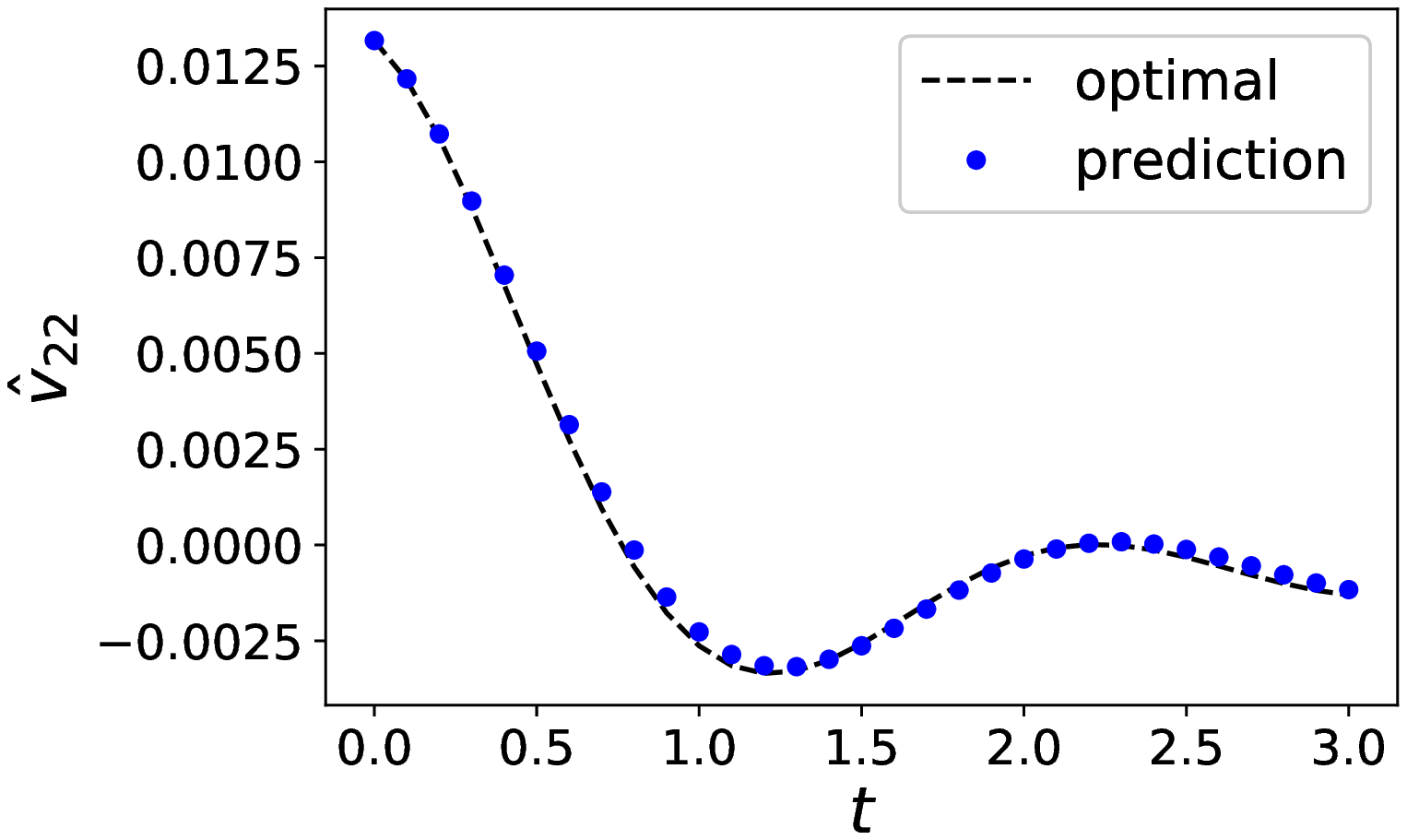}}
	{\includegraphics[width=0.24\textwidth]{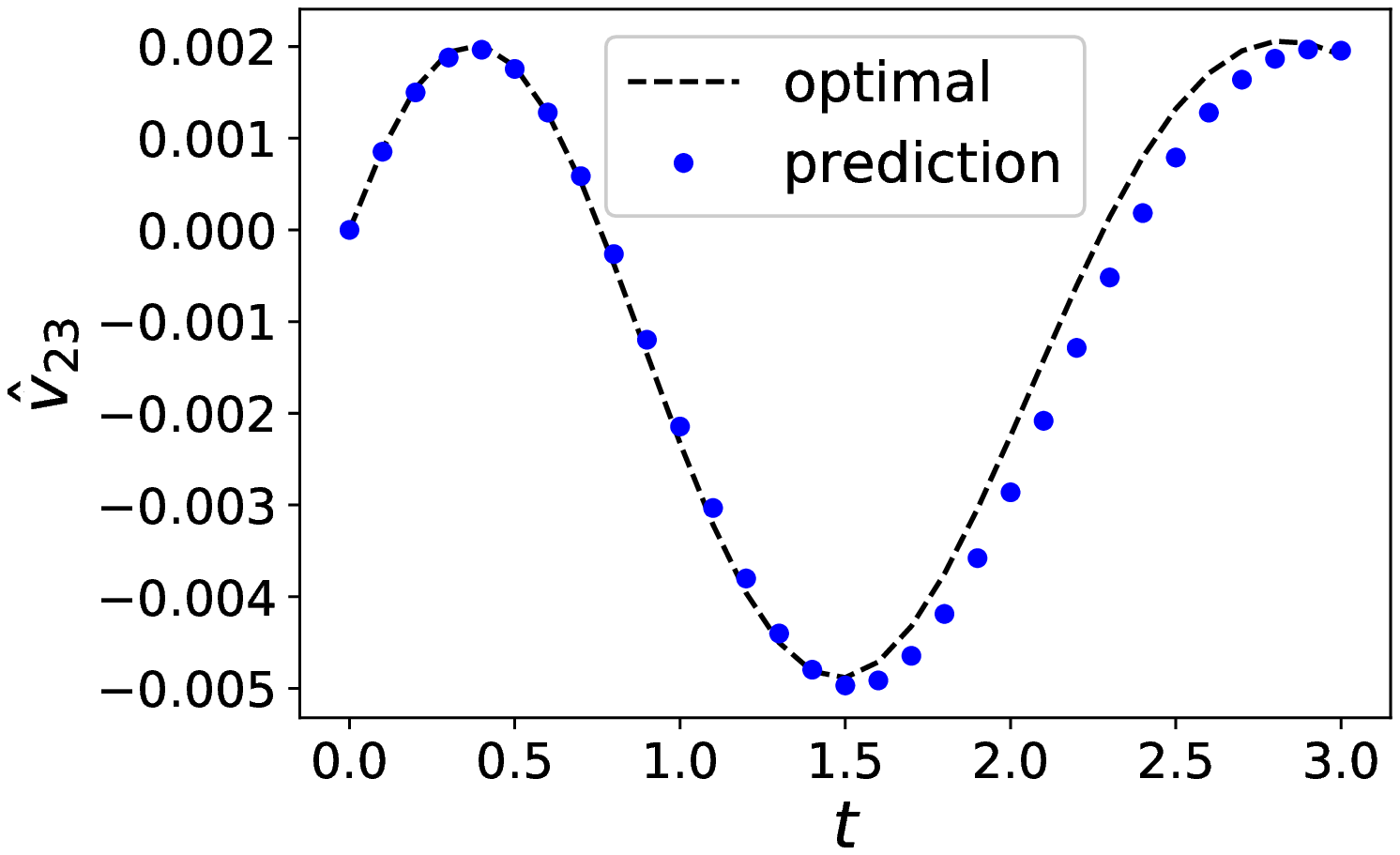}}
	{\includegraphics[width=0.24\textwidth]{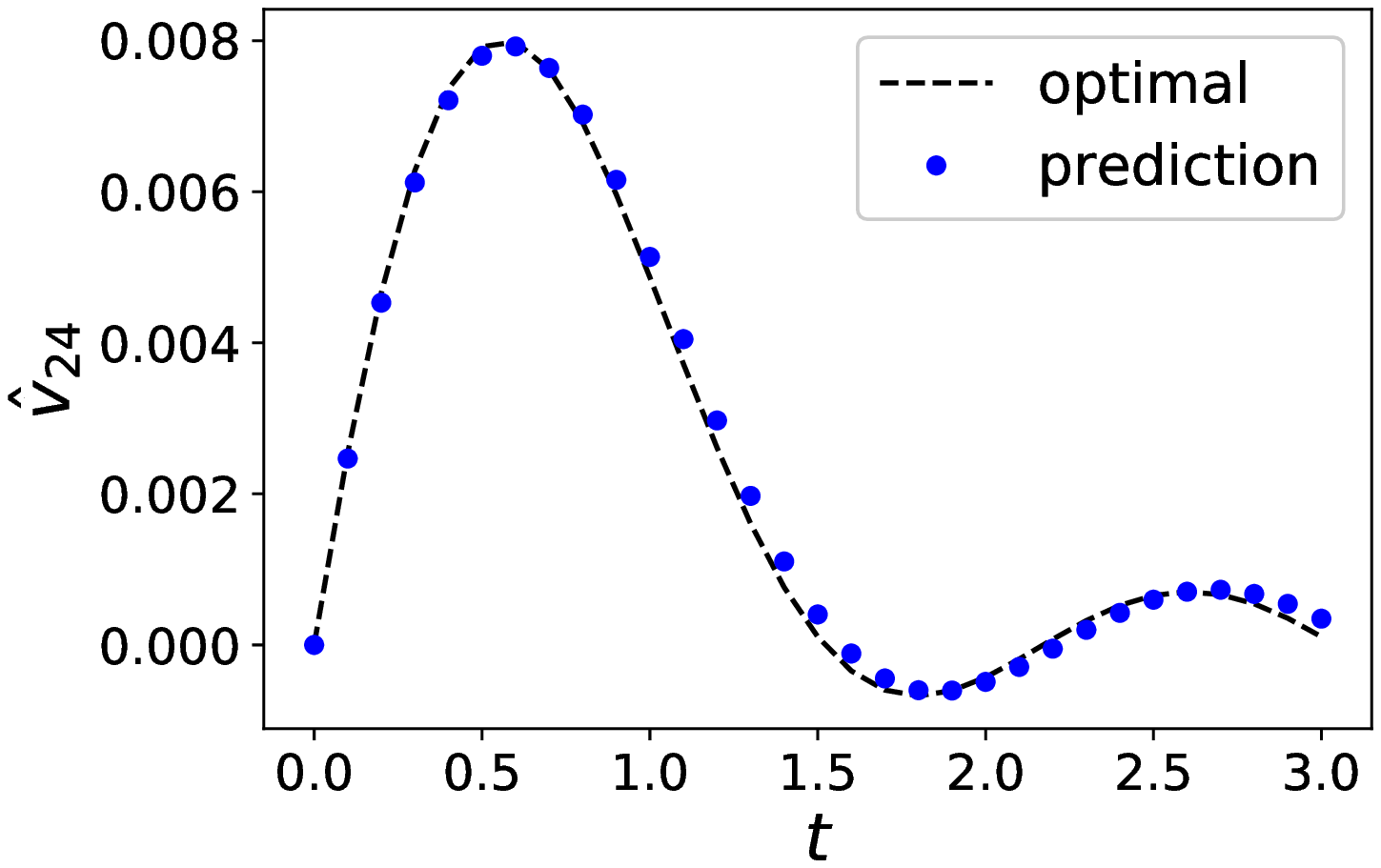}}
	{\includegraphics[width=0.24\textwidth]{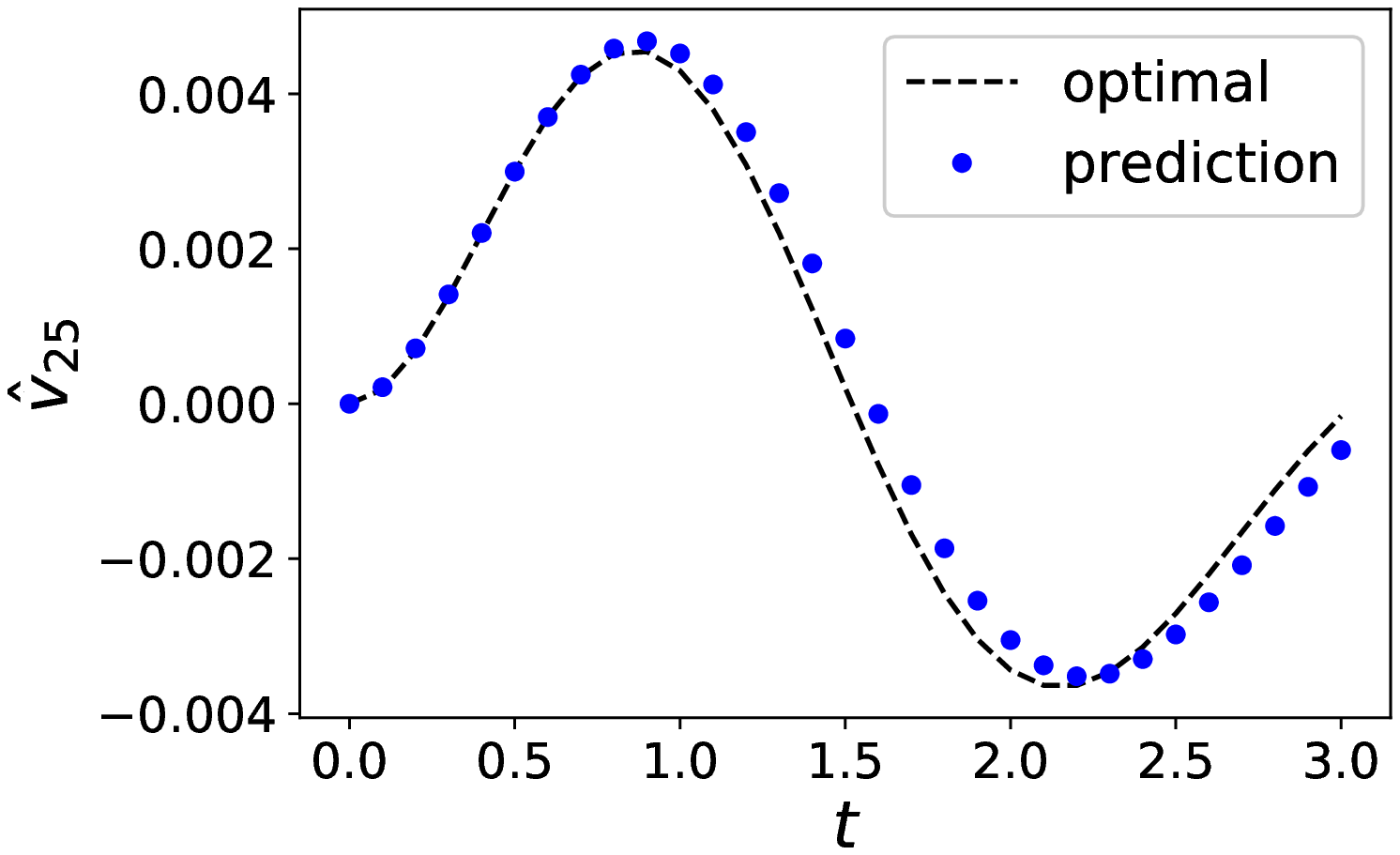}}
	\caption{\small
		Example 5: Evolution of the expansion coefficients for the learned model solution and the projection of the true solution.
	}\label{fig:ex5_coef}
\end{figure}

The comparison between the network model prediction and the exact
solution is shown in Figs.~\ref{fig:ex5_soluC} as solution contours,
and in Fig.~\ref{fig:ex5_soluA} as solution slice profiles at certain locations.
The relative error in the network prediction is shown in
Fig.~\ref{fig:ex5_error}. The time evolution of the modal expansion
coefficients is shown in
Fig.~\ref{fig:ex5_coef}, along with the evolution of the orthogonal
projection coefficients of the exact solution.
Again, we observe good accuracy in the network prediction.
%%%%%%

\section{Conclusion} \label{sec:conclusions}

In this paper, we presented a data-driven framework for learning unknown time-dependent autonomous PDEs, 
based on training of deep neural networks, particularly, those based on residual networks. 
Instead of identifying the exact terms in the underlying PDEs  forms of the unknown PDEs, we proposed to approximately recover the 
evolution operator of the underlying PDEs. 
Since the evolution operator completely characterizes 
the solution evolution,
its recovery allows us to conduct accurate system prediction by recursive use of the operator. 
%However, 
%the evolution operator of PDEs is defined on infinite-dimensional space. 
%In order to reduce the learning problem to finite dimensions, we 
%constrained the evolution operator 
%in a properly
%defined modal space, i.e., Fourier space.
The key to the successful learning of the operator is to reduce the
problem into finite dimension. To this end, we proposed an approach in
modal space, i.e., generalized Fourier space.
Error analysis was conducted to quantify the 
prediction accuracy of the proposed data-driven approach.
We presented a variety of test problems to 
demonstrate the applicability and potential of the method.
More detailed study of its properties, as well as its
applications to more complex problems, will be pursued in future
study.

\bibliographystyle{siamplain}
\bibliography{neural,LearningEqs}

\end{document}